\documentclass{article}
\usepackage{graphicx, amssymb, amsmath, fullpage}
\usepackage[utf8]{inputenc}
\usepackage[LGR,T1]{fontenc}
\usepackage{multicol}
\usepackage{setspace}
\usepackage[para]{footmisc}
\renewcommand{\quote}[1]{\begin{quotation}\noindent #1\end{quotation}}
\newcommand{\Dialogue}{{\it Dialogue}}
\newcommand{\Discourse}{{\it Discourse}}
\newcommand{\SN}{{\it Sidereus Nuncius}}
\usepackage{pgf,tikz,mathrsfs,xcolor,tkz-euclide}
\usetikzlibrary{arrows}
\usetikzlibrary{matrix}
\usetikzlibrary{decorations.pathreplacing,decorations.text,calc}
\providecommand{\tikzfig}[2]{\begin{center}
\begin{tikzpicture}[#2]
#1
\end{tikzpicture}
\end{center}
}

\begin{document}
\title{Galileo, Ignoramus: Mathematics versus Philosophy in the Scientific Revolution}
\author{Viktor Bl{\aa}sj{\"o}}
\date{}
\sloppy

\maketitle

\begin{abstract}\noindent I offer a revisionist interpretation of Galileo's role in the history of science. My overarching thesis is that Galileo lacked technical ability in mathematics, and that this can be seen as directly explaining numerous aspects of his life's work. I suggest that it is precisely {\em because} he was bad at mathematics that Galileo was keen on experiment and empiricism, and eagerly adopted the telescope. His reliance on these hands-on modes of research was not a pioneering contribution to scientific method, but a last resort of a mind ill equipped to make a contribution on mathematical grounds. Likewise, it is precisely {\em because} he was bad at mathematics that Galileo expounded at length about basic principles of scientific method. ``Those who can't do, teach.'' The vision of science articulated by Galileo was less original than is commonly assumed. It had long been taken for granted by mathematicians, who, however, did not stop to pontificate about such things in philosophical prose because they were too busy doing advanced scientific work.\end{abstract}

\begin{multicols}{2}
{\small \tableofcontents}
\end{multicols}

\section{Introduction}

Galileo is overrated, I maintain. I shall systematically go through all his major achievements and all standard arguments as to his alleged greatness, and offer a critical counter-assessment. I go after Galileo with everything but the kitchen sink, but I hasten to add that assembling the case against him is ultimately a means to an end. The stakes are much higher than Galileo's name alone. Galileo is at the heart of fundamental questions: What is the relation between science, mathematics, and philosophy? Between ancient and modern thought? What is the history of our scientific worldview, of scientific method? Galileo is right in the thick of the action on all of these issues. Consequently, I shall use my analysis of Galileo as a fulcrum to articulate a revisionist interpretation of the history of early modern science more broadly. But let's start small, with a simple snapshot of Galileo at work.

\section{Mathematics}

\subsection{Cycloid}
\label{cycloid}

\begin{figure}[pt]\centering
\includegraphics[width=0.65\textwidth]{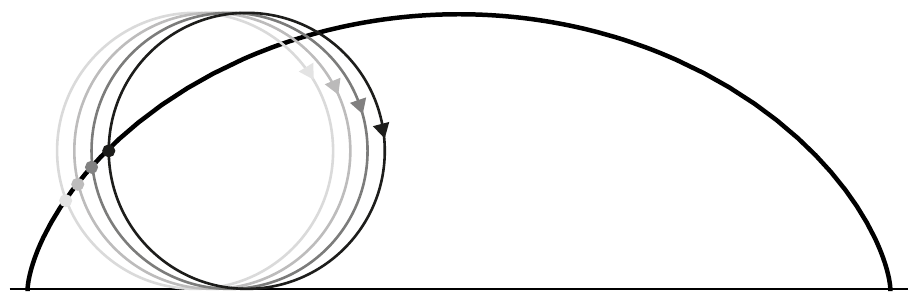}
\caption{The cycloid.}
\label{cycloiddeffig}
\end{figure}

The cycloid is the curve traced by a point on a rolling circle, like a piece of chalk attached to a bicycle wheel (Figure \ref{cycloiddeffig}). Many mathematicians were interested in the cycloid in the early 17th century, including Galileo. What is the area under one arch of the cycloid? That was a natural question in Galileo's time. Finding areas of shapes like that is what geometers had been doing for thousands of years. Archimedes for instance found the area of any section of a parabola, and the area of a spiral, and so on. Galileo wanted nothing more than to join their ranks. The cycloid was a suitable showcase. It was a natural next step following upon the Greek corpus, and hence a chance to prove oneself a ``new Archimedes.''

There was only one problem: Galileo wasn't any good at mathematics. Try as he might, he could not for the life of him come up with one of those clever geometrical arguments for which the Greek mathematicians were universally admired. All those brilliant feats of ingenuity that Archimedes and his friends had blessed us with, it just wasn't happening for Galileo.

Perhaps out of frustration, Galileo turned to the failed mathematician's last resort since time immemorial: trial and error. Unable to crack the cycloid with his intellect, he attacked it with his hands. He cut the shape out of thick paper and got his scales out to have this instrument do his thinking for him. As best as he could gather from these measurements, Galileo believed the area under the cycloid was somewhere near, but not exactly, three times the area of the generating circle.\footnote{\cite[19, 406]{drakeGatwork}.}

This was no way to audition for the pantheon of geometers. Galileo was left red-faced when mathematically competent contemporaries solved the problem with aplomb while he was fumbling with his cutouts. These actual mathematicians proved that the cycloid area was in fact exactly three times the area of the generating circle,\footnote{Roberval in 1634. \cite[232--238]{StruikSB}, \cite{WhitmanCycloidHist}, \cite[350--351]{KlineHistMath}.} even though Galileo had explicitly concluded the contrary on the basis of his little cardboard diorama.

When Galileo heard of others working on the cycloid challenge, he sought help on this ``very difficult'' problem from his countryman Bonaventura Cavalieri, a competent mathematician. ``I worked on it fruitlessly,'' lamented Galileo. ``It needs the mind of a Cavalieri and no other,''\footnote{Galileo to Cavalieri, 24 February 1640, \cite[406]{drakeGatwork}. Cavalieri did not take up the problem---``I too left it aside'' \cite[34]{earlycalcvar}---but Torricelli solved it soon thereafter.} he pleads, tacitly acknowledging his own unmistakably inferior mathematical abilities. It is interesting to contrast this with the very different reaction to the same problem by Galileo's contemporary Ren\'e Descartes, the famous philosopher who was also a vastly better mathematician than Galileo. When Descartes heard of the problem he immediately wrote back to his correspondent that ``I do not see why you attribute such importance to something so simple, that anyone who knows even a little geometry could not fail to observe, were he simply to look.''\footnote{Descartes to Mersenne, 27 May 1638, AT.II.135, \cite[171]{IndivisiblesRevisited}.} He then immediately goes on to give his own proof of the result composed on the spot. Descartes is not famous for his humility, but the fact of the matter is that a number of mathematicians solved the cycloid problem with relative ease, while Galileo was fumbling about with scissors and glue.

In the case of the cycloid, it is an unequivocal fact that Galileo used an experimental approach because he lacked the ability to tackle the problem as a mathematician. If Galileo could have used a more mathematical approach he would unquestionably have done so. It is my contention that what is so glaringly obvious in this case holds for Galileo's science generally. Galileo's celebrated use of experiments in science is not a brilliant methodological innovation but a reluctant recourse necessitated by his pitiful lack of mathematical ability.

The cycloid case also makes it clear {\em why} the mathematically able prefer geometrical proofs to experiments: the latter are hopelessly unreliable. By relying on experiments unchecked by proper mathematics, Galileo got the answer wrong, and not for the first time nor the last. ``Do not think that I am relying on experiments, because I know they are deceitful,''\footnote{\cite[189]{PalmerinoThijssen}. Oeuvres.XI.115.} said Huygens, and all other mathematicians with him. It had always been obvious that mathematics and science can be explored using experiment and observation. As Galileo says: ``You may be sure that Pythagoras, long before he discovered the proof …, had satisfied himself that the square on the side opposite the right angle in a right triangle was equal to the squares on the other two sides''\footnote{Galileo, {\Dialogue}, OGG.VII.75, \cite[85]{wootton}.}---presumably by making numerical measurements on various concretely drawn triangles. But able mathematicians had always known that haphazard trial and error had to be superseded by rigorous demonstration for a treatise to be worth the parchment it is written on. This and only this is why you don't see experimental and numerical data defiling the pages of masterpieces of ancient mathematics and science such as those of Archimedes.\footnote{\S\ref{mathandnature}.}

In this case as in so many others, Galileo's inglorious contribution to the history of thought is to cut off mathematical reasoning at the training-wheels stage, and to air in public what true mathematicians considered unworthy scratch work. He experiments because his intellect comes up short. He cannot reach insights by reason, so he turns to more simplistic, hands-on methods. In physics this ignominious shortcoming has been mistaken for methodological innovation. But in the case of the cycloid its true colours are unmistakable. We see that it is a sign of failure rather than genius. Galileo's empiricism is not science being born; it is science being dumbed down.

\subsection{Mathematicians versus philosophers}
\label{mathvsidiots}

There are mathematicians and then there are other people. Mathematicians have always felt strongly about this. Profound respect for other mathematicians and borderline contempt for anyone else is a persistent theme in their writings. ``I have no doubt that talented and learned mathematicians will agree with me,''\footnote{\cite[6]{CopRevTransl}.} says Copernicus with complete confidence when introducing his new sun-centered astronomical system in 1543, decades before Galileo was born. Copernicus promises to make everything ``clearer than day---at least for those who are not ignorant of the art of mathematics.''\footnote{\cite[24]{CopRevTransl}.} As for non-mathematicians, who cares what those fools think? ``If perchance there are certain idle talkers …, although wholly ignorant of mathematics, …\ [who] dare to reprehend and to attack my work; they worry me so little that I shall …\ scorn their judgements.''\footnote{\cite[7]{CopRevTransl}.} These are the kinds of people who ``on account of their natural stupidity hold the position among philosophers that drones hold among bees.''\footnote{\cite[4]{CopRevTransl}.} ``The studious need not be surprised if people like that laugh at us. Mathematics is written for mathematicians.''\footnote{\cite[7]{CopRevTransl}.}

Johannes Kepler---the best mathematical astronomer of Galileo's age---says the same thing. ``Let all the skilled mathematicians of Europe come forward,''\footnote{\cite[157]{KeplerApologia}.} he implores at one point, confident that mathematical reason compels them to speak with one voice. Those who are not ``skilled mathematicians,'' on the other hand, might as well stay where they are, for they have no credibility as witnesses in scientific matters. Kepler says as much in a remark he labels ``advice for idiots'': ``But whoever is too stupid to understand astronomical science, …\ I advise him that, having dismissed astronomical studies and having damned whatever philosophical opinions he pleases, he mind his own business and betake himself home to scratch in his own dirt patch.''\footnote{Kepler, {\it Astronomia Nova} (1609), 6r, \cite[33]{KeplerAN}.}

This is and always has been the worldview of mathematicians. My thesis in this essay consists in little more than taking them at their word. What does the history of science look like if we adopt the worldview of the mathematicians as a historiographic perspective? If we ``let all the skilled mathematicians come forward'' and leave the philosophers in their ``dirt patch''? Many things do indeed become ``clearer than day---at least for those who are not ignorant of the art of mathematics.''

If we adopt this perspective, virtually Galileo's entire claim to fame evaporates before our eyes. For one thing, mathematicians certainly did not need Galileo to tell them that Aristotle was wrong. Decades before Galileo's first book, Aristotle was already ridiculed by mathematicians. As a contemporary source notes, ``these modern mathematicians solemnly declare that Aristotle's divine mind failed to understand [mathematics], and that as a result he made ridiculous mistakes.''\footnote{Colombe, c.~1611, \cite[223]{GalileoDiscOp}. OGG.III.253--254.} ``And they are right in saying so,'' noted Galileo in the margin next to this passage, ``for he committed many and serious mathematical blunders.''\footnote{Galileo, \cite[223]{GalileoDiscOp}. OGG.III.253--254.}

Hence Descartes's judgement of Galileo's science: ``He is eloquent to refute Aristotle, but that is not hard.''\footnote{Descartes to Mersenne, 11 October 1638, \cite[390]{drakeGatwork}.} This is essentially the thesis of my essay in so many words. It was also, I claim, the standard opinion of Galileo among mathematicians.

Galileo's supposed greatness can only be maintained by ignoring this unequivocal sentiment in the mathematical tradition. Galileo’s claim to fame rests on the assumption that everyone but him was a raving Aristotelian. Galileo himself went out of his way to ensure this framing. His two big books are both dialogues in which Galileo spends hundreds and hundreds of pages in back-and-forth squabbles with a fictive Aristotelian opponent. This is the contrast class Galileo wants us to use when evaluating his achievements. And no wonder. Refuting Aristotle is not hard, so of course Galileo can score some zingers against this feeble opposition. But that’s fish in a barrel, not a scientific revolution.

No other mathematician, ancient or modern, ever engaged with Aristotelian science to anywhere near the same extent as Galileo, let alone devoted the entirety of their main works to refuting it. As Kepler put it: ``The foolish studies of humans have come to such a pitch of vanity that no one's work becomes famous unless he …\ either fortifies himself with the authority of Aristotle, or takes a stand in the battle against him, seeking to show off.''\footnote{Kepler, {\it Astronomiae Pars Optica} (1604), 29, \cite[43]{KeplerO}.} To the mathematicians this was obviously ``foolish vanity'' either way, because they considered it self-evident that Aristotelians were so ridiculous that one should simply dismiss them as so many ``idiots'' in one sentence of the introduction and leave it at that, just as Copernicus and Kepler do in their major books. Descartes did the same, calling followers of Aristotle ``less knowledgable than if they had abstained from study.''\footnote{\cite[57]{DescartesMethod}.} Tycho Brahe---a prominent astronomer in the generation between Copernicus and Galileo---likewise lamented ``the oppressive authority of Aristotle'': ``Aristotle's individual words are worshipped as though they were those of the Delphic Oracle.''\footnote{Tycho Brahe, 1589, \cite[124]{westmanastrrole}.} Kepler once castigated a colleague with the ultimate insult: ``your theory will attract lecturers and philosophers.'' Lest anyone mistake this for a compliment, Kepler explained what he meant: ``it will offer a way out to the enemies of the physics of the sky, the patrons of ignorance.''\footnote{Kepler to David Fabricius, 10 November 1608, \cite[81]{BaumgardtKepler}.}

Such was always the universal opinion of philosophers among mathematicians. Mathematically competent people were united and had nothing but complete contempt for Aristotelian philosophy and the like. By the time Galileo comes along and belabours this point it has been old news for hundreds of years. As Galileo himself says: ``the philosophers of our times …\ philosophize …\ as men of no intellect and little better than absolute fools.''\footnote{Galileo, \cite[10]{SantillanaCrime}.} Precisely. Which is why there is little value in writing several thick books hammering home this point and little else, which is what Galileo did. That’s just beating a dead horse as far as the mathematicians are concerned.

The standard view among modern scholars is very different: ``Galileo's significance for the formation of modern science lies partly in his discoveries and opinions in physics and astronomy, but much more in his refusal to allow science to be guided any longer by philosophy.''\footnote{\cite[xvii]{VShortGalileo}.} But in this refusal Galileo is just following the mathematicians, who had always maintained that ``only mathematics can provide sure and unshakeable knowledge'' and that other ``divisions of theoretical philosophy should rather be called guesswork than knowledge.''\footnote{Ptolemy, {\it Almagest}, I.1, \cite[36]{AlmagestToomer}.} The only difference is that, while competent mathematicians took this for granted and got on with real science, Galileo stopped to dwell on the matter for hundreds of pages.

``It is evident that Galileo admired Archimedes as his mentor,'' while ``no philosopher had ever taken the slightest interest in Archimedes.''\footnote{\cite[10]{drakeGS}.} This was the one thing that set Galileo apart from the philosophical hordes of his day. ``What was novel in Galileo's approach, when contrasted with that of his Aristotelian opponents, was …\ his faith in the relevance of mathematics.''\footnote{\cite[11]{sheaGrev}.} It is true that Galileo got this much right. He had a crush on the right person. Insofar as Galileo ever accomplished anything, it stems from his infatuation with Archimedes. His only regret was that his own mathematical abilities were so vanishingly small compared to those of his hero, as he frankly admitted: ``Those who read his works realise only too clearly how inferior are all other minds compared with Archimedes', and what small hope is left of ever discovering things similar to those he discovered.''\footnote{Galileo, \cite[1]{sheaGrev}. OGG.I.215--216.}

This is Galileo's unique position in the history of science. He is the worst of the mathematicians and the best of the philosophers. He is too incompetent to do mathematics but competent enough to admire it. Relative to the mathematicians he is a dud who has to spend his time spelling out obvious philosophical points because he lacks the ability to do real mathematics. Relative to the philosophers he is a revelation and a watershed, articulating for the first time in philosophical prose many central tenets of the modern scientific worldview.

Evaluating Galileo thus comes down to understanding the historical relation between mathematics and philosophy. If you believe that mathematicians from Greek times onwards considered it obvious that mathematics encompassed an expansive, daringly creative, empirically informed, interconnected study of all quantifiable aspects of the world, and that philosophers talking about mathematics and science are basically just derivative commentators, then Galileo is nobody. If you believe that mathematics before Galileo was a technocratic enterprise that was intrinsically compartmentalised and limited to traditionally defined specific tasks, and that the empirical-mathematical method of modern science was a great ``conceptual'' revelation coming from philosophy, then you may indeed be inclined to view Galileo as the ``father of modern science.'' It is my contention that the former is the right way to understand the history of scientific thought.

\begin{figure}[p!]\centering
\includegraphics[width=0.75\textwidth]{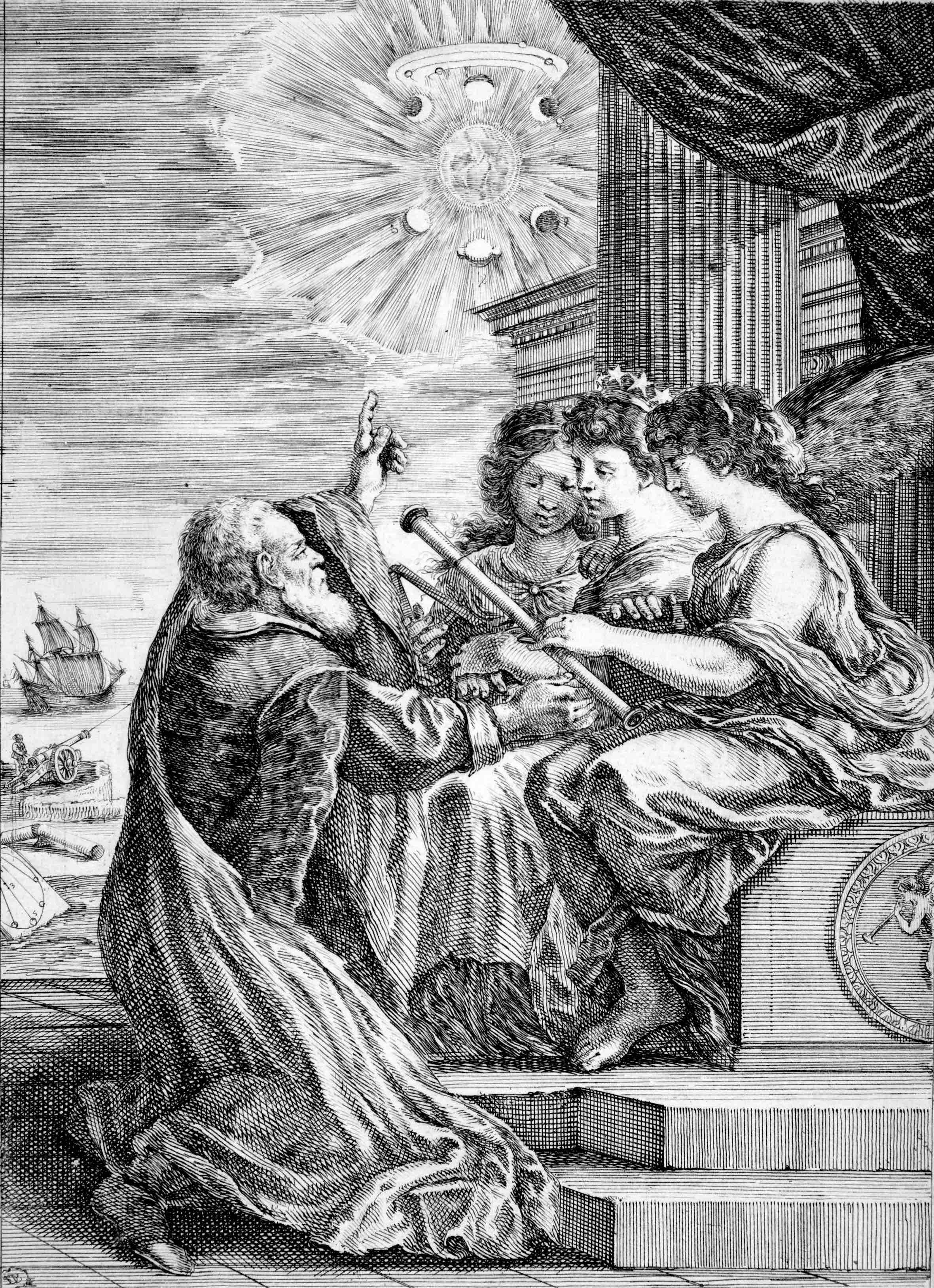}
\caption[Galileo kneeling before muses.]{Galileo kneeling before personifications of Astronomy, Perspective, and Mathematics (far side, with compasses, hiding her face in shadows and averting her eyes, perhaps symbolising Galileo's uneasy relationship with mathematics). Frontispiece for {\it Opere di Galileo Galilei}, Bologna, 1656.
}
\label{G3musesfig}
\end{figure}

\subsection{Professor}

Galileo was a ``Professor of Mathematics'' for some twenty years. But we must not let the title fool us. The position had nothing to do with the research frontier in the field. In modern terms it is more comparable to that of high school teacher. Galileo taught basic and practical courses,\footnote{\cite[286--290]{vallerianisymp}, \cite[34--35, 45--46]{drakeGatwork}.} and his lectures were unoriginal and usually cribbed from standard sources.\footnote{\cite[12, 23, 34--35, 51]{drakeGatwork}.} His mathematical lectures went no further than elementary Euclidean geometry. He also had to teach a basic astronomy course ``mainly for medical students, who had to be able to cast horoscopes.''\footnote{\cite[35--36]{drakeGatwork}.} They ``needed it to determine when not to bleed a patient'' and the like.\footnote{\cite[46]{heilbron}.} Perhaps Galileo didn't mind, for he seems to have been quite open to astrology judging by the fact that he cast horoscopes for his own family members and friends without renumeration.\footnote{\cite{RutkinJHA}, \cite[55]{drakeGatwork}, \cite[68]{heilbron}, \cite{kollerstromsymp}.} Alas, he did not enjoy much success as an astrologer: “In 1609, Galileo …\ cast a horoscope for the Grand Duke Ferdinand I, foretelling a long and happy life. The Duke died a few months later.”\footnote{\cite[19]{camcomp}.}

As soon as Galileo had made a name for himself he eagerly resigned his lowly university post in favour of a court appointment. This freed him from teaching duties and boosted his finances, but Galileo also had an additional demand:
\quote{I desire that in addition to the title of ``mathematician'' His Highness will annex that of ``philosopher''; for I may claim to have studied more years in philosophy than months in pure mathematics.\footnote{Galileo, negotiating his Tuscan court appointment in 1610, \cite[64]{GalileoDiscOp}.}}
This is traditionally taken as a request for a kind of promotion: {\em in addition} to being a great mathematician, Galileo also wanted recognition in philosophy, which in some circles was considered more prestigious and in any case included what today is called science (then ``natural philosophy''). But I think a more literal reading is in order. Galileo is not only declaring himself a philosopher; he is also confessing his ignorance in mathematics. Taken literally, his statement implies that he could not have spent more than two or three years on pure mathematics---which indeed sounds about right considering his documented mediocrity in this field.

\subsection{Sector}\label{sector}

\begin{figure}[pt]\centering
\includegraphics[width=0.25\textwidth]{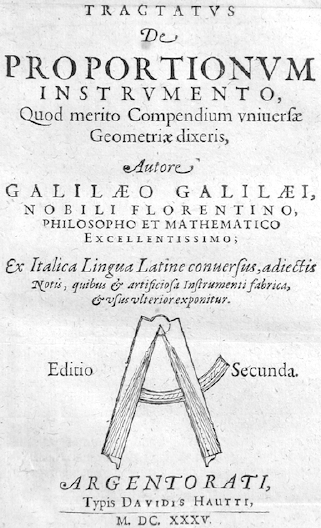}  \includegraphics[width=0.18\textwidth]{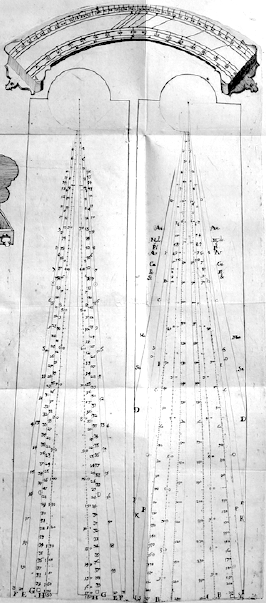}
\caption[Galileo's sector.]{Galileo's ``geometric and military compass'': a slide-rule-style calculation device.
}
\label{Gsectorfig}
\end{figure}

For his teaching, Galileo developed an instrument---his ``geometric and military compass'' or sector (Figure \ref{Gsectorfig})---that could be used to facilitate certain computations, somewhat like a slide rule. Is this a great testament to his mathematical genius? Should Galileo be credited with inventing ``the first mechanical computing device,''\footnote{\cite[104]{DrakeSciAmSector}.} as one sensationalistic headline has it? Obviously not, since far more sophisticated mechanical computing devices were developed in antiquity (such as the astrolabe and the Antikythera mechanism, although they were devised for different problems). Nor is it clear to what extent, if any, Galileo's instrument was original in its time. ``It has always been supposed that Galileo took over a calculating device already in service,''\footnote{\cite[106--107]{DrakeSciAmSector}.} and various at least partial precursors are documented,\footnote{\cite{DrakeSciAmSector}, \cite{Gcompass}, \cite[vi--x]{RobertsonMathInstr}.} but it is hard to say anything definitive since very few instruments from this era have survived. One recent study declines to pronounce on Galileo's degree of originality, noting only that, at any rate, ``students did not necessarily buy the compass because it was better or more original than other similar devices one could buy in Germany, France, or England. …\ They trusted Galileo the way the residents of a certain neighborhood trust the local baker.''\footnote{\cite[12]{BiagioliInstr}.} In any case it is undisputed that a ``very similar instrument,'' developed in England ``entirely independently of Galileo's work,'' was extensively described in print eight years before Galileo published anything about his device.\footnote{Drake, \cite[10]{Gcompass}. Thomas Hood, {\it The Making and Use of the Geometricall Instrument, Called a Sector}, London, 1598.}

Galileo is upfront about the target audience of his sector: it is ``suited to …\ first youthful studies''\footnote{\cite[39]{Gcompass}.} and an aid ``for those who cannot easily manage numbers.''\footnote{\cite[21]{Gcompass}.} The fact that the latter group arguably includes Galileo himself is perhaps not coincidental. Ancient Greek folklore has it that some king found geometrical study too onerous and demanded to be brought to its fruits by an easier route, only to be reprimanded that ``there is no royal road to geometry.'' Many a mathematician has shared this beautifully egalitarian sentiment. But not Galileo, who leapt at any opportunity to kiss the feet of princes. ``Not entirely improper was the request of that royal pupil who sought …\ an easier and more open road,''\footnote{\cite[41]{Gcompass}.} Galileo declares in his instrument manual. While Galileo encouraged this shortcut mindset, more committed mathematicians took a principled stance to the contrary: “The true way of Art is not by Instruments, but by demonstration: …\ it is a preposterous course of vulgar Teachers, to beginne with Instruments, and not with the Sciences, and so in stead of Artists, to make their Schollers onely doers of tricks, and as it were jugglers.”\footnote{William Oughtred, {\it To the English gentrie, and all others studious of the mathematicks} (1634), 27--28.} The divide remains as current as ever. Anyone familiar with modern mathematics classrooms have observed this unmistakeable relationship: the weaker the student (or teacher), the greater their reliance on calculators. From this point of view, Galileo's instrument can be seen as one more indication of that aversion to serious mathematical work that characterises his career as a whole.

\subsection{Book of nature}
\label{bookofnature}

Here is one of Galileo's most quoted passages:
\quote{This grand book, the universe, …\ stands continually open to our gaze. But the book cannot be understood unless one first learns to comprehend the language and read the letters in which it is composed. It is written in the language of mathematics, and its characters are triangles, circles, and other geometric figures without which it is humanly impossible to understand a single word of it.\footnote{Galileo, {\it Assayer} (1623), \cite[237--238]{GalileoDiscOp}.}}
Like so much else Galileo gets credit for, this idea had obviously been a truism among mathematicians for thousands of years. Pythagoras had put it succinctly more than two millennia before: ``All is number.''\footnote{\cite[Chapter 4]{Cornelli}.} Plato agreed, saying ``God is always doing geometry.''\footnote{\cite[I.387]{LoebMath}, slightly paraphrased.} Mathematicians have lived by these words ever since.\footnote{\S\ref{mathandnature}.} It was obvious even to for instance Grosseteste, a medieval bishop, that ``the usefulness of considering lines, angles, and figures is very great, since it is impossible to understand natural philosophy without them.''\footnote{Grosseteste, first half of 13th century, \cite[IX.59]{Grosseteste1912}, \cite[385]{GrantSourceBook}.}

Galileo's version of this platitude is better read as a frank confession of his reactionary restriction to simplistic mathematics. ``Triangles and circles'' are hardly the cutting edge. Ellipses, for instance, are passed over in silence by Galileo, not just in this quote but throughout all his works, even though Kepler had already demonstrated many years before that planets move in ellipses about the sun.\footnote{Kepler, {\it Astronomia Nova} (1609).}
\quote{Unfortunately, for reasons that remain unclear or controversial, Galileo never did pay the proper attention to Kepler's writings and so ignored or neglected the elliptical nature of planetary orbits.\footnote{\cite[36]{findefending}}}
There is nothing ``unclear'' about it. The reason Galileo ignored Kepler's writings is abundantly obvious: He didn't have anywhere near the mathematical ability to understand them. Kepler had already anticipated as much:
\quote{It is extremely hard these days to write mathematical books. …\ There are very few suitably prepared readers. …\ How many mathematicians are there who put up with the trouble of working through the Conics of Apollonius of Perga?\footnote{Kepler, {\it Astronomia Nova} (1609), 2r, \cite[17]{KeplerAN}.}}
Not many indeed, and certainly not Galileo. Yet Kepler saw no other option than to write mathematics all the same and ``let it await its reader for a hundred years.''\footnote{\cite[170]{KeplerEpiHar}.}

Compared to Kepler's visionary science, Galileo's book of nature is a children's book sticking to baby mathematics like ``triangles and circles.'' The traditional way to read Galileo's phrase that without these it is ``humanly impossible to understand'' nature is to see it as stipulating a minimum knowledge requirement: if you don't know math, you're screwed. But it can also be read as a confession of an upper limit on Galileo's own understanding: anything beyond ``triangles and circles'' and I'm lost. This would also explain Galileo's equally reactionary rejection of ``infinities and indivisibles, the former incomprehensible to our understanding by reason of their largeness, and the latter by their smallness.''\footnote{\cite[34]{galileo2newsci2ndedWT}, OGG.VIII.73.} Virtually all competent mathematicians in Galileo's time were eagerly exploring infinitesimal methods, prefiguring what we call calculus today. But Galileo found all of this ``incomprehensible.'' When he did try to dabble with such issues his arguments were ``rejected immediately by Descartes as sophistic, and then on similar grounds by Bonaventura Cavalieri, Galileo's own disciple, perhaps the greatest expert on indivisibles in the early seventeenth century.''\footnote{\cite[VII.76]{WallaceColl}. See also \cite[153--154]{buttspitt}.} Such shortcomings undermine core parts of his physics, as we shall see in the next section.

\subsection{Calculus of speeds}
\label{Calculusofspeeds}
\label{Ginfintesimals}

\begin{figure}\centering
\includegraphics[width=\textwidth]{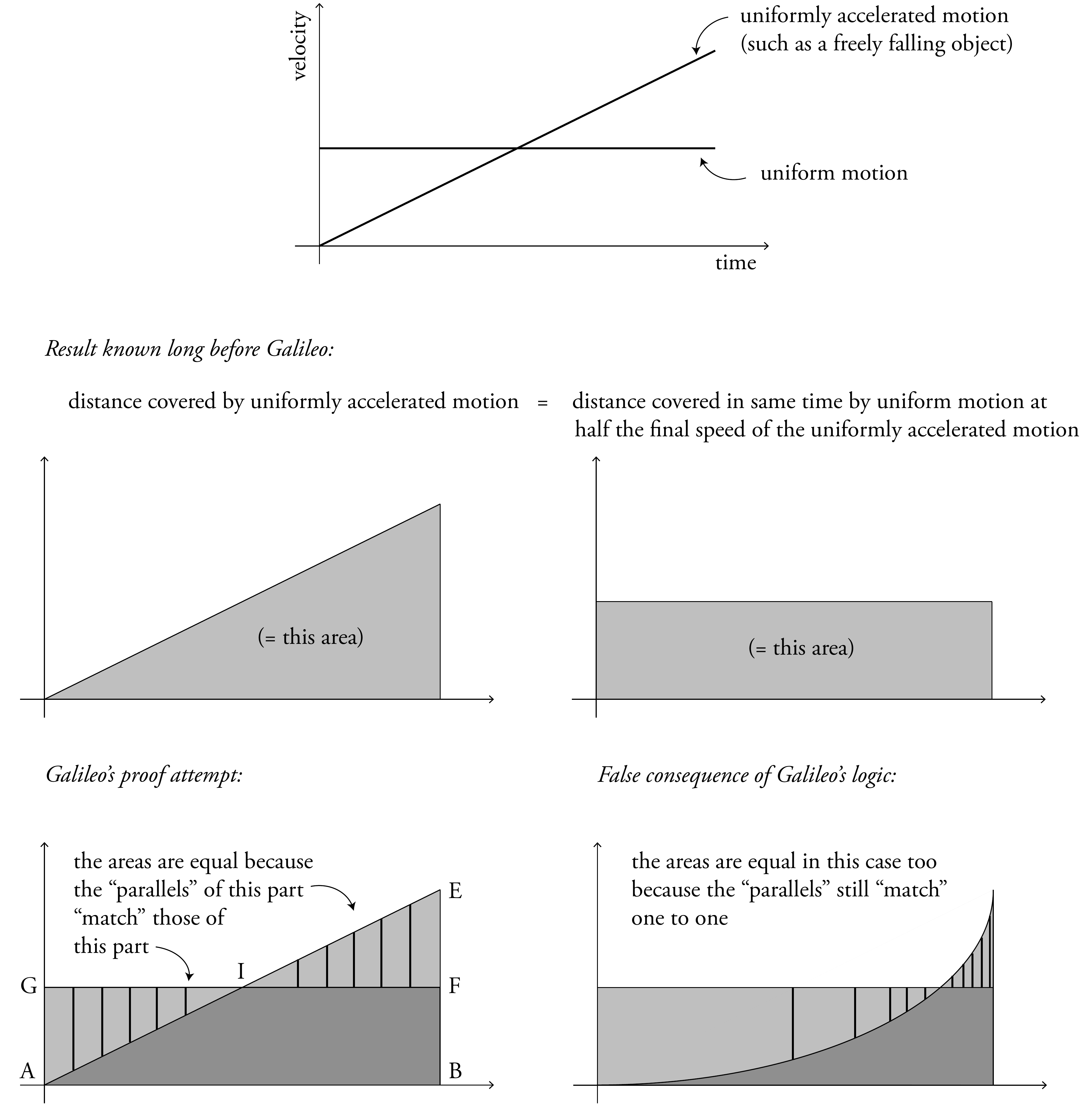}
\caption{Galileo's dubious way of comparing two motions.}
\label{velocitytimeparadox}
\end{figure}

Galileo's ``technique for dealing with problems involving infinitesimals …\ is incompatible with basic features of integration …\ and even leads, if explicated systematically, to absurd consequences.''\footnote{\cite[237]{LimitsPreclasMech}.} The core issue is illustrated in Figure \ref{velocitytimeparadox}. Galileo tries to prove areas equal by ``matching'' the lines that make them up one by one:
\quote{The aggregate of all parallels contained in the quadrilateral [ABFG is] equal to the aggregate of those included in triangle AEB. For those in triangle IEF are matched by those contained in triangle GIA, while those which are in the trapezium AIFB are common. Since each instant and all instants of time AB correspond to each point and all points of time AB, from which points the parallels drawn and included within triangle AEB represent increasing degrees of the increased speed, while the parallels contained within the parallelogram represent in the same way just as many degrees of speed not increased but equable, it appears that there are just as many momenta of speed consumed in the accelerated motion according to the increasing parallels of triangle AEB, as in the equable motion according to the parallels of the parallelogram GB. For the deficit of momenta in the first half of the accelerated motion (the momenta represented by the parallels in triangle AGI falling short) is made up by the momenta represented by the parallels of triangle IEF.\footnote{\cite[165--166]{galileo2newsci2ndedWT}, OGG.VIII.208--209.}}
The problem, as shown in Figure \ref{velocitytimeparadox}, is that this reasoning invites false results because ``parallels'' can be ``matched'' also for shapes that do not have the same area. ``Neither Galileo's manuscripts nor his published texts indicate that he noticed this clearly nonsensical implication of the proof technique.''\footnote{\cite[250--251]{LimitsPreclasMech}.} Many of his contemporaries were more alert, including Cavalieri, who ``criticizes
Galileo for not having emphasized that the indivisibles have to be taken as equidistant.''\footnote{\cite[251]{LimitsPreclasMech}. Bonaventura Cavalieri to Galileo, June 28, 1639. OGG.XVIII.67.} Cavalieri's reaction is the most charitable option: Galileo neglected to ``emphasise'' (or even mention, in fact) some implicit assumptions that rule out the false applications of his method of proof.

A better way of thinking about these kinds of problems is to view areas as composed of thin rectangles rather than lines. The correct conclusion of Figure \ref{velocitytimeparadox} works because the two areas can be sliced into thin rectangles that correspond to an equal rectangle in the other figure. The false conclusion of Figure \ref{velocitytimeparadox}, on the other hand, doesn't go through this way, because while each rectangle ``fits'' in a natural place in the other area (just as the lines do, as illustrated), the rectangles positioned this way do not cover the full figure, but rather leave gaps between them. Galileo's proof ought to have included some clarifications or precautions to specify in what way his method avoids false conclusions on these kinds of grounds.

If we look only at the examples in Figure \ref{velocitytimeparadox}, it may seem that we are interpreting Galileo quite maliciously. It is quite evident that the two triangles are equal, while the areas under the curved velocity diagram are not. Surely, we may be tempted to say, it is clear what Galileo must have meant, even if he may not have expressed it with great precision.

\begin{figure}\centering
\includegraphics[width=0.9\textwidth]{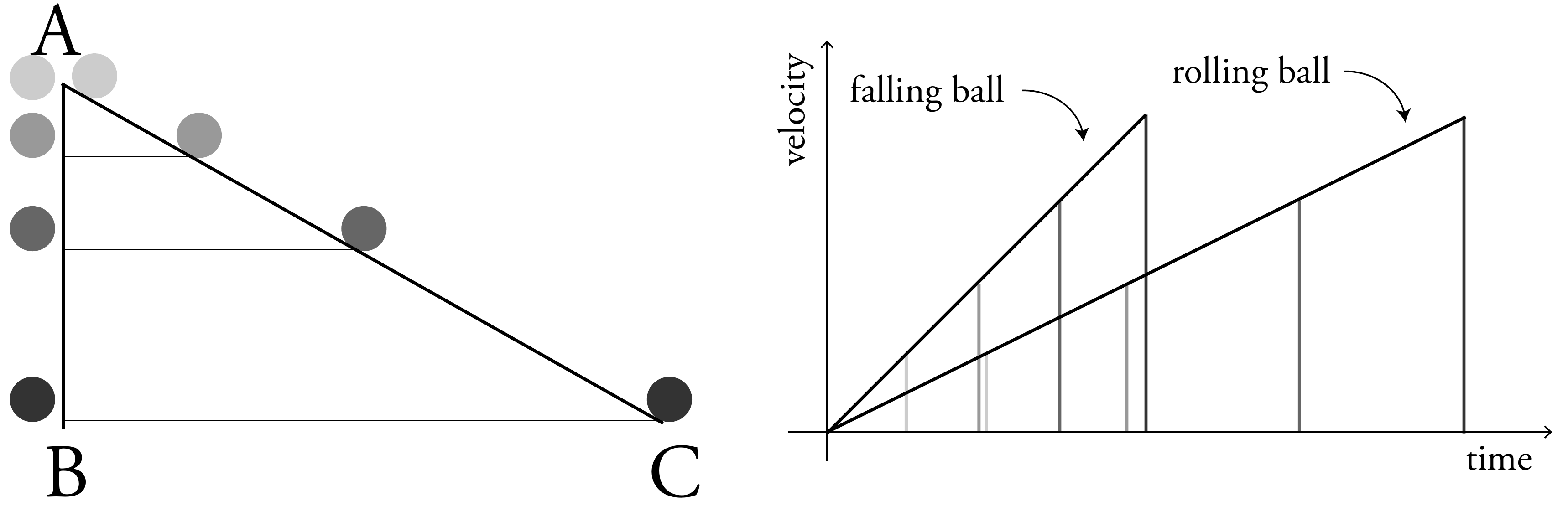}
\caption{A falling ball compared to a ball rolling down an inclined plane.}
\label{velocitytimeparadoxplane}
\end{figure}

\begin{figure}\centering
\includegraphics[width=0.9\textwidth]{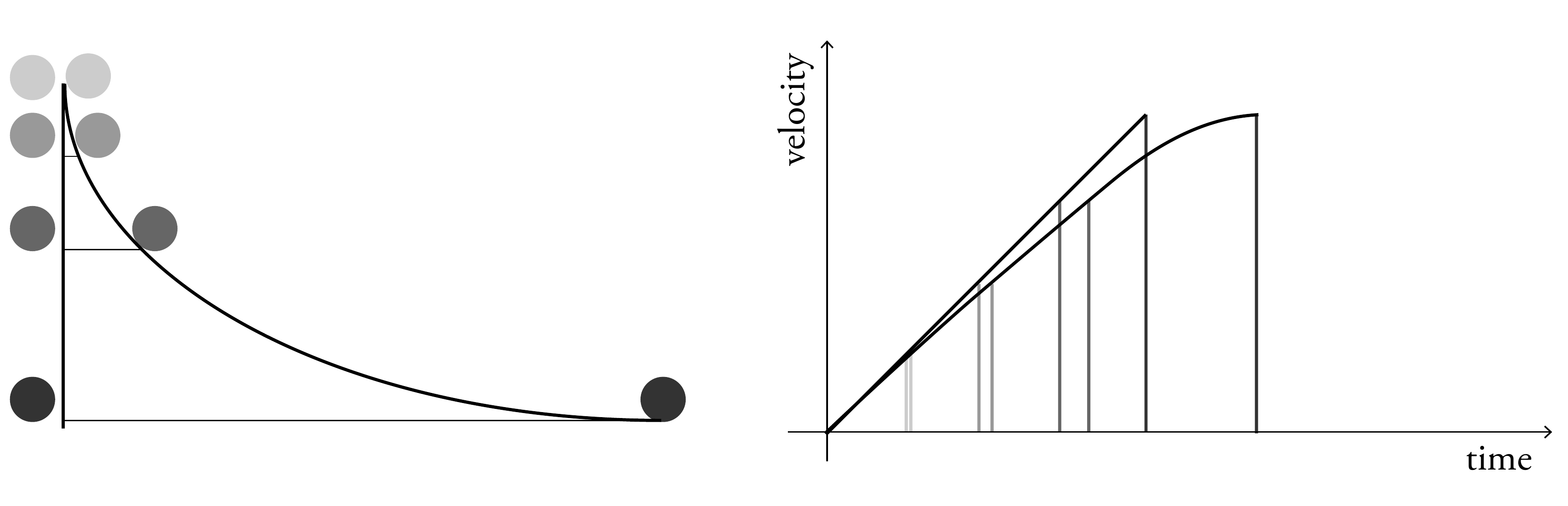}
\caption{Variant of Figure \ref{velocitytimeparadoxplane}. The ramp along which the ball is rolling is no longer straight.}
\label{velocitytimeparadoxcurve}
\end{figure}

But it gets worse. Galileo applies his ``matching'' technique again later, in a way that very much complicates the picture. He applies this method of proof when comparing the motion of a falling body to that of a ball rolling down an inclined plane (Figure \ref{velocitytimeparadoxplane}). Galileo wants to show that time it takes the ball to complete the motion is proportional to the distance it has to travel: for example, if $AC=2AB$, then the rolling ball needs twice as long as the falling ball to reach the ground. Galileo tries to derive this result from the fact that the two balls have the same velocity when they are the same vertical distance from the starting point. This fact establishes a one-to-one correspondence or ``matching'' between the speeds of the falling and the rolling body. Galileo claims that the result follows from this:
\quote{If …\ parallels from all points of the line AB are supposed drawn as far as line AC, the momenta or degrees of speed at both ends of each parallel are always matched with each other. Thus the two spaces AC and AB are traversed at the same degrees of speed. But …\ if two spaces are traversed by a moveable which is carried at the same degrees of speed, then whatever ratio those spaces have, the times of motion have the same [ratio]. Therefore the time of motion through AC is to the time through AB as the length of plane AC is to the length of vertical AB, which was to be demonstrated.\footnote{\cite[176]{galileo2newsci2ndedWT}, OGG.VIII.216--217.}}
But again the rules of this method is unclear. From Galileo's words, there seems to be nothing stopping us from applying the same proof, word for word, to the situation of a non-straight path, as in Figure \ref{velocitytimeparadoxcurve}. But this clearly produces false results, for in such cases the time it takes to complete the motion is no longer proportional to the distance travelled. So Galileo's proof technique can be used to prove false results again.

The example of Figure \ref{velocitytimeparadoxplane} makes it much harder to argue that Galileo's method can be fixed by adding some common-sense implicit assumptions or restrictions. Galileo clearly thought his method could be applied even in cases where the time-velocity diagrams have different shape and different spacing of the matched lines, as in Figure \ref{velocitytimeparadoxplane}. In terms of the area interpretation mentioned above, this reasoning is dependent on matching rectangles of one area to thinner rectangles of another area, but for Galileo's reasoning to work the ratio of the thicknesses of these rectangles must remain the same throughout. This is what sets the permissible conclusion of Figure \ref{velocitytimeparadoxplane} apart from the impermissible conclusion of Figure \ref{velocitytimeparadoxcurve}. But note that this is a {\em global} property: it has to do with how {\em all} the slices relate to one another---the ratio must be the same for {\em all} pairs of rectangles. This seems to be directly incompatible with Galileo's text. Galileo quite clearly highlights only the {\em local, pairwise} matching of velocities and draws his conclusion from there. There is no indication in Galileo's text that he considers a constant-ratio requirement essential to the legitimate application of his proof technique.

Hence we see why scholars have argued that ``closer inspection of [Galileo's] proof [of Figure \ref{velocitytimeparadox}] in the light of [his use of this technique in the case of Figure \ref{velocitytimeparadoxplane}] indeed makes it clear that this proof also works equally well for nonuniformly accelerated motion,'' and therefore leads to false results.\footnote{\cite[250]{LimitsPreclasMech}.}

\section{Physics}

\subsection{Fall and weight}
\label{fallandweight}

Galileo often portrays himself as defeating obstinate philosophers who would rather cling to the words of Aristotle than believe empirical evidence made clear through  experimental demonstration. Thus he imagines his enemies to say things like: ``You have made me see this matter so plainly and palpably that if Aristotle's text were not contrary to it …\ I should be forced to admit it to be true.''\footnote{\cite[125]{galileodialogueML}.} It makes life easy for Galileo to pretend that his enemies are fools blind to facts. In reality, Galileo's ``anti-Aristotelian polemics …\ were directed only at straw men.'' Galileo concocted these caricature Aristotelians in order to ``have fun with them, let them play the buffoon in his dialogues, and thus enhance his own image in the eyes of his readers.''\footnote{\cite[I.377]{WallaceColl}.}
Galileo's ploy was well calculated. It tricks many of his readers to this day into believing the fairytale of Galileo the valiant knight singlehandedly fighting for truth in world beset by dogmatism. But ``excessive claims for Galileo the dragon slayer have to be muted,''\footnote{\cite[142]{PittBG}.} to say the least.

This is clearly seen in the case of a famous question: do heavy objects fall faster than lighter ones? Aristotle had answered: yes. Twice as heavy means twice as fast, according to Aristotle.\footnote{E.g.\ {\it De Caelo}, I.6: ``A given weight moves a given distance in a given time; a weight which is as great and more moves the same distance in a less time, the times being in inverse proportion to the weights. For instance, if one weight is twice another, it will take half as long over a given movement.''} Legend has it that Galileo shocked the world when he dropped some balls of differing weights from the tower of Pisa and revealed them to fall at the same speed. But the notion that it required some kind of radical conceptual innovation by a scientific genius like Galileo to realise that one could test the matter by experiment is ludicrous. Of course one can drop some rocks and see if it works: this much has been obvious to anyone since time immemorial. Indeed, Philoponus---an unoriginal commentator---had clearly and explicitly rejected Aristotle's law of fall by precisely such an experiment more than a thousand years before Galileo.\footnote{\cite[172--173]{ClagettGS}.} So if you want to believe that experimental science and empirical verification was a radical new insight then this puts you in quite a pickle. If that was a revolution, then why is it found for the first time in this mediocre commentator from the 6th century? If people like Philoponus had the key to science, then why did they sit around and write commentary upon commentary on Aristotle? Their contributions to mathematics and science is otherwise zero. How likely is it that such elementary scientific principles eluded generations of the best mathematical minds the world has ever seen, in the age of Archimedes, only to be then discovered by derivative and subservient thinkers in an age where the pinnacle of mathematical expertise extended little further than the ability to multiply three-digit numbers?

Long before Philoponus, for that matter, Lucretius had clearly stated that, in the absence of air resistance, all objects fall at the same speed regardless of weight,\footnote{Lucretius, {\it De rerum natura}, II:225--239.} a result many still believe Galileo was the first to discover.\footnote{E.g.\ \cite[4]{Gowers}.} As ever, Galileo gets credit for elementary ideas that are thousands of years old.

Nor was Galileo original in his own age. Have you ever heard anyone calling Benedetto Varchi ``the father of modern science''? Yet here is his statement of ``Galileo's'' great insight, expressed two decades before Galileo was even born:
\quote{The custom of modern philosophers is always to believe and never to test that which they find written in good authors, especially Aristotle. …\ [But it would be] both safer and more delightful to otherwise and sometimes descend to experience in some cases, as for example in the motion of heavy bodies, in which both Aristotle and all other philosophers without ever doubting the fact, have believed and affirmed that according to how much a body is more heavy, by so much more [speed] will it descend, the test of which shows it not to be true.\footnote{Benedetto Varchi, 1544, \cite[9]{SettleExpr}.}}
The fact that run-of-the-mill humanists like Varchi were arguing this way long before Galileo shows what an obvious idea it was.

It is not clear whether Galileo did in fact carry out such an experiment from the tower of Pisa when he was teaching at the university there, as legend would have it.\footnote{\cite[19--20]{drakeGatwork}.} Personally I find the story plausible inasmuch as a field day throwing rocks surely had great appeal to a professor whose strong suit certainly wasn't pure mathematical theory. But be that as it may. It is in any case perfectly clear that many people carried out such experiments around that time, independently of Galileo. Simon Stevin, for example, certainly did, and published his results years before Galileo made his experiment, if indeed the latter ever took place. Stevin dropped lead balls of different weights from a height of 30 feet onto a surface that would make much noise, and noted that they banged to the ground in unison.\footnote{\cite[511]{StevinWorksI}.} Similarly, Mersenne was already dropping weights out out Parisian chamber windows before he heard of Galileo's work.\footnote{\cite[III.274--275]{MersenneCorr}, \cite[249]{VergaraCaffarelli}.}

``To crown the comedy, it was an Aristotelian, Coresio, who in 1612 claimed that previous experiments had been carried out from too low an altitude. In a work published in that year he described how he had improved on all previous attempts---he had not merely dropped the bodies from a high window, he had gone to the very top of the tower of Pisa. The larger body had fallen more quickly than the smaller one …, and the experiment, he claimed, had proved Aristotle to have been right all the time.''\footnote{\cite[83]{Butterfield}.} Girolamo Borro, another Aristotelian philosopher, similarly investigated whether lead falls faster than wood. ``We took refuge in experience, the teacher of all things,'' he says, and hence ``threw these two pieces of equal weight from a rather high window of our house at the same time. …\ The lead descended more slowly. …\ Not only once but many times we tried it with the same results.''\footnote{Girolamo Borro, {\it De motu et levium} (1575), \cite[7]{SettleExpr}.} These philosophers were united in their emphasis on actually testing their theories practically and even in taking care to mitigate possible sources of error in the measurements. So much for Aristotelians hating experiments. On the contrary, appeal to experiments were commonplace long before Galileo.

But these people got the result wrong, you say. Maybe they didn’t experiment at all, or if they did they messed it up somehow. Well, so did Galileo. Galileo's own early experiments also produced the wrong results.\footnote{Just as they did for the cycloid, \S\ref{cycloid}.} ``If an observation is made, the lighter body will, at the beginning of the motion, move ahead of the heavier and will be swifter,''\footnote{Galileo, \cite[38]{DrabkinDrake}.} Galileo reports, but if the fall is long enough the heavier body will eventually overtake it. Galileo devotes a full chapter to following this up, ``in which the cause is given why, at the beginning of their natural motion, bodies that are less heavy move more swiftly than heavier ones.''\footnote{Galileo, \cite[106]{DrabkinDrake}.} So when Galileo started experimenting on this he got the wrong result, and he also believed himself to a have a good theory “explaining” those false results. We now know that the true cause for the errant results was not a theoretical one like Galileo imagined, but more likely the pedestrian circumstance that we are not good at dropping one object from each hand at the same time. Modern experiments show that people tend to drop the lighter object sooner, ``even though the subject `feels' that she has released both simultaneously.''\footnote{\cite[14]{SettleExpr}.} This is why Galileo and others ended up thinking that lead balls started out slower than lighter bodies and only then picked up speed.

With experiments being so inconclusive, it is no wonder that Galileo relied more on a theoretical argument in his published account. His supporters would have us believe that “Galileo showed that Aristotle’s rule could be refuted by logic alone.”\footnote{\cite[35]{drakeessays1}.} Let us listen to his argument, which is supposedly a ``splendid and incontrovertible'' model example of ``cast-iron {\it reductio ad absurdum} reasoning.''\footnote{\cite[41, 46]{thusspoke}.}
\quote{By a short and conclusive demonstration, we can prove clearly that it is not true that a heavier moveable is moved more swiftly than another, less heavy …\ If we had two moveables whose natural speeds were unequal, it is evident that were we to connect the slower to the faster, the latter would be partially retarded by the slower …\ But the two stones joined together make a larger stone …;  therefore this greater stone is moved less swiftly than the lesser one. But this is contrary to your assumption. So you see how, from the supposition that the heavier body is moved more swiftly than the less heavy, I conclude that the heavier moves less swiftly.\footnote{\cite[66--67]{galileo2newsci}, OGG.VIII.107--108.}}
``From this we conclude that both great and small bodies, of the same specific gravity, are moved with like speeds.''\footnote{\cite[68]{galileo2newsci}, OGG.VIII.109.} Furthermore, ``if one were to remove entirely the resistance of the medium, all materials would descend with equal speed.''\footnote{\cite[75]{galileo2newsci}, OGG.VIII.116.}

With this argument Galileo allegedly exposes a fundamental logical inconsistency in the Aristotelian theory of fall. But he doesn't. Aristotle is perfectly clear: heavier objects fall faster. So when you put the heavy and the light together they will fall faster. The inconsistency arises only when one inserts the additional assumption that when you put two bodies together the lighter will retard the heavier and slow it down. But there is no basis for this latter assumption in Aristotle. It is a fiction that Galileo has made up. Only by dishonestly misrepresenting the view he is trying to refute in this way is he able to draw his triumphant conclusion.

A more honest form of the argument, which doesn't depend on misrepresenting Aristotle, is the following, which Galileo knew but didn't include in his published treatises:
\quote{Two identical bricks would fall side by side; no doubt about that. If a piece of string was tied to them they still would. Shortening the string could not change that. Hence two bricks tied together end to end would fall at the same speed as either brick alone. Now throw away the string and glue the bricks together; no reason appears why this double brick of double weight should fall faster than two bricks tied together---or either one alone.\footnote{\cite[35]{drakeessays1}, paraphrasing notes Galileo wrote at Pisa. This argument was known before Galileo. For example, a very similar argument appears in Benedetti, {\it Diversarum speculationum mathematicarum et physicarum liber} (1585), Chapter 10 \cite[206]{DrakeDrabkin}.}}
So the real crux of the argument is the claim that two bricks held side by side should behave the same way in terms of fall whether they are glued together or not. This is not a bad argument, but it falls short of being a matter of ``logic alone.'' For instance, imagine you are taking a basketball free throw. You can choose between trying to hit the hoop with either two bricks glued together or two bricks merely held side by side as you throw them. Would you really say that ``no reason appears'' why nature should treat the two cases differently, so you might as well go with the loose bricks? I don’t think so. Then why should this assumption be accepted in the case Galileo describes? Of course these two cases are quite different, but my point is that even though the assumption about falling bricks may feel ``obvious,'' it is really dependent on experience. It makes little sense to claim that we know purely by {\it a priori} intuition how the bricks behave in the case of simple fall, yet that we do not have such {\it a priori} knowledge of other scenarios like that bricks thrown through a hoop. If one case is ``logic alone'' then why isn't the other? So, if we are being honest, we are back to having to rely on experiment after all. It is not ``logic'' that guarantees that bricks behave the same way whether they are tied together or not; it is experience.

Galileo indeed discusses the experimental side of the matter too in his treatise. He admits that actual experiments do not come out in accordance with his law because of air resistance. But, he says, the fit is much better than for Aristotle’s law. To make this point Galileo provides specific measurements of how much the slower ball lags behind the heavier one:
\quote{Aristotle says, ``a hundred-pound iron ball falling from the height of a hundred braccia hits the ground before one of just one pound has descended a single braccio.'' I say that they arrive at the same time. You find, on making the experiment, that the larger anticipates the smaller by two inches; that is, when the larger one strikes the ground, the other is two inches behind it. And now you want to hide, behind those two inches, the ninety-nine braccia of Aristotle, and speaking only of my tiny error, remain silent about his enormous one.\footnote{\cite[68]{galileo2newsci}, OGG.VIII.109. Of course the statement attributed to Aristotle is not a literal quote but rather an inference from his law of fall.}}
But this is fake data. Galileo cooks the numbers to sound much more convincing in favour of his theory. The actual lag or difference between the two bodies is more than 20 times greater than the fake data Galileo reports in his published so-called masterpiece. ``In no case could Galileo have consistently achieved the results he reported.''\footnote{\cite[132]{hillprojection}.}

Nevertheless it remains true that Galileo's law doesn't fare as poorly as that of Aristotle in this experiment. So Galileo has indeed managed to improve on a two-thousand-year-old claim, made passingly by a non-mathematician. Aristotle himself obviously did not intend his so-called law as quantitative science. He only introduces it very parenthetically as a stepping-stone toward making the philosophical point that there can be no such thing as an object of infinite weight. Aristotle has no interest in this ``law'' beyond drawing this isolated metaphysical conclusion from it.

Galileo's entire case rests on his readers considering Aristotle to be a great authority. No wonder he clings to this framing and uses it as the trope of his dialogues. If we admit the truth---that Aristotle's law had been refuted more than a thousand years before, and that the notion of relying on Aristotle for quantitative science would never have entered the mind of any mathematically competent person in Galileo's time---then what does Galileo have to show for himself? An unproven claim that doesn't even fit the fake data Galileo has specifically concocted for the purpose, let alone the many experiments that proved the opposite, including his own before he knew which way he was supposed to fudge the data. Galileo likes to portray himself as doing the world a great service by defeating the rampant Aristotelianism all around him. The truth is that he is rather doing himself a great service by pretending that these Aristotelian opinions are ever so ubiquitous, so that he can inflate the importance of his own contributions, the feebleness of which would be all too evident if he addressed actual scientists instead of straw men Aristotelians.

\subsection{Law of fall}
\label{lawoffall}

Here is ``Galileo's'' law of fall, familiar to everyone from high school physics: ignoring air resistance, the acceleration of a freely falling object is constant. Equivalently, the velocity it acquires during its fall is proportional to time, and the distance fallen is proportional to time squared. The equivalence of these three formulations of the law of fall is elementary and immediately suggests itself to any mathematically competent person. In modern terms we can express the three version of the law in calculus terms as $y''=g$, $y'=gt$, and $y=gt^{2}/2$, respectively, but the straightforward equivalence of these three relationships long predates the calculus.\footnote{Cf.\ \S\ref{continuitythesis}.} Galileo set out to prove as much, but since he botched his treatment of calculus-style deductions,\footnote{\S\ref{Ginfintesimals}.} ``it is not Galileo but Descartes who was the first to conceptualize the relations of time, space, and velocity in a manner consistent with the conceptual framework of classical mechanics.''\footnote{\cite[281]{LimitsPreclasMech}. Descartes to Huygens, 18 February 1643, AT.III.807--808.}

Motion on an inclined plane (without friction) is closely related to this, because a ball rolling down a plane will acquire the same speed as it would have in free fall through the same vertical distance.\footnote{As we would say today, in anachronistic terms, since all balls covering the same vertical distance trade in the same amount of potential energy, they get the same amount of kinetic energy out of it.} All of these laws for free fall and for the inclined plane were indeed clearly stated by Galileo.\footnote{In various places in both the {\Dialogue} and the {\Discourse}, conveniently compiled and quoted in \cite[206--207]{VergaraCaffarelli}. OGG.VII: 47, 50--52, 225, 248, 252. OGG.VIII: 198, 202--205, 208--209.} The equivalence with the inclined plane is important since it makes empirical verification of the law much easier. A ball rolling down a slope is basically a slow-motion version of falling.
\quote{In a way it is surprising that the …\ law for the spontaneous descent of heavy bodies had not been recognized long before the 17th century. Measurements sufficient to put the law within someone's grasp are quite simple. Equipment for making them had not been lacking---a gently sloping plane, a heavy ball, and the sense of rhythm with which everyone is born.\footnote{\cite[1]{drakeHistFall}.}}
Relying on the sense of rhythm means that no time-keeping device is necessary. One places small grooves or bells along the ramp at intervals that the ball ought to cover in equal times according to the law of fall to be tested. The ball will make a click every time it rolls across one of the markings, so hearing whether the law holds is as easy as telling when a musician is off beat.\footnote{\cite[88--90]{drakeGatwork}.} Nevertheless, it is ``rather unlikely'' that Galileo is truthful in his exorbitant claim “that he had obtained precise correspondence in hundreds of trials on inclined planes.”\footnote{\cite[26]{drakeessays1}.}

Descartes seems to have arrived at the correct law of fall independently of Galileo. In fact, Descartes once held the mistaken belief that velocity is proportional to distance rather than to time, and he initially retained this belief even after reading Galileo's work, which he disliked.\footnote{\cite[142]{buttspitt}, \cite[10]{LimitsPreclasMech}. \S\ref{Descartes}.} The question of whether velocity is proportional to distance or to time was indeed a tricky and elusive one.\footnote{\cite[24--69]{LimitsPreclasMech}, \cite[100--102]{drakeGatwork}, \cite[12--13]{buttspitt}.} Galileo too got it wrong for a long time,\footnote{For example in a letter to Sarpi, 1604, OGG.X.115.} and at one point still gets himself confused and in effect uses the wrong form even in his mature {\Dialogue}.\footnote{In the course of the argument regarding centrifugal whirling discussed in \S\ref{centrifugal}. \cite[321]{chalmersnicholas} and \cite[121, 124]{hillprojection} make the case that Galileo's error in effect amounts to assuming that speed in free fall is proportional to distance rather than to time.}

In any case, Beeckman independently discovered the correct law of fall,\footnote{\cite[10]{LimitsPreclasMech}, \cite[417--418]{ClagettMA}.} as did Harriot, both at about the same time as Galileo and completely independently of him. ``In his mathematical analysis of [uniformly accelerated] motion, Harriot discovered that this actually includes two distinct possibilities: either the velocity of fall increases in proportion to the time elapsed or it increases in proportion to the space traversed. Harriot performed free fall experiments in order to decide between the two possibilities. …\ Harriot concluded that velocity increases in proportion to time and thus arrived at the law of time proportionality, which is correct also according to classical mechanics.''\footnote{\cite[238]{SchemmelHarriot}.} Of course he also immediately derived the equivalent form of the law that  distance fallen is proportional to time squared, in keeping with our observation that this equivalence is trivial and obvious.

\subsection{Gravitational constant}
\label{gravitationalconstant}

All objects fall with the same acceleration, we have now concluded. But how fast is that exactly? What is this same universal acceleration that every object shares? The answer is well-known to any student of high school physics: the constant of acceleration is $g\approx 9.8$ m/sec$^2$. But Galileo messed this up. He gives values equivalent to less than half of the true value.\footnote{\cite[223]{galileo2sys1sted}.} According to his defenders, ``clearly, round figures were taken here in order to make the ensuing calculation simple.''\footnote{\cite[69]{drakeHistFall}.} In other words, Galileo ``used arbitrary data.''\footnote{\cite[70]{drakeHistFall}.} And that’s what the people trying to {\em defend} Galileo are saying. Isn’t the law of fall supposed to be one of Galileo’s greatest discoveries? Why did he use fake data? Why not use real data? It was readily doable. His contemporaries did it. Why not do a little work to get the details right when you are publishing your supposed key results in your mature treatise? Is that really too much to ask? Instead of reporting make-believe evidence with a straight face, as Galileo does.

Competent and serious readers were in disbelief at Galileo’s inaccurate data. They certainly did not think it was fine to “use arbitrary data” in order to get ``round'' numbers for simplicity. Nor did they think it was “clear” that this is what Galileo was doing, contrary to what Galileo’s defenders are forced to argue when they try to excuse his inexcusable behaviour. Mersenne put it clearly: ``I doubt whether Mr Galileo has performed the experiment on free fall on a plane, since …\ the intervals of time he gives often contradict experiment.''\footnote{Mersenne, {\it Harmonie universelle}, I.122, \cite[261]{lewismersenne}.} Being a serious and diligent scientist, Mersenne did the work to find the correct value, unlike Galileo. The same goes for Riccioli, another contemporary of Galileo, who likewise did his own detailed experiments and found that they did not at all agree with the data Galileo had made up.\footnote{\cite[204]{Graney}.}

Once again the historical record fails to conform to Galileo's narrative of himself as a novel master experimenter refuting dogmas everyone believes in. Instead, contemporaries like Mersenne were already doing free fall experiments of their own before hearing of Galileo's,\footnote{\S\ref{fallandweight}.} and immediately recognised Galileo's data as fraudulent. Galileo is not addressing actual scientists. To them, much of what he has to say is disappointingly shallow and lacks serious scientific follow-through. Galileo's claim to fame relies on us not knowing this, and buying into the fiction that everyone at the time was an Aristotelian who had never heard of such a thing as experiment and quantitative science.

\subsection{Planetary speeds}
\label{planetaryspeeds}

\begin{figure}[pt]
\tikzfig{
\tkzDefPoints{0/0/O};
\draw[gray] (0,0) circle (5.2cm);
\draw[gray] (0,0) circle (4.4cm);
\draw[gray] (0,0) circle (3.6cm);
\draw[gray] (0,0) circle (2.6cm);
\draw[gray] (0,0) circle (1.8cm);
\draw[gray] (0,0) circle (1.1cm);

\tikzset{compass style/.append style={<-}};
\tkzDrawArc[R,color=black,thick](O,5.2)(75,84);
\tkzDrawArc[R,color=black,thick](O,4.4)(68,86);
\tkzDrawArc[R,color=black,thick](O,3.6)(56,88);
\tkzDrawArc[R,color=black,thick](O,2.6)(40,90);
\tkzDrawArc[R,color=black,thick](O,1.8)(5,92);
\tkzDrawArc[R,color=black,thick](O,1.1)(-110,94);

\node[inner sep=0pt] (sun) at (0,0)
    {\includegraphics[width=0.5cm]{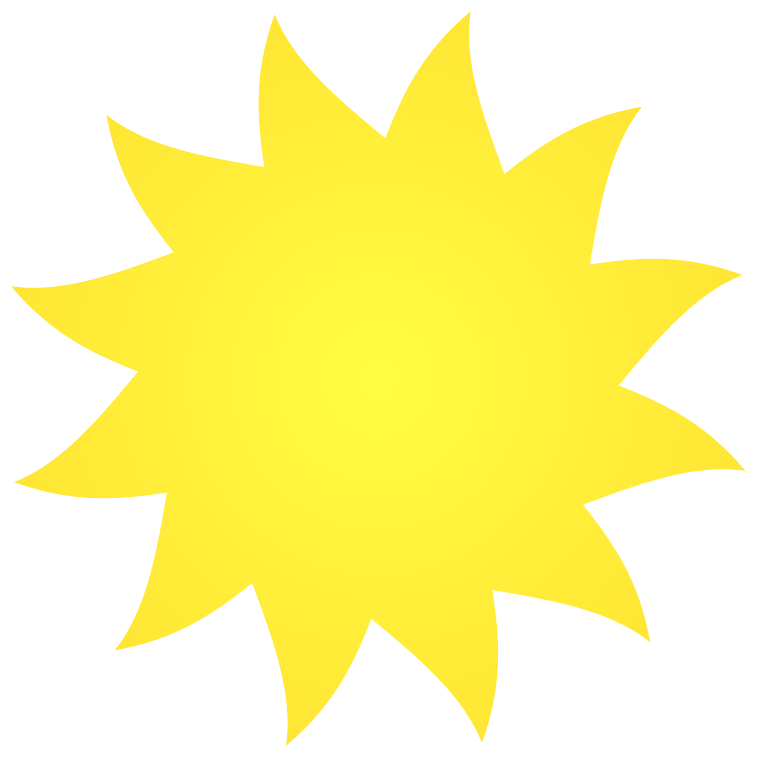}};

\node[inner sep=0pt] (mercury) at ({7.5*cos(94)},{7.5*sin(94)})
    {\includegraphics[width=0.15cm]{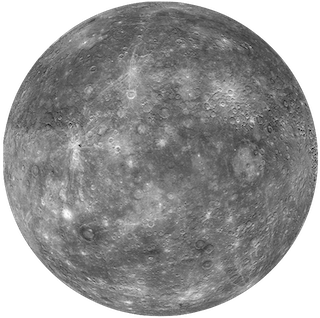}};
\node[inner sep=0pt] (venus) at ({7.5*cos(92)},{7.5*sin(92)})
    {\includegraphics[width=0.2cm]{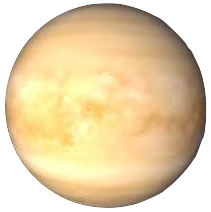}};
\node[inner sep=0pt] (earth) at ({7.5*cos(90)},{7.5*sin(90)})
    {\includegraphics[width=0.25cm]{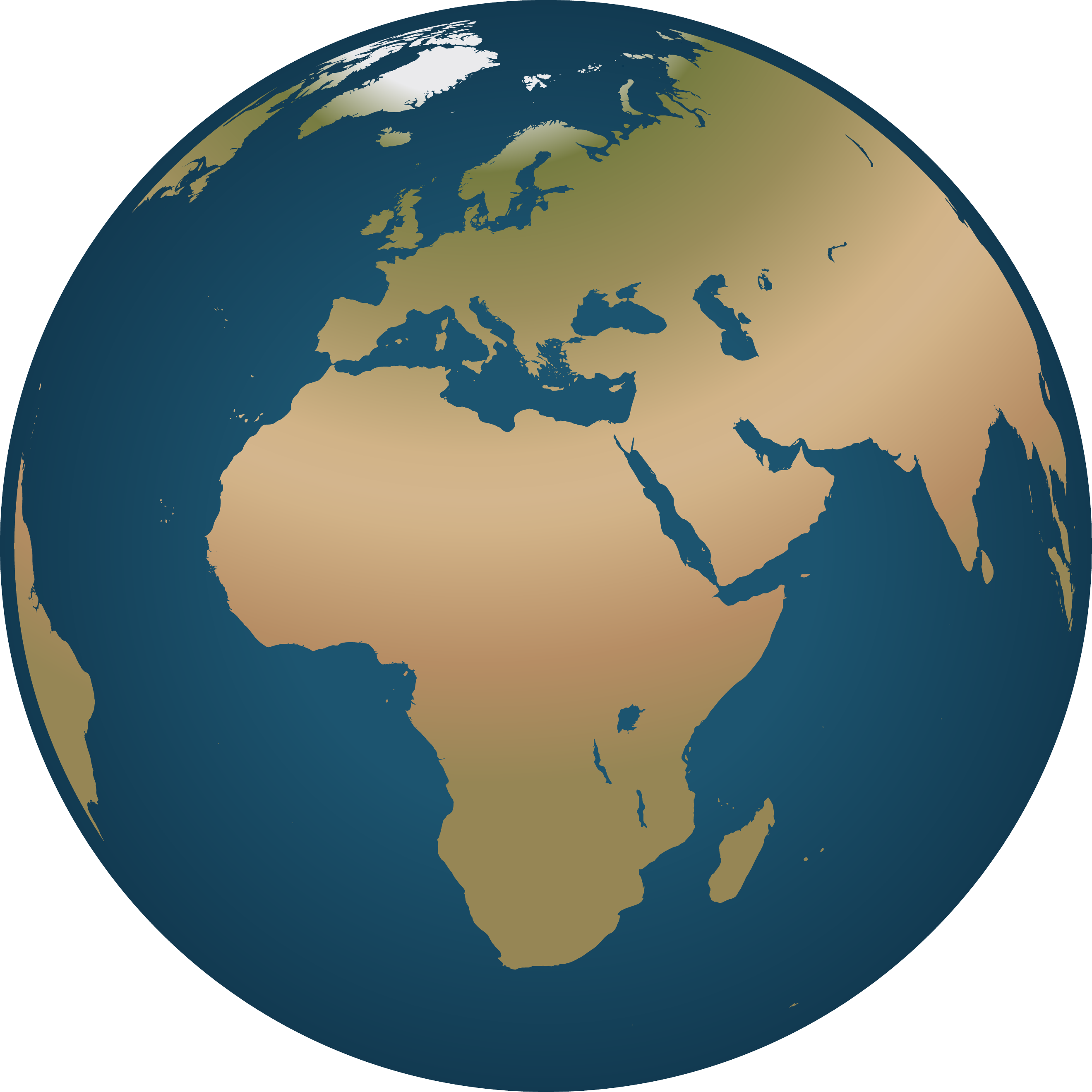}};
\node[inner sep=0pt] (mars) at ({7.5*cos(88)},{7.5*sin(88)})
    {\includegraphics[width=0.25cm]{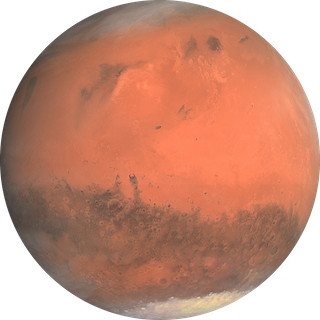}};
\node[inner sep=0pt] (jupiter) at ({7.5*cos(86)},{7.5*sin(86)})
    {\includegraphics[width=0.35cm]{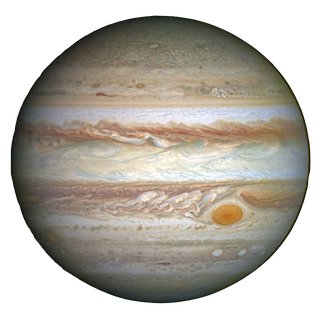}};
\node[inner sep=0pt] (saturn) at ({7.5*cos(84)},{7.5*sin(84)})
    {\includegraphics[width=0.4cm]{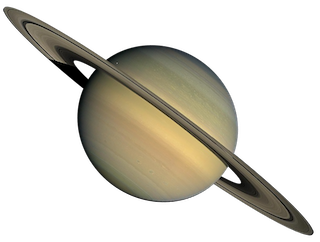}};

\draw ({7*cos(94)},{7*sin(94)}) -- ({1.1*cos(94)},{1.1*sin(94)});
\draw ({7*cos(92)},{7*sin(92)}) -- ({1.8*cos(92)},{1.8*sin(92)});
\draw ({7*cos(90)},{7*sin(90)}) -- ({2.6*cos(90)},{2.6*sin(90)});
\draw ({7*cos(88)},{7*sin(88)}) -- ({3.6*cos(88)},{3.6*sin(88)});
\draw ({7*cos(86)},{7*sin(86)}) -- ({4.4*cos(86)},{4.4*sin(86)});
\draw ({7*cos(84)},{7*sin(84)}) -- ({5.2*cos(84)},{5.2*sin(84)});

}{scale=0.6}
\caption[Planetary speeds.]{Galileo's erroneous theory that the orbital speeds of the planets are equal to the speeds they would have acquired through free fall if dropped from a common height.}
\label{pisandropfig}
\end{figure}
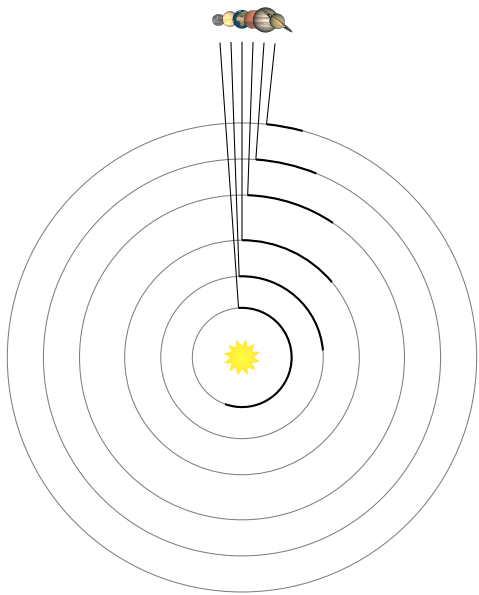

Before praising Galileo for stating the correct law of fall, it is sobering to consider the foolish uses he put it to. One is what has been called his ``Pisan Drop'' theory of planetary speeds.\footnote{\cite[116]{heilbron}.} The planets orbit the sun at different speeds. Mercury has a small orbit and zips around it quickly. Saturn goes the long way around in a big orbit and it is also moving very slowly. Galileo imagines that these speeds were obtained by the planets falling from some faraway point toward the sun, and then being somehow deflected into their circular orbits at some stage during this fall (Figure \ref{pisandropfig}). That supposedly explains why the planets have the speeds they do.

Galileo describes it as follows in the {\Dialogue}. “Suppose all the [planets] to have been created in the same place …\ descending toward the [sun] until they had acquired those degrees of velocity which originally seemed good to the Divine mind. These velocities being acquired, …\ suppose that the globes were set in rotation [around the sun], each retaining in its orbit its predetermined velocity. Now, at what altitude and distance from the sun would have been the place where the said globes were first created, and could they all have been created in the same place? To make this investigation, we must take from the most skilful astronomers the sizes of the orbits in which the planets revolve, and likewise the times of their revolutions.” Using this data and “the natural ratio of acceleration of natural motion” (that is, the constant $g$), one can compute “at what altitude and distance form the center of their revolutions must have been the place from which they departed.” According to Galileo this shows that indeed all the planets were dropped from a single point and their orbital data “agree so closely with those given by the computations that the matter is truly wonderful.”\footnote{\cite[29]{galileo2sys1sted}.}

Galileo was so proud of this erroneous argument that he repeated it in his second major work, the {\Discourse}, as well:
\quote{Our Author [Galileo] …\ may at some time have had the curiosity to try whether he could assign a determinate sublimity [i.e., distance from the sun] from which the bodies of the planets left from a state of rest, and were moved through certain distances in straight and naturally accelerated motion, and were then turned at the acquired speeds into equable motions. This might be found to correspond with the sizes of their orbits and the times of their revolutions. …\ Indeed, I seem to remember that he told me he had once made the computation, and also that he found it to answer very closely to the observations.\footnote{\cite[233]{galileo2newsci}, OGG.VIII.284.}}
Galileo omits the details though. He has one of the characters in his dialogue say that “making these calculations …\ would be a long and painful task, and perhaps one too difficult for me to understand,” whereupon Galileo's mouthpiece in the dialogue confirms that “the procedure is indeed long and difficult.”\footnote{\cite[30]{galileo2sys1sted}.}

Actually there is nothing ``difficult'' about it, at least not to mathematically competent people. Mersenne immediately ran the calculations and found that Galileo must have messed his up, because his scheme doesn't work.\footnote{Marin Mersenne, {\it Harmonie Universelle}, II.103--107, \cite[233, note 22]{galileo2newsci}. Later Newton made the same observation, \cite[144]{newtonprincipiacohened}.} There is no such point from which the planets can fall and obtain their respective speeds. Galileo's precious idea is so much nonsense, which evidently must have been based on an elementary mathematical error in calculation.

\subsection{Moon fall}
\label{moonfall}

Galileo tried to compute how long it would take for the moon to fall to the earth, if it was robbed of its orbital speed. “Making the computation exactly,” according to himself, he finds the answer: 3 hours, 22 minutes, and 4 seconds.\footnote{\cite[224]{galileo2sys1sted}.} This is way off the mark because Galileo assumes that his law of fall---that is, constant gravitational acceleration---extends all the way to the moon. This assumption is erroneous; the force of gravity diminishes with distance according to Newton’s inverse-square law. Ironically, Galileo’s purpose with this calculation was to refute the claim of another scholar that the fall would take about six days, which is a much better value: in fact it would take the moon almost five days to fall to the earth.\footnote{\cite[480]{galileo2sys1sted}.} That’s Galileo, the great hero of quantitative science, in action for you: bombastically claiming to refute others with his “exact calculations,” only to make fundamental mistakes and err orders of magnitude worse than his opponents did.

\subsection{Resistance of the medium}
\label{resistance}

As we have seen, Aristotelians were often as inclined to experiment as Galileo---a point obscured by Galileo's pretences to the contrary when it suited his purposes. Elsewhere it suited Galileo better to feign other straw men. ``In one of the dialogues of Galileo, it is Simplicius, the spokesman of the Aristotelians---the butt of the whole piece---who defends the experimental method of Aristotle against what is described as the mathematical method of Galileo.''\footnote{\cite[80]{Butterfield}.}

Consider for example the question of the resistance of the medium (such as air or water) on a moving object. Aristotle stated a law regarding how a body moves faster in a rarer medium than in a dense one. Galileo, in an early text, criticises Aristotle for accepting this ``for no other reason than experience''; instead one must ``employ reasoning at all times rather than examples,'' ``for we seek the causes of effects, and these are not revealed by experience.''\footnote{\cite[7]{sheaGrev}. OGG.I.260, 263.} Alas, Galileo's own law as to how resistance depends on density of the medium is also ``incompatible with classical mechanics'' as one study puts it---a polite, scholarly way of saying it's wrong.\footnote{\cite[270]{LimitsPreclasMech}.}

Employing some more ``reasoning'' along the same lines, Galileo decided that air resistance doesn't really increase appreciably with speed: ``The impediment received from the air by the same moveable when moved with great speed is not very much more than that with which the air opposes it in slow motion.''\footnote{\cite[226]{galileo2newsci}, OGG.VIII.277.} A surprising conclusion to modern bicyclists, among others. Yet ``experiment gives firm assurance of this,''\footnote{\cite[226]{galileo2newsci}, OGG.VIII.277.} Galileo promises. Alas, once again ``the statement is false, and the experiment adduced in its support is fictitious.''\footnote{Drake, \cite[226]{galileo2newsci}.} More fake data, in other words. This is quickly becoming a pattern.

\subsection{Inertia}
\label{inertia}

Newton's law of inertia says: ``Every body perseveres in its state of being at rest or of moving uniformly straight forward, except insofar as it is compelled to change its state by forces impressed.''\footnote{Newton, {\it Principia} (1687), Law 1. \cite[416]{newtonprincipiacohened}.} Is this important conception, so crucial to classical mechanics, due to Galileo? No. Even the most ardent defenders of Galileo ``freely grant that Galileo formulated only a restricted law of inertia''\footnote{\cite[608]{drakeinertia}.} and that ``he neglected to state explicitly the general inertial principle'' that everyone knows today, which was instead correctly ``formulated two years after his death by Pierre Gassendi and René Descartes.''\footnote{\cite[601]{drakeinertia}. Descartes, {\it Principia philosophiae} (1644), II.37, 39. On Gassendi, see \cite[150]{PalmerinoThijssen}.}

The charitable interpretation of trying to attribute to Galileo some kind of ``restricted law of inertia'' is a murky business. According to one author who tries to do so, ``in my opinion the essential core of the inertial concept lies in the ideas …\ of a body's indifference to motion or to rest and its continuance in the state it is once given. This idea is, to the best of my knowledge, original with Galileo.''\footnote{\cite[606]{drakeinertia}.} You could very well argue that that’s not really inertia at all because it doesn’t involve the straightness of the direction of the motion, nor does it explicitly say that the motion keeps going at a perpetual uniform speed. It only focusses on indifference of motion versus rest and preservation of the state of motion. But it is a dangerous game to cherry-pick some aspects of inertia. With this logic, one could equally well argue that that Aristotle had ``the essential core of the inertial concept,'' for he already highlighted indifference of motion versus rest and preservation of the state of motion two thousand years before Galileo:
\quote{No one could say why a thing once set in motion should stop anywhere; for why should it stop here rather than there? So that a thing will either be at rest or must be moved ad infinitum, unless something more powerful get in its way.\footnote{Aristotle, {\it Physics}, IV.8, 215a.}}
If this is inertia, then Aristotle was the pioneering near-Newtonian who conceived it, not Galileo.

This claim is rather isolated in Aristotle and didn’t really form part of a sustained and coherent physical treatment of motion comparable to how we use inertia today. Aristotle as usual is focussed on much more philosophical purposes. So you might say: that’s a one-off quote taken out of context which sounds much more modern than it really is.

Indeed. But then again the same could be said for Aristotle’s so-called law of fall that Galileo refuted with so much fanfare.\footnote{\S\ref{fallandweight}.} This too is only mentioned in passing very briefly and plays no systematic role in Aristotle’s thought. Yet Galileo takes great pride in defeating this incidental remark, and his modern fans applaud him greatly for it. So if we want to dismiss Aristotle’s inertia-like statement as insignificant, then, by the same logic, we ought to likewise dismiss all of Galileo’s exertions to refute his law of fall as completely inconsequential as well. If we argue that statements such as those of Aristotle don’t count as scientific principles unless they are systematically applied to explain various natural phenomena, then we would have to conclude that there was no Aristotelian science of mechanics at all. This, of course, would be a disastrous concession to make for advocates of Galileo’s greatness, since so much of Galileo’s claim to fame is based on contrasting his view with so-called “Aristotelian” science.

What about the {\em rectilinear} character of inertia? The thing keeps going {\em straight}. Was this key point grasped by Galileo? The following passage may appear to suggest as much:
\quote{A projectile, rapidly rotated by someone who throws it, upon being separated from him retains an impetus to continue its motion along the straight line touching the circle described by the motion of the projectile at the point of separation. …\ The projectile would continue to move along that line if it were not inclined downward by its own weight. …\ The impressed impetus, I say, is undoubtedly in a straight line.\footnote{\cite[191, 193]{galileo2sys1sted}.}}
In Newtonian terms this would be interpreted in terms of the law of inertia. But it hardly seems Galileo is thinking about it that way. He calls it “impetus.” What is “impetus” and why should we equate it with inertia? Will “impetus” run out? Is the motion caused by the “impetus” perpetual and uniform? In many other sources at the time, ``loss of impetus by projectiles was likened to …\ the diminution of sound in a bell after it is struck, or heat in a kettle after it is removed from the fire.''\footnote{\cite[244]{drakeGS}.} This conception is perfectly compatible, to say the least, with what Galileo writes. In fact, ``neither in the {\it Discourses} nor in the {\it Dialogue} does Galileo anywhere assert the eternal conservation of rectilinear motion.''\footnote{\cite[175]{KoyreGS}.} On the contrary, he explicitly rejects it: “Straight motion cannot be naturally perpetual.”\footnote{\cite[32]{galileo2sys1sted}.} “It is impossible that anything should have by nature the principle of moving in a straight line.”\footnote{\cite[19]{galileo2sys1sted}.} It is easy to understand, then, why Galileo's defenders are so eager to insist on characterising ``the essential core of the inertial concept'' in a way that does not involve its rectilinear character, since Galileo clearly and explicitly {\em rejected} rectilinear inertia.

If there is any inertia in Galileo, it is horizontal rather than rectilinear inertia. For example:
\quote{To some movements [bodies] are indifferent, as are heavy bodies to horizontal motion, to which they have neither inclination …\ or repugnance. And therefore, all external impediments being removed, a heavy body on a spherical surface concentric with the earth will be indifferent to rest or to movement toward any part of the horizon. And it will remain in that state in which it has once been placed; that is, if placed in a state of rest, it will conserve that; and if placed in movement toward the west, for example, it will maintain itself in that movement. Thus a ship, for instance, having once received some impetus through the tranquil sea, would move continually around our globe without ever stopping.\footnote{\cite[113--114]{GalileoDiscOp}.}}
\quote{Motion in a horizontal line which is tilted neither up nor down is circular motion about the center; …\ once acquired, it will continue perpetually with uniform velocity.\footnote{\cite[28]{galileo2sys1sted}. See also \cite[265--269]{coffainertia}.}}
Again, as with the sling and the projectile, one can debate whether this is inertia per se. In Newtonian mechanics too a hockey puck on a spherical ice earth would glide forever in a great circle, even though this is not inertial motion. But this agreement with Newtonian mechanics only holds if the object is prevented from moving downward, as the puck is by the ice, or the ship by the water. Galileo seems to have believed horizontal inertia to hold also for objects travelling freely through the air, which is not compatible with Newtonian mechanics. For example:
\quote{I think it very probable that a stone dropped from the top of the tower …\ will move, with a motion composed of the general circular movement and its own straight one.\footnote{\cite[165]{galileo2sys1sted}. See also \cite[270]{coffainertia}.}}
Once again it is not entirely clear that this is supposed to represent inertia at all. It is conceivable that, in Galileo's conception, the circular motion itself is not a force-free, default motion, but rather a motion caused or contaminated by some force or other. Who knows? Galileo just isn't clear about these kinds of things.

In his discussion of projectile motion too Galileo seems to rely on a horizontal conception of inertia. He speaks unambiguously of “the horizontal line …\ which the projectile would continue to follow with uniform motion if its weight did not bend it downward.”\footnote{\cite[199]{galileo2sys1sted}.} But he does {\em not} make the same claim for projectiles fired in non-horizontal directions. It rather seems as if he studiously avoided committing himself on that point.

An argument against attributing circular or horizontal inertia to Galileo, on the other hand, is his discussion attempting to prove that the the earth's rotational speed, no matter how much it increased, could not cause objects on its surface to be whirled off into space.\footnote{As will be discussed further in \S\ref{centrifugal}.} Galileo uses some intricate arguments to try to prove this, which he could have avoided by a simple appeal to circular inertia if he truly believed in the latter.

In conclusion, ``there is no general principle of inertia in Galileo's work.''\footnote{\cite[329]{chalmersnicholas}.} He definitely never stated the correct law of inertia. He often spoke of what appears to be a kind of circular inertia, in particular for horizontal motion (that it to say, circular motion around the center of the earth). But altogether it is impossible to attribute to him any one consistent view on the matter. Newton and Descartes, like the good mathematicians that they are, state concisely and explicitly what the exact fundamental assumption of their theory of mechanics are. Their laws of inertia are crystal clear and specifically announced to be basic principles upon which the rest of the theory is built. Galileo never comes close to anything of this sort. He uses the casual dialogue format of his books to hide behind ambiguities. One moment he seems to be saying one thing, then soon thereafter something else, like an opportunist who doesn't have a systematically worked out theory but rather adopts whatever assumptions are most conducive to his goals in any given situation.

\subsection{Relativity of motion}
\label{relativity}

When teaching basic astronomy at Padua, Galileo explained to his students that Copernicus was undoubtedly wrong about the earth's motion. The earth doesn’t move, Galileo explained. Because, if the earth moved, a rock dropped from a tower would strike the ground not at its foot but some distance away, since the earth would have moved during the fall. In support of this claim, ``Galileo observed that a rock let go from the top of a mast of a moving ship hits the deck in the stern.''\footnote{\cite[71]{heilbron}.} This had indeed been reported as an experimental fact by people who actually carried it out.\footnote{\cite[277--278]{wootton}.}

Of course this is completely backwards and the opposite of Galileo’s later views that he is famous for. To be sure, these lectures do not necessarily say anything about Galileo’s personal beliefs. In all likelihood he simply taught the party line because it was the easiest way to pay the bills. But at least the episode does show that the simplistic narrative that “the experimental method” forced the transition from ancient to modern physics is certainly wrong. On the contrary, experimental evidence was among the standard arguments for the conservative view well before Galileo got into the game.

In his later works Galileo of course affirms the opposite of what he said in those lectures: the rock will fall the same way relative to the ship regardless of whether the ship is standing still or travelling with a constant velocity. He gives a very vivid and elaborate description of this principle:
\quote{Shut yourself up with some friend in the main cabin below decks on some large ship, and have with you there some flies, butterflies, and other small flying animals. Have a large bowl of water with some fish in it; hang up a bottle that empties drop by drop into a wide vessel beneath it. With the ship standing still, observe carefully how the little animals fly with equal speed to all sides of the cabin. The fish swim indifferently in all directions; the drops fall into the vessel beneath; and, in throwing something to your friend, you need throw it no more strongly in one direction than another, the distances being equal; jumping with your feet together, you pass equal spaces in every direction. When you have observed all these things carefully (though doubtless when the ship is standing still everything must happen in this way), have the ship proceed with any speed you like, so long as the motion is uniform and not fluctuating this way and that. You will discover not the least change in all the effects named, nor could you tell from any of them whether the ship was moving or standing still. In jumping, you will pass on the floor the same spaces as before, nor will you make larger jumps toward the stern than toward the prow even though the ship is moving quite rapidly, despite the fact that during the time that you are in the air the floor under you will be going in a direction opposite to your jump. In throwing something to your companion, you will need no more force to get it to him whether he is in the direction of the bow or the stern, with yourself situated opposite. The droplets will fall as before into the vessel beneath without dropping toward the stern, although while the drops are in the air the ship runs many spans. The fish in their water will swim toward the front of their bowl with no more effort than toward the back, and will go with equal ease to bait placed anywhere around the edges of the bowl. Finally the butterflies and flies will continue their flights indifferently toward every side, nor will it ever happen that they are concentrated toward the stern, as if tired out from keeping up with the course of the ship, from which they will have been separated during long intervals by keeping themselves in the air. And if smoke is made by burning some incense, it will be seen going up in the form of a little cloud, remaining still and moving no more toward one side than the other.\footnote{\cite[186--187]{galileo2sys1sted}.}}
Galileo's prose is as embellished with fineries as this little curiosity-cabinet of a laboratory that he envisions. But is it any good of an argument? Insofar as it is, the credit is perhaps due to Copernicus himself, who had already made much the same point a hundred years before.\footnote{``When a ship floats on over a tranquil sea, all the things outside seem to the voyagers to be moving in a movement which is an image of their own, and they think on the contrary that they themselves and all the things with them are at rest. So it can easily happen in the case of the movement of the Earth that the whole world should be believed to be moving in a circle. Then what would we say about the clouds and the other things floating in the air or falling or rising up, except that not only the Earth and the watery element with which it is conjoined are moved in this way but also no small part of the air?'' \cite[17]{CopRevTransl}.}
So Galileo's relativity argument, like so much else he says, is old news. The primary contribution of his version is literary ornamentation---adding some butterflies and whatnot, while saying nothing new in substance.

Today the ``Galilean'' principle of relativity says that the phenomena in the cabin cannot be used to distinguish between the ship being at rest or moving with constant velocity in a straight line. But Galileo clearly has another scenario in mind: he sees the ship as travelling along a great circle around the globe. This is the kind of motion he believes cannot be distinguished from rest, in keeping with his misconceived idea of horizontal inertia.\footnote{\S\ref{inertia}.} This principle of relativity---the actually ``Galilean'' one---is of course false.\footnote{\cite[355]{DijksterhuisMech}.}

In fact it’s even worse than that. Galileo’s purpose with this passage about the ship is to argue, erroneously, that his ship cabin experiment ``alone shows the nullity of all those [arguments based on the motions of objects such as bird and clouds and falling objects relative to the earth] ad­duced against the motion of the earth.''\footnote{\cite[186]{galileo2sys1sted}.} That is to say, the rotation of the earth has no detectible mechanical effects, according to Galileo. But he is wrong. The earth's rotation is in fact mechanically detectible. The Foucault pendulum is a device that does behave differently due to the rotation of the earth than if the earth had been stationary, contrary to Galileo's claim that all local mechanical phenomena would be identical regardless of whether the earth was fixed or uniformly rotating, just as there is no way to detect the uniform motion of the ship from inside the cabin.

So the attribution of the principle of relativity of motion to Galileo in modern textbooks is doubly mistaken. First of all, relativity of motion and the idea of an inertial frame had been noted long before and was invoked by Copernicus to much the same end as Galileo. Moreover, Galileo’s principle is wrong in itself, because it’s about motion in a great circle, not in a straight line. And furthermore Galileo's purpose in introducing it is to draw another false conclusion from it, namely that the earth’s rotation is undetectable.

\subsection{Centrifugal force}
\label{centrifugal}

\begin{figure}[pt]\centering
\includegraphics[width=0.6\textwidth]{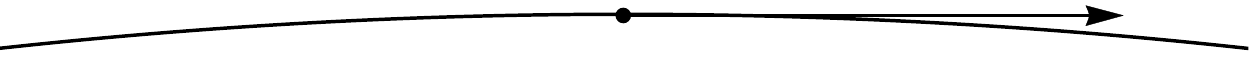}
\caption{Tangential motion of an object if gravity stopped acting on it. The object was located on the equator of the earth. The arrow shows how far it has gone in ten minutes.}
\label{Earthtangentfig}
\end{figure}

Galileo wished to refute this ancient argument: ``The earth does not move, because beasts and men and buildings'' would be thrown off.\footnote{\cite[220]{galileo2newsci2ndedWT}.} Picture an object placed at the equator of the earth, such as a rock lying on the African savanna. 
Imagine this little rock being “thrown off” by the earth’s rotation. In other words, the rock takes the speed it has due to the rotation of the earth, and shoots off with this speed in the direction tangential to its motion. The spectacle will be rather underwhelming at first: since the earth is so large, the tangent line is almost parallel to the ground, and since the speed of the rock and of the earth are the same they will keep moving in tandem. So rather than shooting off into the air like a canon ball, the rock will slowly begin to hover above the ground, a few centimeters at a time. For instance, Figure \ref{Earthtangentfig} shows how far the rock has gotten ten minutes after being thrown off the earth. As we can see, the deviation of the tangential path from the curved surface of the earth is only now becoming noticeable in relation to the size of the earth.

Of course this is not what happens to an actual rock, because gravity is pulling it back down again. The rock stays on the ground since gravity pulls it down faster than it rises due to the tangential motion. How can we compare these two forces quantitatively? Since we know the size and rotational speed of the earth, it is a simple task (suitable for a high school physics test) to calculate how much the rock has risen after, say, one second. This comes out as about 1.7 centimeters. We need to compare this with how far the rock would fall in one second due to gravity. Again, this is a standard high school exercise (equivalent to knowing constant of gravitational acceleration $g$). The answer is about 4.9 meters. This is why the rock never actually begins to levitate due to being “thrown off”: gravity easily overpowers this slow ascent many times over.

But of course this conclusion depended on the particular size and speed and mass of the earth. We could make the rock fly by spinning the earth fast enough. For example, if we run the above calculations again assuming that the earth rotates 100 times faster, we find that, instead of rising a measly 1.7 centimeters above the ground in one second, the rock now soars to 170 meters in the same time. The fall of 4.9 meters due to gravity doesn’t put much of dent in this, so indeed the rock flies away.

These things were calculated correctly in Galileo's time.\footnote{By Mersenne. \cite[113]{ThinkingWithObjects}.} But Galileo, alas, gets all of this horribly wrong. Even though we are supposed to celebrate Galileo as the discoverer of the law of fall, it is apparently too much to ask that he work out this very basic application of it. As noted,\footnote{\S\ref{gravitationalconstant}.} Galileo did not offer a serious estimate for the constant of gravitational acceleration $g$, unlike his contemporaries. Therefore he did not have the quantitative foundations to carry out the above analysis, which high school students today can do in five minutes.

Worse yet, Galileo maintains that no such analysis is needed in the first place, because he can “prove” that the rock will never be thrown off regardless of the rotational velocity. “There is no danger,” Galileo assures us, “however fast the whirling and however slow the downward motion, that the feather (or even something lighter) will begin to rise up. For the tendency downward always exceeds the speed of projection.”\footnote{\cite[197]{galileo2sys1sted}.} He even offers us “a geometrical demonstration to prove the impossibility of extrusion by terrestrial whirling.”\footnote{\cite[198]{galileo2sys1sted}.}

\begin{figure}[pt]\centering
\includegraphics[width=0.75\textwidth]{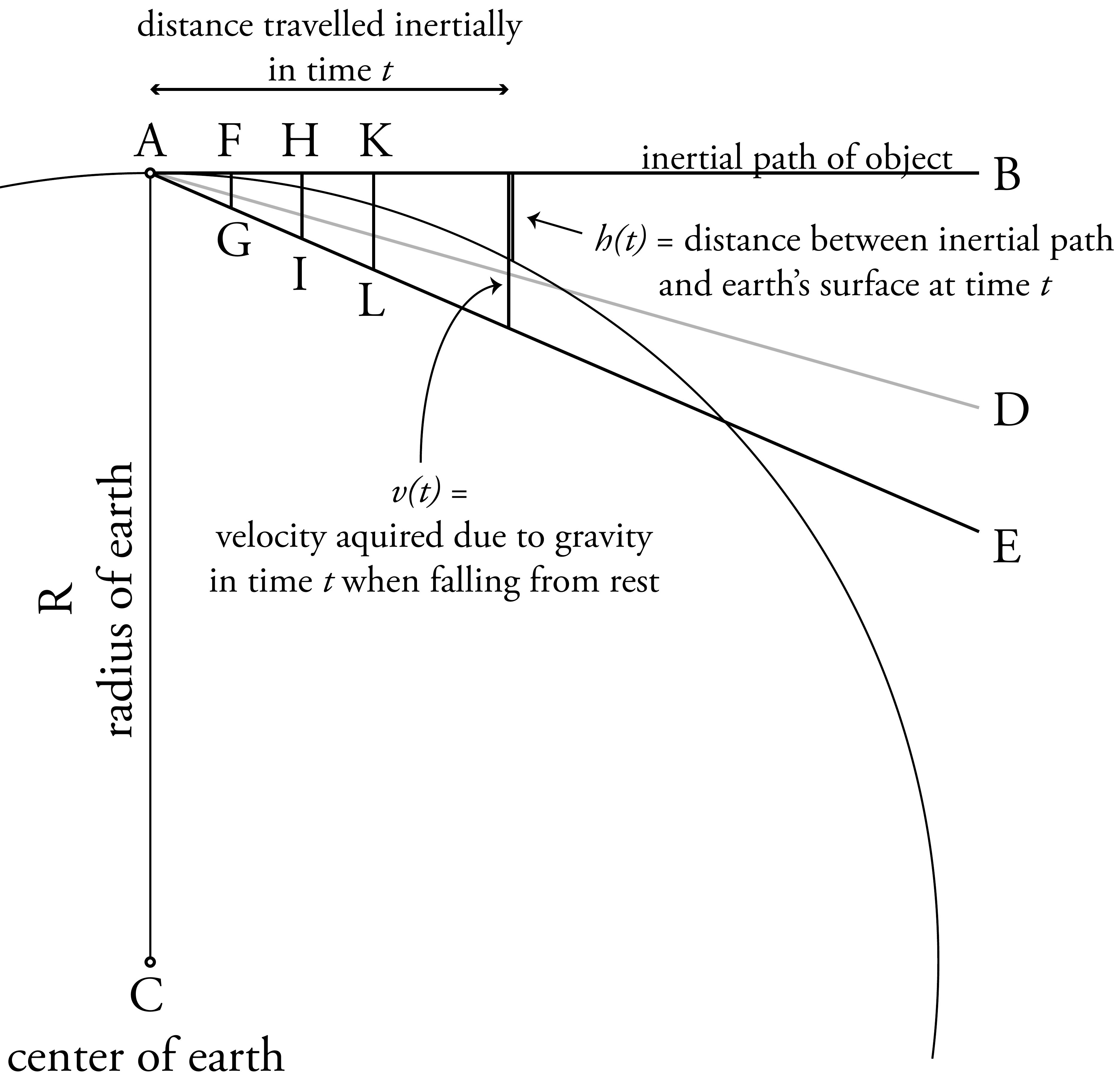}
\caption[Impossibility of centrifugal projection.]{Galileo's ``proof'' that centrifugal projection can never hurl objects off the earth. If gravity stops acting on an object at $A$, it would move inertially in the tangential direction $AB$. Since inertial motion has uniform speed, it would reach the equally spaced points $AFHK$ in equal time intervals. If the object had instead been dropped from rest, it would have acquired a certain downward speed in those same time intervals. These speeds are represented in the diagram by $FG$, $HI$, $KL$. Since the velocity acquired in free fall is proportional to time, $AGILE$ is a straight line. The slope of the line depends on the magnitude of gravitational acceleration, but for the purposes of this argument this value does not matter; in other words, we could just as well consider the speeds to be determined by some other line $AD$. The impossibility of centrifugal projection follows, according to Galileo, from the fact that as we consider smaller and smaller time intervals (that is to say, as we zoom it at $A$), the distance $h(t)$ required to catch up with the earth shrinks very rapidly to zero, while the speed $v(t)$ acquired from fall shrinks only linearly to zero. Therefore, says Galileo, the speed of fall $v(t)$ will, for some small enough $t$, be more than enough to cover the distance $h(t)$ and then some. In other words, the object will never get off the ground.}
\label{centrifugfig}
\end{figure}

Galileo's so-called ``demonstration'' is shown in Figure \ref{centrifugfig}.\footnote{\cite[231--234]{galileodialogueML}. The errors in Galileo's argument have been analysed by \cite{chalmersnicholas}, \cite{hillprojection}.} It is indeed a qualitative argument that rules out all possible cases of centrifugal projection, regardless of the rotational speed of the earth $V$, the radius of the earth $R$, or the magnitude of gravitational acceleration $g$. It is true, as Galileo says, that the ratio $v(t)/h(t)$ goes to $\infty$ as $t$ goes to zero. But this is obviously comparing apples to oranges, namely a velocity with a distance. The relevant comparison is between $h(t)$ and the distance $d(t)$ covered by free fall in this time. Galileo evidently felt that since in small time intervals $v(t)$ is overwhelmingly larger than $h(t)$, then $d(t)$ must surely be larger than $h(t)$ as well. But this is false. Instead,
\[ \lim_{t\rightarrow 0} \frac{d(t)}{h(t)} = \lim_{t\rightarrow 0} \frac{gt^{2}/2}{R-\sqrt{R^{2} -V^{2}t^{2}}} = \frac{gR}{V^{2}}.\]
In other words, $d(t)$ does not always overpower $h(t)$, as Galileo mistakenly believes. Rather, whether $d(t)$ is greater or smaller than $h(t)$ for small $t$ depends on the specific parameters of the situation in question. A strong gravitational acceleration $g$, or a big radius of the rotational path $R$, makes it easier for the object to ``catch up'' with the surface of the earth, while a big rotational speed $V$ makes it harder. Whether the object catches up with the surface or flies away depends on the relation between these parameters.

Galileo's claim to fame as a “mathematiser of nature” is certainly done no favours by this episode. He doesn’t know how to quantify his own law of fall, and doesn't understand basic implications of it. His physical intuition is categorically wrong on a qualitative level, and worse than that of the ancients he is trying to refute (whose stance was quite reasonable and would be accurate if the earth was spinning faster). He even offers a completely wrongheaded geometrical “proof” that the ancients' conception is impossible, even though so-called “Galilean” physics leads to the opposite conclusion in an elementary way.

\subsection{Circular path of fall}
\label{pathoffallabs}

\begin{figure}[pt]\centering
\includegraphics[width=0.7\textwidth]{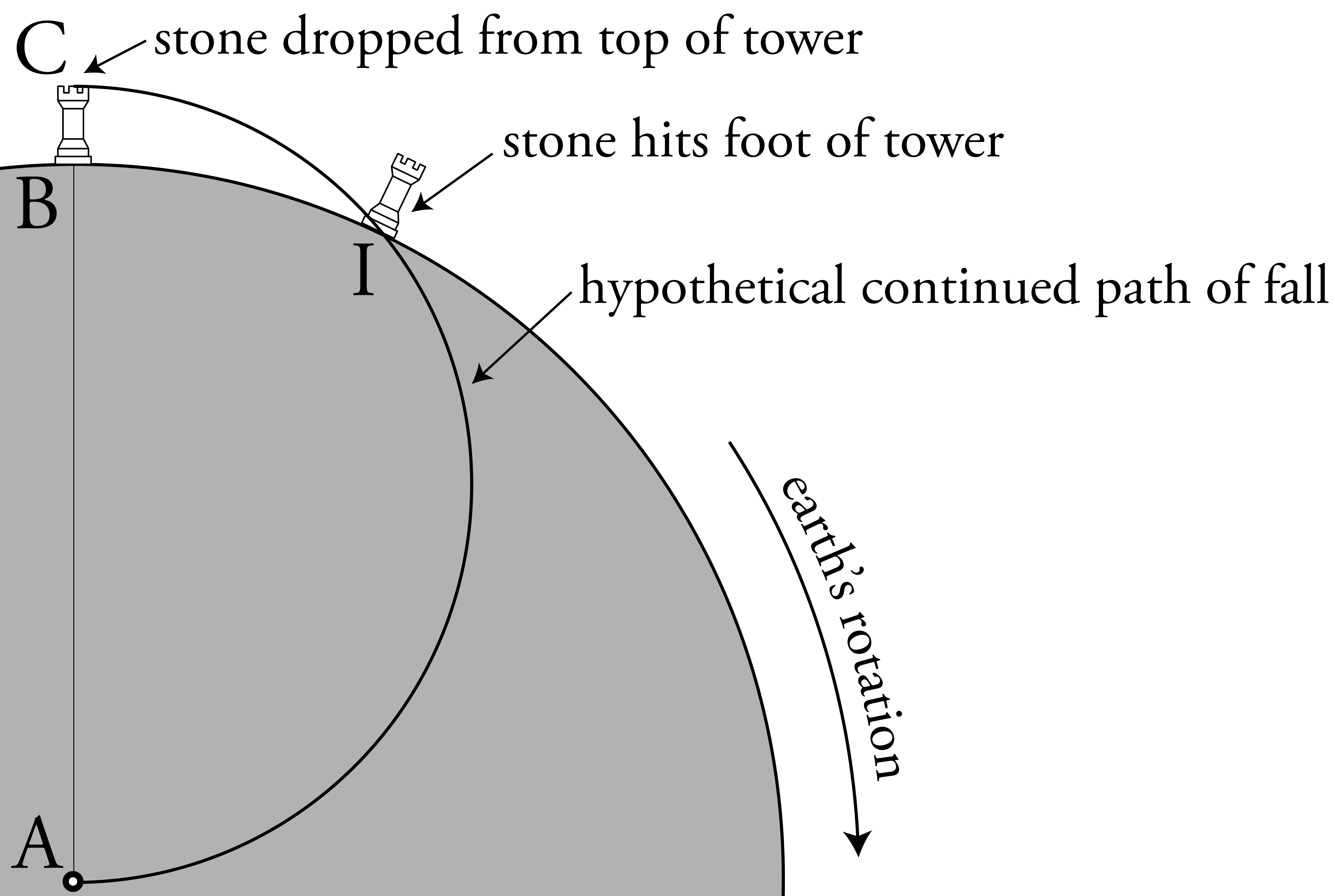}
\caption[Circular path of fall.]{Circular path of fall of a rock dropped from a tower, according to Galileo. ``$AB$ [is the radius of] the terrestrial globe. Next, prolonging $AB$ to $C$, the height of the tower $BC$ is drawn. …\ The semicircle $CIA$ …, along which I think it very probable that a stone dropped from the top of the tower $C$ will move, with a motion composed of the general circular one [due to the rotation of the earth] and its own straight one [due to gravity].'' \cite[192]{galileodialogueML}, OGG.VII.191.}
\label{circularfallfig}
\end{figure}

\begin{figure}[pt]\centering
\includegraphics[width=0.7\textwidth]{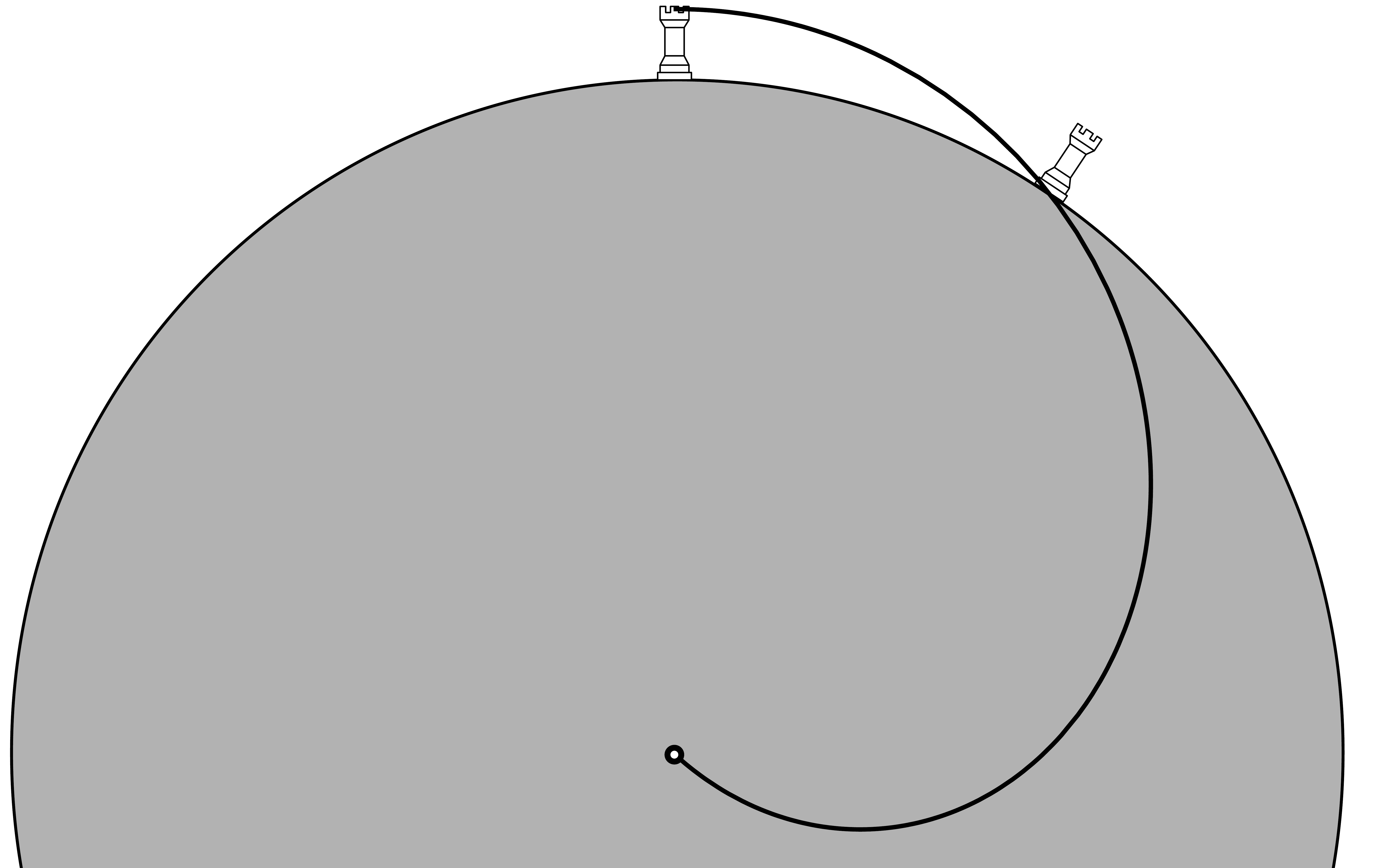}
\caption[Spiral path of fall.]{Corrected version of the path of fall of Figure \ref{spiralfallfig}. From Galileo's assumptions it follows that the path should be a spiral rather than a semicircle.}
\label{spiralfallfig}
\end{figure}

A rock dropped from the top of a tower falls in a straight line to the foot of the tower. But its path of fall is not actually straight if we take into account the earth's rotation. Seen from this point of view---that is to say, from a vantage point that doesn't move with the rotation of the earth---what kind of path does the rock trace? Galileo answers, erroneously, that it will be a semicircle going from the top of the tower to the center of the earth (Figure \ref{circularfallfig}):
\quote{If we consider the matter carefully, the body really moves in nothing other than a simple circular motion, just as when it rested on the tower it moved with a simple circular motion. …\ I understand the whole thing perfectly, and I cannot think that …\ the falling body follows any other line but one such as this. …\ I do not believe that there is any other way in which these things can happen. I sincerely wish that all proofs by philosophers had half the probability of this one.
\footnote{\cite[192--193]{galileodialogueML}, OGG.VII.191.}}
This is ``so obviously false (and besides incompatible with his own theory of uniformly accelerated motion of falling bodies) that one may wonder that Galileo did not see it himself.''\footnote{\cite[335]{KoyreDoc}.} Once again Galileo doesn't understand basic implications of his own law. Mersenne readily spotted Galileo's error, whereupon Fermat observed that the path should be a spiral (Figure \ref{spiralfallfig}), not a semicircle.\footnote{\cite[336, 342]{KoyreDoc}, \cite{HMmersennespiral}, \cite[556]{galileodialogueML}.} This would be the right answer given Galileo's assumptions, namely that the path is generated by composing uniform angular motion with uniformly accelerated radial motion toward the center of the earth.\footnote{As stated in \cite[192]{galileodialogueML} and again later when he admitted Fermat's correction \cite[343]{KoyreDoc}.} This implies that the path of fall is $r=r_{0}-a\theta^{2}$ in polar coordinates, which is indeed a spiral. This is still not the true path of fall, since Galileo's assumption that his law of fall remains unchanged in the interior of the earth is itself false. But I am not concerned here with criticising Galileo on such anachronistic grounds. Much worse is the fact that he got the wrong answer even if we grant his own assumptions.

When his embarrassing error was pointed out to him, Galileo replied that ``this was said as a jest, as is clearly manifest, since it is called a caprice and a curiosity.''\footnote{Galileo to Pierre Carcavy, 5 June 1637, OGG.XVII.89, \cite[135]{sheaGrev}.} Some defence this is! Far from offering exonerating testimony, Galileo actually openly pleads guilty to the main charge: namely, that his science is a joke. But in reality ``it is hard to believe that Galileo had really meant his solution of the trajectory of the falling body to be merely a joke.''\footnote{\cite[343]{KoyreDoc}.} If Galileo truly meant his argument to be taken merely in jest, then why did he say that he “considered the matter carefully” and “sincerely wished that all proofs by philosophers had half the probability of this one” and so on? Many of Galileo's errors come with these kinds of bombastic claims where Galileo is editorialising about how remarkably convincing his own arguments are. Isn’t that convenient? Throw out a bunch of half-baked guesses, and when they turn out right you can claim credit for stating it with such confidence while a more responsible scientist may have been exercising prudent caution. And when the guesses turn out wrong, you can apparently just write it off as a “joke” and pretend that that was what you intended all along, even though you published it with all those extremely assertive phrases right in the middle of your big definitive book on the subject. It’s easy to be “the father of science” if you can count on posterity to play along with this double standard.

\subsection{Projectile motion}
\label{projectilemotion}

Pick up a rock and throw it in front of you. The path of its motion makes a parabola. That’s Galileo’s great discovery, right? Well, not really. Galileo does claim this but he doesn’t prove it. Even Galileo’s own follower Torricelli acknowledged this. The result is “more desired than proven,”\footnote{Torricelli, 1644, \cite[275]{LimitsPreclasMech}.} as he says, very diplomatically. And the reason why Galileo doesn’t prove it is a revealing one. It is due to a basic physical misunderstanding.

The right way to understand the parabolic motion of projectiles like this is to analyse it in terms of two independent components: the inertial motion and the gravitational motion. If we disregard gravity, the rock would keep going along a straight line forever at exactly the same speed. That’s the law of inertia. But gravity pulls it down in accordance with the law of fall.\footnote{\S\ref{lawoffall}.} The rock therefore drops below the inertial line by the same distance it would have fallen below its starting point in that amount of time if you had simply let it fall straight down instead of throwing it. A staple fact of elementary physics is that the resulting path composed of these two motions has the shape of a parabola.

Galileo does not understand the law of inertia, and that is why he fails on this point. If the projectile is fired horizontally, such as for instance a ball rolling off a table, then Galileo does prove that it makes a parabola. He proves it the right way, the way just outlined, by composition of inertial and gravitational motion.\footnote{\cite[217, 221--222]{galileo2newsci2ndedWT}, OGG.VIII.269, 272--273.}

But if you throw the rock at some other angle, not horizontally, then Galileo doesn’t dare to give such an analysis. ``Although [Galileo's] {\it Discorsi} takes it for granted that the trajectory for oblique projection is a parabola, no derivation of this proposition is presented.''\footnote{\cite[237]{LimitsPreclasMech}.} ``At the point in the systematic treatment of projectile motion in the {\it Discorsi} where oblique projection is actually dealt with and correctly stated to yield a parabolic trajectory, there is simply a gap in the argumentation, and no derivation is offered for this claim.''\footnote{\cite[237]{LimitsPreclasMech}.}

\begin{figure}[pt]\centering
\includegraphics[width=0.8\textwidth]{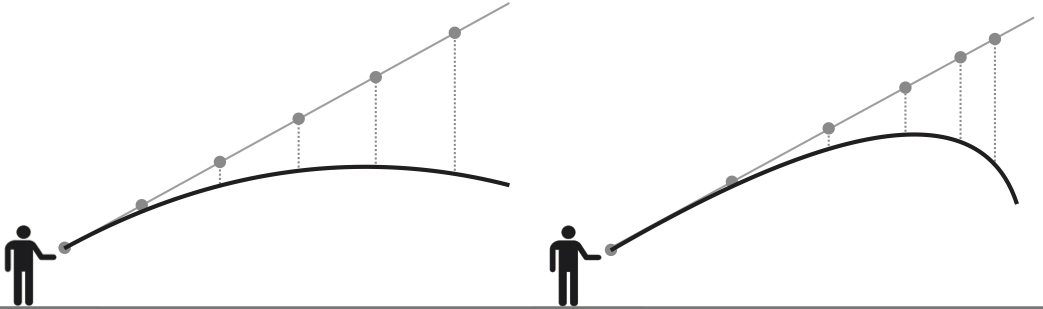}
\caption[Projectile motion.]{Left: Correct conception of projectile motion. The dots indicate uniform inertial motion in the firing direction. Right: Erroneous conception of projectile motion, as drawn by Galileo in unpublished manuscripts. The dots indicate a decelerating motion in the firing direction, as if the projectile was struggling to ascend the incline. In both cases the rectilinear motion is composed with an independent vertical motion according to the law of fall. Based on \cite[94, 96]{schemmelfox}.}
\label{projectilefig}
\end{figure}

Galileo's failure is quite clearly due to his not daring to believe in uniform inertial motion in any other direction than along the horizontal. He seems to fear that the law of inertia is perhaps not true for such motions. He equivocates and never takes a clear stand, because he is unsure whether it is the case or not. He is worried that the rectilinear component of the projectile's motion should be seen not as uniform but rather as gradually slowing down, like a ball struggling to roll up a hill or an inclined plane. In the latter case the trajectory is still a parabola, though not an ``upright'' one. See Figure \ref{projectilefig}.

In his final account, Galileo correctly ``postulated upright parabolas for all angles of projection. Galileo's reasoning for this shape is, however, untenable in classical mechanics. What is more, Galileo was unable to derive it from the consideration of two component motions.''\footnote{\cite[234]{SchemmelHarriot}.} 
``Galileo was …\ confronted with a contradiction between the inclined-plane conception of projectile motion and his claim that the trajectory is an upright parabola for all angles of projection, a contradiction he was never able to resolve.''\footnote{\cite[234]{SchemmelHarriot}.} Since he only trusted the horizontal case, Galileo tried to analyse other trajectories in terms of this case. To this end he assumed, without justification, that a parabola traced by an object rolling off a table would also be the parabola of an object fired back up again in the same direction.\footnote{\cite[245]{galileo2newsci2ndedWT}, OGG.VIII.296. \cite[234]{SchemmelHarriot}, \cite[227, 236]{LimitsPreclasMech}.} In other words, ``he takes the converse of his proposition without proving or explaining it,''\footnote{Descartes to Mersenne, 11 October 1638, AT.II.387. \cite[391]{drakeGatwork}.} as Descartes---a mathematically competent reader---immediately pointed out.

Instead, ``it was Galileo's disciples who first derived the parabolic trajectory for oblique projection, although they present it merely as an explication of Galileo's {\it Discorsi},'' which it is not.\footnote{\cite[7]{LimitsPreclasMech}.} Indeed, ``even before Galileo's {\it Discorsi} appeared in print, Bonaventura Cavalieri published a derivation of the parabolic trajectory that is consistent with classical mechanics and is not restricted to horizontal projection.''\footnote{\cite[284]{LimitsPreclasMech}.} Cavalieri was a good mathematician, unlike Galileo. He was also Galileo's countryman and in some sense ``disciple,'' and was very generous in deferring credit to Galileo.

The failures of Galileo's treatment of projectile motion confirms his misconception that inertia is limited to horizontal motion, which, as we have seen, was already independently suggested by other passages.\footnote{See \S\ref{inertia}.} Some have tried to argue that ``if Galileo never stated the law [of inertia] in its general form, it was implicit in his derivation of the parabolic trajectory of a projectile.''\footnote{\cite[602]{drakeinertia}.} This would have been a very good argument if Galileo had treated parabolic trajectories correctly. But he didn't, so the evidence goes the other way: Galileo's bungled treatment of parabolic motion is yet more proof that he did not understand inertia.

Even apart from the above errors and omissions, the mathematical details of Galileo's presentation of projectile motion are very clumsy. Galileo's ``calculations are unnecessarily complicated, and were greatly simplified by Torricelli in …\ 1644, a complete revision and enlargement …\ which …\ makes Galileo's demonstrations and procedures obsolete.''\footnote{\cite[53]{OxHBHistPhys}.} Once again Galileo's text bears the marks of an amateur mathematician, in other words. And once again his followers almost immediately cleaned up his mess in more mathematically able works that were full of deference to Galileo. ``While …\ inspired by veneration of Galileo, Torricelli is more logical in his treatise.''\footnote{\cite[91]{HallBall}.} Hence later mathematicians who used Torricelli's better but reverential account rather than Galileo's original for the mathematical details could easily be left with a much more flattering impression of the mathematical quality of ``Galileo's'' theory than if they had studied Galileo's own treatise in detail. Perhaps it is not so strange, then, that posterity got a bit confused and attributed much more to Galileo than he actually earned.

\subsection{Practical ballistics}
\label{ballistics}

Galileo made no theoretical use of his theory of projectile motion. In particular, he made no connection between this theory and the motion of the planets, the moon, or comets---a huge missed opportunity. Instead Galileo erroneously claimed that his theory would be practically useful for people who were firing cannons. ``He seems to have written [this theory] only to explain the force of cannon shots fired at different elevations,''\footnote{Descartes to Mersenne, 11 October 1638, \cite[391]{drakeGatwork}.} as Descartes put it. Descartes correctly denied that Galileo's idealised theory of projectile motion would correspond to practice. Galileo was less prudent. ``In many passages Galileo remarks that the theory of projectiles is of great importance to gunners. He made little or no distinction between his theory and useful ballistics; he believed---though without experiment---that he had discovered methods sufficiently accurate within the limitations of military weapons to be capable of direct application in the handling of artillery.''\footnote{\cite[91]{HallBall}.} This belief, however, was completely wrong. ``Galilean ballistics were not and could not be the fruit of experimental method; experiment …\ at once reveals that projectiles do not move in a parabola.''\footnote{\cite[96]{HallBall}.} A contemporary of Galileo put the matter to experimental test:
\quote{I was astonished that such a well-founded theory responded so poorly in practice. …\ If the authority of Galileo, to which I must be partial, did not support me, I should not fail to have some doubts about the motion of projectiles, [and] whether it is parabolical or not.\footnote{Giovan Battista Renieri to Torricelli, 2 August 1647, \cite[97--98]{HallBall}.}}
In reply, Torricelli ``stated outright that his book was written for philosophers, not gunners.'' ``According to Torricelli [the Galilean theory] had no connection with practical gunnery or with real projectiles.''\footnote{\cite[98]{HallBall}.} He also attributed this view to Galileo, but ``Torricelli's assertion that neither he nor Galileo ever believed that the science of motion had anything to do with practical affairs …\ is contradicted by Galileo's letters.''\footnote{\cite[100]{HallBall}.} ``Here Galileo himself must be charged with confusion.''\footnote{\cite[100]{HallBall}.} This is evident for example from the extensive tables that Galileo included in his big book: ballistic range tables based on his theory. These long tables make no sense other than as a practical guide for firing cannons.\footnote{\cite[90]{HallBall}.} So clearly Galileo thought his theory was practically viable, which it is absolutely not.

\subsection{Catenary}\label{catenary}

\begin{figure}[pt]\centering
\includegraphics[width=0.4\textwidth]{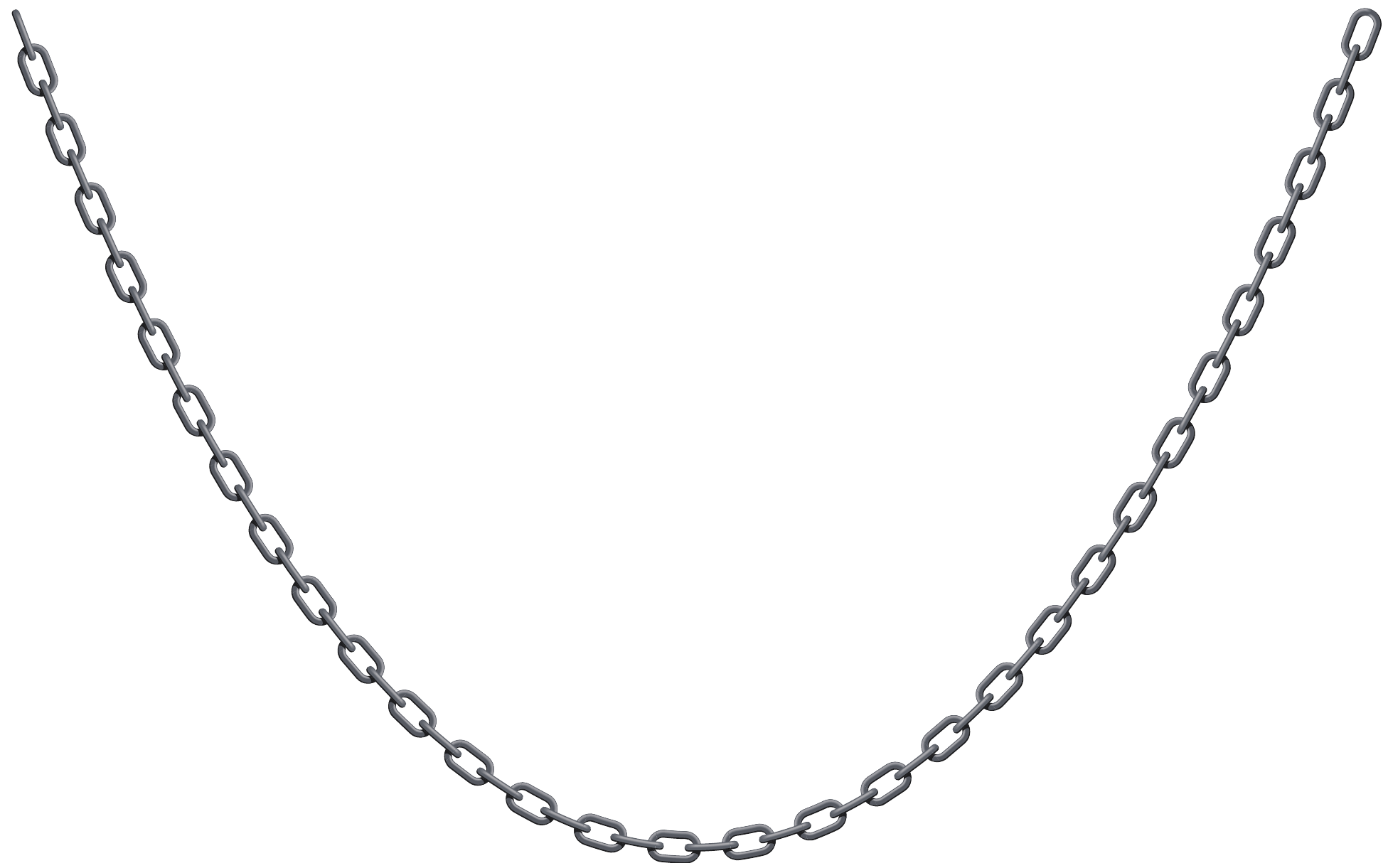} \hspace{.5cm} \includegraphics[width=0.4\textwidth]{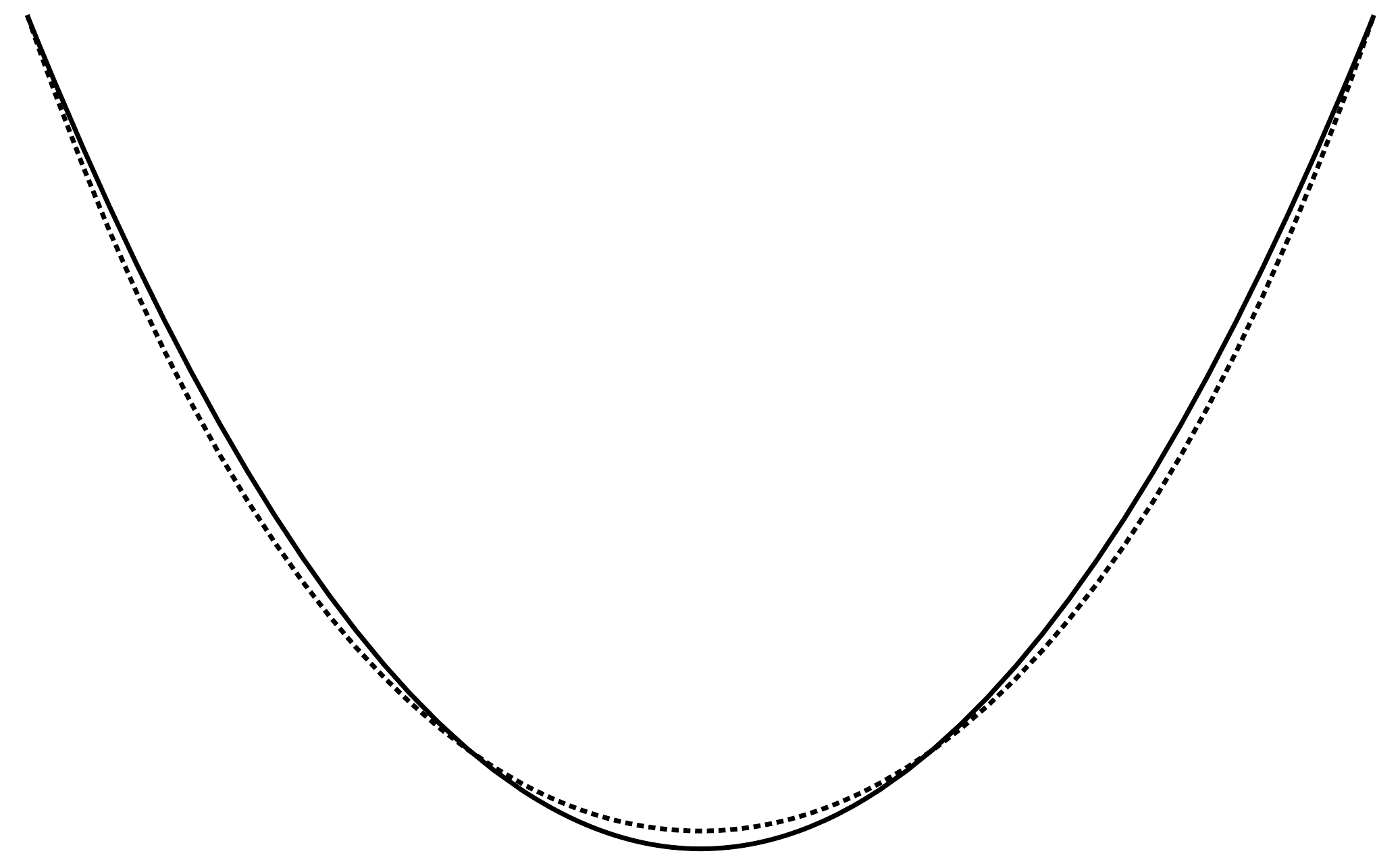}
\caption{Left: The catenary, or shape of a hanging chain, which Galileo erroneously believed to be a parabola. Right: The catenary (dotted) compared to a parabola (solid) of equal arc length between the same endpoints.}
\label{catenaryfig}
\end{figure}

The shape of a hanging chain (Figure \ref{catenaryfig}) looks deceptively like a parabola. It is not, but Galileo fell for the ruse: ``Fix two nails in a wall in a horizontal line …\ From these two nails hang a fine chain …\ This chain curves in a parabolic shape.''\footnote{\cite[143]{galileo2newsci}, OGG.VIII.186.} More competent mathematicians proved him wrong: Huygens demonstrated that the shape was not in fact parabolic.\footnote{\cite{BukowskiChain}.} Admittedly, his proof is from 1646, which is four years after Galileo's death. So one may consider Galileo saved by the bell on this occasion, since he was proved wrong not by his contemporaries but only by posterity. It is not fair to judge scientists by anachronistic standards. On the other hand Huygens was only seventeen years old when he proved Galileo wrong. So another way of looking at it is that a prominent claim in Galileo's supposed masterpiece of physics was debunked by a mere boy less than a decade after its publication.

In any case, Galileo thus ascribed to the catenary the same kind of shape as the trajectory of a projectile. He considered this to be no coincidence but rather due to a physical equivalence of the forces involved in either case.\footnote{\cite[256]{galileo2newsci2ndedWT}, OGG.VIII.309.} Indeed, Galileo made much of this supposed equivalence and “intended to introduce the chain as an instrument by which gunners could determine how to shoot in order to hit a given target.”\footnote{\cite[118]{whiteelephant}.}

Galileo also tried to test experimentally whether the catenary is indeed parabolic. To this end he drew a parabola on a sheet of paper and tried to fit a hanging chain to it. His note sheets are preserved and still show the holes where he nailed the endpoints of his chain.\footnote{\cite[39]{whiteelephant}.} The fit was not perfect, but Galileo did not reject his cherished hypothesis. Instead of questioning his theory, he evidently reasoned that the error was due merely to a secondary practical aspect, namely the links of the chain being too large in relation to the measurements. Therefore he tried it with a longer chain, and found the fit to be better. In this way he evidently convinced himself that he was right after all.\footnote{\cite[92--104]{whiteelephant}.}

The catenary case thus undermines two of Galileo's main claims to fame. First it brings his work on projectile motion into disrepute. The composition of vertical and horizontal motions that we are supposed to admire in that case looks less penetrating and perceptive when we consider that Galileo erroneously believed it to be equivalent to the vertical and horizontal force components acting on a catenary. Secondly, Galileo's reputation as an experimental scientists par excellence is not helped by the fact that his experiments in this case led him to the wrong conclusion, apparently because his love of his pet hypothesis led him to a biased interpretation of the data and a sweeping under the rug of an experimental falsification.

One may add a further complication to this account. While Galileo clearly and unequivocally stated that the catenary is a parabola in the passage quoted above, later in the same work he returns to the topic again and now speaks in less definitive terms, seemingly saying that the fit is approximate only:
\quote{The cord thus hung, whether much or little stretched, bends in a line that is very close to parabolic. …\ The similarity is so great that if you draw a parabolic line in a vertical plane surface …\ and then hang a little chain from the extremities of the base of the parabola thus drawn, you will see by slackening the little chain now more and now less, that it curves and adapts itself to the parabola and the agreement will be the closer, the less curved and the more extended the parabola drawn shall be.\footnote{\cite[256--257]{galileo2newsci2ndedWT}, OGG.VIII.310.}}
So the catenary is expressly {\em not} a parabola, only close to it, which is the correct view. But then again Galileo immediately goes on to add that ``then with a chain wrought very fine, one might speedily mark out many parabolic lines on a plane surface.''\footnote{\cite[257]{galileo2newsci2ndedWT}, OGG.VIII.310.} So now we are back to the erroneous view again. The insertion of the qualifier that the chain be ``very fine'' might be taken to suggest that while a typical ``cord'' will not be perfectly parabolic, an ideal chain with infinitesimally small links would be, which is not true.

Galileo may well have left the matter ambiguous on purpose. We can be absolutely sure that, if the catenary had in fact turned out to be equal to the parabola, then Galileo would certainly have been widely praised for having announced this fact. The passages about approximate fit could then be written off as having to do with practical limitations only, while the unequivocal statements that the catenary is parabolic could be taken as Galileo's core theoretical claim. On the other hand, now that we know that the catenary is not in fact a parabola, Galileo's defenders can twist his words to say that that is what he meant all along too: it is the parts about the deviation from a perfect fit that are the prescient and wise ones worthy of a pioneer of science, while the clear statements suggesting perfect fit are just simplifications for pedagogical purposes. Galileo left his text so ambiguous that he could always be construed as being ``right'' regardless of whether the facts of the matter turned out one way or the other.

\subsection{Pendulum}
\label{pendulum}

“With regard to the period of oscillation of a given pendulum, [Galileo] asserted that the size of the arc [i.e., the height of the starting position of the pendulum] did not matter, whereas in fact it does.”\footnote{\cite[26]{drakeessays1}. \cite[97]{galileo2newsci2ndedWT}, OGG.VIII.139. \cite[451]{galileo2sys1sted}.} Galileo's allegedly experimental report on pendulums in the {\Discourse} is clearly fabricated---or “exaggerated,” to use the diplomatic term.\footnote{\cite[139--140]{whiteelephant}.} Mersenne did the experiment and rejected Galileo's claim.\footnote{Mersenne, {\it Les nouvelles pensées de Galilée} (1639), 72--73, \cite[209]{OxHBHistPhys}.} Galileo's friend Guidobaldo del Monte did the same, but when he told Galileo of his error Galileo rejected the experiment and insisted in his claim.\footnote{\cite[86]{whiteelephant}, OGG.X.97--100.} Instead of admitting what experiments made by sympathetic and serious scientists showed, Galileo preferred to defend his false theory with ``conscious deception''\footnote{\cite[131]{hillprojection}.} (more commonly known as lying).

\subsection{Brachistochrone} \label{brachistochrone}

\begin{figure}[t]\centering
\includegraphics[width=0.75\textwidth]{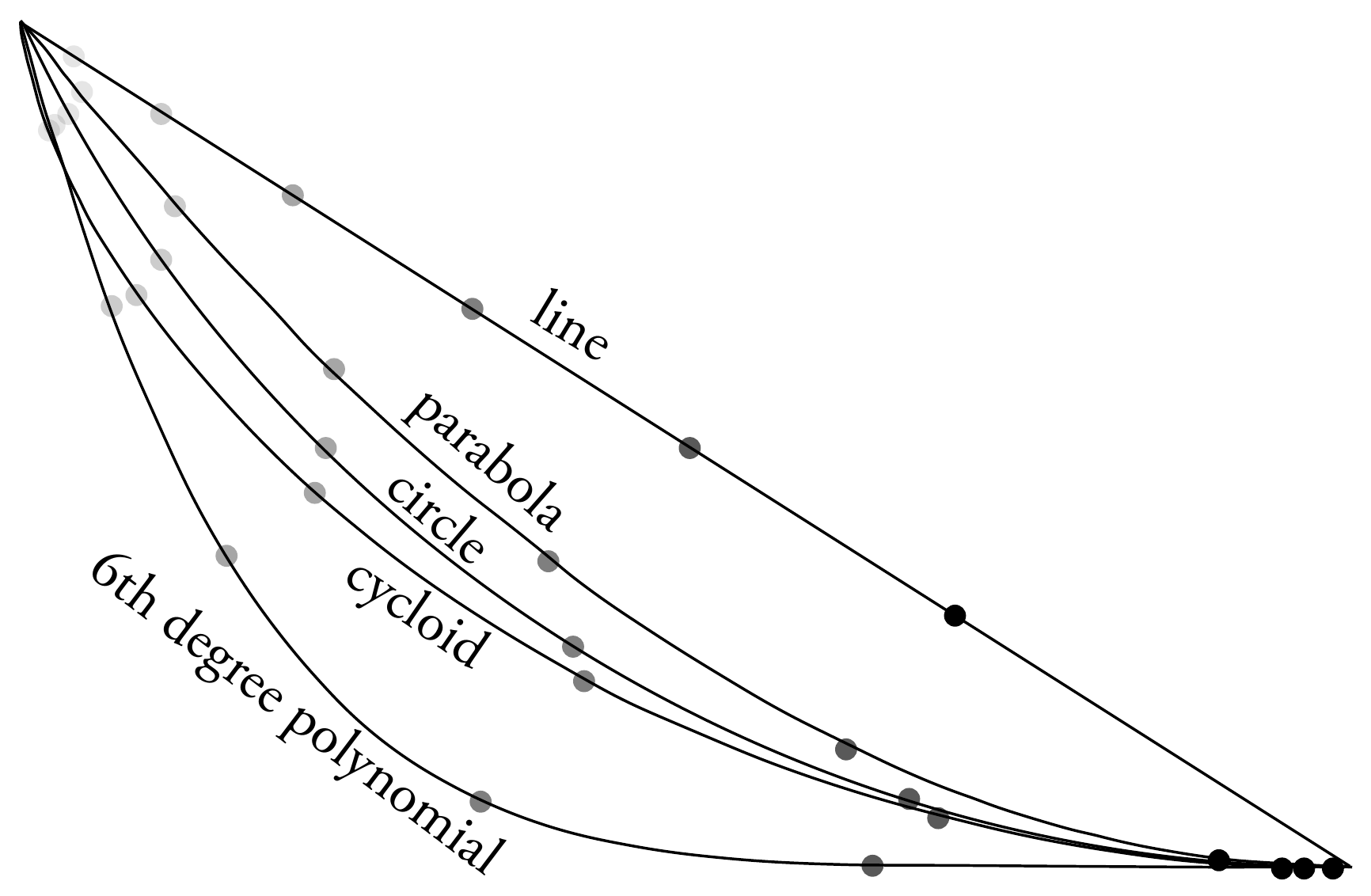}
\caption{The brachistochrone problem. The ball rolls down the fastest on a cycloidal ramp, not a circular ramp as Galileo thought. Points of the same shade correspond to the same moment in time. Based on \cite{Sanchis}.}
\label{brachistochronefig}
\end{figure}

The brachistochrone problem asks for the path along which a ball rolls down the quickest from one given point to another. Galileo believed himself to have proved that the optimal curve was a circular arc: ``From the things demonstrated, it appears that one can deduce that the swiftest movement of all from the terminus to the other is not through the shortest line of all, which is the straight line, but through the circular arc.''\footnote{\cite[213]{galileo2newsci}, OGG.VIII.263. \cite[451]{galileo2sys1sted}.} Actually the fastest curve is not a circle but a cycloid. This was only proved in the 1690s, using quite sophisticated calculus methods.\footnote{\cite[184]{BlasjoTrBook}.} We cannot blame Galileo for not possessing advanced mathematical tools developed only half a century after his death. Nevertheless it is one more addition to his sobering pile of erroneous assertions about various physical problems. We are supposed to celebrate him for being the first to discover the parabolic path of projectile motion, and conveniently forget that at the same time he was wrong on the brachistochrone, wrong on the catenary, wrong on the shape of the isochrone pendulum, etc. With all these errors stacking up, one starts to wonder whether the one thing he got right was much more than luck. Throw up enough half-baked guesses and something's bound to stick.

\subsection{Floating bodies}

\begin{figure}[tp]\centering
\includegraphics[width=\textwidth]{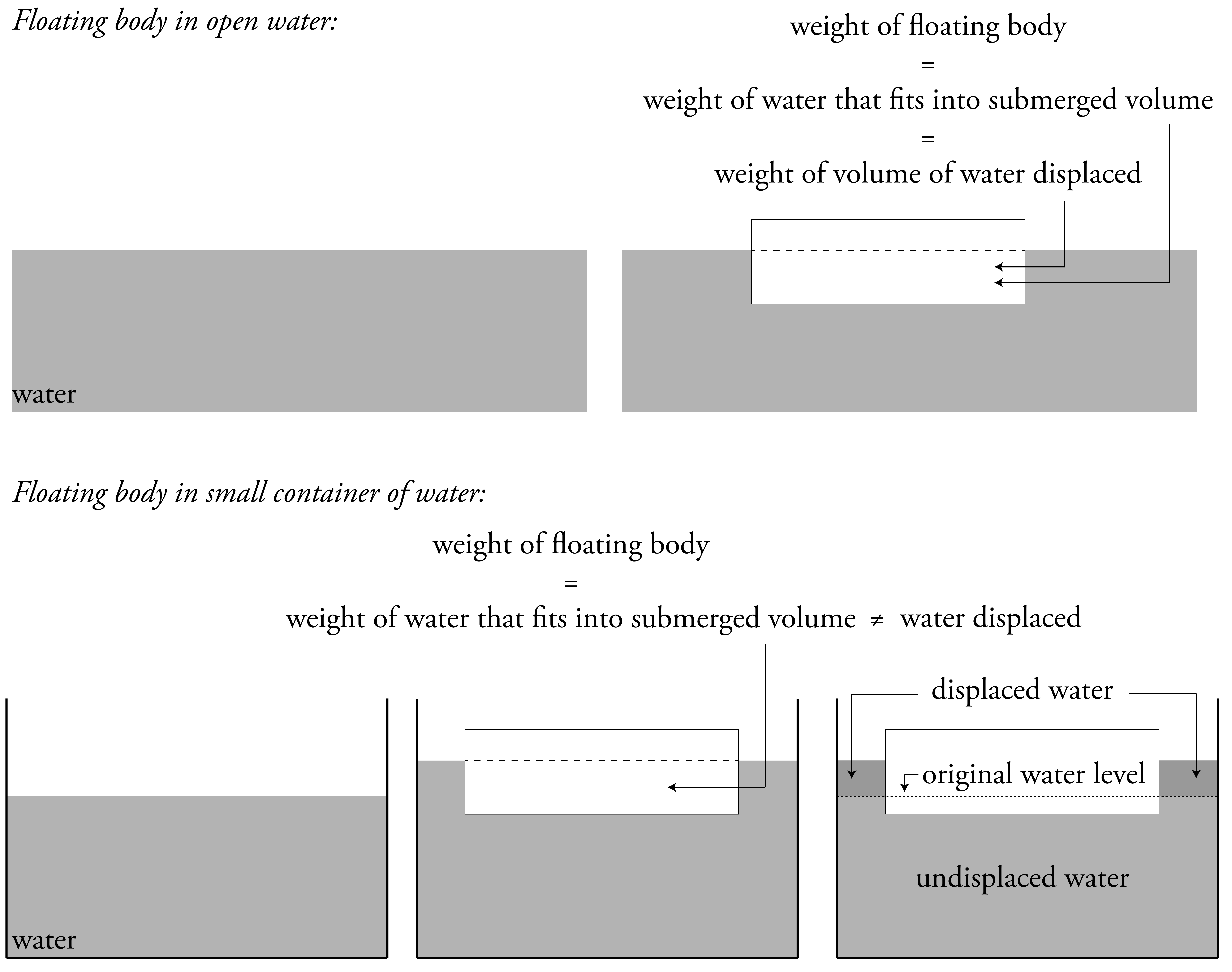}
\caption[Archimedes' principle.]{Archimedes' principle. In open water, such as an ocean, the rise in water level when a body is submerged is negligible. Therefore the submerged volume of the body is equal to the volume of water displaced. In a closed container, however, this equality does not hold, since the water level changes appreciably when a body is submerged.}
\label{Afloat1fig}
\end{figure}

Archimedes studied the science of floating bodies almost two thousand years before Galileo. Archimedes's work is excellent and spot on correct. Galileo had trouble understanding the most basic aspects of it, however. Archimedes's fundamental finding is that the weight of a floating body is equal to the weight of however much water would fit into the space occupied by the part of the floating body that is below the water level. On the open sea, an equivalent way of phrasing this is that the weight of a floating body equals the weight of the water it displaces. In a small container such as a bucket or a bowl, however, Archimedes's principle cannot be phrased in terms of displaced water in this way. Figure \ref{Afloat1fig} explains why.

Galileo, in his early work on floating bodies, got himself confused on precisely this point. He made the elementary blunder of using the formulation of the law in terms of displaced volume of water for containers. Thus his work was ``rooted in the false assumption that the volume of water that is displaced when a solid is immersed is equal to the total volume of the solid immersed whereas it is only equal to the volume of that part of the solid under the initial level of the water.'' It is baffling how ``the `Father of Experimental Science' could make such a mistake after some twenty years of University teaching.''\footnote{\cite[19]{sheaGrev}.}

Some scholars, however, have raved that ``the whole of the Scientific Revolution is encapsulated'' in how Galileo handled the matter, for in this very case we can ``watch our modern conception of common sense (as applied to science …) being born.''\footnote{\cite[§8.2]{woottoninvention}.} On this reading, Galileo didn't make a dumb mistake, as I claim. Instead he ``discovered'' an ``anomaly'':\footnote{\cite[§8.2]{woottoninvention}.}
\quote{A block of wood floating in a tank displaces less than its own weight in water. …\ Archimedes’ principle did not apply. …\ Indeed, the actual water in the tank could weigh less than the block of wood it was lifting---which, according to Archimedes’ principle as traditionally understood, was impossible. (You can do this yourself by putting a small amount of water inside a wine cooler and then floating a bottle of wine in it.)\footnote{\cite[§8.2]{woottoninvention}.}}
As a result of this ``discovery'':
\quote{Boldly, [Galileo] had gone back, reanalysed Archimedes’ principle, and revised it. He had then tested his new theory with a quite different experiment. This back and forth movement between theory and evidence, hypothesis and experiment has come to seem so familiar that it is hard for us to grasp that Galileo was doing something fundamentally new. Where his predecessors had been doing mathematics or philosophy, Galileo was doing what we call science.\footnote{\cite[§15.7]{woottoninvention}.}}
I beg to differ. Galileo may have ``discovered'' an ``anomaly'' in ``Archimedes’ principle as traditionally understood,'' but only because people were much too stupid to read Archimedes. The erroneous formulation of Archimedes's principle is obviously not due to Archimedes, and the foolish confusion about it would never occur to anyone who had actually read Archimedes's work with understanding. It is true that Galileo used experiment to uncover a mistake, but that mistake would not have occurred to a mathematically competent person in the first place. Galileo's use of experiment in this case was not a new innovation leading to progress in the theory of hydrostatics; rather it was a way for the mathematically inept to correct a misconception that arose from not having understood the mathematical theory.

It is self-evident that a theory of how floating bodies behave can be tested with a bucket of water and some blocks of wood. This was not a ``fundamentally new'' idea ``born'' with Galileo. It is ``common sense'' not only to us but to anyone who ever lived. It was just as obvious to Archimedes and the other ancient Greeks as it is to us. Archimedes did not neglect such experiments because he lacked the modern scientific method. Rather, he did not discuss such experiments because his sophisticated mathematical treatment went miles beyond the kinds of baby steps that occupied Galileo. It is a fact that Archimedes's theory of hydrostatics is excellent and correct, even when it goes into highly specific and nontrivial claims (see \S\ref{ArchHydrostatics}). It makes no sense to imagine that he somehow achieved this without scientific and experimental methods. It makes perfect sense, however, that he was more interested in the advanced mathematical theory he was developing than in explaining how childishly simple experiments can be used to avoid blockheaded misunderstandings such as that committed by Galileo.

\subsection{Tides}
\label{tides}

At the end of his famous {\it Dialogue}, Galileo lists what he considers to be his three best arguments for proving that the earth moves around the sun. One of these is his argument “from the ebbing and flowing of the ocean tides,”\footnote{\cite[536]{galileodialogueML}.} or high and low tidal water. Galileo believed the tides were caused by the motion of the earth. This is truly one of his very worst theories, even though he was ever so proud of it.

\begin{figure}[pt]\centering
\includegraphics[width=0.3\textwidth]{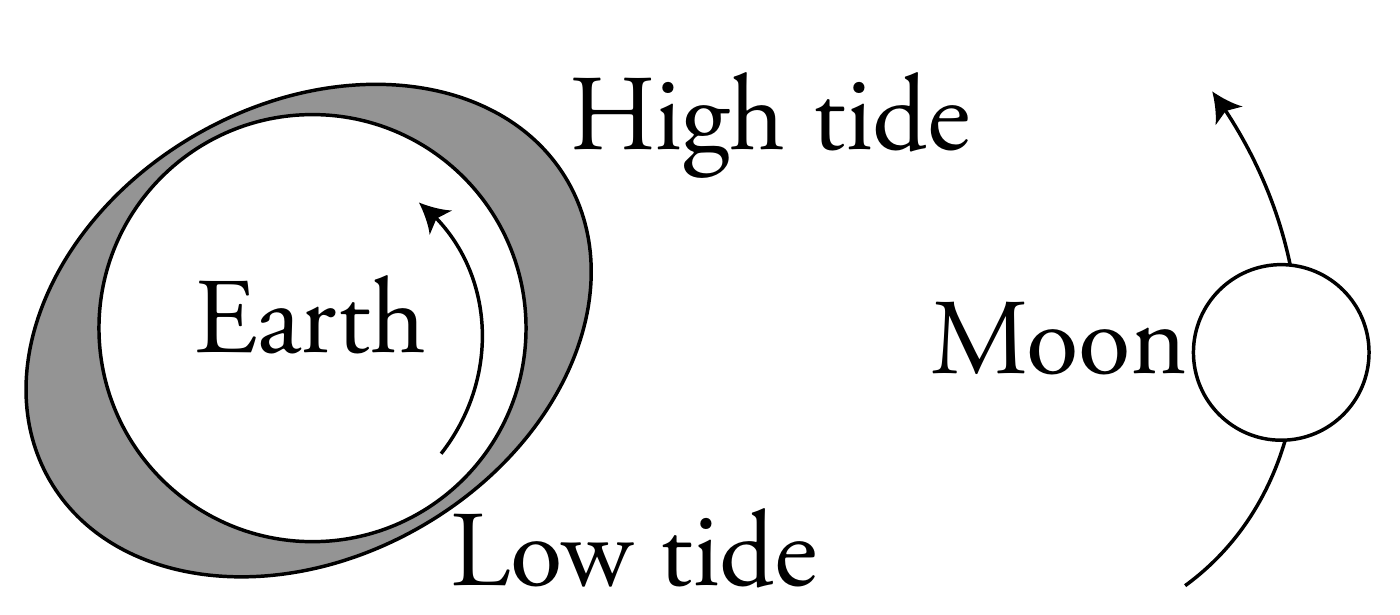}
\caption{The correct theory of tides.}
\label{correcttidesfig}
\end{figure}

First things first. How do the tides work? As we know today, ``the ebb and flow of the sea arise from the action of the sun and the moon,'' as Newton proved.\footnote{\cite[835]{newtonprincipiacohened}.} The moon, and to a lesser extent the sun, pull water towards them, causing our oceans to bulge now in one direction, then the other (Figure \ref{correcttidesfig}). This was clearly understood already in Galileo's time. Kepler explained it perfectly, and many others proposed lunar-attraction theories of the tides.\footnote{\cite[200]{PalmerinoThijssen}, \cite[172]{sheaGrev}.} In fact, the lunisolar theory of tides is found already in ancient sources, including the causal role of the sun and moon, and descriptions of the effects in extensive and accurate detail.\footnote{\cite[308]{RussoSciRev}.}

Galileo, however, got all of this completely wrong. Why should ``the tides of the seas follow the movements of the fireballs in the skies,''\footnote{Kepler to Herwart, 9--10 April 1599, \cite[51]{BaumgardtKepler}.} as Kepler had put it? Galileo considered the very notion ``childish'' and ``occult,'' and declared himself ``astonished'' that ``Kepler, enlightened and acute thinker as he was, …\ listened and assented to the notion of the Moon's influence on the water.''\footnote{\cite[356]{sheaOxSW}, \cite[462]{galileo2sys1sted}.}
\quote{There are many who refer the tides to the moon, saying that this has a particular dominion over the waters …\ [and] that the moon, wandering through the sky, attracts and draws up toward itself a heap of water which goes along following it.\footnote{\cite[419]{galileo2sys1sted}.}}
Yes, many indeed believed such things. And they were right. But Galileo would have none of it. This theory is not ``one which we can duplicate for ourselves by means of appropriate devices.'' How indeed could we ever ``make the water contained in a motionless vessel run to and fro, or rise and fall''? Certainly not by moving some heavy rock located thousands of miles away. ``But if, by simply setting the vessel in motion, I can represent for you without any artifice at all precisely those changes which are perceived in the waters of the sea, why should you reject this cause and take refuge in miracles?''\footnote{\cite[421]{galileo2sys1sted}.} That’s Galileo’s objection to the lunar theory of tides: It’s hocus-pocus. It assumes the existence of mysterious forces that we cannot otherwise observe or test. Proper science should be based on stuff we can do in a laboratory, like shaking a bowl of water.

\begin{figure}[pt]\centering
\includegraphics[width=0.23\textwidth]{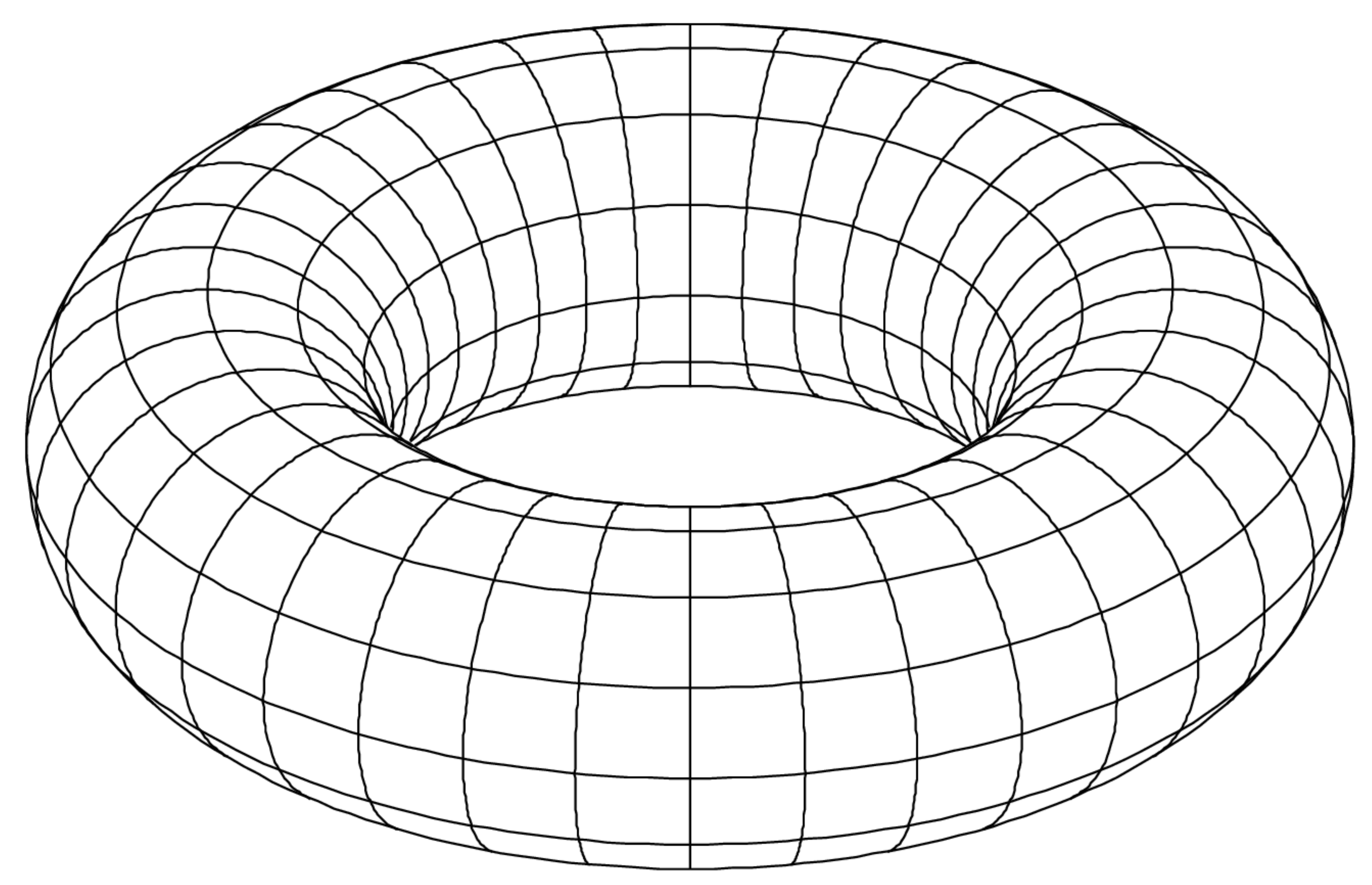} \hspace{0.5cm}\includegraphics[width=0.2\textwidth]{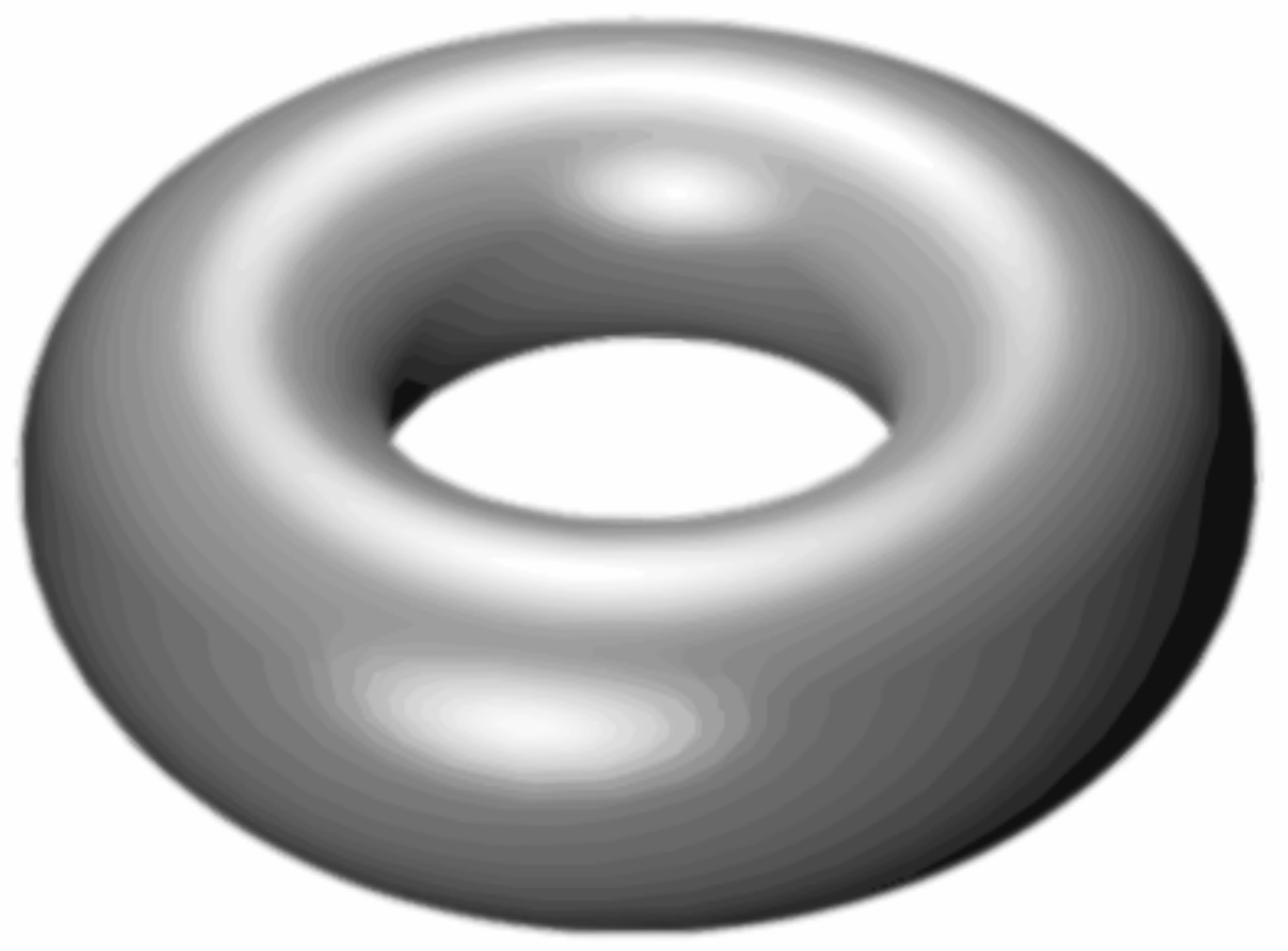}
\caption[Torus.]{Galileo's theory of the tides can be understood in terms of a torus.}
\label{torusfig}
\end{figure}

\begin{figure}[pt]\centering
\includegraphics[width=0.4\textwidth]{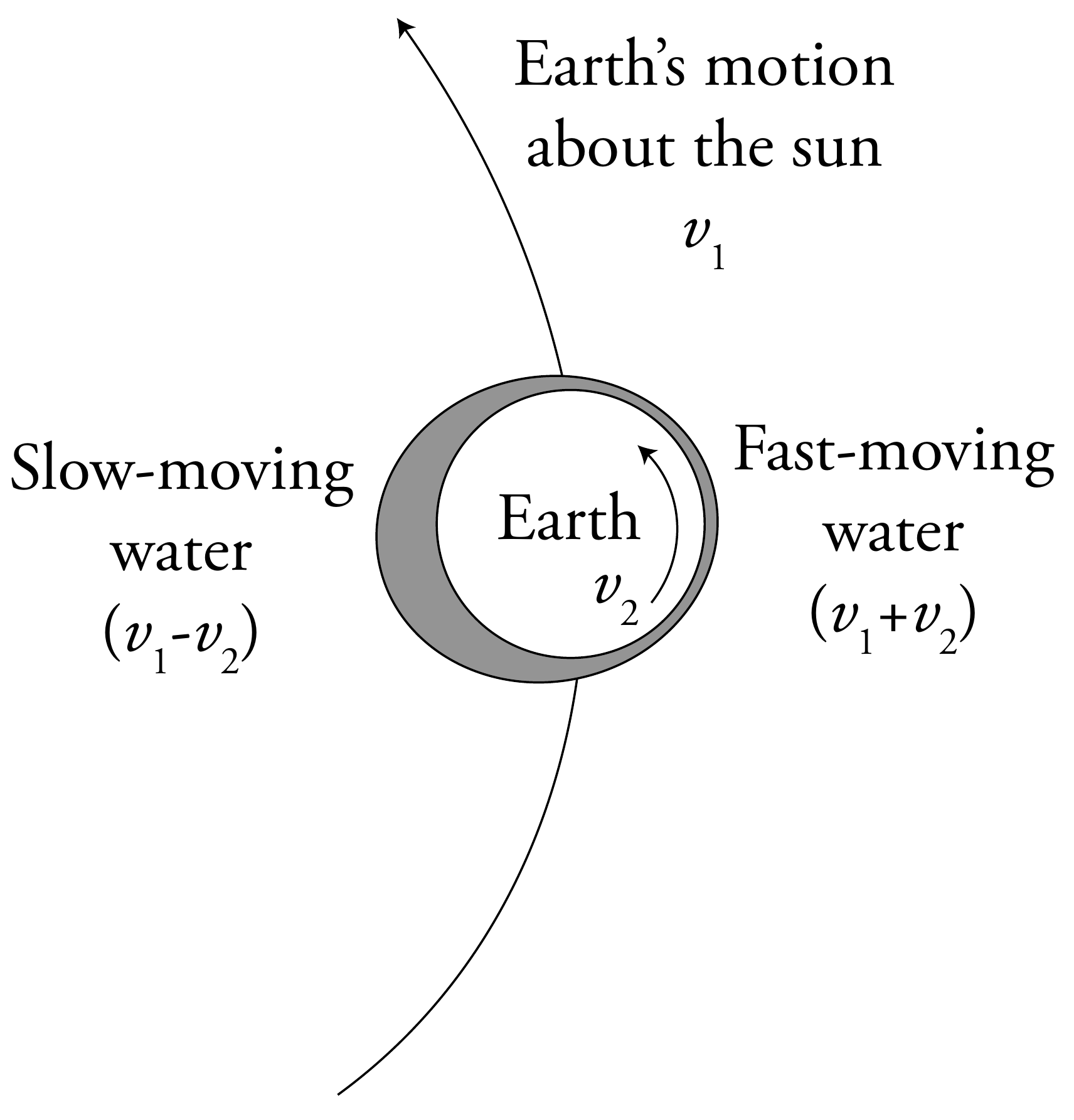}
\caption{Galileo's theory of tides.}
\label{Gtidesfig}
\end{figure}

Galileo's own theory of the tides can be described as follows (Figure \ref{Gtidesfig}). Picture a torus (Figure \ref{torusfig}), laying flat on the ground, filled halfway with water. This represents the water encircling the globe of the earth. Now spin the torus in place, around its midpoint, like a steering wheel. This represents the rotation of the earth. For symmetry reasons it is clear that the water will remain equally distributed around the full circle. Because there is no reason for the surface of the water to become higher or lower in one place of the tube than another. Every part is rotating equally, everything is symmetrical, and therefore no asymmetrical distribution of water could arise from this process.

But now picture the torus sitting on a merry-go-round. The represents the earth orbiting the sun. And also at the same time it is still spinning around its own midpoint like before. Now there is asymmetry in the configuration. The torus has one part facing inward toward the middle of the merry-go-round, and one part facing outward. As far as the rotation of the torus around its midpoint is concerned, these two parts of the torus are moving in opposite directions. But as far as the rotation of the merry-go-round is concerned they are both going the same way. So one part of the torus spins along with the orbital motion, and one part against the orbital motion. Therefore one part of the water moves faster than the other. It is boosted by the merry-go-round rotation helping it along in the direction it was already going, while the other part of the torus is slowed down; the merry-go-round cancelling its efforts by going in the opposite direction.

So you have fast-moving water and slow-moving water. The fast-moving water will catch up with the slow water and pile up on top of it, creating a high tide. And the space it vacated will not be replenished because the slow water behind it isn’t keeping up. Hence the low tide. In Galileo's words:
\quote{Mixture of the annual and diurnal motions causes the unevenness of motion in the parts of the terrestrial globe. …\ Upon these two motions being mixed together there results in the parts of the globe this uneven motion, now accelerated and now retarded by the additions and subtractions of the diurnal rotation upon the annual revolution.\footnote{\cite[427]{galileo2sys1sted}.}}
Thus, in fact, ``the flow and ebb of the seas endorse the mobility of the earth.''\footnote{\cite[416]{galileo2sys1sted}.}

Unfortunately, Galileo's theory is completely out of touch with even the most rudimentary observational facts about tidal waters. High and low water occur six hours apart. In the lunisolar theory this is explained very naturally as an immediate consequence of its basic principles. The rotation of the earth takes 24 hours. There’s a wave of high water pointing toward the moon, basically (give or take a bit of lag in the system), and then another high water mark diametrically opposite to it on the other side of the earth. So that’s two highs and two lows in 24 hours, so 6 hours apiece.

Galileo's theory, on the other hand, implies that high and low water should be twelve hours apart rather than six, as he himself says: ``there resides in the primary principle no cause of moving the waters except from one twelve-hour period to another.''\footnote{\cite[432]{galileo2sys1sted}.} So Galileo immediately find himself on the back foot, having to somehow talk himself out of this obvious flaw of his theory. To this end he alleges that “the particular events observed [regarding tides] at different times and places are many and varied; these must depend upon diverse concomitant causes,”\footnote{\cite[428]{galileo2sys1sted}.} such as the size, depth and shape of the sea basin, and the internal forces of the water trying to level itself out.\footnote{\cite[428--432]{galileo2sys1sted}.} The fact that everyone could observe two high and two low tides per day Galileo thus wrote off as purely coincidental:
\quote{Six hours …\ is not a more proper or natural period for these reciprocations than any other interval of time, though perhaps it has been the one most generally observed because it is that of our Mediterranean.\footnote{\cite[432--433]{galileo2sys1sted}.}}
Galileo even mistakenly believed that he had data to prove his erroneous point, namely that tides twelve hours apart are ``daily observed in Lisbon.''\footnote{\cite[177]{sheaGrev}, \cite[128]{GA}. OGG.V.388--389.} That is not in fact how tides behave in Lisbon. Others corrected Galileo on this point, which is why he did not include this claim in the published {\it Dialogue}. But that did not stop him from publishing the rest of his tidal theory without change.

The fact that Mediterranean and Atlantic tidal periods are identical is obviously a huge problem for Galileo's theory. The lunar theory of the tides immediately and correctly explains why these two periods are equal, and why they are both six hours specifically. According to Galileo's theory, both the equality and the six-hour period are purely coincidental. Any reasonable scientist interested in objectively evaluating these two hypotheses would find these facts significant. And Galileo clearly understood this perfectly well, which is why he was so keen to seize upon the false Lisbon data. Yet, in the published {\it Dialogue}, Galileo simply suppresses these inconvenient facts. Ever the opportunist, Galileo is very happy to highlight the importance of certain data when he believes it proves his point, only to then pretend that that data doesn't exist when it turns out it proved the opposite of what he wanted.

There is a further complication involved in Galileo's theory, which ``caused embarrassment to his more competent readers.''\footnote{\cite[182]{sheaGrev}.} The inclination of the earth's axis implies that the effects Galileo describe should be strongest in summer and winter. Unfortunately for Galileo’s theory, the reverse is the case. Actually the tides are most extreme in spring and fall because they receive the maximum effects of the sun’s gravitational pull.\footnote{Drake, \cite[490]{galileo2sys1sted}.} Galileo got himself confused on this point because he was again relying on false data. Galileo---the self-declared enemy of relying on textual authority, who often mocked his opponents for believing things simply because it said so in some book---was the one in this case who found in some old book the claim that tides are greatest in summer and winter, took this for fact and derived this supposed effect from his own theory.\footnote{\cite[296]{drakeGatwork}.}

The embarrassing mismatch between Galileo's theory and basic facts is on display in another episode as well:
\quote{Galileo …\ attacked [those who] postulated that an attractive force acted from the Moon on the ocean for failing to realize that water rises and falls only at the extremities and not at the center of the Mediterranean. [But his opponents] can hardly be blamed for failing to detect this [so-called] phenomenon: it only exists as a consequence of Galileo’s own theory.\footnote{\cite[420]{sheaOxSW}.}}
In other words, Galileo was so biased by his wrongheaded theory that he used its erroneous predictions as ``facts'' with which to attack those who were actually right.

But all of the above is not yet the worst of it. There is an even more fundamental flaw in Galileo's theory: it is inconsistent with the principle of relativity he himself espoused. Think back to his scenario of the scientist locked in a cabin below deck of a ship that may or may not be moving.\footnote{\S\ref{relativity}.} Galileo's conclusion on that occasion was that no physical experiment could detect whether the ship was moving or not. But the torus tidal simulation, if it really worked as Galileo claimed, certainly could detect such a motion. If we put the torus on the floor of the cabin and spun it, one part would be spinning with the direction of motion of the ship, and another part against it. Hence high and low water should arise, by the same logic as in Galileo's tidal theory. If the ship stood still, on the other hand, no such effect would be observed, of course. So we have a way of determining whether the ship is moving, which is supposed to be impossible. And indeed it is impossible. But if that’s so then Galileo’s theory of the tides cannot possibly work because it is inconsistent with this principle.

This objection against Galileo’s theory was in fact raised immediately already by contemporary readers:
\quote{They draw attention to a difficulty raised by several members about the proposition you make that the tides are caused by the unevenness of the motion of the different parts of the earth. They admit that these parts move with greater speed when they [go] along [with] the annual motion than when they move in the opposite direction. But this acceleration is only relative to the annual motion; relative to the body of the earth as well as to the water, the parts always move with the same speed. They say, therefore, that it is hard to understand how the parts of the earth, which always move in the same way relative to themselves and the water, can impress varying motions to the water.\footnote{Jean-Jaques Bouchard to Galileo, 1633, \cite[176]{sheaGrev}.}}
That is to say, picture the earth moving along its orbit and also rotating around its axis at the same time. Hit pause on this animation and mark two diametrically opposite spots on the equator. Then hit play, let it run for a second or two, and then pause it again. Now, compare the new positions of the two marked spots with their original position. One will have moved further than the other. But that’s in a coordinate system that doesn’t move with the earth. That type of inequality of speed is irrelevant. What is needed to create tides is a different kind of inequality of speed of the water. Namely, a difference in speed relative to the earth itself and to the other water. Tides arise when a fast-moving part of the water catches up with a slow-moving part of the water. But that is to say, these waters are fast and slow in their speed of rotation about the earth. So inequality of speed in a coordinate system centred on the earth. But no inequality of this type arises from the motion of the earth about the sun. Galileo had no solution to this accurate objection.

So, to sum up, Galileo small-mindedly rejected the correct theory of the tides, based on the sun and the moon, even though this was widely understood by his contemporaries. He then proposed a completely wrongheaded theory of his own, which is based on elementary errors of physical reasoning that are inconsistent with his own principles. These flaws were readily spotted by his contemporaries. Furthermore, his theory is fundamentally at odds with the most basic phenomena, which he tried to explain away by attributing them to untestable, ad hoc secondary effects. He also adduced several false observational so-called “facts” in support of his theory.

No wonder many have felt that Galileo’s “ill-fated theory of the tides is a skeleton in the cupboard of the scientific revolution.”\footnote{\cite[186]{sheaGrev}.} But this is a problem only if one assumes that Galileo was science personified. If we accept instead that Galileo was an exceptionally mediocre mind, who constantly got wrong what mathematically competent people like Kepler got right, then we see that Galileo’s skeletons belong only to himself, not to the scientific revolution. It’s not that the scientific revolution was flawed. It’s just that Galileo was. If we restrict ourselves to mathematically competent people then we don’t have to deal with this kind of nonsense.

\section{Astronomy}

\subsection{Adoption of Copernicanism}
\label{adoption}

Does the earth move around the sun, or vice versa? Copernicus worked out the right answer long before Galileo was even born, as did the best Greek mathematicians thousands of years earlier.\footnote{\S\ref{Greekastro}.} Yet somehow Galileo has ended up with much of the credit:
\quote{If one wonders why the Copernican theory, with almost no adherents at the beginning of the seventeenth century, had pretty much swept the field by the middle, the answer …\ is above all [Galileo's] {\it Dialogue}.\footnote{Swerdlow, \cite[267]{camcomp}.}}
\quote{Galileo wrote the book that won the war [and] made belief in a moving earth intellectually respectable.\footnote{\cite[68]{GingerichCopVS}.}}
This may be right in a limited sense: maybe indeed the ignorant masses needed a book like Galileo's to dumb it down for them before they could finally come to their senses. But mathematically competent people were already convinced long before and had no use for Galileo telling them the ABCs.

When Copernicus made the earth go around the sun, he confidently declared that ``I have no doubt that talented and learned mathematicians will agree with me.''\footnote{\cite[6]{CopRevTransl}.} He was right. In 1600, long before Galileo had published a single word on the matter, there were already over a dozen committed Copernicans. One historian counted ten, besides Copernicus himself: ``Thomas Digges and Thomas Harriot in England; Giordano Bruno and Galileo Galilei in Italy; Diego de Zuniga in Spain; Simon Stevin in the Low Countries; and, in Germany, the largest group---Georg Joachim Rheticus, Michael Maestlin, Christopher Rothmann, and Johannes Kepler.''\footnote{\cite[136]{westmanastrrole}.} This list is not complete. Gemma Frisius publicly ``endorsed physical Copernicanism,''\footnote{\cite[146]{BarkerCop}.} and John Feild ``apparently accepted the views of Copernicus without qualifications.''\footnote{\cite[192]{Dobrzycki}.} ``[Paolo] Sarpi's cosmology was Copernican and heliocentric'' as well: although he did not explicitly declare this publicly, ``the way in which Sarpi referred to Copernicus proves his acceptance of the latter's theory.''\footnote{cite[98]{Kainulainen}.}

Another early supporter was Achilles Gasser, a contemporary of Copernicus. He praises Copernicus for having brought about “the restoration--or rather, the rebirth---of a true system of astronomy.” Much like Copernicus himself, Gasser considers it a foregone conclusion that competent mathematicians will see the truth at once, while the ignorant should be dismissed:
\quote{[Copernicus] has not only demonstratively proven his theory among the mathematicians, $\ldots$ but has also immediately been regarded as having perpetrated a heresy, and indeed---by many others incapable of understanding this matter---is already being condemned.}
But one must “set truth free $\ldots$ in spite of the critical gaze of the plebs.” Ultimately, “there is no doubt that this new thing [i.e., heliocentrism] will one day be accepted without bitterness by all educated people as something both agreeable and useful.”\footnote{Gasser, writing in the 1540s. \cite[461--462, 465]{Danielson}.}

That makes fifteen in total. Fifteen professed believers in the new astronomy, before Galileo has published anything. Such a number already at this early stage is not ``almost no adherents.'' What do you expect? How many ``talented and learned mathematicians'' do you think there were in pest-ridden, blood-letting, witch-burning Europe of 1600? And how many of them were interested in the Copernican question and formed an opinion on it, even though that was a philosophical question beyond the scope of the official computational task of the astronomer? And, among those in turn, how many were prepared to declare allegiance to a flagrantly heretical opinion in an age where the religious thought-police routinely employed vicious torture and burnt dissenters alive? Including one person on the very list just mentioned, in fact.\footnote{Giordano Bruno. \cite{Bruno2018}.} So fifteen avowed Copernicans may well be regarded as quite a crowd considering the circumstances.

Galileo's own assessment of the situation in 1597 is apt indeed:
\quote{I have preferred not to publish, intimidated by the fortune of our teacher Copernicus, who though he will be of immortal fame to some, is yet by an infinite number (for such is the multitude of fools) laughed at and rejected.\footnote{Galileo to Kepler, 4 August 1597, \cite[41]{drakeGatwork}.}}
This confirms that social pressures to avoid the issue were very real---enough to ``intimidate'' Galileo and surely many others. More importantly, Galileo is also making my main point for me: even at this very early stage---long before Galileo has published a single word on the matter, and long before the invention of the telescope---every serious astronomer has rejected the old astronomy already. By Galileo's own reckoning, there were only ``fools'' left to convince. On this point he is exactly right.

There were most likely a seizable number of closet Copernicans who figured ``don't ask, don't tell'' was the best policy to avoid needless conflicts with the intolerant. ``Scriptural and theological concerns, attempts at censorship, apologies and forms of self-censorship accompanied Copernicus’s ideas from the beginning.''\footnote{\cite[271]{OmodeoCop}.} “The wars of religion may well have inhibited outright advocacy of heliocentrism by late sixteenth-century Frenchmen,” for example: it is “suggestive of a causal relationship” that initially positive attitudes toward Copernicus gave way to negative ones in step with religious strictures.\footnote{\cite[83, 85]{BaumgartnerScep}.}

Kepler faced theological pushback for the heliocentrism of his {\it Mysterium cosmographicum} (1597). Regarding one of his critics, Kepler says: ``I really cannot believe that he is opposed to this doctrine [Copernicanism]. He pretends [to be so] in order to appease his colleagues.'' Indeed, Kepler himself recommended a similar attitude to his friends: ``If anybody approaches us privately, let us tell him our opinion frankly. Publicly let us be silent.''\footnote{\cite[63, 66]{VoelkelComp}.}

By some indications, William Gilbert was ``a probable Copernican,'' even though he stopped short of openly endorsing it in print.\footnote{\cite[144, 155]{BarkerCop}.} The same goes for Gemma Frisius's student Johannes Stadius, who ``was inclined to accept the Copernican system.''\footnote{\cite[156]{OmodeoCop}.} Similarly, Harriot had a number of followers in England who were enthusiastic believers in heliocentrism, but ``the local political and intellectual milieu …\ forced them into what we can term preventive self-censorship.''\footnote{\cite[146]{Bucciantini}.} Another example is Mersenne, who does not go on the official list because ``at no time during his life did he find any proof so overwhelming that he felt like challenging the Church on the matter.''\footnote{\cite[20]{HineMersenneCop}.} And this despite the fact that he was one of the most enthusiastic readers of Galileo's {\Dialogue}. He remained uncommitted for political reasons, it seems.

Another indication that many were silently receptive to Copernicanism is the fact that most of the leading astronomers of the 16th century owned Copernicus's book, and many of them wrote extensive notes in the margins,\footnote{\cite[xvi]{GingerichCencus}.} as was the habit at the time. Books were printed with generous margins because everyone was expected to take detailed notes as they read. And indeed they did. Owen Gingerich conducted a thorough census of all surviving copies of Copernicus's book. He looked at all this marginalia that this large group of serious, competent readers of Copernicus's book has written. That's a group far larger than those fifteen mentioned above. Some of them were probably secretly convinced that Copernicus was right; others studied the work because they saw it as they duty to keep up with the best technical mathematical astronomy of the day whether they agreed with it or not. Either way, they took meticulous notes as they painstakingly worked through this long and technical treatise.

This large group of serious astronomers did not include Galileo, however. Galileo's dilettantism is so blatant and shameless that Gingerich could hardly believe his eyes:
\quote{I had long supposed that Galileo was not the sort of astronomer who would have read Copernicus' book to the very end. …\ Still, when I saw the copy in Florence, my reaction was one of scepticism that it was actually Galileo's copy, since there were so few annotations in it. …\ This copy had no technical marginalia, in fact, no penned evidence that Galileo had actually read any substantial part of it. …\ Eventually, …\ I realized that my scepticism was unfounded and that it really was Galileo's copy.\footnote{\cite[142, 200]{GingerichNobody}.}}
But there is no need for surprise. Galileo was a poor mathematician. He did not have the patience or ability to understand serious mathematical astronomy, let alone make any contribution to it.

\subsection{Pre-telescopic heliocentrism}
\label{pretelescopehelio}

\begin{figure}[pt]\centering

\tikzfig{

\draw[gray,dotted,rotate=90,postaction={decorate,decoration={text along path, text align=center,text={|\footnotesize\color{gray}|fixed stars}}}] (0,0) circle (6.4cm);

\draw[gray,rotate=90,postaction={decorate,decoration={text along path,text align=center,text={|\footnotesize|Saturn}}}] (0,0) circle (5.6cm);
\draw[gray,rotate=90,postaction={decorate,decoration={text along path,text align=center,text={|\footnotesize|Jupiter}}}] (0,0) circle (4.8cm);
\draw[gray,rotate=90,postaction={decorate,decoration={text along path,text align=center,text={|\footnotesize|Mars}}}] (0,0) circle (4.0cm);
\draw[gray,rotate=90,postaction={decorate,decoration={text along path,text align=center,text={|\footnotesize|Sun}}}] (0,0) circle (3.2cm);
\draw[gray,rotate=90,postaction={decorate,decoration={text along path,text align=center,text={|\footnotesize|Venus}}}] (0,0) circle (2.4cm);
\draw[gray,rotate=90,postaction={decorate,decoration={text along path,text align=center,text={|\footnotesize|Mercury}}}] (0,0) circle (1.6cm);
\draw[gray,rotate=90,postaction={decorate,decoration={text along path,text align=center,text={|\footnotesize|Moon}}}] (0,0) circle (0.8cm);

\node[inner sep=0pt] (earth) at (0,0)
    {\includegraphics[width=0.5cm]{pic/earth.pdf}};
\node[inner sep=0pt] (earth) at (-0.4,0.693)
    {\includegraphics[width=0.2cm]{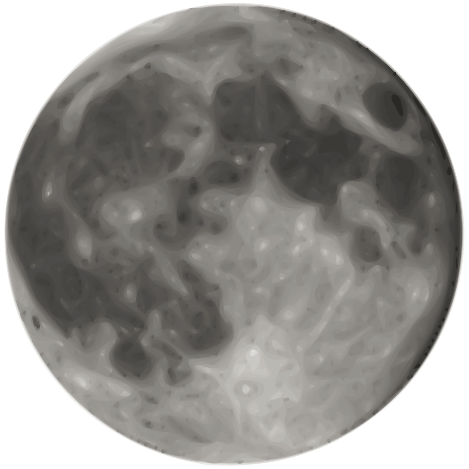}};
\node[inner sep=0pt] (sun) at (3.151,0.556)
    {\includegraphics[width=0.5cm]{pic/sun.pdf}};

\node[inner sep=0pt] (mercury) at (1.576,-0.278)
    {\includegraphics[width=0.15cm]{pic/mercury.png}};
\node[inner sep=0pt] (venus) at (2.175,1.014)
    {\includegraphics[width=0.2cm]{pic/venus.png}};
\node[inner sep=0pt] (mars) at (-0.695,3.939)
    {\includegraphics[width=0.2cm]{pic/mars.png}};
\node[inner sep=0pt] (jupiter) at (-2.4,-4.157)
    {\includegraphics[width=0.35cm]{pic/jupiter.png}};
\node[inner sep=0pt] (saturn) at (-4.850,2.8)
    {\includegraphics[width=0.4cm]{pic/saturn.png}};

}{scale=0.8}

\caption[Ptolemy's system.]{Ptolemy's system. All heavenly bodies orbit the stationary earth.}
\label{ptolemyfig}
\end{figure}

Why were people convinced Copernicans before Galileo? Arguably the most compelling reason was that the Copernican system explained complex phenomena in a simple, unified way. Here are some basics observational facts: The planets sometimes move forwards and sometimes backwards. Also, Mercury and Venus always remain close to the sun. In a geocentric system, with the earth in the center, there is no inherent reason why those things should be like that. Nevertheless these are the most basic and prominent facts of observational astronomy. In the Ptolemaic system (Figure \ref{ptolemyfig}), the geocentric system of antiquity, these facts are accounted for by introducing complicated secondary effects and coordinations beyond the basic model of simple circles. So planetary orbits are not just simple circles but combinations of circles in complicated ways that also happen to be coordinated with one another in particular patterns. In the Ptolemaic system there is no particular reason why these complicated constructions should be just so and not otherwise. We have to accept that it just happens to be that way. So Ptolemy could account for, or accommodate, the phenomena, but he can hardly have been said to have explained them. The basic idea, that planets move in circles around the earth, is on the back foot from the outset. It is inconsistent with the most basic data and is therefore forced to add individual quick-fixes for these phenomena.

\begin{figure}[pt]
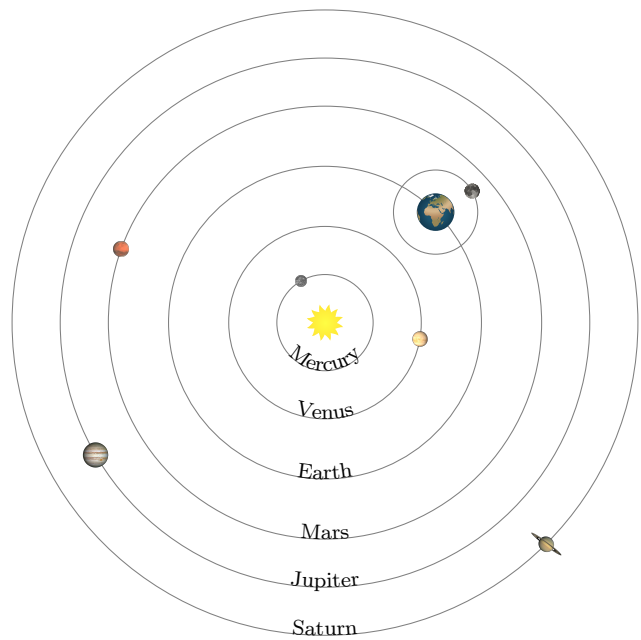

\tikzfig{

\draw[gray,rotate=90,postaction={decorate,decoration={text along path,text align=center,text={|\footnotesize|Saturn}}}] (0,0) circle (5.2cm);
\draw[gray,rotate=90,postaction={decorate,decoration={text along path,text align=center,text={|\footnotesize|Jupiter}}}] (0,0) circle (4.4cm);
\draw[gray,rotate=90,postaction={decorate,decoration={text along path,text align=center,text={|\footnotesize|Mars}}}] (0,0) circle (3.6cm);
\draw[gray,rotate=90,postaction={decorate,decoration={text along path,text align=center,text={|\footnotesize|Earth}}}] (0,0) circle (2.6cm);
\draw[gray,rotate=90,postaction={decorate,decoration={text along path,text align=center,text={|\footnotesize|Venus}}}] (0,0) circle (1.6cm);
\draw[gray,rotate=90,postaction={decorate,decoration={text along path,text align=center,text={|\footnotesize|Mercury}}}] (0,0) circle (0.8cm);
\draw[gray] (1.838,1.838) circle (0.7cm);

\node[inner sep=0pt] (earth) at (1.838,1.838)
    {\includegraphics[width=0.5cm]{pic/earth.pdf}};
\node[inner sep=0pt] (earth) at (2.444,2.188)
    {\includegraphics[width=0.2cm]{pic/moon.pdf}};
\node[inner sep=0pt] (sun) at (0,0)
    {\includegraphics[width=0.5cm]{pic/sun.pdf}};

\node[inner sep=0pt] (mercury) at (-0.4,0.693)
    {\includegraphics[width=0.15cm]{pic/mercury.png}};
\node[inner sep=0pt] (venus) at (1.576,-0.278)
    {\includegraphics[width=0.2cm]{pic/venus.png}};
\node[inner sep=0pt] (mars) at (-3.383,1.231)
    {\includegraphics[width=0.2cm]{pic/mars.png}};
\node[inner sep=0pt] (jupiter) at (-3.811,-2.2)
    {\includegraphics[width=0.35cm]{pic/jupiter.png}};
\node[inner sep=0pt] (saturn) at (3.677,-3.677)
    {\includegraphics[width=0.4cm]{pic/saturn.png}};

}{scale=0.8}
\caption[Copernicus's system.]{Copernicus's system. The planets orbit the stationary sun.}
\label{copernicusfig}
\end{figure}

In the Copernican system (Figure \ref{copernicusfig}), it’s the opposite. The phenomena are here instead an immediate, natural consequence of the motion of the earth. It becomes obvious and unavoidable that outer planets appear to stop and go backwards as the earth is speeding past them in its quicker orbit. It becomes obvious and unavoidable that Mercury and Venus are never seen far from the sun, since, like the sun, they are always on the “inside” of the earth’s orbit. Ad hoc secondary causes and just-so numerical coincidences in parameter values are no longer needed to accommodate these facts; instead they follow at once from the most basic assumptions of the system.

Galileo makes this point in the {\Dialogue}:
\quote{What are we to say of the apparent movement of a planet, so uneven that it not only goes fast at one time and slow at another, but even goes backward a long way after doing so? To save these appearances, Ptolemy introduces vast epicycles, adapting them one by one to each planet, with certain rules for incongruous motions---all of which can be done away with by one very simple motion of the earth.\footnote{\cite[397]{galileo2newsci2ndedWT}. It is not true, of course, that all epicycles can be done away with. Galileo is exaggerating and oversimplifying, as ever.}}
\quote{You see, gentlemen, with what ease and simplicity the annual motion---if made by the earth---lends itself to supplying reasons for the apparent anomalies which are observed in the movements of the five planets. …\ It removes them all. …\ It was Nicholas Copernicus who first clarified for us the reasons for this marvelous effect.\footnote{\cite[400]{galileo2newsci2ndedWT}.}}
``This alone ought to be enough to gain assent for the rest of the [Copernican] doctrine from anyone who is neither stubborn nor unteachable.''\footnote{\cite[398]{galileo2newsci2ndedWT}.}

This is all good and well. Galileo is absolutely right. But this was a core point of Copernicanism from the outset, that was already a hundred years old and common knowledge by the time Galileo repeated it. Copernicus himself made the point clearly. He challenged those who deny heliocentrism to explain why Mercury and Venus always stay close to the sun: “what cause will they allege why these planets do not also make longitudinal circuits separate and independent of the sun, like the other planets”?\footnote{\cite[21]{CopRevTransl}.} Similarly, retrograde motion of the outer planets arises because “the movement of the Earth is speedier, so that it outruns the movement of the planet.”\footnote{\cite[240]{CopRevTransl}.} Thus the motion of the earth gives a unified explanation for what the deniers of heliocentrism must attribute to incidental properties of individual planetary models: “And so it is once more manifest that all these apparent movements---which the ancients were looking into by means of the epicycles of the individual planets---occur on account of the movement of the Earth.”\footnote{\cite[240]{CopRevTransl}.} Thus is it striking “what power and effect the assumption of the revolution of the Earth has in the case of the apparent movement in longitude of the wandering stars [i.e.\ planets] and in what sure and necessary order it places all the appearances. …\ Accordingly by means of the assumption of the mobility of the Earth we shall do with perhaps greater compactness and more becomingly what the ancient mathematicians thought to have demonstrated by means of the immobility of the Earth.”\footnote{\cite[311]{CopRevTransl}. The force of this argument was acknowledged by Tycho Brahe, among others; see \S\ref{tycho}.}

It is instructive, however, to compare this with Galileo's theory of tides.\footnote{\S\ref{tides}.} The correct, lunisolar theory of the tides explains the basic phenomena in a simple and natural way as immediate consequences of the first principles of the theory. That’s exactly the same point that we made about the Copernican system. The correct theory of the tides thus has the same kind of credibility as the Copernican system. So, by Galileo’s logic, this “ought to be enough to gain assent.” But in the case of the tides it seems Galileo was the “stubborn and unteachable” one. He insisted on a theory which---like that of Ptolemy---could only account for basic facts by invoking arbitrary and unnatural secondary causes unrelated to the primary principles of the theory.

It’s a sign that your theory has poor foundations if the foundations themselves are good for nothing and all the actual explanatory work is being done by emergency extras duct-taped on later to specifically fix obvious problems with the foundations. Intelligent people realised this, which is why they turned to the sun-centered view of the universe. Galileo paid lip service to the same principle when he wanted to ride on the coattails of their insights. But, if he had been consistent in his application of this principle, he should have used it to reject his foolish theory of the tides.

\subsection{Tycho Brahe's system}
\label{parallax}\label{tycho}

If you move from one side of room to another, your view of everything on the walls will change. The wall you are approaching will appear to ``grow,'' while the wall behind you will shrink and occupy a smaller part of your field of vision. This is called parallax. 

A major problem with Copernicus's theory is the absence of stellar parallax in the course of a year. If the earth moves in an enormous circle around the sun, we should be at one moment close to some constellation of stars, and then half a year later much further away from them. Hence we should see them sometimes big and ``up close,'' and sometimes shrunk into a small area, like a faraway wall at the end of a long corridor. But this does not happen. The night sky is immutable. As far as 17th-century astronomers could detect, the constellations all look exactly the same throughout the year, just as if we never move an inch.

To maintain the earth's motion in spite of this, it is therefore necessary to postulate, as Copernicus does, that ``the fixed stars …\ are at an immense height away.''\footnote{\cite[27]{CopRevTransl}.} The diameter of the earth's orbit is so small in relation to such an astronomical distance that our feeble little motion is all but tantamount to standing still. That is why no parallax can be detected. This is the correct explanation, as we now know, but in the 16th century it didn't sound too convincing.

Tycho Brahe was one of the sceptics. He was the most exacting astronomical observer in the pre-telescopic era, but even he, with his very advanced and precise observations, could find no parallactic effect. He calculated that, if Copernicus was right, the stars would have to be at least 700 times more distant than Saturn for this to happen. The universe would not have been designed with so much wasted space, he reasoned.\footnote{\cite[253]{Siebert}.} 

\begin{figure}[pt]
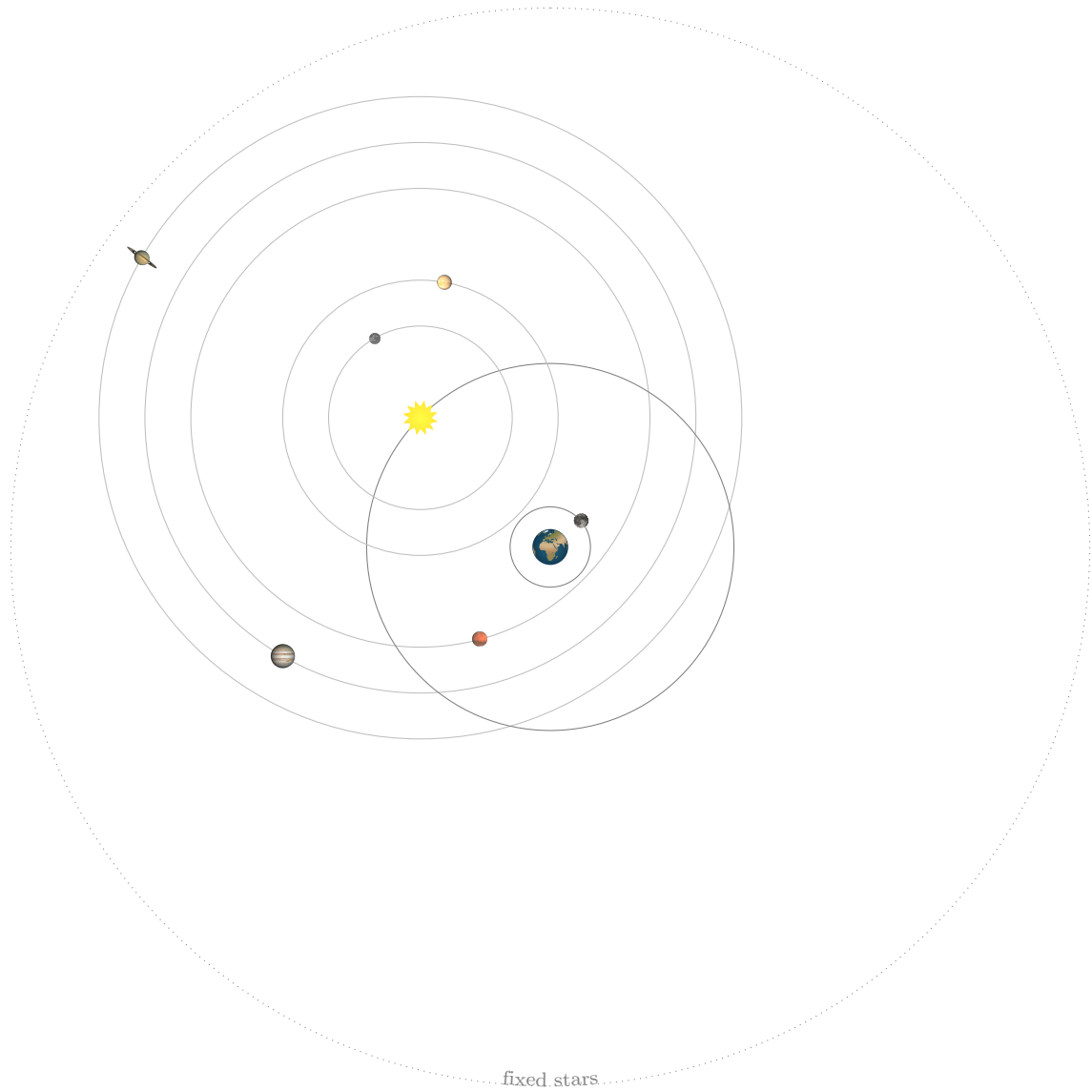
\centering
\tikzfig{

\coordinate (S) at (-2.263,2.253);

\draw[gray,dotted,rotate=90,postaction={decorate,decoration={text along path, text align=center,text={|\footnotesize\color{gray}|fixed stars}}}] (0,0) circle (9.4cm);

\draw[lightgray] (S) circle (5.6cm);
\draw[lightgray] (S) circle (4.8cm);
\draw[lightgray] (S) circle (4.0cm);
\draw[gray] (0,0) circle (3.2cm);
\draw[lightgray] (S) circle (2.4cm);
\draw[lightgray] (S) circle (1.6cm);
\draw[gray] (0,0) circle (0.7cm);

\node[inner sep=0pt] (earth) at (0,0)
    {\includegraphics[width=0.5cm]{pic/earth.pdf}};
\node[inner sep=0pt] (earth) at (0.536,0.459)
    {\includegraphics[width=0.2cm]{pic/moon.pdf}};
\node[inner sep=0pt] (sun) at (S)
    {\includegraphics[width=0.5cm]{pic/sun.pdf}};

\node[inner sep=0pt] (mercury) at (-3.063,3.639)
    {\includegraphics[width=0.15cm]{pic/mercury.png}};
\node[inner sep=0pt] (venus) at (-1.846,4.612)
    {\includegraphics[width=0.2cm]{pic/venus.png}};
\node[inner sep=0pt] (mars) at (-1.227,-1.610)
    {\includegraphics[width=0.2cm]{pic/mars.png}};
\node[inner sep=0pt] (jupiter) at (-4.663,-1.904)
    {\includegraphics[width=0.35cm]{pic/jupiter.png}};
\node[inner sep=0pt] (saturn) at (-7.112,5.053)
    {\includegraphics[width=0.4cm]{pic/saturn.png}};

}{scale=0.8}
\caption[Tycho Brahe's system.]{Tycho Brahe's system. The planets orbit the sun, while the sun orbits the stationary earth.}
\label{tychofig}
\end{figure}

Tycho Brahe therefore devised a system of his own (Figure \ref{tychofig}), in which the earth remained the center of the universe, while the planets orbit the sun. This removed the problem of parallax. Moreover, the arguments for the Copernican system based on explanatory simplicity (\S\ref{pretelescopehelio}) apply equally well to the Tychonic system. Tycho indeed acknowledged that Copernicus ``expertly and completely circumvents all that is superfluous or discordant in the system of Ptolemy.''\footnote{Tycho, \cite[96]{DuhemSavePhen}.} As we have seen, Galileo repeated this old argument in his \Dialogue as proof for the Copernican hypothesis, while conveniently neglecting to acknowledge that serious mathematical astronomers from Tycho onwards had already reconciled the full force of this argument with geocentrism.

In fact, Tycho's system is equivalent to the Copernican one as far as the relative positions of the heavenly bodies are concerned. Tycho and Copernicus describe the same planetary motions, but they choose a different reference point in terms of which to describe them. Kepler illustrates the point with an analogy: the same circle can be traced on a piece of paper by either rotating the pen arm of a compass around the fixed leg, or by keeping the compass fixed while rotating the paper underneath it.\footnote{Kepler, {\it Harmonices Mundi} (1619), V, \cite[175]{KeplerEpiHar}.} Because of this equivalence, the Copernican and Tychonic systems are by necessity on equal footing as far as the simplicity arguments of \S\ref{pretelescopehelio} are concerned.

One might feel that the Tychonic system is less physically plausible than those of Ptolemy or Copernicus. Indeed, traditional conceptions had it that planets were enclosed in translucent crystalline spheres, like the layers of an onion. Both the Ptolemaic and Copernican systems are basically compatible with such an ``onion'' conception of the cosmos. The Tychonic system clearly is not: planets are crossing each other's ``orbs'' all over the place. But Tycho had some good counterarguments to this.\footnote{\cite[16]{KeplerEpiHar}.} By a careful study of the paths of comets, he proved that they evidently passed through the alleged crystalline spheres with ease. Furthermore, he pointed out that these alleged spheres did not refract light, as glass or other materials had been known to do since antiquity.

All in all, Tycho's system was a serious scientific theory with good arguments to its credit. This is another reason why our headcount of Copernicans above is misleading.\footnote{\S\ref{adoption}.} The number of people who rejected the Ptolemaic system was certainly greater than the number of outright Copernicans, and the middle road of Tycho was by no means blind conservatism but rather a viable system based on the latest mathematical astronomy.

Galileo, however, liked to pretend otherwise. The full title of his famous book reads: {\it Dialogue Concerning the Two Chief World Systems: Ptolemaic and Copernican}. It certainly made Galileo's life a lot easier to frame his fictional debate with fictional opponents in those antiquated terms. That way he could battle two-thousand-year-old ideas instead of having to engage with the latest mathematical astronomy.

More serious and mathematically competent people had a very different view of which were ``the two chief world systems.'' Around 1600, long before Galileo enters the fray, Kepler considers it obvious that the Ptolemaic system is obsolete:
\quote{Today there is practically no one who would doubt what is common to the Copernican and Tychonic hypotheses, namely, that the sun is at the centre of motion of the five planets, and that this is the way things are in the heavens themselves---though in the meantime there is doubt from all sides about the motion or stability of the sun.\footnote{\cite[147]{KeplerApologia}.}}
Later, after the telescope has brought its new evidence, not much has changed. Kepler is a bit more assured that ``today it is absolutely certain among all astronomers that all the planets revolve around the sun.''\footnote{Kepler, {\it Harmonices Mundi} (1619), V, \cite[175]{KeplerEpiHar}.} But the battle between Copernicus and Tycho remained far from settled: ``…\ the hypothesis not only of Copernicus but also of Tycho Brahe, whereof either hypotheses are today publicly accepted as most true, and the Ptolemaic as outmoded.''\footnote{Kepler, {\it Harmonices Mundi} (1619), V, \cite[169]{KeplerEpiHar}.} ``The theologians may decide which of the two hypotheses …\ ---that of Copernicus or that of Brahe---should henceforth be regarded as valid [for] the old Ptolemaic is surely wrong.''\footnote{Kepler, 1619, \cite[136]{BaumgardtKepler}.}

Nor was this merely Kepler's opinion. Historians who study much more minor figures have also concluded that, indeed, ``the Ptolemaic system already had been set aside, at least among mathematical astronomers,''\footnote{\cite[208]{Magruder}. See also \cite[210--211]{Bucciantini}.} well before Galileo wrote his {\Dialogue}. But Galileo, in his great book, like a schoolyard bully secretly too scared to pick on someone his own size, preferred to pretend that the old Ptolemaic system was still the enemy of the day. To be sure, there were still a ``multitude of fools'' left to convince, and perhaps indeed Galileo did so more effectively than anyone else. But that proves at most that Galileo should be praised as a populariser, not as a scientist. To mathematically competent astronomers he was beating a dead horse.

\subsection{Against Tycho}
\label{againsttycho}

Some scholars have made too much of the above arguments in favour of the Tychonic system, however. They have concluded that ``it is fair to say that, contrary to [the standard view], science backed geocentrism.''\footnote{\cite[132]{Graney}.} If you take the works of Galileo to be the extent of ``science'' then this conclusion is indeed defensible. But the conclusion is false if you include genuinely talented scientists like Kepler. Unlike Galileo, Kepler dared to take on the Tychonic system and he gave a long list of compelling arguments against it.\footnote{{\it Astronomia Nova} (1609), 3r--3v, \cite[20--23]{KeplerAN}, {\it Epitome Astronomiae Copernicanae} (1620), IV.5, \cite[71--76]{KeplerEpiHar}.} Tycho's system is equivalent to the Copernican one in terms of relative position of the planets, but ``that Copernicus is better able than Brahe to deal with celestial physics is proven in many ways.''\footnote{Kepler, {\it Astronomia Nova} (1609), 3r, \cite[21]{KeplerAN}.} Let us sample just two of them. One is based on ``the magnitude of the moveable bodies'':
\quote{For just as Saturn, Jupiter, Mars, Venus, and Mercury are all smaller bodies than the solar body around which they revolve; so the moon is smaller than the Earth …\ [and] so the four [moons] of Jupiter are smaller than the body of Jupiter itself, around which they revolve. But if the sun moves, the sun which is the greatest …\ will revolve around the Earth which is smaller. Therefore it is more believable that the Earth, a small body, should revolve around the great body of the sun.\footnote{Kepler, {\it Epitome Astronomiae Copernicanae} (1620), IV.5, \cite[73]{KeplerEpiHar}.}}
This is of course completely correct also from a modern point of view. It is also verified by elementary experience. For instance, take a lead ball and a ping-pong ball, and tie them together with a string. If you flick the ping-pong ball it will start spinning around the lead ball. But if you flick the lead ball it will roll straight ahead without any regard for the ping-pong ball, which will simply be dragged along behind it. So the lighter object adapts its motion to the heaver one but not conversely, just as Kepler says. The planets are not tied to the sun with a string but the point generalises. You can observe the same principle with a big and a small magnet for example: the little one is moved by the bigger, not the other way around.

More generally, the Copernican picture is more readily susceptible to a mechanical explanation of planetary motion than the Tychonian picture. Such mechanical explanations were eagerly sought in the 17th century. Slings and magnets were powerful analogies for trying to conceptualise the kinds of forces involved in heavenly motions.\footnote{Kepler's {\it Astronomia Nova} (1609) gives a quasi-magnetic account of planetary motion in great detail.} Another popular conception was that the planets were carried along in their orbits by some sort of vortex made up of invisible particles.\footnote{Descartes, {\it Traité du monde et de la lumière} (1633).} Whichever mechanical model one favours, Copernican cosmology is much more naturally explained in mechanical terms than a Tychonic cosmos.

Another of Kepler's arguments against Tycho has to do with the motion of the planets in latitude,\footnote{{\it Epitome Astronomiae Copernicanae} (1620), IV.5, \cite[72]{KeplerEpiHar}.} that is to say ``up and down.'' Suppose you are in a dark room, and someone has place a glowing object in a sling and is whirling it around you. If you are asked to describe the motion you will notice first of all that it is a circular motion. But it may also be the case that the orbit of the object is not exactly parallel to the floor. You will therefore notice that the object is higher is some part of its orbit and lower in another part. This is the motion in latitude. The orbit of the object in the sling lies in a single plane, but it is inclined with respect to your plane of reference, which is a plane parallel to the floor through your eyes. If you walked around and observed the object being slung around in its orbit from various vantage points, the apparent ``height'' of the object at any given moment would depend both on its position in its orbit, and your position. Nevertheless, it would always be geometrically determined from the simple hypotheses that as the object orbits it remains in a single plane, and as you walk around your eye remains in a single plane, and these two planes are not the same but are rather inclined somewhat wit respect to one another.

This is exactly what happens with the planets. In the Copernican system, latitude motions can be explained very simply: each planet moves in a single plane, but these orbital planes are not perfectly aligned. All of the planes pass through the sun, but they differ by a few degrees in their inclinations. Seen from the earth, this will give the impression that the planets are moving a bit ``up and down'' in addition to their predominant motion in longitude along the ecliptic circle in the sky determined by the earth's orbit.

In the earth-static systems of Ptolemy and Tycho Brahe, such an explanation of planetary motion in latitude is not possible. Somehow, longitude and latitude motion must {\em inherently} belong the motion of the each planet. This is physically less plausible than the simple Copernican explanation. It is easy to imagine physical processes that produce circular motions that take place in a single plane. This happens for example with an object in a sling or a point on a rotating solid. But it is much more difficult to imagine a physical process that could move an object in a circle and also make it bob up and down in a systematic way at the same time.

This is just two of Kepler's many arguments against Tycho. We need not discuss them all. My point here is only that from the low quality of Galileo's case for Copernicus one cannot infer that ``big surprise: in 1651 the geocentric cosmology had science on its side.''\footnote{Owen Gingerich, back cover blurb of \cite{Graney}.} The mistake in this argument is the implicit assumption that Galileo represents ``science.'' Galileo was an opportunistic populariser. If you want to know the state of serious science at the time you have to read instead mathematicians like Kepler.

\subsection{The telescope}

The year is 1609. Galileo is well into his middle age; a frail man not infrequently bedridden with rheumatic or arthritic pains. Had he died from his many ailments, in this year, at the age of 45, he would have been all but forgotten today. He would have been an insignificant footnote in the history of science, no more memorable than a hundred of his contemporaries. It has often been said that mathematics is a young man's game. Newton had his {\it annus mirabilis} in his early twenties---``the prime of my age for invention,'' as he later said. Kepler was the same age when he finished his first masterpiece, the {\it Mysterium Cosmographicum} (1596). Galileo was already nearly twice this age, and he had nothing to show for it but some confused piles of notes of highly dubious value. In short, as a mathematician the ageing Galileo had proved little except his own mediocrity.

It is this run-of-the-mill nobody that first hears of a new invention: the telescope. Here was his chance. He only had to point this contraption to the skies and record what he saw. No need anymore for mathematical talent or painstaking scientific investigations. For twenty years he had tried and failed to gain scientific fame the hard way, but now a bounty of it lay ripe for the plunder. All you needed was eyes and being first.

Galileo first heard about the mysterious new ``optical tube for seeing things close'' in July 1609. A week or two later a traveller offered one for sale in Padua and Venice at an outrageous price---about twice Galileo's yearly salary.\footnote{\cite[139--140]{drakeGatwork}.} The enterprising salesman found no takers, but the sense of opportunity remained in the air. And it was an opportunity tailor-made for Galileo: finally a path to scientific fame that required only handiwork and none of that tiresome thinking in which he was so deficient.

The design of telescopes was still a trade secret among the Dutch spectacle-makers who had stumbled upon the discovery. But acting fast was of the essence. Making a basic telescope is not rocket science. Soon many people figured out how to make their own. ``It took no special talent or unique inventiveness to come up with the idea that combining two different lenses …\ would create a device allowing people to see faraway objects enlarged.''\footnote{\cite[22]{Bucciantini}.} Reading glasses and magnifying glasses were already in common use: they obviously made text appear bigger, so it was not a far-fetched idea to use them to magnify more distant objects as well. And the external shape of the telescopes people reported seeing suggested that at least two lenses were combined in a long cylinder.\footnote{It is generally believed that Galileo had not seen a telescope before he made his own, but it is possible that he might have. \cite{biagiolitelescopecopied}. If not, he was relying on eye-witness accounts of others.} It didn't take a genius, therefore, to soon strike upon the simple recipe Galileo found: take one convex and one concave lens and stick them in a tube, and look through the concave end. That's it. No theoretical knowledge of optics played any part in this; it was purely a matter of hands-on craftsmanship and trial-and-error.\footnote{As Galileo himself in effect says in the {\it Assayer} (1623), \cite[245]{GalileoDiscOp}.}

About a month after first hearing of the telescope, Galileo has managed to build his own, with $8\times$ magnification, and begins a campaign to leverage it into a more lucrative appointment for himself. He gives demonstrations to various important dignitaries---``to the infinite amazement of all,''\footnote{Galileo, letter to his brother-in-law, \cite[141]{drakeGatwork}.} according to himself---and enters multiple negotiations about improved career prospects. Between hands-on optical trials and lens grinding, showmanship demonstrations, shrewd self-marketing, and juggling potential job offers, Galileo must have had a busy couple of months indeed. On top of this his regular teaching duties at the university were just starting again in the fall.

We can easily understand, therefore, why the scientific importance of the new instrument for astronomy was not realised right away. At first neither Galileo nor anyone else thought of the telescope as primarily an astronomical instrument. Galileo instead tried to market it as ``a thing of inestimable value in all business and every undertaking at sea or on land,''\footnote{Galileo to Donà, 24 August 1609, \cite[149]{heilbron}.} such as spotting a ship early on the horizon. But the moon makes an obvious object of observation, especially at night when there is little else to look at. Perhaps indeed moon observations were part of Galileo's sales pitch routine more or less from the outset,\footnote{\cite[142]{drakeGatwork}.} though as a gimmick rather than science.

But this was soon to change. In the dark of winter, the black night sky is less bashful with its secrets than in summer. It monopolises the visible world from dinner to breakfast; it seems so eager to be seen that it would be rude not to look. In January, Galileo takes up the invitation and spots moons around Jupiter.\footnote{\cite[146]{drakeGatwork}.} This changes everything. Suddenly it is clear that the telescope is the key to a revolution in astronomy. Eternal scientific fame is there for the taking for whoever is the first to plant his flag on the shores of this terra incognita. For the next two months Galileo goes on a frenzied race against the clock. He writes during the day and raids of the heavens for one precious secret after another at night. In early March he has cobbled together enough to claim the main pearls of the heavens for himself. He rushes his little booklet into print with the greatest haste: the last observation entry is dated March 2, and only ten days later the book is coming off the presses\footnote{\cite[77]{Bucciantini}, \cite[157]{drakeGatwork}.}---a turnaround time modern academic publishers can only dream of, even though they do not have to work with hand-set metal type and copper engravings for the illustrations.

``I thank God from the bottom of my heart that he has pleased to make me the sole initial observer of so many astounding things, concealed for all the ages.''\footnote{Galileo to Vinta, 30 January 1610, \cite[153]{heilbron}.} So wrote Galileo, and his palpable relief is fully justified. Little more than dumb luck---or, as he would have it, the grace of God---separated Galileo from numerous other telescopic pioneers who also produced telescopes and made the same discoveries independently of Galileo, such as Simon Marius who discovered the moons of Jupiter one single day later than Galileo. ``A delay of only three or four months would have set him behind several of his rivals and undercut his claim to priority regarding several key discoveries with the telescope.''\footnote{\cite[171]{BaumgartnerFCorr}.} Perhaps it was not the grace of God, but Galileo's desperation, born of decades of impotence as a mathematician, that drove him to publish first. Being incapable of making any contribution to the mathematical side of science and astronomy, he perhaps needed and craved this shortcut to stardom more than anyone else.

Accordingly, Galileo greedily sought to milk every last drop of fame he could from the telescope. ``I do not wish to show the proper method of making them to anyone,'' he admitted; rather ``I hope to win some fame.''\footnote{\cite[61]{GalileoDiscOp}.} His competitors quickly realised that ``we must resign ourselves to obtaining the invention without his help.''\footnote{Peiresc to Pace, 7 November 1610, \cite[156]{Bucciantini}.} Still six years after his booklet of discoveries, people who thought science should be a shared and egalitarian enterprise were rightly upset by Galileo's selfish quest for personal glory: ``How long will you keep us on the tenterhooks? You promised in your {\it Sidereal Message} to let us know how to make a telescope so that we could see all the things that are invisible to the naked eye, and you haven't done it to the present day.''\footnote{Malatesta Porta to Galileo, 13 September 1616, \cite[398]{sheaOxSW}, OGG.XII.281.} Galileo never missed a chance to mock stuffy Aristotelian professors for thinking ``that truth is to be discovered, not in the world or in nature, but by {\em comparing texts} (I use their own words).''\footnote{\cite[419]{sheaOxSW}. OGG.X.422.} But if he genuinely wanted them to turn to nature he could have shared his techniques for telescope construction. In truth it served his own interests very well that these people were left with no choice but ``comparing texts'' while he claimed the novelties of the heavens for himself.

\subsection{Optics}
\label{optics}

Some believed ``the telescope carries spectres to the eyes and deludes the mind with various images $\ldots$\ bewitched and deformed.''\footnote{\cite[80]{DrakeOMalley}.} Perhaps these peculiar ``Dutch glasses'' were but a cousin of the gypsy soothsayer's crystal ball? The “transmigration into heaven, even whil’st we remain here upon earth in the flesh,”\footnote{Robert Hooke, {\it Micrographia} (1665), 234.} may indeed seem like so much black magic. Add to this the very numerous imperfections of early telescopes, which often made it very difficult even for sympathetic friends to confirm observations, not to mention gave ample ammunition to outright sceptics.

Indeed, we find Galileo on the defensive right from the outset, just a few pages into his first booklet. Seen through the telescope, the moon appeared to have enormous mountains and craters but its boundary was still perfectly smooth. ``I am told that many have serious reservations on this point'': for if the surface of the moon is ``full of chasms, that is countless bumps and depressions,'' then ``why is the whole periphery of the full Moon not seen to be uneven, rough and sinuous?''\footnote{Galileo, {\it Sidereus Nuncius} (1610), \cite[64]{sheaSN}.} Galileo replies that this is because the Moon has an atmosphere, which ``stop[s] our sight from penetrating to the actual body of the Moon'' at the edge only, since there ``our visual rays cut it obliquely.''\footnote{Galileo, {\it Sidereus Nuncius} (1610), \cite[65]{sheaSN}.} So when we look at the edge of the moon our line of sight spends more time passing through the atmosphere of the moon and that's why it's blurred. It is therefore ``obvious,'' says Galileo, that ``not only the Earth but the Moon also is surrounded by a vaporous sphere.''\footnote{Galileo, {\it Sidereus Nuncius} (1610), \cite[93]{sheaSN}.} This is of course completely wrong.

Another puzzling fact was that the planets were magnified by the telescope, but not the stars: they remained the same point-sized light spots no matter what the strength the telescope. Some even mistook this for a ``law that the enlargement appears less and less the farther away [the observed objects] are removed from the eye.''\footnote{Grassi, \cite[268]{drakeGatwork}.} Galileo tried to explain the matter, but once again he gets it completely wrong.\footnote{\cite[110]{sheaSN}. A correct explanation was given in 1665 (ibid.).}

Clearly, in light of all these challenges to the reliability and consistency of the telescope, it was important to understand its basis in theoretical optics. That is why, presumably, Galileo felt obliged to swear at the outset that ``on some other occasion we shall explain the entire theory of this instrument.''\footnote{Galileo, {\it Sidereus Nuncius} (1610), \cite[31]{GalileoDiscOp}.} To those aware of his mathematical shortcomings, it will come as no surprise that Galileo never delivered on this promise. Kepler---a competent mathematician---took up the task instead, and in the process came up with a fundamentally new telescope design better than that of Galileo.\footnote{Kepler, {\it Dioptrice} (1611). Kepler's telescope uses two convex lenses instead of Galileo's pair of one convex and one concave lens.} ``Galileo kept silent [on Kepler's work on optics]. Apparently, this self-styled mathematical philosopher was not interested in the mathematical properties of the instrument that had brought him fame.''\footnote{\cite[35]{FJDLenses}.} Galileo's excuse for ignoring Kepler's excellent work of the optics of lenses was that it was allegedly ``so obscure that it would seem that the author did not understand it himself.'' A modern scholar comments that ``this is a curious statement since the {\it Dioptrice}, unlike other works by Kepler, is remarkably straightforward.''\footnote{\cite[11]{sheaSN}.} Still not straightforward enough for Galileo, evidently. Indeed, Galileo's naive conception of optics was still mired in the old notion that seeing involved rays of sight spreading outward from the eye rather than conversely.\footnote{Galileo made repeated statements to this effect, collected in \cite[415--416]{sheaOxSW}.}

\subsection{Mountains on the moon}
\label{moonmountains}

``The moon is not robed in a smooth and polished surface but is in fact rough and uneven, covered everywhere, just like the earth's surface, with huge prominences, deep valleys, and chasms.''\footnote{Galileo, {\SN} (1610), \cite[28]{GalileoDiscOp}.} It is all too easy to cast this report by Galileo as a revolutionary discovery. The ``Aristotelian'' worldview rested on a sharp division between the sublunar and heavenly realm. Our pedestrian world is one of constant change---a bustling soup of the four elements (earth, water, air, fire) mixing and matching in fleeting configurations. The heavens, by contrast, were a pristine realm of perfection and immutability, governed by its very own fifth element entirely different from the physical stuff of our everyday world. If we are predisposed to view Galileo as the father of modern science, a pleasing narrative readily suggests itself: With his revolutionary discovery of mountains on the moon, Galileo disproved what ``everybody'' believed.
\quote{Every educated person in the sixteenth century took as well-established fact …\ that the Moon was a very different sort of place from the Earth. …\ The lunar surface, according to the common wisdom, was supposed to be as smooth as the shaven head of a monk.\footnote{\cite[30]{GNewUniverse}.}}
\quote{In those years virtually no one questioned the ontological difference between heaven and earth. …\ The difference between Earth and the heavenly bodies was an absolute truth for astrologers and astronomers, theologians and philosophers of every ilk and school. …\ If the Moon turned out to be covered with mountains, just like Earth, a millenary representation of the sky would be shattered.\footnote{\cite[7]{Bucciantini}.}}
Using data and hard facts to expose its prejudices, this narrative goes, Galileo sent an entire worldview crashing down. Furthermore, by revealing the similarity of heaven and earth, Galileo opened the door to a unification of terrestrial and celestial physics\footnote{Galileo was ``staking a claim to enemy territory. Celestial physics would no longer be an Aristotelian preserve, with astronomers confined to the computations of motions. …\ [Galileo's observations regarding] the Moon was to mark a new chapter in physics and to subvert not only the traditional departmental division of academic labour, but also the world-picture that sustained it.'' \cite[156]{sharratt}.}---in other words, led us to the brink of Newton's insight that a moon and an apple are governed by the very same gravitational force.

The problem with this narrative lies in one word: ``everybody.'' The Aristotelian worldview is not what ``everybody'' believed; it is what one particular sect of philosophers believed. As ever, Galileo's claim to fame rests on conflating the two. If we compare Galileo to this Aristotelian sect---as Galileo wants us to do---then indeed he comes out looking pretty good. Members of this sect did indeed try to deny the mountainous character of the moon in back-pedalling desperation, for instance by arbitrarily postulating that the mountains were not on the surface of the moon at all but rather enclosed in a perfectly round, clear crystal ball.\footnote{\cite[177]{thusspoke}. OGG.XI.118.} If we mistake this for the state of science of the day, then indeed Galileo will appear a great revolutionary hero.

But to anyone outside of that particular sect blinded by dogma, the idea of a mountainous moon had been perfectly natural for thousands of years. It is obvious to anyone who has ever looked at the moon that its surface is far from uniform. Clearly it has dark spots and light spots. If one wanted to maintain the Aristotelian theory one could try to argue, as many people did, that this is some kind of reflection or marbling effect in a still perfectly spherical moon.\footnote{\cite[115]{Bucciantini}. \cite[176]{thusspoke}. OGG.XI.92.} Whatever one thinks of the plausibility of such arguments, they are certainly defensive in nature: the Aristotelian theory is on the back foot trying to explain away even the most rudimentary phenomena that any child is familiar with. The idea of an irregular moon is an obvious and natural alternative explanation. Which is why we find for instance in Plutarch, a millennium and a half before Galileo, the suggestion that ``the Moon is very uneven and rugged.''\footnote{\cite[4]{sheaSN}.}

If we look to actual scientists and mathematically competent people instead of dogmatic Aristotelians, we find that ``Galileo's'' discovery of mountains on the moon was already accepted as fact long before. Kepler had already ``deduced that the body of the moon is dense …\ and with a rough surface,'' or in other words ``the kind of body that the earth is, uneven and mountainous.''\footnote{Kepler, {\it Astronomiae Pars Optica} (1604), 228, 251, \cite[243, 262]{KeplerO}.} This was also the opinion of his teacher Maestlin before him, who, according to Kepler, ``proves by many inferences …\ that [the moon] also got many of the features of the terrestrial globe, such as continents, seas, mountains, and air, or what somehow corresponds to them.''\footnote{Kepler, {\it Mysterium Cosmographicum} (1596), \cite[165]{KeplerMysterium}.} In a later edition Kepler added the note that ``Galileo has at last throughly confirmed this belief with the Belgian telescope,'' thereby vindicating ``the consensus of many philosophers on this point throughout the ages, who have dared to be wise above the common herd.''\footnote{\cite[169]{KeplerMysterium}.} Indeed, Galileo himself says his observations are reason to ``revive the old Pythagorean opinion that the moon is like another earth.''\footnote{Galileo, {\SN} (1610), \cite[342]{GalileoDiscOp}.} So Galileo's discovery of mountains on the moon was not a revolutionary refutation of what ``everybody'' thought they knew, but rather a vindication of what many thinking people had seen for thousands of years.

Another example along the same lines is ``Earth shine'' (like moon shine, but in the reverse direction): Galileo discusses it in the {\SN} as one of the novelties made clear by the telescope, but in reality it had been correctly explained previously.\footnote{\cite[107]{sheaSN}, Kepler, {\it Astronomiae Pars Optica} (1604), 6.10.}

A similar reality check is in order regarding the ludicrous idea that Galileo's discoveries regarding the moon instigated celestial physics. Kepler had already written an excellent book on celestial physics before the telescope,\footnote{Kepler, {\it Astronomia Nova} (1609).} while Galileo's bumbling attempts on this subject are very poor.\footnote{\S\S\ref{planetaryspeeds}, \ref{moonfall}, \ref{joshuaargument}.}

\subsection{Double-star parallax}
\label{doublestar}

\begin{figure}[pt]\centering
\includegraphics[width=0.6\textwidth]{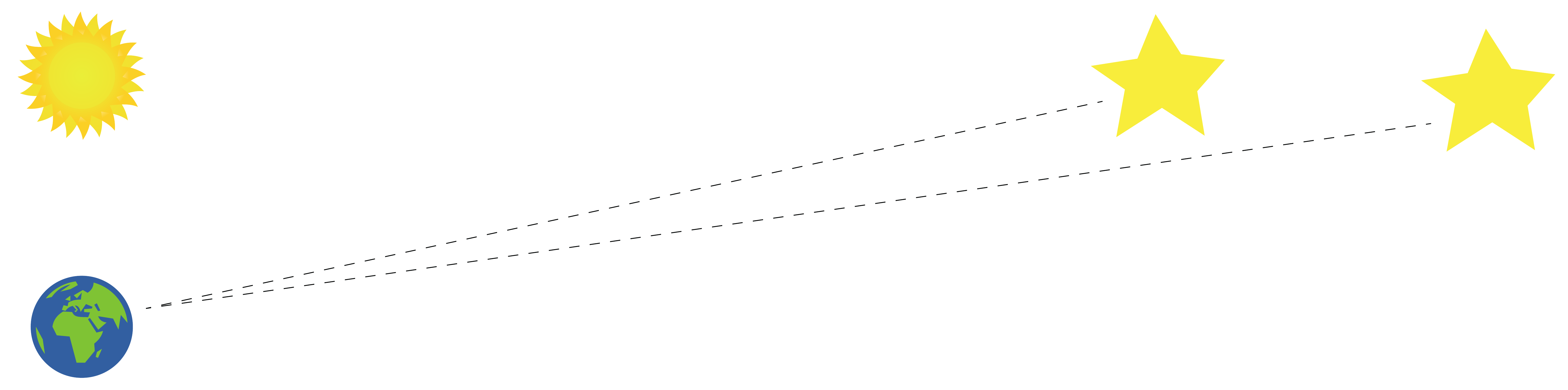}

\vspace{1cm}

\includegraphics[width=0.6\textwidth]{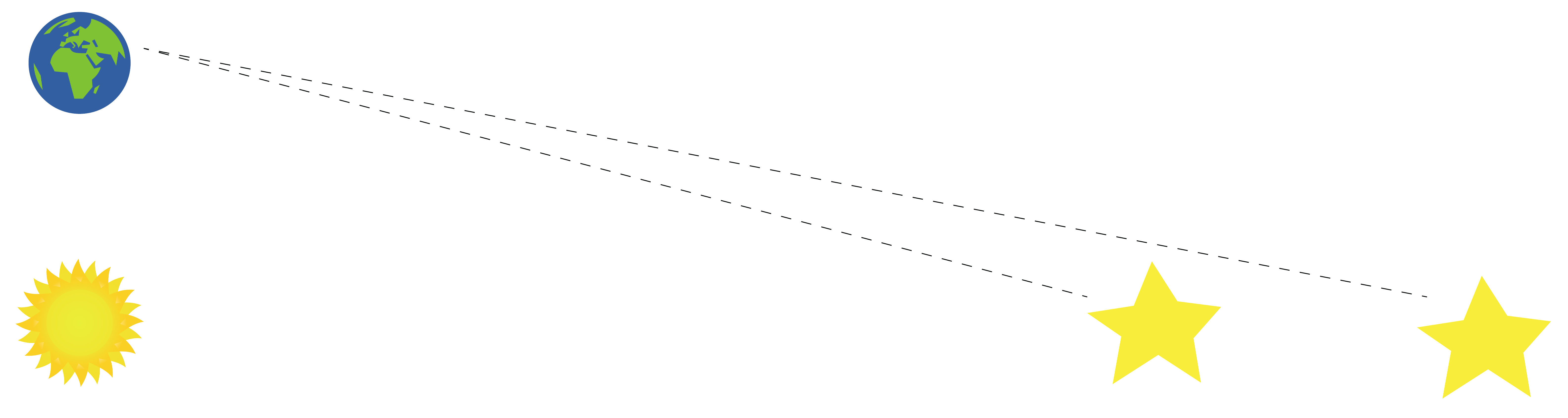}
\caption[Double-star parallax.]{Double-star parallax. Top: To the observer at the earth, the closer star appears bigger and to the left of the further star, which is seen as smaller. Bottom: Half a year later, with the earth in the opposite point in its orbit around the sun, the big star is seen on the other side of the small one.}
\label{dstarparafig}
\end{figure}

The telescope revealed the existence of ``double stars.'' Some stars that had appeared as just a single point of light to the naked eye turned out, when studied with sufficient magnification, to consist of two separate stars. Galileo's former student Castelli was excited about such stars, because he hoped they could be used to prove that the earth moves around the sun. It seemed reasonable to assume that all stars were pretty much alike---they were all just so many suns. Then the apparent size differences in stars would be due not to actual size differences, but only to differences in how far away they are from the observer. A double star would therefore correspond to a geometrical configuration like that in Figure \ref{dstarparafig}. If the earth moves around the sun, it should be possible to see the two stars interchange position in the course of a year. This would certainly not happen if the earth was stationary, so we have striking and undeniable evidence for the motion of the earth.

This is a parallax effect, and as we have seen above astronomers had failed to detect parallax in the past.\footnote{\S\ref{parallax}.} But the traditional method to look for parallax was based on trying to detect subtle shifts in the relative position of stars using tricky precision measurements of angles. The double star case would prove the matter in a much more immediate way, without the need for technical measurements: anyone would be able to see with their own eyes the undeniable fact the the two stars ``switched places'' in the course of a year. And since with this method everything takes place within the field of view of the telescope, there was reason to hope that this new technology could enable success where conventional naked-eye astronomy had failed.

Castelli urged Galileo to make observations of double stars for this purpose, as indeed Galileo did in 1617, when he made detailed observations of the double star Mizar. By assuming that Mizar A and B were in reality as bright as our sun, Galileo estimated from their apparent brightness that Mizar A and B were 300 and 450 times further away than the sun, respectively. This means the above effect should easily be noticeable: ``Mizar A and Mizar B should have swung around each other dramatically as Galileo observed them over time.''\footnote{\cite[263]{GraneyMizar}. See also \cite[260]{Siebert}.} But that didn't happen. They didn't change position at all. Everything remained exactly stationary, as if the earth did not move. Since this failed to confirm his preferred conclusion, Galileo did not publish these results.

Today we know that all the stars in the night sky are much further away than Galileo estimated, and much too far away for any effects of the above sort to be detectable with the telescopes of his time. Galileo's distance estimates were way off because of certain optical effects that make it impossible to judge the distances of stars in the manner outlined above.\footnote{\cite{GraneyGrayson}.} It would be anachronistic to blame Galileo for not knowing these things, which were only understood much later.

But Galileo's way of discussing the matter in the {\Dialogue} is not above reproach. He describes the above procedure but frames it hypothetically: ``if some tiny star were found by the telescope quite close to some of the larger ones,'' they would, if the above effect could be observed, ``appear in court to give witness to such motion in favour of the earth.''\footnote{\cite[382--383]{galileo2sys}. The method of distance estimates via apparent sizes is mentioned at \cite[359--360]{galileo2sys}.} ``This is the very idea that later won Galileo renown and for which he was to be remembered by parallax hunters in the centuries that followed. While it is generally thought that Galileo never tried to detect stellar parallax himself, he is credited with this legacy to future generations.''\footnote{\cite[254]{Siebert}.} In reality he deserves no renown, because the idea was not his own. It had already been explained to him in detail by Ramponi in 1611,\footnote{\cite[254]{Siebert}.} and again by Castelli, who discovered the double star Mizar and explained its importance for parallax to Galileo.\footnote{\cite[259]{Siebert}.} There is no indication that Galileo knew about these things before his friends explained them to him.

Furthermore, Galileo's discussion in the {\Dialogue} is deceitfully framed as a test he has not himself attempted. That way he can pretend that the absence of parallax is not really a problem, by implying that if one just carried out this test one would surely find it. If Galileo was an honest scientist concerned with objectively evaluating the evidence he should have admitted that he had in fact carried out the test, which had failed completely to yield the desired results. But that would have forced him to engage seriously with actual current astronomy like the system of Tycho, which he did not want to do. It was much easier for him to suppress his data and disingenuously insinuate that the outcome of the observation would be the opposite of what he knew it to be.

\subsection{Moons of Jupiter}
\label{jupiter}

The moons of Jupiter were probably the most surprising new discovery made when telescopes were first pointed at the sky. An anecdote related by Kepler conveys some of the excitement: ``My friend the Baron Wakher von Wachenfels drove up to my door and started shouting excitedly from his carriage: `Is it true? Is it really true that he has found stars moving around stars?' I told him that it was indeed so, and only then did he enter the house.''\footnote{Kepler to Galileo, 1610, \cite[10]{SantillanaCrime}.}

It seems Galileo was indeed the first to observe the moons of Jupiter, but only by the smallest possible margin. Simon Marius observed them the very next day.\footnote{\cite[Chapter 5]{GaabMarius}, \cite{Pasachoff}.} In any case one hardly qualifies as the ``Father of Modern Science'' just by looking. Nor does Galileo's account stand out for its scientific excellence. For instance, he tries to ``correct'' Marius on the issue of whether the orbits of the moons are tilted or not with respect to Jupiter's orbital plane. Galileo says they are not: ``it is not true that the four orbits of the satellites [of Jupiter] incline from the plane of the ecliptic; rather, they are always parallel to it.''\footnote{Galileo, {\it Assayer} (1623), \cite[166]{DrakeOMalley}.} In reality, Marius was right and Galileo wrong.

Galileo's mathematical ineptitude is on display in this case as well. “Galileo’s first calculations [of the orbital periods of Jupiter’s moons] were geocentric, not heliocentric. …\ Galileo was treating Jupiter as if it revolved around the Earth, not the Sun. How he ever came to make such an error …\ is an interesting question.”\footnote{\cite[421]{drakeessays1}. See also \cite[35]{sheaSN}. Galileo eventually realised his error himself when his calculations didn't match observations.} Kepler and Marius, meanwhile, understood the matter perfectly and realised at once that this was another good argument against the Ptolemaic system.\footnote{\cite[422]{drakeessays1}.} One Galileo supporter offers a very charitable interpretation: “this throws in doubt the view that by 1611 Galileo was already a Copernican zealot anxious to find every possible argument for the Earth’s motion.”\footnote{\cite[429]{drakeessays1}.} A more plausible explanation, in my opinion, is that Galileo was simply not very competent as a mathematical astronomer. It was not lack of desire to prove the earth's motion that made Galileo miss the point, it was lack of ability.

\subsection{Rings of Saturn}
\label{saturn}

Saturn is ``made of three stars,'' says Galileo.\footnote{OGG.V.237. \cite[218]{PretendSat}.} He is trying to describe the rings of Saturn, but telescopes at the time were not good enough to show that they were rings rather than two companions stars. This only became clear some twenty years after Galileo's death.\footnote{Christiaan Huygens, {\it Systema Saturnium} (1659).} We cannot blame Galileo for this. We can, however, blame Galileo for his lack of balance in evaluating the evidence. He does not say, as an honest scientist might, that this is the best guess on the available evidence and that we can't know for sure until we have better telescopes. Instead he boldly proclaims it as certainty that Saturn is ``accompanied by two stars on its sides,'' ``as perfect instruments reveal to perfect eyes.''\footnote{Galileo, first letter on sunspots (1612), \cite[102]{onsunspots}.}

\begin{figure}[tp]\centering
\includegraphics[height=2cm]{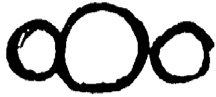} \includegraphics[height=2cm]{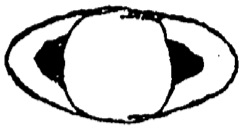} \includegraphics[height=2cm]{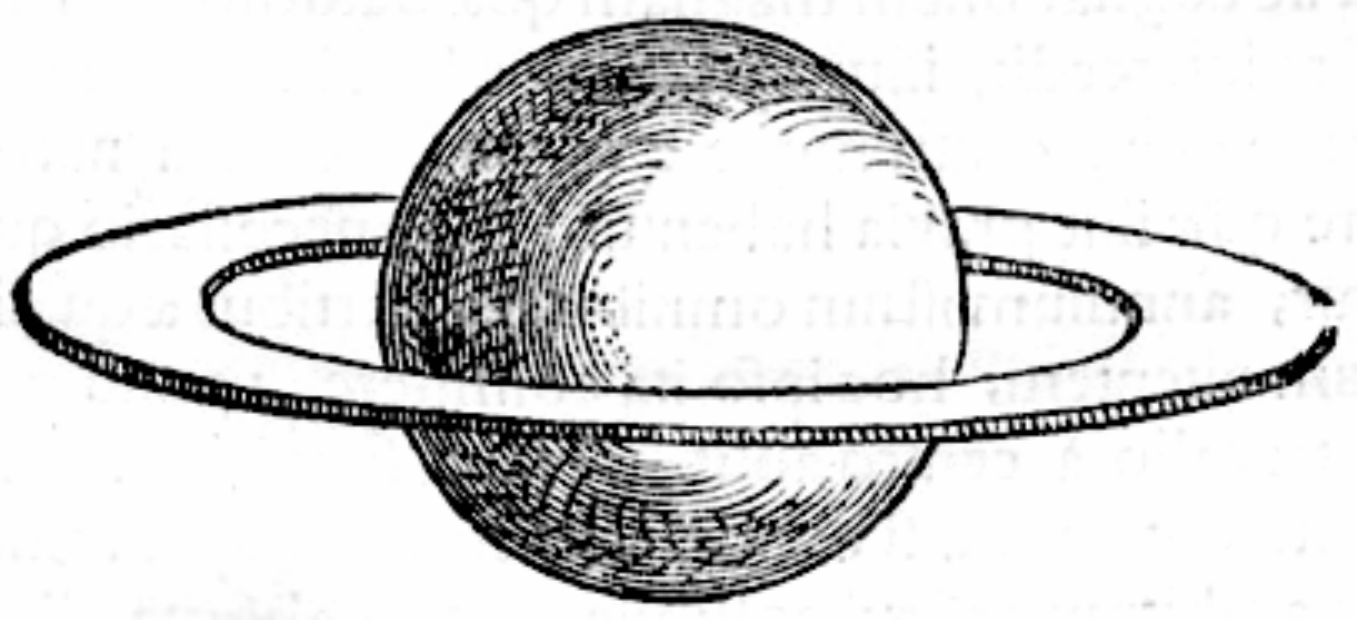}
\caption[Saturn's rings.]{Drawings of Saturn. Left: Galileo, 1610, OGG.X.410. Middle: Galileo, 1616, OGG.XII.276. Right: Huygens, {\it Systema Saturnium} (1659), 47.}
\label{saturnfig}
\end{figure}

In the same vein, Galileo hubristically declared that the appearance of Saturn's ``ears'' would never change:
\quote{I, who have examined [Saturn] a thousand times at different times, with an excellent instrument, can assure you that no change at all is perceived in him: and the same reason, based on the experience which we have of all the other movements of the stars, can render us certain that, likewise, there will be none.\footnote{Galileo, May 1612, OGG.V.110--111, \cite[106]{HeldenSatAnses}.}}
All the more embarrassing then when in fact the appearances changed radically soon thereafter:
\quote{I found him solitary without the assistance of the supporting stars. …\ Now what is to be said about such a strange metamorphosis? Perhaps the two smaller stars …\ have vanished and fled suddenly? Perhaps Saturn has devoured his own children?\footnote{Galileo, 1612, OGG.V.237, \cite[107]{HeldenSatAnses}. ``In classical mythology, Saturn devoured his newborn children to forestall a prophecy that he would be overthrown by one of his sons.'' \cite[403]{sheaOxSW}.}}
So the very thing Galileo said ``thousands'' of observations meant we could be ``certain'' would never happen actually took place almost right away. That was bad publicity at a time when many doubted the reliability of his telescope.

The ``disappearance'' of Saturn's ring was due to the earth passing through the plane of the ring, so that a line of sight from earth was parallel to the plane of the ring. This made the ring invisible, just like a sheet of paper becomes vanishingly thin if looked at exactly sideways. But Galileo did not interpret it that way. Instead, he proposed what he considered to be some ``probable conjectures'' about the future appearance of Saturn's companion stars, apparently based on a model attributing to them a slow revolution. Later he praised himself for ``thinking in my own special way'' and marvelled at how ``I took the courage'' to make such brave conjectures.\footnote{OGG.XVIII.238--239. \cite[111]{HeldenSatAnses}.} Indeed, Galileo liked his model so much that he also ``took the courage'' to lie about having made an observation verifying it. He claims that he ``saw Saturn triple-bodied this year [1612], at about the time of the summer solstice.''\footnote{OGG.V.237. \cite[215]{PretendSat}.} But modern calculations show that the ring of Saturn would have been vanishingly thin at this time. ``Clearly …\ [Galileo] could not have observed the ring at the summer solstice of 1612. …\ Yet the picture of the Saturnian system that was accepted by Galileo implied that the ring should have been visible, so much so that he made a claim to this effect that we know must have been untrue.''\footnote{\cite[218]{PretendSat}.}

Some years later Galileo found that the satellites had undergone another ``strange metamorphosis.'' This time they had ceased to be independent globes at all and instead appeared as half-ellipses attached to the planet. Galileo made a prescient-looking sketch (Figure \ref{saturnfig}), which has led some people to triumph that he discovered the rings of Saturn.\footnote{\cite{BianchiSat}.} But, alas, for all his supposedly admirable courage and precious ``special way of thinking'' he failed to interpret them correctly. Instead he reverted back to the erroneous idea of the satellites as ``collateral globes'' in the {\Dialogue}.\footnote{OGG.VII.287. \cite[110]{HeldenSatAnses}.}

\subsection{Sunspots}
\label{sunspots}

When astronomers first turned telescopes to the heavens, one of the most disturbing sights was ``filth on the cheeks of the Sun.''\footnote{Alessandro Allegri, {\it Lettere di Ser Poi Pedante} (Bologna: Vittorio Benacci, 1613), 14, \cite[77]{onsunspots}.} ``I saw the sonne in this manner:''\footnote{\cite[26]{onsunspots}. This is the earliest recorded telescopic observation of sunspots.}
\begin{center}\includegraphics[height=3cm]{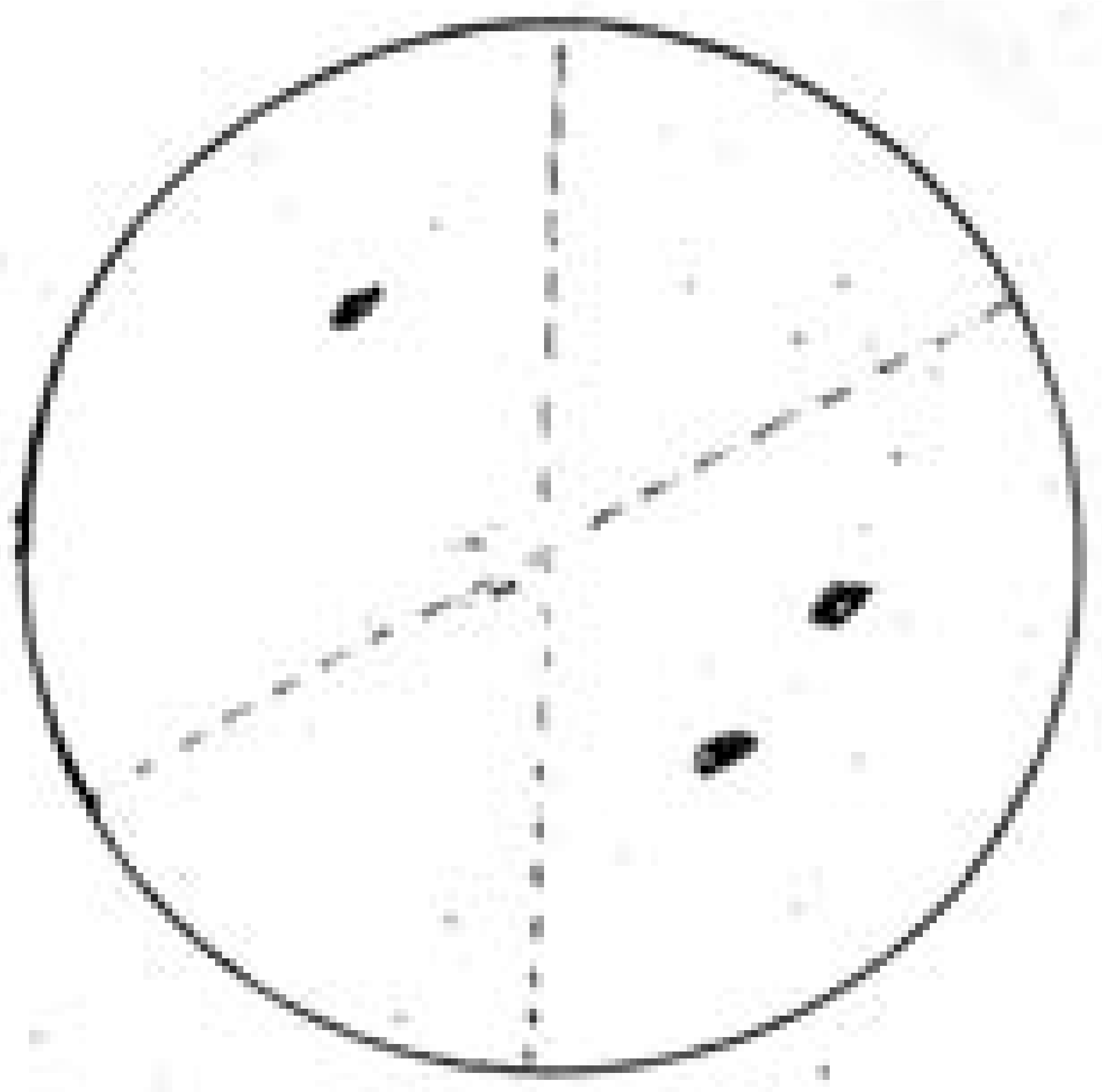}\end{center}
So records Thomas Harriot in his observational diary. Such spots on the sun soon drew much attention across Europe. While some ``neglected to observe them, being afraid, $\ldots$ that the image might burn my eye,''\footnote{Maelcote, \cite[47]{onsunspots}.} others figured God had given them two eyes for a reason and ``burned'' them alternately in the name of science. Thus Harriot ``saw it twise or thrise, once with the right ey \& other time with the left'' before ``the Sonne was to cleare.''\footnote{\cite[26]{onsunspots}.} Soon a method was developed for projecting the image of the sun onto a piece of paper so that no burning of the eyes was needed.\footnote{Galileo learned this method from Castelli. \cite[402]{sheaOxSW}.} This was convenient enough, though there would have been no shortage of martyrs of science willing to pay with their eyes for wisdom.

``Galileo insisted to his dying day that he was the first to have seen'' sunspots but in reality he was ``probably preceded by the Dutch astronomer Johann Fabricius, who was the first to publish information about them.''\footnote{\cite[402]{sheaOxSW}.} To boost his priority case, Galileo later claimed he had seen sunspots already in 1610,\footnote{\cite[401]{galileodialogueML}.} rather than in 1611 as documented, but this seems to be a lie since ``if he had first observed them …\ in 1610 it is extremely probable that he would have mentioned that fact in his {\it Letters} printed in 1613, when the priority issue was first hotly debated.''\footnote{\cite[333]{drakeGatwork}.} At that time of the {\it Letters}, Galileo was keen to establish his priority over Scheiner, who published on sunspots in 1612. At that time Galileo was evidently unaware of the more obscure publication by Fabricius, which had appeared in 1611.\footnote{{\it De maculis in sole observatis narratio}, Wittenberg, 1611.} Hence it later became important for Galileo to push his discovery back even earlier, whereupon he conveniently asserts that he saw them in 1610. For that matter, even pre-telescopic astronomers had noticed the phenomenon of sunspots: it ``was certainly mentioned in the time of Charlemagne, and possibly was referred to by Virgil.''\footnote{Drake, \cite[82]{GalileoDiscOp}.}

In any case, the game was now on to explain the nature of the spots. Scheiner was concerned ``to liberate the Sun's body entirely from the insult of spots,''\footnote{Third letter of {\it Tres Epistolae} (1611), \cite[73]{onsunspots}.} for ``who would dare call the Sun false?''\footnote{Third letter of {\it Tres Epistolae} (1611), \cite[67]{onsunspots}.} He found a way of accomplishing this by arguing that the sunspots were ``many miniature moons,''\footnote{Third letter of {\it Tres Epistolae} (1611), \cite[72]{onsunspots}.} rather than blemishes on the sun itself.

Galileo on the other hand eagerly embraced sunspots as an opportunity to stick it to his Aristotelian enemies, ``for this novelty appears to be the final judgement of their philosophy.''\footnote{Galileo to Barberini, 2 June 1612, OGG 11:311, \cite[83]{onsunspots}.} Thus Galileo placed the spots on the sun itself, arguing that ``clouds about the Sun'' was the most plausible explanation, for one ``would not find anything known to us that resembles them more.''\footnote{First letter on sunspots (1612), \cite[99]{onsunspots}.} It is true that sunspots are on the surface of the sun---a conclusion, incidentally, which Kepler had already reached before reading Galileo.\footnote{\cite[213]{drakeGatwork}.} However, sunspots are not clouds above the surface of the sun, as Galileo believed, but rather dark pits or depressions in the solar surface. ``Scheiner …\ entertained the possibility of this hypothesis, while Galileo resolutely discarded it as unworthy of serious consideration.''\footnote{\cite[67]{sheaGrev}.}

Galileo loved claiming new discoveries as his own and using them as ammunition in his philosophical disputes. But he soon lost interest when it came to the detailed work of actual science. Again and again he makes careless errors and jumps to conclusions with premature confidence, while his competitor Scheiner does the meticulous observational work that Galileo had no patience for. For instance, Galileo erroneously claimed---supposedly based on ``a great number of most diligent observations of this particular''\footnote{Third letter on sunspots (1612), \cite[268]{onsunspots}.}---that all sunspots had the same orbital period, regardless of latitude. In fact, sunspots near the sun's equator orbit quicker than those near the poles by a few days. Galileo was corrected by Scheiner.\footnote{\cite[67]{sheaGrev}, \cite[268]{onsunspots}.} Similarly, ``the sun's disc, as we normally see it in the sky, appears to be uniformly bright, but this impression is dispelled by even the most perfunctory telescopic observation which reveals that the brightness decreases from the centre towards the limb. Scheiner made this discovery but Galileo dismissed it.''\footnote{\cite[66]{sheaGrev}.}

Galileo also did not miss the opportunity to make some mathematical errors as usual. He tried to compute the perspective aspect of the sunspots' motion: how does their apparent speed along the sun's disc vary, given that their actual direction of motion turns more away from us the closer they get to the edge? Galileo's attempted demonstration covers three pages and contains at least as many errors.\footnote{\cite[357--359]{onsunspots}, \cite[57]{sheaGrev}.}

\subsection{Sunspots and heliocentrism}
\label{sunspotshelio}

Sunspots can be used as evidence that the earth moves around the sun. In his {\Dialogue}, Galileo considered this one of his three best arguments in favour of Copernicus.\footnote{\cite[462]{galileo2sys1sted}.} “For the sun has shown itself unwilling to stand alone in evading the confirmation of so important a conclusion [i.e., that the earth orbits the sun], and instead wants to be the greatest witness of all to this.”\footnote{\cite[345]{galileo2sys1sted}.}

\begin{figure}[pt]\centering
\includegraphics[width=0.45\textwidth]{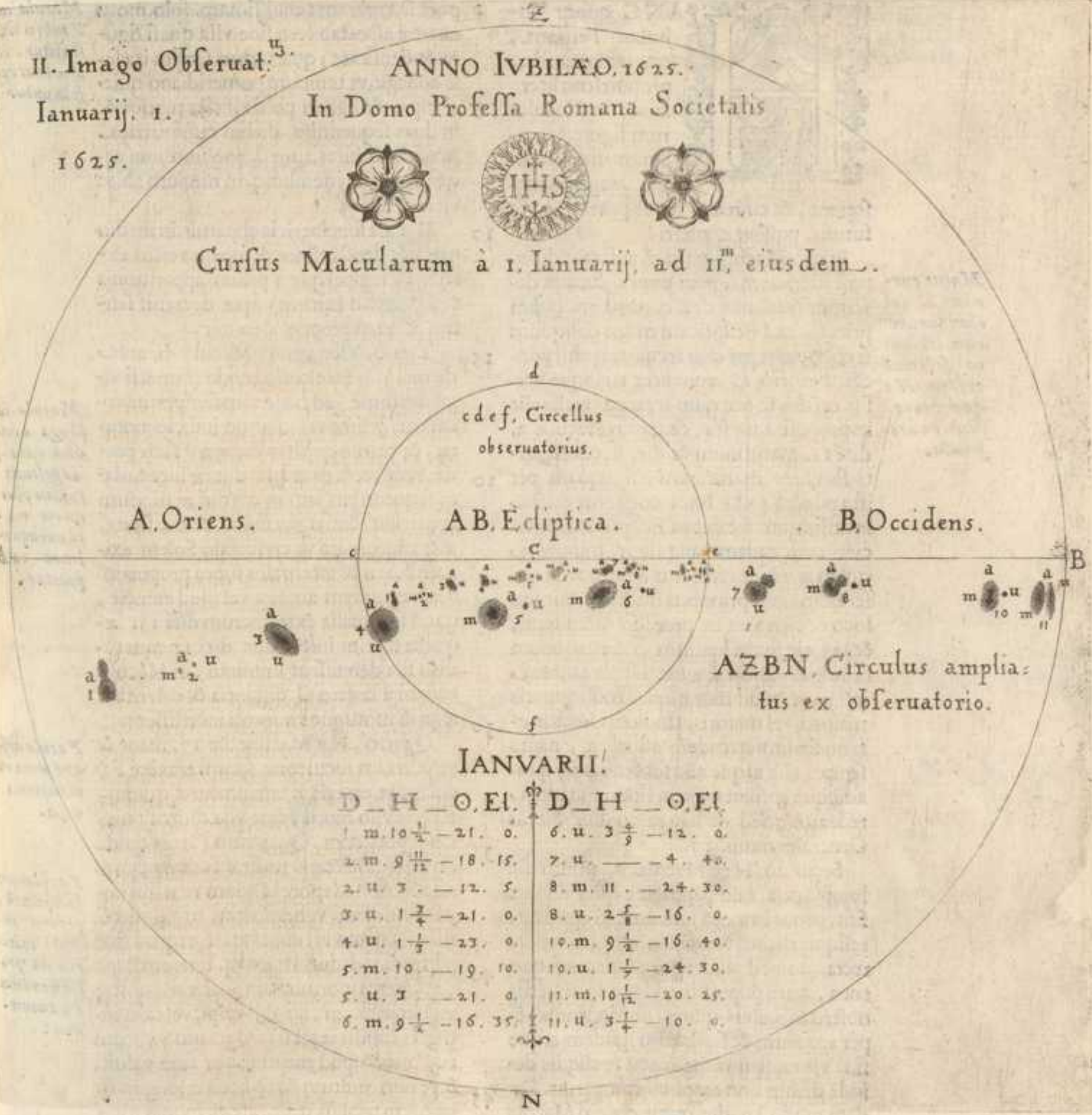} \includegraphics[width=0.45\textwidth]{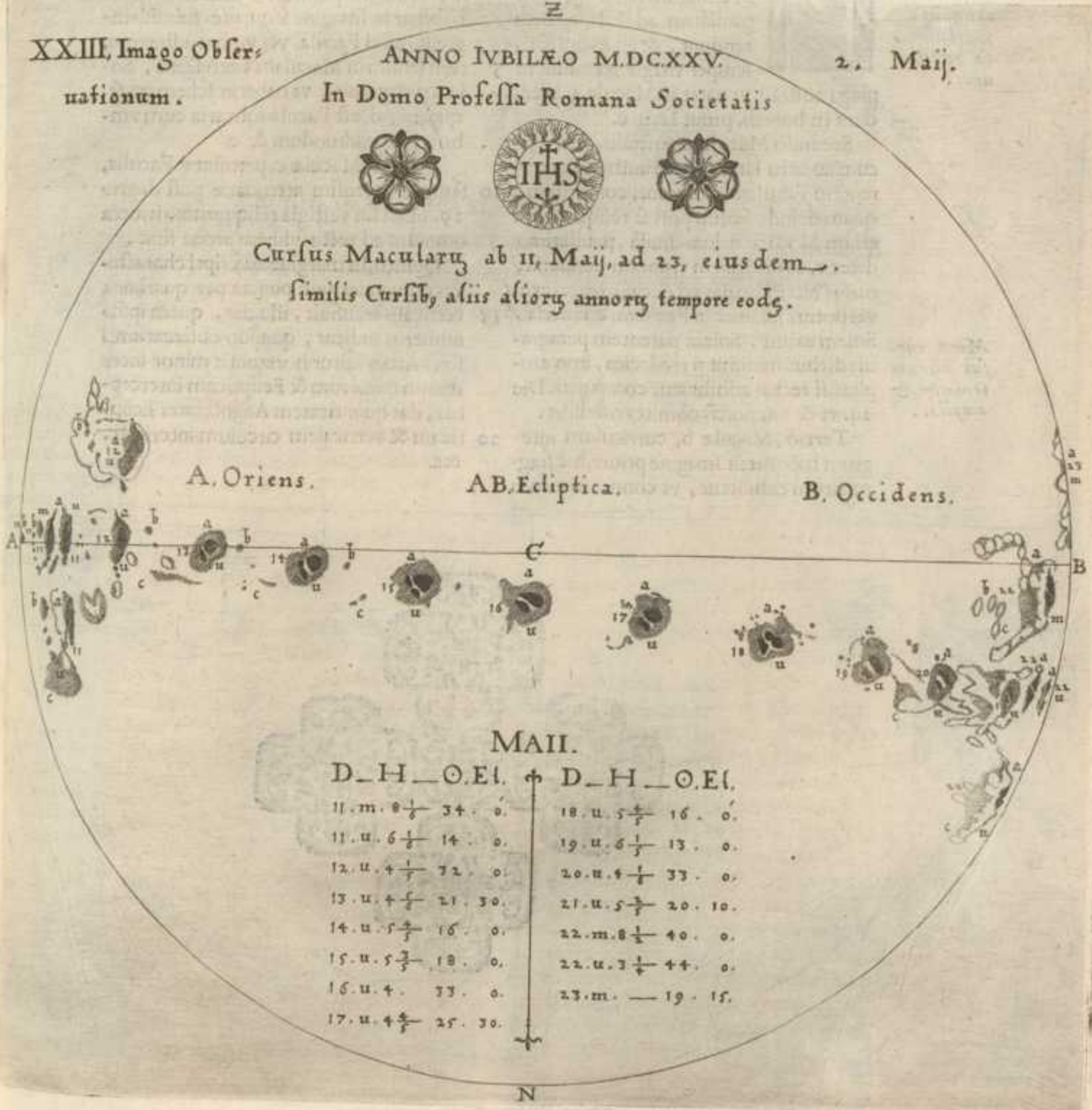}

\includegraphics[width=0.45\textwidth]{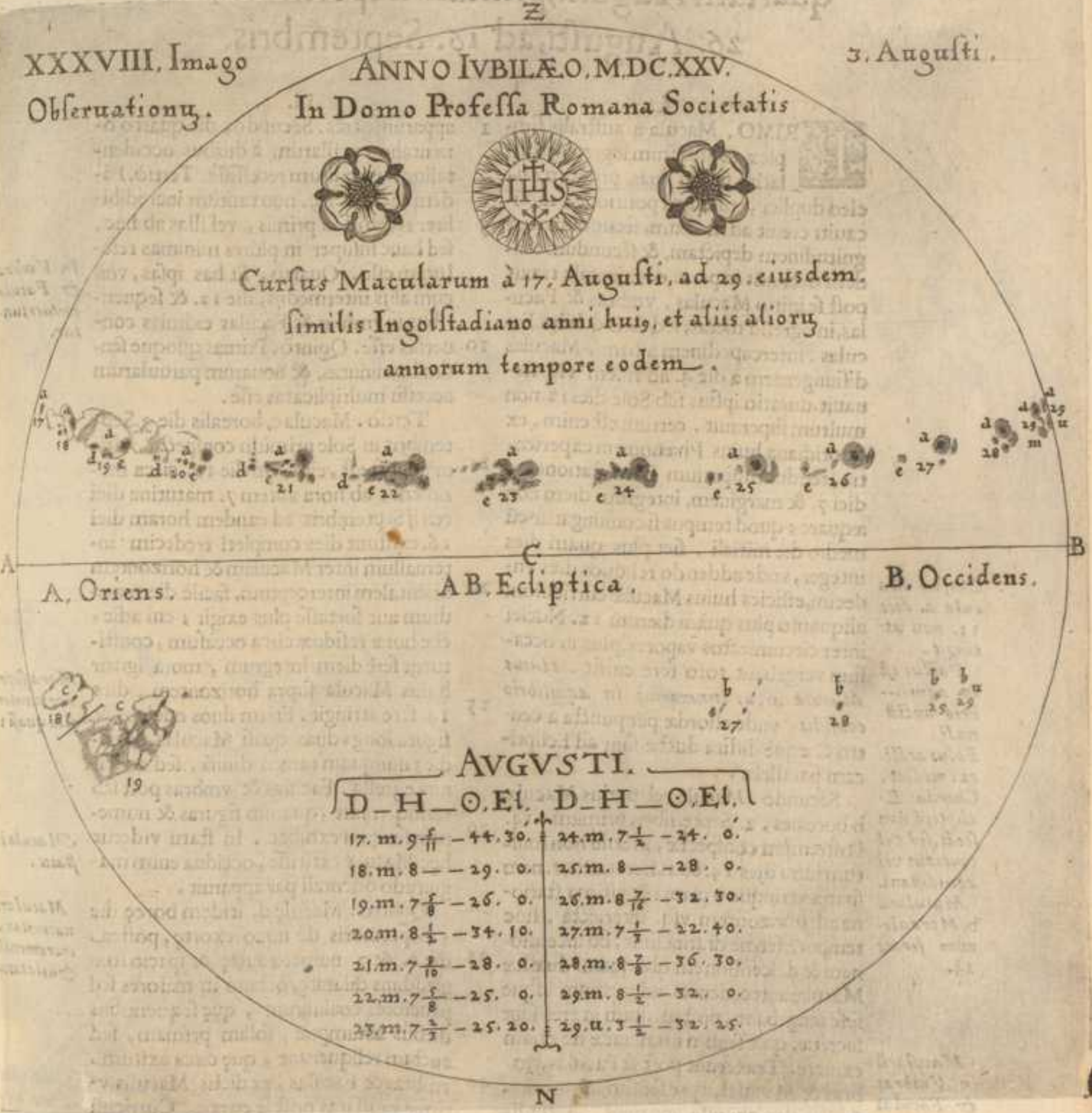}
\includegraphics[width=0.45\textwidth]{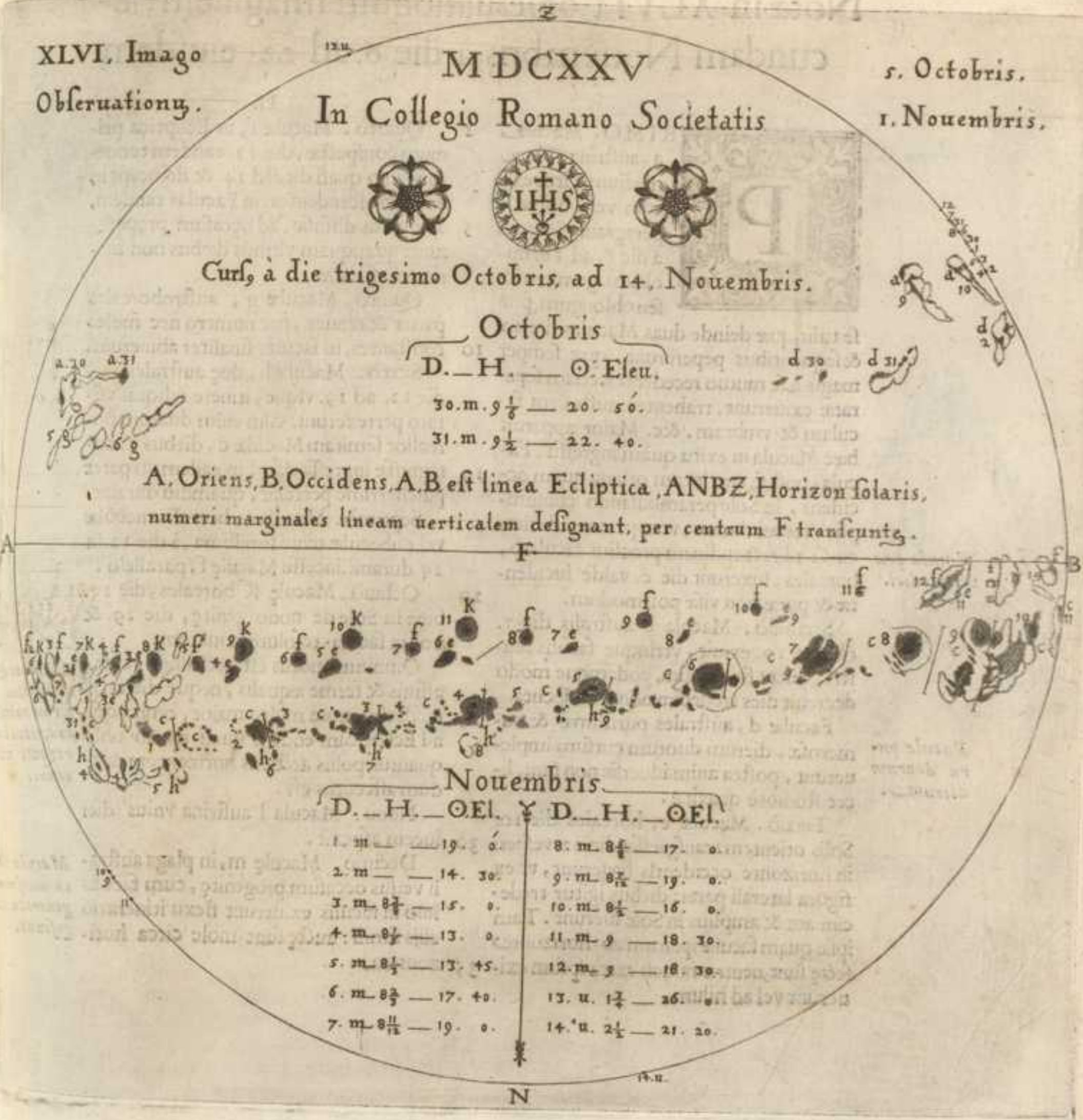}
\caption[Sunspots.]{Paths traced by sunspots over a series of days at different times of year. From Scheiner, {\it Rosa Ursina} (1630).}
\label{sunspotspathfig}
\end{figure}

\begin{figure}[pt]\centering
\includegraphics[width=0.22\textwidth]{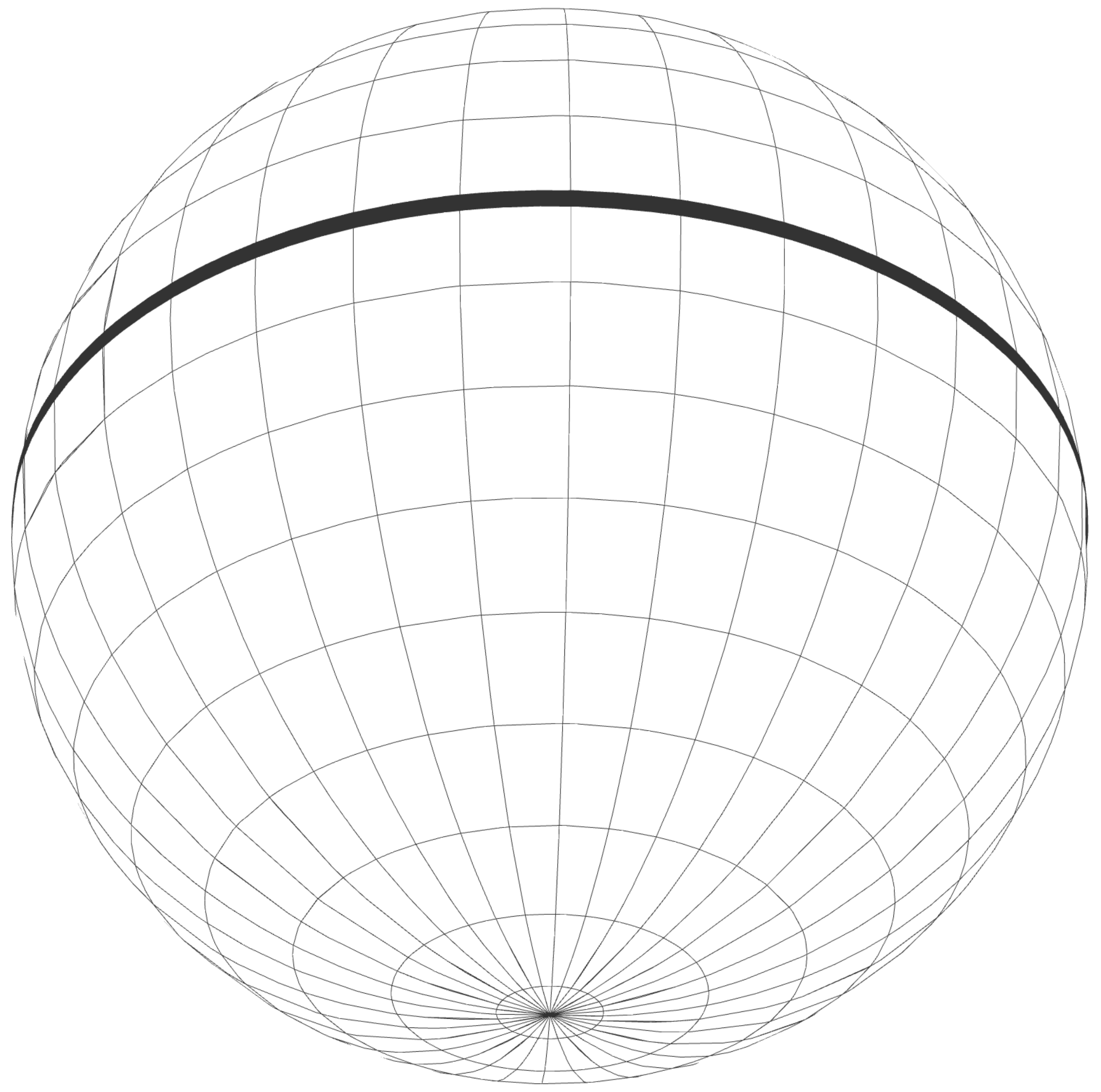} \includegraphics[width=0.22\textwidth]{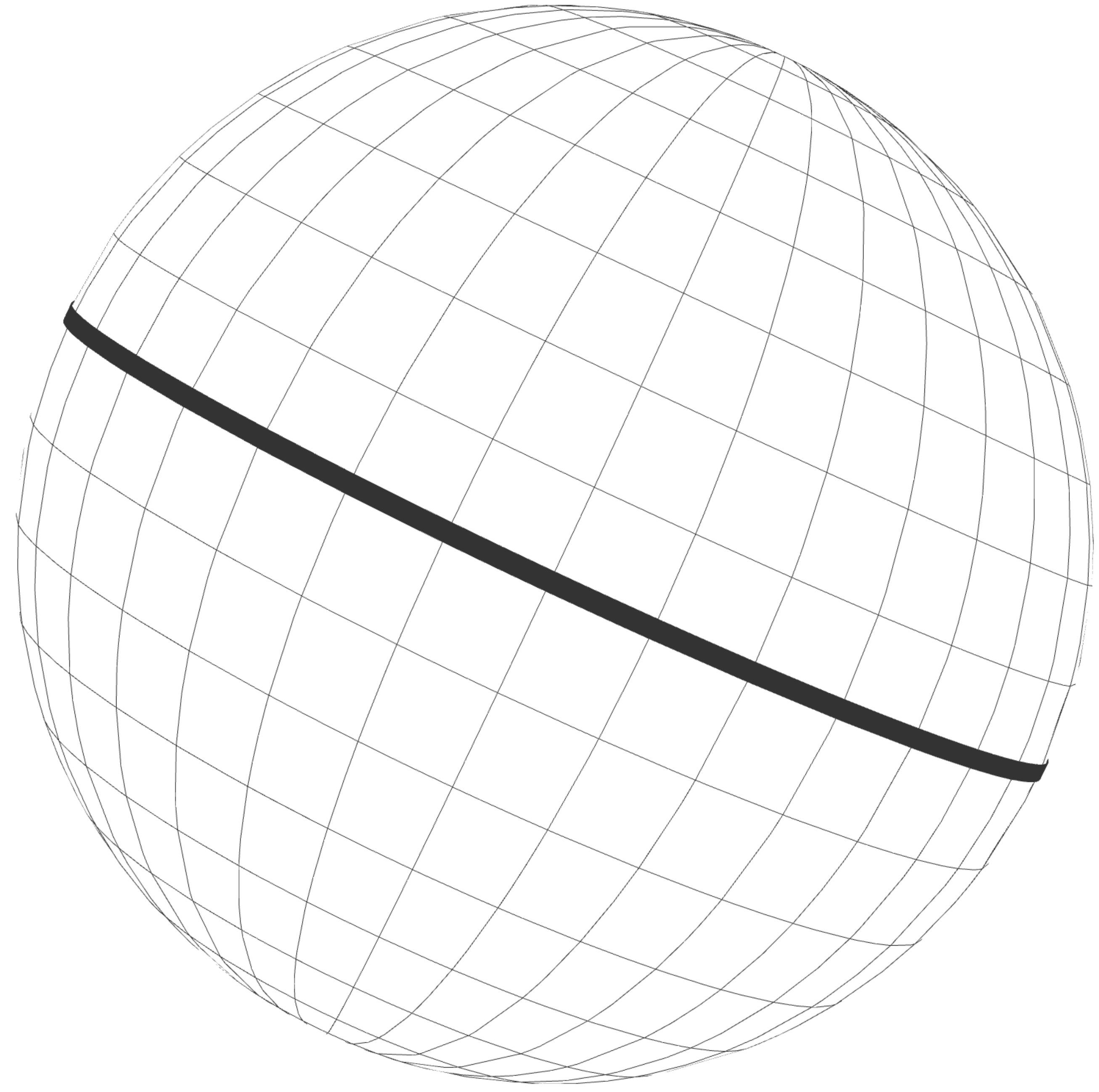} \includegraphics[width=0.22\textwidth]{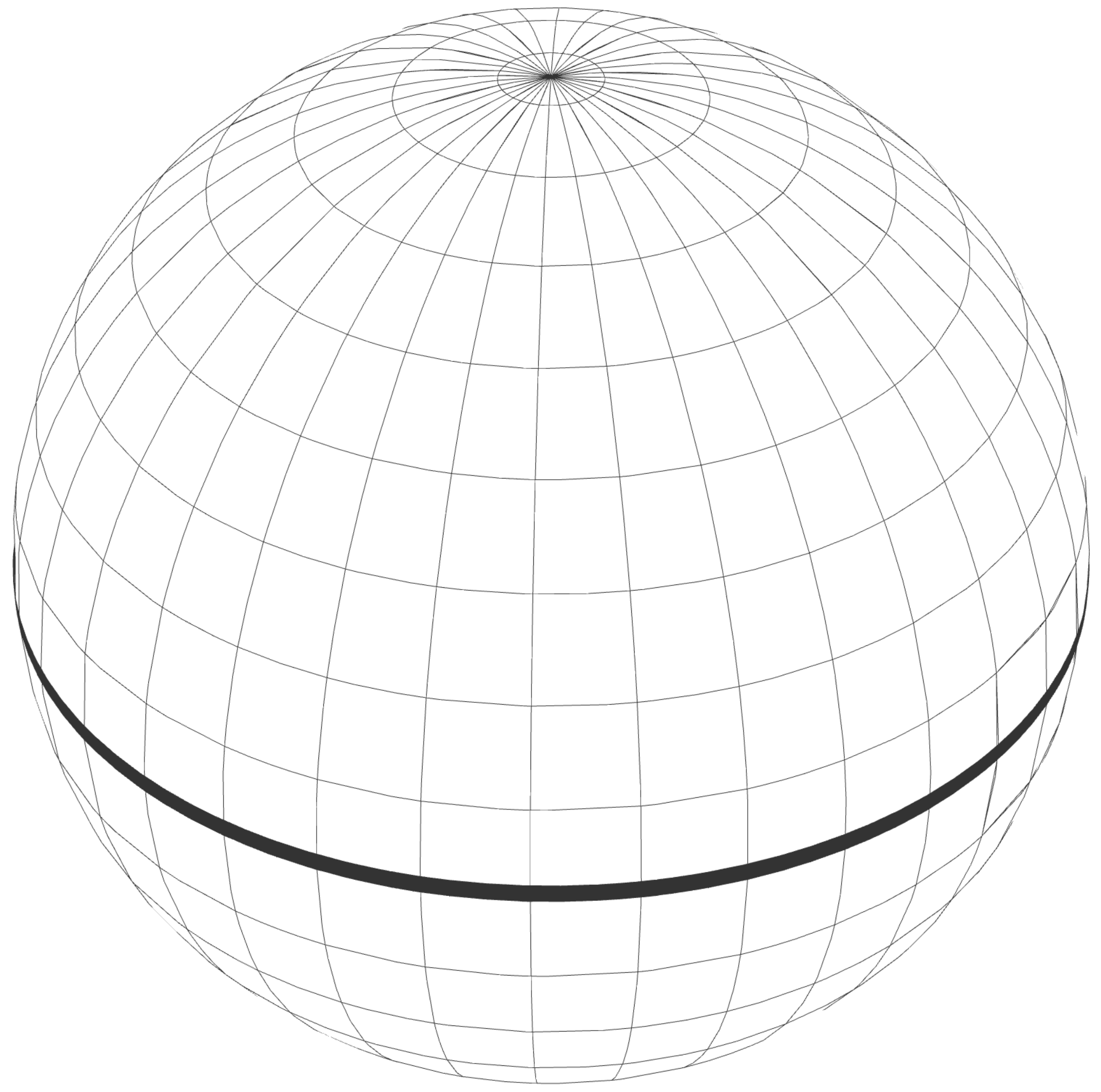} \scalebox{-1}[1]{\includegraphics[width=0.22\textwidth]{pic/sphereviewside}} 
\caption{Equatorial circle of a sphere viewed from different vantage points.}
\label{sphereviewfig}
\end{figure}

The Copernican argument from sunspots goes as follows.\footnote{Cf.\ \cite[348--351]{galileo2sys1sted} and Figures \ref{sunspotspathfig} and \ref{sphereviewfig}.} Imagine a standard globe of the earth standing on a table. Its axis is a bit tilted, of course---the north pole is not pointing straight up. Have a seat at one side of the table and face the globe. What do you see? Focus on the equator. What kind of shape is it? If the north pole is facing in your direction, the equator will make a ``happy mouth'' or U shape. If you move to the opposite side of the table, where you see mostly the southern hemisphere, the equator is instead a sad mouth shape. From the sides, the equator looks like a diagonal line.

Now, suppose the sun had its equator marked on it. And suppose that in the course of a year we would see it as alternately as a happy mouth, straight diagonal, sad mouth, straight diagonal, etc. That would correspond exactly to us moving around the table, looking at a stationary globe from different vantage points. In the same way, if the sun's equator exhibited those appearances, the most natural explanation would be that the earth is moving around it and we are looking at its equator slightly from above, from the side, slightly from below, etc.

The sun does not have the equator conveniently marked on its surface, but not far from it. The sun is spinning rather quickly, making a full turn in less than a month. As it spins, a point on its surface traces out an equatorial or at least latitude circle. So by tracking the paths of sunspots over the course of a few weeks, we in effect see equatorial and other latitude circles being marked on the surface of the sun.

So the shapes of the paths show that we are looking at the sun from alternating vantage points. But this does not necessarily mean that we are moving around the sun. The same phenomena could be accounted for from a geostatic or Ptolemaic point of view by saying that the sun is so to speak wobbling, showing us different sides of itself in the course of a year. You can see this with your globe on the table. Instead of moving around the table and looking at the globe from different sides, you can have a friend tilt the globe, pointing its axis now this way and now that. If you let the axis spin around in a conical motion, this will produce the exact same visual impressions for you as if you had moved around the table.

In order to use the sunspot paths as evidence for Copernicus, then, Galileo needed to dismiss this alternative explanation. He did so by attacking it as physically implausible. To account for the sunspots phenomena from a Ptolemaic point of view, the sun had to orbit the earth, and spin on its own axis, and have its axis wobble in a conical motion. These diverse motions, says Galileo are “so incongruous with each other and yet necessarily all attributable to the single body of the sun.”\footnote{\cite[355]{galileo2sys1sted}.} Surely this is a geometrical fiction that would never happen in an actual physical body.

Actually, such an ``incongruous'' combination of motions is not only possible but a plain fact. The earth, in fact, has exactly such a combination of motions, as had been known since Copernicus.\footnote{\cite{Topper}.} The earth has a conical wobbling motion which means that the north pole is pointing to a slightly different spot among the stars from year to year, returning to its original spot after 26,000 years. This is the explanation for the so-called precession of the equinoxes, an important technical aspect of classical astronomy. So if Galileo's argument about ``incongruous'' motions disproves the Ptolemaic explanation of sunspots, it also disproves Copernicus's correct explanation of the precession of the equinoxes.

Galileo conveniently neglects to bring up this rather obvious problem with his argument. Whether he did so out of ignorance or dishonesty is hard to say, but either way is none too flattering. Any serious mathematical astronomer was well acquainted with the precession of the equinoxes and of course considered it an essential requirement that any serious astronomical system account for this phenomenon. Galileo, though, is not a serious mathematical astronomer. He is a simplistic populariser who simply ignores technical aspects like these. And it is only because of this oversimplification that he is able to maintain his argument against the Ptolemaic interpretation of sunspots.

In any case, in the 1610s, when he was first studying sunspots, Galileo completely missed all of this. He lacked the disposition to do painstaking scientific research like Scheiner. Instead, with premature hubris, he soon imagined that he had “looked into and demonstrated everything that human reason could attain to” regarding sunspots.\footnote{\cite[346]{galileo2sys1sted}, referring to the 1610s.} Many years later he was still convinced that his was the last word on the matter: “writing …\ apropos of recent news that Scheiner would soon publish a thick folio volume on sunspots, he remarked that any such book would surely be filled with irrelevancies, as there was no more to be said on the subject than he had already published in his Letters on Sunspots.”\footnote{\cite[185]{drakeGS}.}

When Scheiner's much better work on sunspots came out, Galileo realised he had to completely reverse his earlier proclamations, made with such arrogant confidence. With unwarranted pomposity, he had claimed to offer ``observations and drawings of the solar spots, ones of absolute precision, in their shapes as well as in their daily changes in position, without a hairsbreadth of error.''\footnote{First letter on sunspots (1612), \cite[104]{onsunspots}. The observations themselves were sent with the second letter (1612).} 
According to Galileo, the sunspots were ``describing lines on the face of the sun'':\footnote{First letter on sunspots (1612), \cite[90]{onsunspots}.} ``they travel across the body of the sun $\ldots$ in parallel lines.''\footnote{Second letter on sunspots (1612), \cite[109]{onsunspots}.} In fact, ``I do not judge that the revolution of the spots is oblique to the plane of the ecliptic, in which the earth lies.''\footnote{Second letter on sunspots (1612), \cite[113]{onsunspots}. The same claim is repeated in the third letter (1612), \cite[255]{onsunspots}.} In other words, every sunspot path is straight as an arrow, just as the equator of a globe would be from every side if its axis was perfectly vertical.

But Scheiner showed that the sun's axis has a $7^{\circ}15'$ inclination and that the paths exhibit exactly the alternating diagonal and U shapes described above. He published this result in his {\it Rosa Ursina} (1630), the folio Galileo had mocked as bound to be superfluous, but from whence he now realised his error.\footnote{\cite[315, 327]{onsunspots}.}

When Galileo finally realised that inclined sunspot paths spoke in favour of heliocentrism, he immediately threw all his old observations “without a hairsbreadth of error” out the window and rushed the pro-Copernican argument into print without making any new observations, as is clear from the fact that the published argument “displays entire ignorance or complete neglect of the observational data,” his vague descriptions being “utterly wrong” and “almost the exact opposite” of the careful data published by Scheiner.\footnote{\cite[186--187]{drakeGS}. See also \cite[335]{drakeGatwork}.}

Galileo did not want to admit his debt to Scheiner, however, so he pretended that he had come upon this discovery independently,\footnote{\cite[346]{galileo2sys1sted}.} and lied that he had made “very careful observations for many, many months, and noting with consummate accuracy the paths of various spots at different times of the year, we found the results to accord exactly with the predictions.”\footnote{\cite[352]{galileo2sys1sted}.} In reality, “the evidence is unequivocal: Galileo …\ must have had a copy of Scheiner's book in front of him as he wrote this section.” By pretending otherwise, “Galileo has deliberately set out to efface Scheiner from the historical record and to deny his debt to him. It is impossible to find any excuse for this behaviour.”\footnote{\cite[208--209]{wootton}.}

\subsection{Phases of Venus}
\label{venus}

\begin{figure}[pt]\centering
\includegraphics[width=0.9\textwidth]{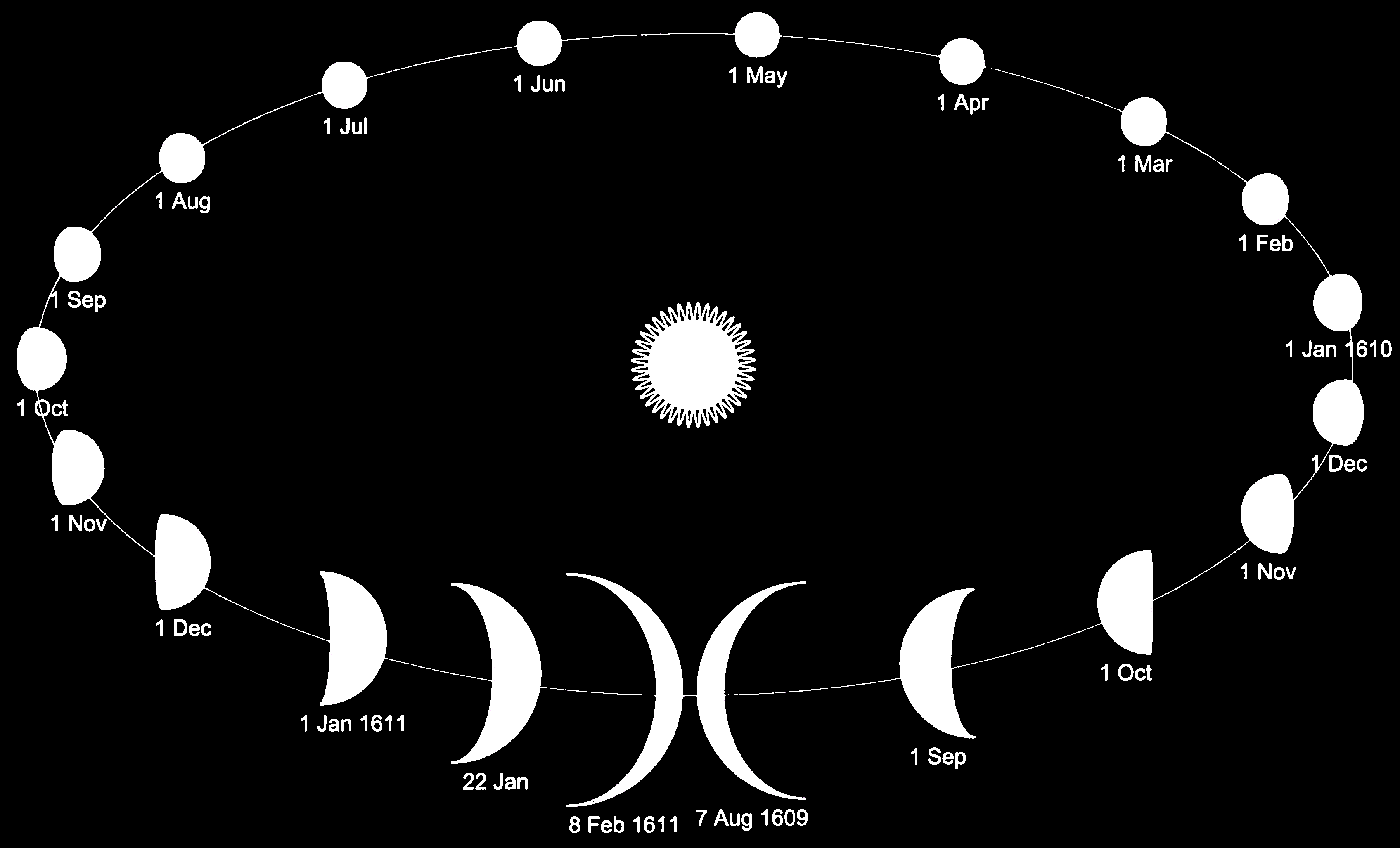}
\caption[Phases of Venus.]{Actual appearance of phases of Venus around the time of Galileo's observations. Based on calculations by Rob van Gent.}
\label{venusphasesillustrationfig}
\end{figure}

To the naked eye, Venus is just a dot of light. But telescopic magnification reveals it as a sphere, only half of which is bright, namely the half facing the sun. Thus the half of Venus facing us exhibits phases like the moon, being sometimes crescent, sometimes half full, sometimes gibbous, and so on (Figure \ref{venusphasesillustrationfig}). These appearances show that Venus orbits the sun, contrary to Ptolemaic cosmology.

Galileo observed the phases of Venus in 1610. Here is a timeline of events:
\begin{itemize}
\item On December 5, Castelli wrote to Galileo and pointed out that it ought to be possible to confirm the Copernican hypothesis by observing phases of Venus with a telescope.\footnote{\cite[11]{westfalltelescope}.}
\item On December 11, perhaps right after receiving Castelli’s letter,\footnote{It is hard to say when exactly Galileo would have received the letter from Castelli. \cite[24]{westfalltelescope} discusses evidence regarding mail delivery times and finds it “easily possible” that Galileo could have received the letter before December 11. \cite[203]{DrakeVenus}, on the other hand, finds the probability of this “vanishingly small” on the basis of other evidence regarding mail delivery times, as well as arguments from other references in their correspondence.} Galileo announced something “just observed by me which involves the outcome of the most important issue in astronomy and, in particular, contains in itself a strong argument for the …\ Copernican system,” namely that he has observed Venus exhibiting phases like the moon.\footnote{\cite[24]{westfalltelescope}. Galileo's announcement was initially in the form of an anagram, the meaning of which he revealed only later.} Before this there is no record that Galileo knew anything about the phases of Venus.\footnote{\cite[25]{westfalltelescope}.}
\item On December 30, Galileo gave for the first time an account of his observations of Venus. At this point he claims to have observed Venus for about three months, and gives an accurate and fairly detailed description of its appearance during this period.\footnote{\cite[109--110]{PalmieriVenus}.}
\end{itemize}
The timing of Galileo's December 11 letter is certainly very suspicious. Was Galileo ignorant of the phases of Venus and their importance before it was pointed out to him? Did Galileo steal the idea from Castelli? Very possibly.

In fact, in a letter of November 13 of the same year, Galileo seems to state expressly that he had no new planetary discoveries to report, implying that he did not yet know about the phases of Venus, contrary to his later assertions. His defenders claim that Galileo’s phrase should instead be read as saying that he has made no new discoveries “around” the planets, that is, discovered no new moons.\footnote{\cite[200]{DrakeVenus}.} Perhaps so. It would make sense for Galileo to be on the lookout for moons and the like, as he had found for Jupiter and Saturn.\footnote{\S\S\ref{jupiter}, \ref{saturn}.} But that's all the more reason for him to miss Venus' phases. At this time Venus was well over half full, so its shape would not have been very remarkable unless you were specifically paying attention to it. He could easily have missed it if he was too busy moon hunting and looking only ``around'' the planets.

But if Galileo had not observed the phases of Venus before Castelli’s letter, then how could he later give an accurate description of their appearance dating back two months before this letter? Easy: by fabricating data and passing them off as actual observations, as he did on many other occasions.\footnote{\S\S\ref{fallandweight}, \ref{gravitationalconstant}, \ref{resistance}, \ref{pendulum}, \ref{saturn}, \ref{sunspotshelio}.} After all, he was surely concerned to get the important pro-Copernican argument from the phases of Venus on the record as quickly as possible and claim it for himself, and for this purpose it would be important to have observed Venus’ fully gibbous appearance in the fall (the next opportunity to observe it in this form would be months away). So Galileo certainly had a strong motive to fabricate this data. Making observations throughout most of December, after receiving Castelli's letter, would also have been enough to give him great confidence that the heliocentric explanation for the phases of Venus was right. So faking the data was not risky.

Galileo's defenders have a counterargument to this. They claim that Galileo could not have fabricated the data in question even if he had wanted to. According to them, the changes in appearance of Venus during these months were so complex and “non-linear” that Galileo could never have given such an accurate account if he had not if fact made these observations. Specifically, Galileo correctly describes the fact that the transition from a gibbous to semicircular phase is quite rapid, while a roughly semicircular phase lingers for a considerable time.
\quote{Castelli’s letter cannot have been the spark that ignited Galileo’s programme of observation of Venus. It was simply too late. If he only then had started observing Venus, he would have seen it already nearing the exact semicircular phase, thus completely missing the non-linear patterns of change. And he could not possibly have been able to calculate the duration of one month for the “lingering” phenomenon. In other words, Galileo cannot have predicted Venus’s non-linear patterns of behaviour by re-constructing them ‘backwards’. For a Copernican it might have been easy to predict that Venus should display phases. However, it is one thing to predict this type of behaviour qualitatively and quite another to predict the non-linear patterns of change of Venus’s phases. A quantitative analysis would have required of Galileo a sophisticated mathematical theory that he did not have. There remains only one possibility, namely, that Galileo really did observe Venus’s non-linear patterns of behaviour.\footnote{\cite[117]{PalmieriVenus}.}}
I say that, on the contrary, Galileo could easily have reconstructed these phenomena. He would not have needed any sophisticated mathematics as all. All he would have had to do would have been to simulate its appearance by looking at a half-painted sphere representing Venus from vantage points corresponding to the Earth’s position relative to it.

\begin{figure}[tp]\centering
\includegraphics[width=\textwidth]{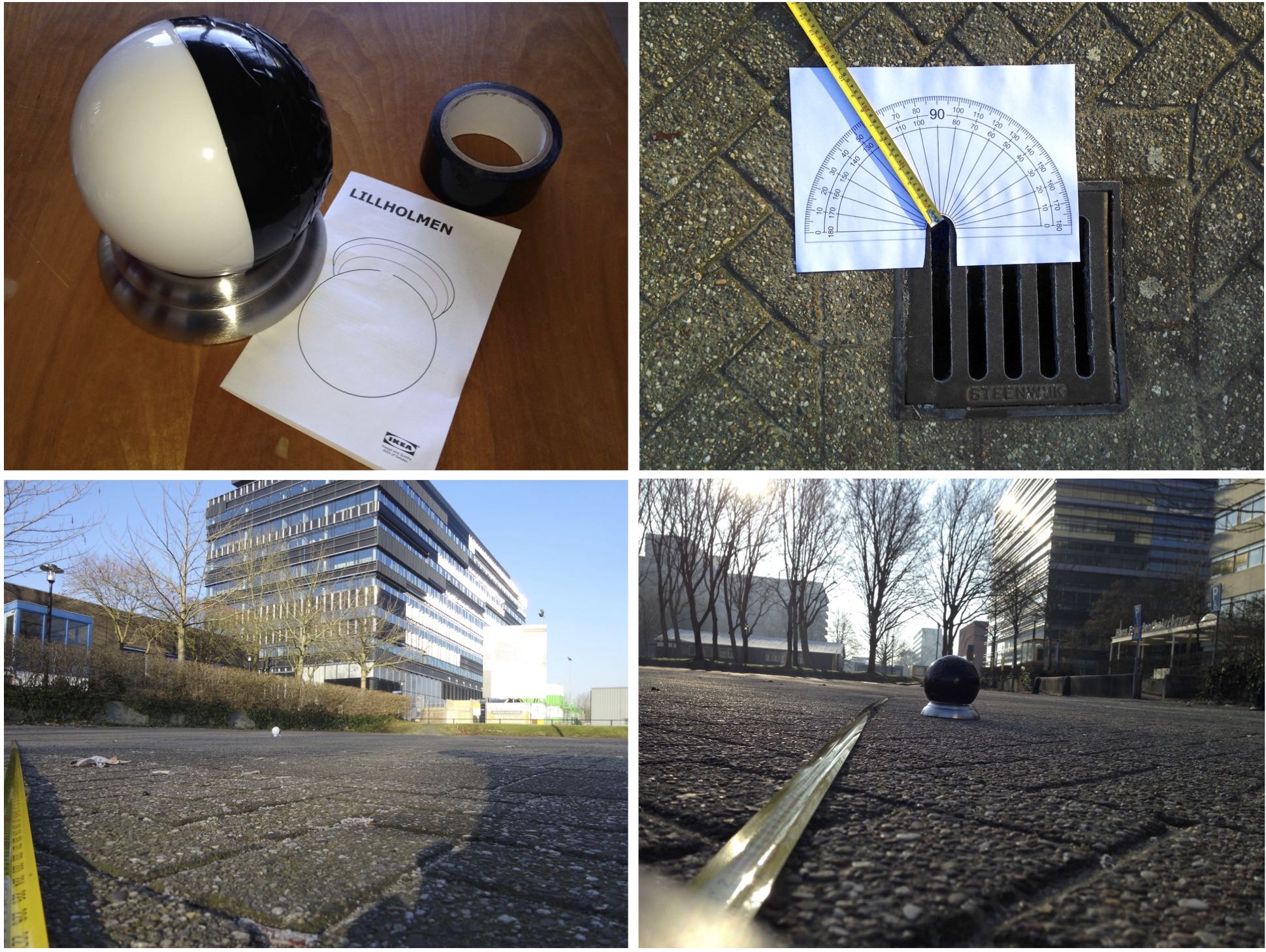}
\caption[Simulating phases of Venus.]{The simple setup I used to simulate the phases of Venus.}
\label{Vsimtools}
\end{figure}

\begin{figure}[tp]\centering
\includegraphics[width=0.8\textwidth]{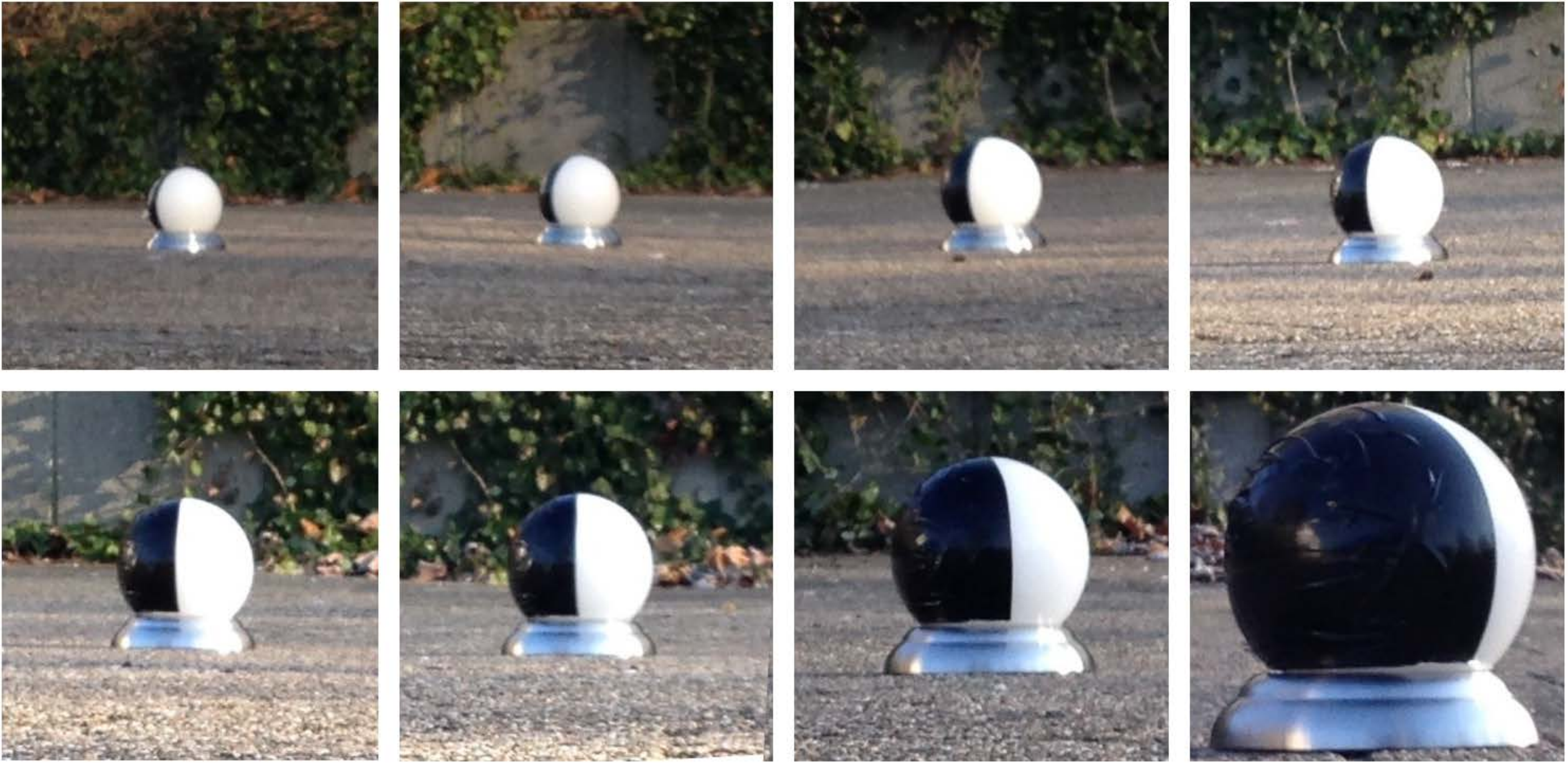}
\caption[Physical simulation of phases of Venus.]{Simulated appearances of Venus as seen from the Earth for Sep 4, Sep 25, Oct 16, Nov 6, Nov 27, Dec 18, 1610, Jan 8, Jan 29, 1611.}
\label{Vsim}
\end{figure}

\begin{figure}[tp]\centering
\includegraphics[width=0.8\textwidth]{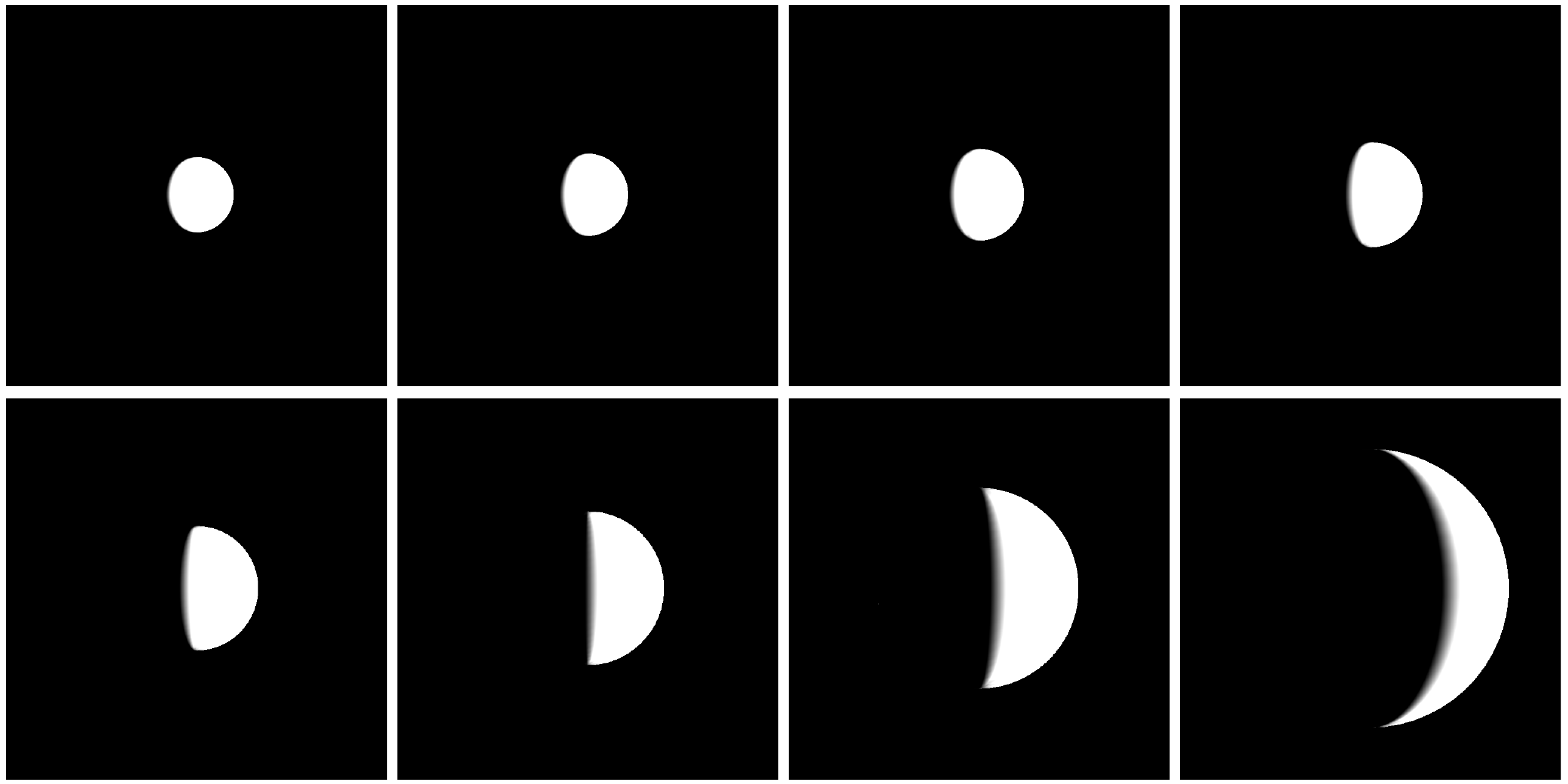}
\caption[Simulation of phases of Venus.]{The same simulation as in Figure \ref{Vsim}, carried out using 3D software.}
\label{VsimMathematica}
\end{figure}

I carried out such a simulation using very simple means (Figure \ref{Vsimtools}). The results are shown in Figure \ref{Vsim}. I used a white spherical lamp as Venus. I covered half of it in black to represent the half not illuminated by the Sun. I pointed the white half toward a center point (the Sun) 4.34 meters away. I then marked off a circle of radius 6 meters with the same center, representing the orbit of the Earth. I used the fact that Venus was seen exactly semicircular on December 18 to find where the Earth must have been it its orbit that day.\footnote{\cite[212]{PetersVenus1610}.} I placed a camera at this position and photographed the Venus sphere. I then used a protractor positioned at the Sun to reconfigure the setup to correspond to other dates, counted forward and backwards from December 18 in 21 day increments using the simplest possible estimation for the motions of these planets (I simplistically assumed uniform circular motions for the Earth and Venus, so there is no advanced mathematical astronomy involved in any way, just basic calculations using the radii and orbital times of these two planets). I again photographed Venus from these positions. I did all of this in a rough-and-ready way in an empty parking lot using crude measurements. I also recreated the exact same setup using 3D software (Figure \ref{VsimMathematica}), which shows the results of this simulation without the accidental imperfections of my physical demonstration.

Galileo could easily have completed such a simulation from start to finish in just a few hours. And of course the idea of illustrating the phases of the moon by an illuminated or half-painted sphere had been commonplace since antiquity, so Galileo would not have needed much imagination to come up with this scheme.

The results of this simple simulation are very close to the true appearances.\footnote{Figure \ref{venusphasesillustrationfig} shows exactly computed actual appearances. Such modern reconstructions of the actual appearances are also given in \cite{PalmieriVenus}, \cite{GingerichVenus1610}, \cite{PetersVenus1610}.} In particular, the simulation is easily sufficient to reproduce the allegedly so unpredictable “non-linear” phenomena that Galileo got right in his December 30 report. So the claim that it would have been impossible for Galileo to recreate these appearances after the fact is definitely false. We may note also that one argument that Galileo’s account has “the ring of a record of visual impressions rather than an account coloured by calculations” in that it “has a highly visual character.”\footnote{\cite[213--214]{PetersVenus1610}.} Obviously this is consistent with my simulation just as well as actual observations.

The hypothesis that Galileo simulated his Venus observations by using such a model is lent some further credibility by its close parallels with his treatment of sunspots.\footnote{\S\ref{sunspotshelio}.} As we have seen, Galileo realised that sunspots constituted an important pro-Copernican argument only quite late, based on the input of others, and needed to act fast in writing something about it without having the time for thorough observations. I suggest that this is an exact parallel of what happened also in the case of the phases of Venus.

This parallel undermines the common assumption that Castelli’s idea must already have been obvious to Galileo. One scholar, for example, thinks it “would be to dignify the idea beyond reasonable measure” to view Castelli’s suggestion as a significant insight; rather, “the thought that Venus might have phases was ‘in the air’” and hence Castelli’s contribution is to be considered quite trifling.\footnote{\cite[92]{AriewVenus}.} Another historian argues along similar lines that Galileo had no need to be spurred to action by Castelli’s letter, only by news of others making advanced telescopic observations. Around a day or two before hearing from Castelli, Galileo had received another letter, reporting that Clavius and his assistants at Rome had observed the moons of Jupiter.\footnote{\cite[200--201]{DrakeVenus}.} So Galileo now had serious competitors in the realm of advanced telescopic observations, or so it would have seemed to him. Presumably they would turn to the other planets next, and perhaps anticipate the discovery of the phases of Venus, whence Galileo’s sudden urgency. For “the problem was to have a good telescope, not to posses reasoning power that astronomers had never lacked.”\footnote{\cite[203]{DrakeVenus}.} The sunspots case is a counterexample to this claim: if there was no shortage of “reasoning power,” Galileo should have realised the potential importance of sunspots much earlier and not let himself be beaten to the punch about their curved appearance by his arch-rival Scheiner. The fact of the matter is that the sunspots argument for heliocentrism eluded Galileo for twenty years, despite the fact the was passionately committed to proving heliocentrism in novel ways, and despite the fact that he himself had written specifically and in detail about the very phenomenon at stake, and despite the fact that the argument is very simple. By analogy, this suggests that Castelli’s idea about Venus could very well have been news to Galileo: if he could somehow miss the sunspots argument for twenty years despite all of this, then he could certainly have failed to think of the Venus argument during his one initial frantic year of telescopic observations, when he had a myriad other novelties and issues to deal with all at once.

But perhaps the most interesting aspect of the parallel between the two cases is the possibility that they both involved the use of physical models to simulate celestial appearances. For, in the {\Dialogue}, one speaker reports regarding the appearances of the paths traced on the surface of the sun by sunspots as seen from the Earth that Galileo “assisted my understanding by representing the facts for me upon a material instrument, which was nothing but an astronomical sphere, making use of some of its circles---though a different use from that which they ordinarily serve.”\footnote{\cite[348]{galileo2sys1sted}.} The same sentiment is repeated later: the appearances of the sunspot paths “will become better fixed in my mind when I examine them by placing a globe at this tilt and then looking at it from various angles.”\footnote{\cite[352]{galileo2sys1sted}.} This is very closely analogous to the Venus simulation I outlined above, suggesting that the latter would have been quite natural to Galileo, and in keeping with his style of reasoning.\footnote{For further examples of Galileo preferring to think with physical objects, see \cite[67--71]{camcomp}.}

The following, then, are generally accepted facts about the sunspots matter: Galileo claimed to have conducted careful observations when he had not; according to his own account, Galileo simulated observations by looking at a physical sphere from a variable vantage point corresponding to the position of the Earth; Galileo failed to see an important pro-Copernican argument for a long time, despite it being simple and very naturally connected to his own work. The fact that these things did happen in the sunspots case suggests that they very well could have happened also in the Venus case.

In conclusion, if Galileo had wanted to fabricate or reconstruct Venus observations he had not made, he could easily have done so. His December 30 account, where he describes appearances of Venus going back to October, is perfectly consistent with the hypothesis that he only started serious observations after receiving Castelli’s letter in December and simulated earlier observations using a simple physical model. There are, as we have seen, furthermore a number of circumstantial indications that this would have been very much in keeping with his character and habits.

\subsection{Comets}
\label{comets}

``Have you seen the fleeting comet with its terrifying tail?''\footnote{Grassi, \cite[4]{DrakeOMalley}.} This was the question on everyone's lips in 1618, following the appearance of a comet ``of such brightness that all eyes and minds were immediately turned toward it.'' ``Suddenly, men had no greater concern than that of observing the sky. …\ Great throngs gathered on mountains and other very high places, with no thought for sleep and no fear of the cold.''\footnote{Grassi, \cite[6]{DrakeOMalley}.} ``That stellar body with its menacing rays …\ was considered as a monstrous thing,''\footnote{Grassi, \cite[4, 6]{DrakeOMalley}.} and, according to many, surely a cosmic omen foretelling imminent disasters.

Some urged a more dispassionate approach, arguing that ``the single role of the mathematician'' is merely to ``explain the position, motion, and magnitude of those fires.''\footnote{Grassi, \cite[6--7]{DrakeOMalley}.} Indeed, ``the mathematician'' had been so engaged for generations. Tycho Brahe, for instance, had worked extensively on comets.\footnote{As we noted in \S\ref{tycho}.}

It would be difficult for Galileo to enter this game, since he was such a poor mathematician. Not coincidentally, he has an argument for why one should ignore the serious mathematical astronomy of comets, namely that such accounts are hopelessly inconsistent:
\quote{Observations made by Tycho and many other reputable astronomers upon the comet's parallax …\ vary among themselves. …\ If …\ complete faith …\ be placed in them, one must conclude either that the comet was simultaneously below the sun and above it, …\ or else that, because it was not a fixed and real object but a vague and empty one, it was not subject to the laws of fixed and real things.\footnote{Galileo, {\it Assayer} (1623), \cite[257--258]{DrakeOMalley}.}}
Kepler is flabbergasted that someone who calls himself a geometer could write such drivel:
\quote{Certainly so far as Galileo is concerned, he, if anyone, is a skilled contributor of geometrical demonstrations and he knows …\ what a difference there is between the incredible observational diligence of Tycho and the indolence common to many others in this most difficult of all activities. Therefore, it is incredible that he would criticize as false the observations of all mathematicians in such a way that even those of Tycho would be included.\footnote{Kepler, appendix to {\it Hyperaspistes} (1625), \cite[351]{DrakeOMalley}.}}
The paradox disappears if one recognises that Galileo is not a skilled geometer after all.

Unlike serious mathematical astronomers (and perhaps precisely in order to avoid having to engage with them), Galileo maintained that comets were not physical bodies travelling through space at all, but rather a chimerical atmospheric phenomena. He believes that ``their material is thinner and more tenuous than fog or smoke.''\footnote{Galileo, {\it Assayer} (1623), \cite[254]{GalileoDiscOp}.} ``In my opinion,'' says Galileo, comets have ``no other origin than that a part of the vapour-laden air surrounding the earth is for some reason unusually rarefied, and …\ is struck by the sun, and made to reflect its splendour.''\footnote{\cite[81]{sheaGrev}. OGG.VI.94.}

Galileo's vapour theory of comets is inconsistent with basic observations, as he himself admits. If comets are nothing but ``rarefied vapour''---that is to say, some kind of pocket of thin gas---then you'd imagine that their natural motion would be straight up, like a helium balloon. Indeed Galileo does propose that comets have such paths. But then he at once admits that this doesn't fit the facts: ``I shall not pretend to ignore that if the material in which the comets takes form had only a straight motion perpendicular to the surface of the earth …, the comet should have seemed to be directed precisely toward the zenith, whereas, in fact, it did not appear so. …\ This compels us either to alter what was stated, …\ or else to retain what has been said, adding some other cause for this apparent deviation. I cannot do the one, nor should I like to do the other.''\footnote{\cite[82--83]{sheaGrev}. OGG.VI.98.} Bummer, it doesn't work. But Galileo sees no way out, so he just leaves it at that.

Galileo's contemporaries were not impressed. ``[Grassi's] criticism of Galileo is on the whole penetrating and to the point. He was quick to spot Galileo's inconsistencies. …\ Grassi produced an impressive array of arguments to show that vapours could not explain the appearance and the motion of the comets [as Galileo had claimed].''\footnote{\cite[84]{sheaGrev}.} For instance, the speeds of comets do not fit Galileo's theory. According to Galileo's theory, the vapours causing the appearance of comets rise uniformly from the surface of the earth straight upwards. Therefore the comet should appear to be moving fast when it is close to the horizon, and then much slower when it is higher in the sky. Just imagine a red helium balloon released by a child at a carnival: it first it shoots off quickly, but soon you can barely tell if it's rising anymore, even though it keep going up at more or less the same speed, because your distance and angle of sight is so different. But comets do not behave like that. Detailed observations of the comet of 1618 showed a much more constant speed than Galileo's hypothesis requires. 
 
Galileo also offered another very poorly considered argument against the correct view of comets as orbiting bodies, namely that their orbits would have to be unrealistically big: ``How many times would the world have to be expanded to make enough room for an entire revolution [of a comet] when one four-hundredth part of its orbit takes up half of our universe?''\footnote{Galileo, \cite[77]{sheaGrev}.} This is a poor argument, because the universe must indeed be very big and then some according to Copernican theory, in order to explain the absence of stellar parallax.\footnote{\S\ref{parallax}.} Since the earth's motion is observationally undetectable, the orbit of the earth must be minuscule in relation to the distance to the stars. That means there is plenty of room for comets. But Galileo conveniently pretends otherwise in his argument against comets. Evidently, Galileo ``was so intent on refusing Tycho that he failed to notice that he was pleading for a universe in which there would be no room for the heliocentric theory'' either.\footnote{\cite[88]{sheaGrev}.}

\section{Salon science}

\subsection{Galileo, populariser}

{\it The Cambridge Companion to Galileo} poses for itself the question: “What did Galileo actually do that made his image so great and so long-standing?” Its answer is not a list of great scientific accomplishments but rather: “Certainly his was the first main effort that fired the vision of science and the world that went well beyond limited intellectual circles.”\footnote{\cite[1--2]{camcomp}.} Galileo was a populariser, in other words. ``It was to the man of general interests that Galileo originally addressed his works.''\footnote{Drake, \cite[3]{GalileoDiscOp}.} Indeed, Galileo embraced this role, praising himself for ``a certain natural talent of mine for explaining by means of simple and obvious things others which are more difficult and abstruse.''\footnote{Galileo, {\it Assayer} (1623), \cite[265]{GalileoDiscOp}.}

I agree with these learned authors that Galileo wrote for the vulgar masses. I must add only one point, which they omit, namely that Galileo was driven to turn to popularisation {\em because he was so bad at mathematics}. ``Galileo scarcely ever got around to writing for physicists.''\footnote{Drake, \cite[2]{GalileoDiscOp}.} Yes, and he was scarcely {\em able} to do so either. The two are not unrelated.

Take for instance the ``new stars''\footnote{Novas or supernovas in modern parlance.} that appeared in Galileo's lifetime. One appeared in 1572. It was studied with great care by Tycho Brahe. Another appeared in 1604, when Galileo was 40 years old and an established professor of mathematics. But Galileo didn't make a contribution based on serious astronomy as Tycho had done. Instead he gave public lectures on the nova to a layman audience totalling more than a thousand people.\footnote{\cite[105]{drakeGatwork}.}

Galileo's little science extravaganzas were a hit at bourgeois dinner parties, as contemporary witnesses describe:
\quote{We have here Signor Galileo who, in gatherings of men of curious mind, often bemuses many concerning the opinion of Copernicus, which he holds for true. …\ He discourses often amid fifteen or twenty guests who make hot assaults upon him. …\ But he is so well buttressed that he laughs them off; and although the novelty of his opinion leaves people unpersuaded, yet he convicts of vanity the greater part of the arguments with which his opponents try to overthrow him. …\ What I liked most was that, before answering the opposing reasons, he amplified them and fortified them himself with new grounds which appeared invincible, so that, in demolishing them subsequently, he made his opponents look all the more ridiculous.\footnote{Querengo, 1616, \cite[452]{sleepwalkers}.}}
Again: Galileo's speciality is burlesque astronomical road shows, not serious science. Galileo's defenders refuse to admit the obvious, and instead go out of their way to try to save face on his behalf:
\quote{The technique …\ is actually a very sound, wise, and proper one; it really amounts to being concerned to avoid the straw-man fallacy; that is, before criticising an opponent, it is a sign of a serious critic to first strengthen the opposing argument as much as possible and interpret it in the most charitable manner; by so doing, one's criticism will really undermine the argument, rather than destroying one's own caricature invented to make one's own task easy.\footnote{\cite[311--312]{finretrying}.}}
The notion that Galileo is concerned with being charitable and avoiding straw men and caricature is preposterous. The truth is the exact opposite. We have already seen this for example in \S\ref{fallandweight}. The following two sections provide further illustrations of the same point.

\subsection{Three dimensions}
\label{3d}

Space is three-dimensional. Why? Because ``three is a perfect number''? Galileo alleges that this is Aristotle's answer, and goes on to refute it triumphantly: ``I do not understand, let alone believe, that with respect to legs, for example, the number three is more perfect than four or two.''\footnote{\cite[11]{galileo2sys}.} You can almost hear across the centuries the laughter erupt among the tipsy dinner-party guests as this clever put-down is delivered. Classic Galileo! Just picture a three-legged beast hobbling about like a limp dog---could anything be more farcical? And yet Galileo's opponents have unwittingly committed themselves to calling this ludicrous spectacle the pinnacle of perfection. So ridiculous are the consequences of Aristotle's metaphysical teachings.

Galileo cannot be beat when it comes to soirée entertainment, but does his argument have any genuine scientific merit? No. Refuting the idea that three is a perfect number is not cutting-edge research by anyone's measure. And even if it was, Galileo's argument still wouldn't be much of contribution to the debates since it is quite clearly more concerned with opportunistic showmanship than serious and balanced engagement with Aristotle's thought. Let us look at the passage from Aristotle on which his critique is based:
\quote{The three dimensions are all that there are. …\ As the Pythagoreans say, the world and all that is in it is determined by the number three, since beginning and middle and end give the number of an `all'. …\ Further, we use the terms in practice in this way. Of two things, or men, we say `both', but not `all': three is the first number to which the term `all' has been appropriated. And in this …\ we do but follow the lead which nature gives. Therefore …\ [three-dimensional] body alone among magnitudes can be complete.\footnote{Aristotle, {\it De Caelo} I.1, trans.~J.~L.~Stocks.}}
When Galileo sets out to attack this point of view, he has Aristotle's mouthpiece Simplicio argue that Aristotle ``has sufficiently proved that there is no passing beyond the three dimensions, …\ and that therefore the body, or solid, which has them all, is perfect.''\footnote{\cite[10]{galileo2sys}.}

We see that there are some differences between Aristotle's account and Galileo's paraphrase. First of all Aristotle hardly pretends to have ``proved'' the matter; rather he merely supports it with some charming historical and etymological considerations. Furthermore, the term ``perfect''---on which Galileo's refutation hinges---is not used by Aristotle in the sense needed for Galileo's argument. The term does not occur at all in the translation I quoted, which instead uses ``complete.'' Regardless of which of these two renditions is more accurate from the point of view of Greek linguistics, it is obvious in any case that Aristotle doesn't mean ``perfect'' in the sense needed for Galileo's rebuttal, which would be something like ``optimal'' or ``ideal.'' Galileo is obviously ``playing on this ambiguity in order to weaken Simplicio's position.''\footnote{\cite[468]{galileo2sys}.} Aristotle is merely making the rather obvious point that three dimensions are the ``whole'' or ``entirety'' of what there is. But it would be no fun for Galileo to admit as much since it would spoil his punchline about the three-legged ogre. So we must conclude that Galileo's argument is not really engaging with Aristotle in an honest way, but rather relies on a cunning equivocation to score an easy point.

Still in modern times though there is no shortage of philosophers who fall for Galileo's argument hook, line, and sinker, and even use it as an exemplary illustration of his argumentative prowess:
\quote{Galileo is quoting Aristotle accurately. Simplicio is not inventing these reasons and merely attributing them to Aristotle. …\ The force of Galileo's rejection …\ is strengthened since Aristotle's account is foolish.\footnote{\cite[29]{pitt}.}}
It is hard to know how to reply to such hogwash. The exact negation of each of those statements would be right on the money.

As for his own explanation of why space has three dimensions, Galileo purports to offer a “geometrical demonstration of threefold dimensionality.”\footnote{\cite[12]{galileo2sys1sted}.} It goes like this. Consider three lines meeting perpendicularly, such as one can effect by, for instance, drawing an ``L'' on a piece of paper, placing it flat on the ground, and letting a plumb line touch the corner of the L. From this configuration Galileo concludes: ``And since clearly no more lines can meet in the said point to make right angles with them …\ then the dimensions are no more than three.''\footnote{\cite[14]{galileo2sys}.} The notion that this is a ``proof'' that there are three dimensions is ridiculous. It's an assertion that the conclusion is ``clear,'' not a proof of it. And yet---perhaps because it comes with lettered geometrical diagrams---it is enough to fool philosophers even today into believing that it is an ``elegant little proof'' and even a ``rigorous demonstration.''\footnote{\cite[32, 30]{pitt}.} As long as readers remain this gullible, Galileo's undeserved reputation as a great geometer is sure to thrive.

\subsection{Babylonian eggs}
\label{Babylonianeggs}

Another prime example of Galileo's rhetorical skills is his debate with Grassi on comets. It is a fact that Galileo is dead wrong on the scientific issues at stake,\footnote{Recall Galileo's very poor work on comets from \S\ref{comets}.} and yet he somehow managed to ``win'' the debate, in the eyes of many. Galileo's rousing mockery of his opponent is so satisfying that many readers are seduced into celebrating it as proof of Galileo's philosophical acumen. You can read Galileo's triumphant put-downs of his opponent and go ``yeah, crush him!'' It's the same kind of pleasure as watching the villain get punched in the face in an action movie. But a little reflection shows that this hero-versus-villain dynamic that Galileo tries to cultivate is a dishonest fiction that has very little to do with reality.

One of Galileo's most celebrated passages concerns eggs. The context is this. Grassi makes the absolutely correct point that comets, if they entered the earth's atmosphere, would quickly heat up to very great temperatures due to friction of the air. In support of this point, Grassi quotes a 10th-century Byzantine author, Suidas, who claimed that ``The Babylonians whirl[ed] about eggs placed in slings …\ [and] by that force they also cooked the raw eggs.''\footnote{Grassi, \cite[119]{DrakeOMalley}.} Grassi also quotes passages describing similar phenomena in Ovid, Lucan, Lucretius, Virgil, and Seneca. ``For who believes that men who were the flower of erudition and speak here of things which were in daily use in military affairs would wish egregiously and impudently to lie? I am not one to cast this stone at those learned men.''\footnote{Grassi, \cite[119--120]{DrakeOMalley}.}

Galileo is unable to answer the substantive point. Indeed, he thinks comets entering the atmosphere would cool down because of the wind rather than heat up because of friction. Galileo is wrong and Grassi is right about the actual scientific issue about comets. But that's nothing Galileo's trademarked sophistry can't work around. Galileo finds a way to ``win'' the debate without actually offering any correct scientific claim regarding the actual subject of comets. He does this by gloatingly attacking Grassi for relying on books rather than experimental evidence:
\quote{If [Grassi] wants me to believe with Suidas that the Babylonians cooked their eggs by whirling them in slings, …\ I reason as follows: ``If we do not achieve an effect which others formerly achieved, then it must be that in our operations we lack something that produced their success. And if there is just one single thing we lack, then that alone can be the true cause. Now we do not lack eggs, nor slings, nor sturdy fellows to whirl them; yet our eggs do not cook, but merely cool down faster if they happen to be hot. And since nothing is lacking to us except being Babylonians, then being Babylonians is the cause of the hardening of eggs, and not friction of the air.'' …\ Is it possible that [Grassi] has never observed the coolness produced on his face by the continual change of air when he is riding post? If he has, then how can he prefer to believe things related by other men as having happened two thousand years ago in Babylon rather than present events which he himself experiences?\footnote{Galileo, {\it Assayer} (1623), \cite[272]{GalileoDiscOp}.}}
Not a few modern philosophers blindly and uncritically fall for the ruse: ``Galileo shot back with …\ a blistering critique in which he pillories [Grassi] and articulates a tough-minded empiricism as an alternative to the mere citation of venerable authority.''\footnote{\cite[135]{McGrewPhilSciAnth}.}

Galileo would no doubt be very pleased that so many readers still to this day come away with the impression that ``tough-minded empiricism'' is what sets him apart from his opponents. That is precisely the intended effect of his ploy. It has very little basis in reality, however. Just a few pages earlier in the same treatise, Grassi describes extensively various laboratory experiments he has carried out himself with regard to another point. ``I decided that no industry or labor ought to be spared in order to prove this by many and very careful experiments,''\footnote{Grassi, \cite[115]{DrakeOMalley}.} says this supposed obstinate enemy of empirical science. So the notion that Galileo is the only one ``tough-minded'' enough to reject authority in favour of experiment is very far off the mark.

Even in the passage criticised, Grassi is clearly not engaged in ``the mere citation of venerable authority.'' Rather he honestly and openly cites sources purporting to truthfully report empirical information, just like any scientist today cites previous works without re-checking all the experiments personally. Grassi does not believe that these authors are automatically right because they are ``venerable authorities.'' Rather he explicitly considers the possibility that they are wrong, but estimates (reasonably but falsely) that they are probably right. For that matter, Galileo himself was not above believing falsehoods on the basis of ``venerable authorities.'' We have seen him make an error of this type in his theory of tides.\footnote{\S\ref{tides}.} He also considers ``credible'' the ancient myth of Archimedes setting fire to enemy ships by means of mirrors focussing the rays of the sun.\footnote{\cite[48]{galileo2newsci2ndedWT}, OGG.VIII.86.} Descartes sensibly took the opposite view.\footnote{Descartes to Mersenne, 11 October 1638, \cite[389]{drakeGatwork}.}

Altogether, Galileo is scoring easy points with his taunts about the eggs, by dishonestly pretending that a simplistic point about empiricism was the crux of the matter. It is worth keeping the context of the passage in mind. Indeed, the pro-Galileo interpretation I quoted above comes with its own origin story: ``In the course of his career [Galileo] engaged in many controversies and made powerful enemies. One of those enemies was the Jesuit Grassi, who published an attack on some of Galileo's works.''\footnote{\cite[135]{McGrewPhilSciAnth}.} This framing goes well with the notion of the ``tough'' Galileo bravely defending himself against ``attacks'' from the ``powerful'' establishment. But the reality is quite different. Grassi was not a ``powerful enemy'': he was a middling college professor just like Galileo. And the conflict did not start with Grassi ``attacking'' Galileo, but precisely the other way around. Grassi published a fine lecture on comets in which he argued, correctly, that the absence of parallax shows that comets are beyond the moon. Galileo is not mentioned in this work. Galileo read the lecture and filled the margins, as one scholar has observed, with an entire vocabulary's worth of savage expletives.\footnote{``There is no misplaced gentleness in the
marginal jottings …\ The expletives alone would make a
vocabulary of good Tuscan abuse: pezzo d'asinaccio, elefantissimo, bufolaccio, villan poltrone, balordone, barattiere, poveraccio, ingratissimo villano, ridicoloso, sfacciato, inurbano.'' \cite[152]{SantillanaCrime}.} Galileo then published an attack on Grassi which was not much more restrained than these marginal notes. Grassi replied to it. It is this reply that is called ``an attack on some of Galileo's works'' in the pro-Galilean quotation above.

In sum, Galileo's celebrated ``pillorying'' of Grassi was not a ``tough'' defence against an ``attack'' on ``some of his works'' by ``powerful enemies.'' The ``enemy'' was not a ``powerful'' arm of ``authority,'' but a conscientious scholar who was right about comets based on good scientific arguments that Galileo rejected. And the enemy was not a cruel aggressor going after ``some works'' by Galileo unprovoked; rather, the ``some works'' in question was an aggressive attack initiated by Galileo in the first place. Furthermore, Galileo's enemy did not favour venerable authority over empiricism, but rather based his analysis of comets on much more thorough empirical work than Galileo did. Altogether, the simplistic contrast between Grassi the credulous believer in authority and Galileo the experimenter has little basis in fact. Galileo is scoring easy points with his taunts about the eggs, by dishonestly pretending that an elementary point about empiricism was the crux of the matter.

\subsection{The Bible}
\label{joshuaargument}

The Bible says next to nothing about astronomical matters. It is more concerned with war, and it is in this context only that it has occasion to speak of the motions of heavenly bodies. Thus, in the Book of Joshua, we find our hero with the upper hand in battle, but alas dusk is drawing close. What a pity if the enemies ``delivered up …\ before the children of Israel'' should be able to get away under the cover of darkness. ``Then spake Joshua to the Lord,'' and he said: ``Sun, stand thou still.'' ``And the sun stood still …\ until the people had avenged themselves upon their enemies.''\footnote{Joshua, X.12--13.} This is the extent of astronomy in the holy book. Nowhere does it say that the earth is in the center of the universe, for example. Just this one passage about how ``the sun stood still'' that one time to ensure that all the infidels could be killed.

Obviously scientists have little reason to engage with such a tangential allusion to cosmology. But Galileo’s philosophical enemies saw an opportunity. By persistently and prominently accusing Galileo of proposing theories contrary to scripture they forced him into a dilemma: either let the argument stand unopposed, and hence let his enemies have the last word, or else get involved with the dangerous matter of scriptural interpretation.\footnote{\cite[Ch.\ 3]{Blackwell}.} Galileo foolishly took the bait. Now all the Aristotelians had to do was to sit back and watch Galileo march to his own ruin in this minefield.

In this context, Galileo offered some common-sense platitudes on the relation between science and religion. Many have found it appealing to see in these writings modern conceptions being born. Was it Galileo who showed how faith and science can coexist? How they need not undermine or conflict with one another since one is about the spiritual and the other about the physical? Galileo indeed makes such a case:
\quote{Far from pretending to teach us the constitution and motions of the heavens and the stars, …\ the authors of the Bible intentionally forbore to speak of these things, though all were quite well known to them. …\ The Holy Spirit has purposely neglected to teach us propositions of this sort as [they are] irrelevant to the highest goal (that is, to our salvation). …\ The intention of the Holy Ghost is to teach us how one goes to heaven, not how heaven goes.\footnote{Galileo, {\it Letter to Duchess Christina} (1615), \cite[184--186]{GalileoDiscOp}.}}
Even a recent Pope praised Galileo for his supposed insight on this subject: ``Galileo, a sincere believer, showed himself to be more perceptive [in regard to the criteria of scriptural interpretation] than the theologians who opposed him.''\footnote{Pope John Paul II, ``Lessons of the Galileo case,'' {\it Origins: Catholic News Service}, November 12, 1992, 22, 372, \cite[291]{mcmullinscripture}.}

I disagree with this papal statement on two grounds. First of all, Galileo was not pioneering a new vision for the roles of science and religion more perceptively than anyone else. Rather, he was merely recapitulating elementary ideas that were virtually as old as organised Christianity itself. ``[Galileo's] exegetical principles were not in any sense novel, as he himself went out of his way to stress. They were all to be found in varying degrees of explicitness in Augustine's {\it De Genesi ad Litteram}''---written twelve centuries before Galileo---``and, separately, they could call on the support of other [even] earlier theologians.''\footnote{\cite[314]{mcmullinscripture}.} Galileo indeed quotes at great length from Augustine and the church fathers. Not that Galileo knew anything about the history of biblical interpretation: ``He had no expertise whatever in that area, so he evidently asked his Benedictine friend, Castelli, to seek out references that would support the exegetical principles he had outlined.''\footnote{\cite[287]{mcmullinscripture}.}

Furthermore, it is highly doubtful whether Galileo was truly ``a sincere believer,'' as he purported to be. Recent scholarship has made a compelling case for “two Galileos, the public Catholic and the private sceptic.”\footnote{\cite[249]{wootton}.}
\quote{The only decisive document we have [is a 1639 letter to Galileo from] Benedetto Castelli, Galileo’s old friend, former pupil and long-time intellectual companion. …\ If anyone was in a position to know if Galileo was or was not a believer it was Castelli. …\ [Castelli writes in his letter that he] has heard news of Galileo that has made him weep with joy, for he has heard that Galileo has given his soul to Christ [in his old age---Galileo was 75 at this point]. Castelli immediately refers to the parable of the labourers in the vineyard (Matthew 20.1–16): even those who were hired in the last hour of the day received payment for the whole day’s work. …\ Then …\ he turns to the crucifixion, and in particular to the two thieves crucified on either side of Christ (Luke 23.39–43). One confessed Christ as his saviour and was saved; the other did not and was damned. …\ Castelli’s [point] is clear and unambiguous. He believes Galileo is coming to Christianity at the last moment, but not too late to save his soul. There is no conceivable interpretation of this letter which is compatible with the generally held view that Galileo was, throughout his career, a believing Catholic.\footnote{\cite[247--248]{wootton}.}}
It is indeed striking that God plays an essential role in the scientific systems of so many other 17th-century scientists---such as Kepler, Descartes, Leibniz, and Newton---but no part whatsoever in Galileo’s. While many of Galileo's contemporaries made great efforts to synthesise scientific and transcendental knowledge that can be nothing but sincere, Galileo brings up religion only when it conveniently serves his purposes, such as when professing in eloquent terms to be more committed to finding the true meaning of scriptural passages than his opponents who used those texts against him.

This sense of opportunism is reinforced by how poorly Galileo's patchwork of trite ideas hold together. ``He says first that saving the appearances is not enough to demonstrate the truth of a hypothesis and ends by remarking that saving the appearances is the most that can be demanded of an hypothesis.''\footnote{\cite[285--286]{mcmullinscripture}.} He argues that ``officials and experts in theology should not arrogate to themselves the authority to issue decrees in professions they neither exercise nor study,''\footnote{\cite[100]{GA}.} but apparently considered himself entitled to issue decrees in a profession that he neither exercised nor studied, namely that of scriptural interpretation. He also says that one should not rely on the Bible in scientific matters, then does precisely this himself, ``violating his own prohibition against using Scripture to support a philosophical thesis about the natural world.''\footnote{\cite[280--281]{mcmullinscripture}.} There is little consistency in his views, except of course that he consistently chooses to espouse whatever principles serve his rhetorical purposes at any given moment.

Galileo's interpretation of the Joshua passage is a prime example of his shameless drive to score rhetorical points at any cost. It is perfectly reasonable to argue that the phrase about the sun ``standing still'' should not be taken too literally. It is commonly accepted, as Galileo observes, that various things in the Bible ``were set down in that manner by the sacred scribes in order to accommodate them to the capacities of the common people, who are rude and unlearned.''\footnote{Galileo, {\it Letter to Duchess Christina} (1615), \cite[181]{GalileoDiscOp}.} Indeed, if the Bible is read literally, ``it would be necessary to assign to God feet, hands, and eyes,'' as Galileo says, but those passages are only figures of speech, according to orthodox Christian understanding. When the Old Testament says that the commandments handed to Moses were ``written with the finger of God,'' the intended point is of course not that God has an actual physical finger and that he needs it to write. It doesn't make a whole lot of sense that he could create the entire universe in under a week, or flood the entire earth at will, yet if he has to write something down he has to painstakingly trace it out in clay with his finger. Perhaps it is the same with the sun ``standing still.'' It's just a phrase adapted to everyday speech, not a scientific account. In fact, even Copernicus himself speaks of ``sunrise'' and ``sunset,'' as Galileo points out,\footnote{Galileo, {\it Letter to Duchess Christina} (1615), \cite[202]{GalileoDiscOp}.} even though the sun doesn't move in his system. So it is hardly unreasonable to think that ``the sacred scribes'' used this kind of common parlance as well, even if they knew that the sun is always stationary.

That's all fine and well.  But Galileo does not stop with this balanced point. He did not become a salon sensation by sensibly looking for middle ground and mutual reconciliation. His audience expects a more triumphant and extravagant finishing blow. Seemingly to this end, Galileo makes the outlandish claim that the Joshua passage in fact literally agrees best with heliocentrism rather than geocentrism:
\quote{If we consider the nobility of the sun …\ I believe that it will not be entirely unphilosophical to say that the sun, as the chief minister of Nature and in a certain sense the heart and soul of the universe, infuses by its own rotation not only light but also motion into other bodies which surround it. …\ So if the rotation of the sun were to stop, the rotations of all the planets would stop too. …\ [Therefore,] when God willed that at Joshua's command the whole system of the world should rest and should remain for many hours in the same state, it sufficed to make the sun stand still. …\ In this manner, by the stopping of the sun, …\ the day could be lengthened on earth---which agrees exquisitely with the literal sense of the sacred text.\footnote{Galileo, {\it Letter to Duchess Christina} (1615), \cite[212--214]{GalileoDiscOp}.}}
This is a terrible argument. It is so unscrupulous that its absurdity can be exposed simply by quoting the words of Galileo himself, written in another context:
\quote{If the terrestrial globe should encounter an obstacle such as to resist completely all its whirling motion and stop it, I believe that at such a time not only beasts, buildings, and cities would be upset, but mountains, lakes, and seas, if indeed the globe itself did not fall apart. …\ This agrees with the effect which is seen every day in a boat travelling briskly which runs aground or strikes some obstacle; everyone aboard, being caught unawares, tumbles and falls suddenly toward the front of the boat.\footnote{\cite[212]{galileo2sys1sted}.}}
So in this manner ``Joshua would have destroyed not only the Philistines, but the whole earth.''\footnote{\cite[439]{sleepwalkers}.} Not to mention that the idea that the sun's rotation on its axis is the only thing moving the planets is completely unsubstantiated in the first place. It seems that Galileo pretended to believe in it on this occasion solely for the sake of being able to make this scriptural argument. Once again the hypocrisy and unbridled opportunism of Galileo's forays into biblical interpretation are plain to see. It is very difficult, if not impossible, to see this interpretation of the Joshua passage as a scientific argument that Galileo genuinely believed. It does make sense, however, as a crowd-pleasing bit of sophistry. If you are an Italian aristocrat who enjoys seeing the learned establishment lose face but don't want to rock the boat yourself, then you can live vicariously through Galileo's snappy comebacks and provocations. To this end it matters little whether they are scientifically sound or not.

\subsection{First Inquisition}

Galileo's famous conflict with the church was entirely unnecessary. It arose precisely because Galileo was a lampooning populariser rather than a mathematical astronomer and scientist. ``[Galileo] was far from standing in the role of a technician of science; had he done so, he would have escaped all trouble.''\footnote{\cite[vii]{SantillanaCrime}.} Copernicus’ book had long been permitted, and Galileo’s own {\it Letters on Sunspots} of 1613 had been censored only where it referred to scripture, not where it asserted heliocentrism. The church establishment had no interest in prosecuting geometers and astronomers.

Today many take for granted that a fundamental rift between science and religion was unavoidable. Some have imagined for instance that Galileo defied the worldview of the church by demoting the earth from its supposedly ``privileged'' position. 20th-century playwright Bertolt Brecht appreciated the dramatic flare of framing the conflict in such terms when he wrote a play about Galileo. He has one of the characters argue the privilege point passionately:
\quote{I am informed that Signor Galilei transfers mankind from the center of the universe to somewhere on the outskirts. Signor Galilei is therefore an enemy of mankind and must be dealt with as such. Is it conceivable that God would trust this most precious fruit of his labor to a minor frolicking star? Would He have sent His Son to such a place? …\ The earth is the center of all things, and I am the center of the earth, and the eye of the Creator is upon me.\footnote{Bertolt Brecht, {\it Galileo}, \cite[72--73]{brecht}, \cite[271]{mcmullinscripture}.}}
But historically this is nonsense, to be sure. Nobody was concerned about this at the time. In fact, classical cosmology clearly stipulated that the Earth was not at all in a privileged position but rather condemned to its very lowly place in the universe. Doesn't everybody know that hell is just below the surface of the earth, while heaven is way up above? Clearly, then, being at the center of the universe is nothing to be proud of.

It was a commonplace argument in Galileo's time ``that the earth is located in the place where all the dregs and excrements of the universe have collected; that hell is located at the centre of this collection of refuse; and that this place is as far as possible from the outermost empyrean heaven where the angels and blessed reside.''\footnote{\cite[377]{Fin2019review}.} Even Galileo himself added to the pile, writing in an Aristotelian mode that ``after the marvellous construction of the vast celestial sphere, the divine Creator pushed the refuse that remained into the center of that very sphere and hid it there lest it be offensive to the sight of the immortal and blessed spirits.''\footnote{OGG.I.344, \cite[109]{sheaSN}.} Contemporaries reasoned alike: considering ``the Vileness of our Earth,'' it ``must be situated at the center, which is the worst place, and at the greatest distance from those Purer and incorruptible Bodies, the Heavens,'' wrote John Wilkins, an Anglican bishop.\footnote{John Wilkins, \cite[110]{sheaSN}.} This is the very opposite of the argument retrospectively imagined by Brecht and other modern minds.

In reality, ``a major part of the Church intellectuals were on the side of Galileo, while the clearest opposition to him came from secular ideas'' and philosophical opponents.\footnote{\cite[xii]{SantillanaCrime}. The same point is made by \cite[288]{drakeGatwork}.} It was only because Galileo got involved with biblical interpretation that he ended up in the crosshairs of the Inquisition. Nobody minded mathematical astronomy, but the question of who has the right to interpret the Bible was the stuff that wars were made of. Luther challenged church authority and emphasised personal understanding of the Bible---``sola scriptura,'' as the motto went. This was the core belief of protestantism, and eradicating protestantism was top of the agenda for the catholic church. Once they baited him into commenting on the Bible, it was all too easy for Galileo's enemies to connect Galileo's otherwise harmless dabbling to this heresy du jour---a matter on which the church could not afford to show any weakness.

There is only one mystery: Why did Galileo walk straight into such an obvious trap? The answer lies, as ever, in his mathematical ineptitude. Galileo was told by church authorities that ``if Galileo spoke only as a mathematician he would have nothing to worry about.''\footnote{\cite[249]{drakeGatwork}, referring to a statement of 1615.} Galileo would presumably have followed this advice if he could. The problem, of course, was that he did not have anything to contribute ``as a mathematician.'' Since a mathematical defence of heliocentrism was beyond his abilities, Galileo was left with no other recourse than to roll the dice and try his luck in the dangerous and unscientific game of scriptural interpretation.

The church was thus reluctantly drawn into these astronomical squabbles and had to do something. The Inquisition settled for a slap on the wrist: in the future, Galileo must not “hold, teach or defend [the Copernican system] in any way whatever,” they decided.\footnote{\cite[147]{GA}.} They also ordered mild censoring of Copernicus’ book, namely the removal of a brief passage concerning the conflict with the Bible and a handful expressions which insinuated the physical truth of the theory.\footnote{\cite[149, 200--202]{GA}.}

\subsection{Second Inquisition}

Galileo did indeed keep quiet for a number of years after being ordered to do so by the Inquisition. But times changed. After waiting for over a decade, Galileo felt it was safe to try the waters again. A new Pope was in power, Urban VIII, who was quite liberal. even said of the 1616 censoring of Copernicus that ``if it had been up to me that decree would never have been issued.''\footnote{\cite[312]{drakeGatwork}.} Galileo had good personal relation with this new open-minded Pope. So Galileo sensed an opening and obtained a permission to publish the Dialogue in 1632. Or rather, as the Inquisition would later put it, he “artfully and cunningly extorted” this permission to publish.\footnote{\cite[290]{GA}.} For when the permission was granted the Pope did not know about the private injunction of 1616 for Galileo to keep off the subject. When this came to light the Pope was outraged and felt, with good cause, that Galileo had been deliberately deceitful and reportedly stated that ``this alone was sufficient to ruin [Galileo] now.''\footnote{\cite[340--341]{drakeGatwork}.}

A special commission was thus appointed. It found many inappropriate things in the {\Dialogue}, but this was not a major issue, they noted, for such things “could be emended if the book were judged to have some utility which would warrant such a favor.”\footnote{\cite[222]{GA}.} The problem was instead that Galileo “overstepped his instructions” not to treat heliocentricism.\footnote{\cite[219]{GA}.} 

The same report also points out that Galileo had disrespected the Pope on another point as well. The Pope had asked Galileo to include the argument that since God is omnipotent he could have created any universe, including a heliocentric one. So even though the church does not agree with Copernicus, their own logic, namely belief in God's omnipotence, can be used to legitimate at least considering the possibility of this hypothesis. So that's a useful argument that Galileo could have used to try to find at least a little bit of common ground with his opponents. But instead of using it for such purposes of reconciliation as intended, Galileo used it to fuel the fires of conflict even more. He made had placed the Pope’s favourite argument “in the mouth of a fool,” the commission observed.\footnote{\cite[221]{GA}.} Galileo made Simplicio, the dumb character in the {\it Dialogue} who constantly expresses the wrong ideas and is proven wrong at every turn, be the one who spoke the Pope's words. He hardly did himself any favours with this disrespectful move.

Following these findings, the second Inquisition proceedings took place in 1633: 17 years after the first Inquisition where Galileo had gotten off easy, and one year after the publication of his inflammatory {\it Dialogue} in defence of Copernicanism. The outcome was a forgone conclusion. Galileo’s defence was transparently dishonest. He pretended that, in the {\it Dialogue}, “I show the contrary of Copernicus’s opinion, and that Copernicus’s reasons are invalid and inconclusive.”\footnote{\cite[262]{GA}.} This is of course pure nonsense. In private correspondence shortly before, Galileo had spoken more honestly, and stated that the book was ``a most ample confirmation of the Copernican system by showing the nullity of all that had been brought by Tycho and others to the contrary.''\footnote{Galileo to Diodati, October 1629, \cite[310]{drakeGatwork}.} But now before the Inquisition he had to pretend otherwise. In light of the accusations, Galileo continued, “it dawned on me to reread my printed {\it Dialogue},” and “I found it almost a new book by another author.”\footnote{\cite[277--278]{GA}.} These transparent lies did little to save him. He was forced to abjure. The {\it Dialogue} was prohibited, but not for its contents but rather, in the words of the Inquisition’s sentence, “so that this serious and pernicious error and transgression of yours does not remain completely unpunished” and as “an example for others to abstain from similar crimes.”\footnote{\cite[291]{GA}.}

There is a popular myth that Galileo muttered ``eppur si muove''---``yet it moves'' (the earth moves, that is)---as he rose from his knees after abjuring before the Inquisition. But this is certainly false.\footnote{\cite[356--357]{drakeGatwork}.} Obviously the Inquisition would not have tolerated such insubordination, especially since the whole point the trial in the first place was to punish Galileo for his defiance. Galileo had been shown the instruments of torture, and such a rebellious exclamation would have been the surest way to have them dusted off for the occasion. Today no historian believes the myth that Galileo mumbled these words before the Inquisition. Yet it remains instructive in warning us of the lengths many Galilean idol worshippers are willing to go to, who do not want to admit the many ignominious historical facts about their hero. The sheer multitude of such myths now universally regarded as busted should leave us open to the distinct possibility that we have not gotten to the end of them yet.

A similar myth, appealing to anti-religion ideologues, is that ``the great Galileo $\ldots$ groaned away his days in the dungeons of the Inquisition, because he had demonstrated $\ldots$ the motion of the earth.''\footnote{\cite[120]{voltaireW13}.} But in reality Galileo was sentenced more for his provocateurism than for his science, and furthermore he was never imprisoned in any ``dungeon.'' He was sentenced to house arrest. A visitor ``reported that [Galileo] was lodged in rooms elegantly decorated with damask and silk tapestries.''\footnote{\cite[166]{findefending}.} Soon thereafter he retired to ``this little villa a mile from Florence,'' where ``nearby …\ I had two daughters whom I much loved''\footnote{Galileo to Diodati, 25 July 1634, \cite[328]{heilbron}.} and where he also received many friends and guests. Many today would pay dearly for such a retirement. Galileo got it as a ``punishment.''

\section{Galileo evaluated}

\subsection{Conclusion of the above}

Galileo's catalogue of errors is extensive. It makes for quite a list even if we restrict ourselves only to central claims and only to matters his mathematically superior contemporaries did better than him: wrong value for cycloid area (\S\ref{cycloid}); erroneous infinitesimal geometry (\S\S\ref{Ginfintesimals}, \ref{fallandweight}); wrong value of gravitational acceleration $g$ (\S\ref{fallandweight}); erroneous theory of planetary speeds (\S\ref{planetaryspeeds}); erroneous claim that path of fall is semicircular in absolute space (\S\ref{pathoffallabs}); lack of proof that projectile motion is parabolic (\S\ref{projectilemotion}); erroneous theory that tides are caused by the motion of earth (\S\ref{tides}); erroneous claim that the simple pendulum is isochronic (\S\ref{pendulum}); erroneous theory that comets are an atmospheric phenomena (\S\ref{comets}). Within a few years of his death he was also corrected on inertia (\S\ref{inertia}) and the shape of a hanging chain (\S\ref{catenary}).

To be sure, other people made mistakes too. Suppose I concede that everyone has an equal comedy of errors to their name. Even so, this would still prove my point that Galileo was a dime a dozen scientist and not at all a singular ``father of modern science.'' But I do not need to concede this much. Galileo's sum of errors are not just par for the course. They are exceptionally poor, and in matters of mathematics altogether astonishing.

The persistent myth of ``Galileo's mathematical genius''\footnote{\cite[I.41]{CostabelLerner}.} must certainly die. Historians will never see Galileo's true colours as long as they keep taking it for granted that Galileo was ``the greatest mathematician in Italy, and perhaps the world'' in his time.\footnote{\cite[303]{heilbron}.} In reality, tell-tale signs of mathematical mediocrity permeate all his works. Like all too many modern mathematics students, he reaches for a calculator or instrument instead of thinking (\S\S\ref{cycloid}; cf.\ \S\S\ref{sector}, \ref{venus}), doesn't do the reading (\S\S\ref{adoption}, \ref{optics}), makes computational mistakes (\S\S\ref{planetaryspeeds}, \ref{sunspots}), methodological mistakes (\S\ref{jupiter}), fundamental mistakes in the entire structure of his theory (\S\ref{projectilemotion}), presents technical material in a clumsy way (\S\S\ref{fallandweight}, \ref{projectilemotion}), doesn't seem to know the difference between proof and assertion (\S\ref{3d}), and ignores the latest research which is much too advanced for him (\S\S\ref{bookofnature}, \ref{tycho}--\ref{againsttycho}, \ref{comets}). Some pages of Galileo would not be out of place somewhere in the middle of the piles of slipshod student homework that some of us grade for a living.

Furthermore, Galileo's errors concern some of his core achievements. I have discussed all of his notable scientific contributions and found much to object in every single case. Notably, Galileo uses ``his'' law of fall erroneously on a number of occasions (\S\S\ref{planetaryspeeds}, \ref{moonfall}, \ref{centrifugal}, \ref{pathoffallabs}). He also not infrequently presents arguments that are demonstrably inconsistent with his core beliefs, such as his tidal theory contradicting his own principle of relativity (\S\ref{tides}), his Joshua argument contradicting his own principle of inertia (\S\ref{joshuaargument}), and his objection to the geocentric explanation of sunspots being inconsistent with his own heliocentrism (\S\ref{sunspotshelio}).

Apollo 15 astronauts performed an experiment on the moon. They dropped a hammer and a feather and found that they fell with the same speed. ``Galileo was correct,'' they concluded in a famous video recording still often shown in science classrooms today. Actually, Lucretius was correct, because he is the one who said this would happen in the absence of air (\S\ref{fallandweight}). Galileo was wrong because he considered it ``obvious'' that the moon had an atmosphere (\S\ref{optics}). If the astronauts wanted to test Galileo’s theory they should not have dropped a hammer and a feather. They should have taken off their helmets and suits and tried to breathe. That would have showed you how “right” Galileo really was.

Posterity have chosen to remember only Galileo's successes while forgetting his numerous errors. It is easy to be a hero of science if you can count on such selective amnesia to always put you in the most flattering light. If there had been air on the moon, the astronauts would have hailed Galileo for this ``discovery'' instead.

The point generalises. Galileo made many claims that would have earned him not a little credit if they had been correct. This includes a number of the errors mentioned above,\footnote{Erroneous theory of planetary speeds (\S\ref{planetaryspeeds}); erroneous theory of tides (\S\ref{tides}); erroneous claim that the simple pendulum is isochronic (\S\ref{pendulum}); erroneous theory that comets are an atmospheric phenomena (\S\ref{comets}); erroneous notion of horizontal inertia (\S\ref{inertia}); erroneous claim the shape of a hanging chain is parabolic (\S\ref{catenary}).} plus his erroneous claims that the motion of the earth is undetectable by experiment (\S\ref{relativity}), that meteors entering the earth's atmosphere would cool down due to wind (\S\ref{Babylonianeggs}), that centrifugal force is independent of radius (\S\ref{centrifugal}), and that the curve of fastest descent is a circular arc (\S\ref{brachistochrone}).

Posterity has a way of ignoring mistakes. It is dangerous to start with what we know and ask of history only who was the first to say it. Such selective retrospection is bound to reward careless scientists who made a hundred wild guesses instead of those who weight evidence carefully before making any rash judgements. Galileo is indeed excessively and erroneously assertive where he should have been much more cautious and aware of the limitations of his evidence in many cases (notably \S\S\ref{saturn}, \ref{sunspotshelio}, in addition to the often over-assertive way in which he made the erroneous claims just outlined). In this way Galileo is undermining his right to claim credit for the things he did get right: his accounts of his correct discoveries may sound very convincing and emphatic, but knowing that he was equally sure of a long list of errors gives us reason to suspect that some of the things he got right are to some extent guesswork propped up with overconfident rhetoric in the hope that readers will mistakenly think his case is stronger than it is.

Galileo's relation to empirical data is also very problematic. He often reports as his own observations results he did not actually obtain (\S\S\ref{fallandweight}, \ref{gravitationalconstant}, \ref{resistance}, \ref{pendulum}, \ref{saturn}, \ref{sunspotshelio}, and maybe \ref{venus}). Elsewhere he relied on false data when others did not (\S\ref{tides}), lied about when he had made his discoveries (\S\ref{sunspots} and maybe \S\ref{venus}), dismissed correct data generally accepted by others (\S\ref{comets}), ignored experimental refutations of his claims (\S\ref{pendulum}), and dishonestly suppressed negative outcomes of his own observations (\S\ref{doublestar}).

We have seen time and time again that virtually all of ``Galileo's'' achievements were either anticipated or at least made independently by others. To name just the most striking case:
\quote{Let us hypothetically assume that a scholar contemporary to Galileo pursued experiments with falling bodies and discovered the law of fall as well as the parabolic shape of the projectile trajectory, that he found the law of the inclined plane, directed the newly invented telescope to the heavens and discovered the mountains on the moon, observed the moons of the planet Jupiter and the sunspots, that he calculated the orbits of heavenly bodies using methods and data of Kepler with whom he corresponded, and that he composed extensive notes dealing with all these issues. In short, let us assume that this man made essentially the same discoveries as Galileo and did his research in precisely the same way with only one qualification: he never in his life published a single line of it. …\ As a matter of fact, the above description refers to a real person, Thomas Harriot.\footnote{\cite[184]{buttneretalsymp}.}}
Actually these discoveries are not identical with those of Galileo but rather go beyond them, because Galileo never ``calculated the orbits of heavenly bodies using methods and data of Kepler,'' as Harriot did, who was a better mathematician.

In the history of science, virtually any given scientific discovery has a complicated genealogy of precursors and independent discoveries and near-discoveries if one only looks closely enough. We must not get carried away and let this kind of thing take away all the honour of discovery, for then no one would ever be credited with anything. It may seem unfair, therefore, to catalogue so many instances in which Galileo was preempted by others, since hardly any scientist would come out of such an examination without having to share their claims to fame with others.

Nevertheless, the list of independent discoveries by others does tell us something. First of all, it shows that we must reject the simplistic idea of Galileo as ``the father of modern science,'' as if his achievement was somehow singular and different in kind from his contemporaries. On the contrary, Galileo's science was business as usual at the time, which is why so many of his discoveries were made independently by others. Furthermore, telescopic astronomy is somewhat of a unique case, as we have noted. If Galileo had never touched a telescope, the march of science would not have been set back more than a matter of months. Many were doing this kind of work, and the fact that Galileo managed to monopolise almost all the fame for it had very little to do with science and much to do with shrewd jostling. But for authors of story-book history of science it is all too tempting to ignore this and focus the narrative on one convenient protagonist.

It is also instructive to compare Galileo to Kepler in these kinds of terms. We can find independent contemporary discoveries for almost everything Galileo did, but not so for Kepler's achievements, even though many of them are still central in modern science. Harriot was a ``second Galileo'' and you could go on to a third or a fourth stand-in without much loss. It would be much harder to find a ``second Kepler.'' In my view it is not hard to see why: Kepler was an excellent mathematician who worked on difficult things, while Galileo didn't know much mathematics and therefore focussed on much easier tasks. The standard story has it that Galileo's insights were more ``conceptual,'' yet at least as deep as technical mathematics.\footnote{Examples of such allegedly deep conceptual developments are for instance the relation between mathematics and the physical world (\cite[Ch.\ 1]{Gorham} and below), and the conception of velocity as an instantaneous quantity (\cite[Ch.\ 3]{Gorham}).} On this account it is imagined that basic conceptions of science that we consider commonsensical today were once far from obvious: we greatly underestimate the magnitude of the conceptual breakthroughs required for these developments because we are biased our modern education and anachronistic perspective. But if this is true, how come that Galileo's ideas---for all their alleged ``conceptual'' avant-gardism---spontaneously sprung up like mushrooms all over Europe? And how come all of those ideas can easily be explained to any high school student today, if they are supposedly so profound and advanced? The same cannot be said for Kepler's ideas. They were neither simultaneously developed by dozens of scientists, nor can they be taught to a modern student without years of specialised training. Perhaps this contrast between Galileo and Kepler says something about what genuine depth in the mathematical sciences looks like.

\subsection{Reception}

Already in his own day, ``Galileo was a celebrity and a hero, receiving endless praise during his lifetime from admirers, friends, and opponents alike.''\footnote{\cite[389]{camcomp}.} Or so {\it The Cambridge Companion to Galileo} tells us. There is some truth to this. Galileo was, after all, a spirited populariser and polemicist, a lightning rod in the conflict with the church, and a successful self-promotor in connection with the telescope. All of which contributed to make him a household name for reasons not primarily based on scientific achievement.

\begin{figure}[tp]\centering
\includegraphics[width=0.9\textwidth]{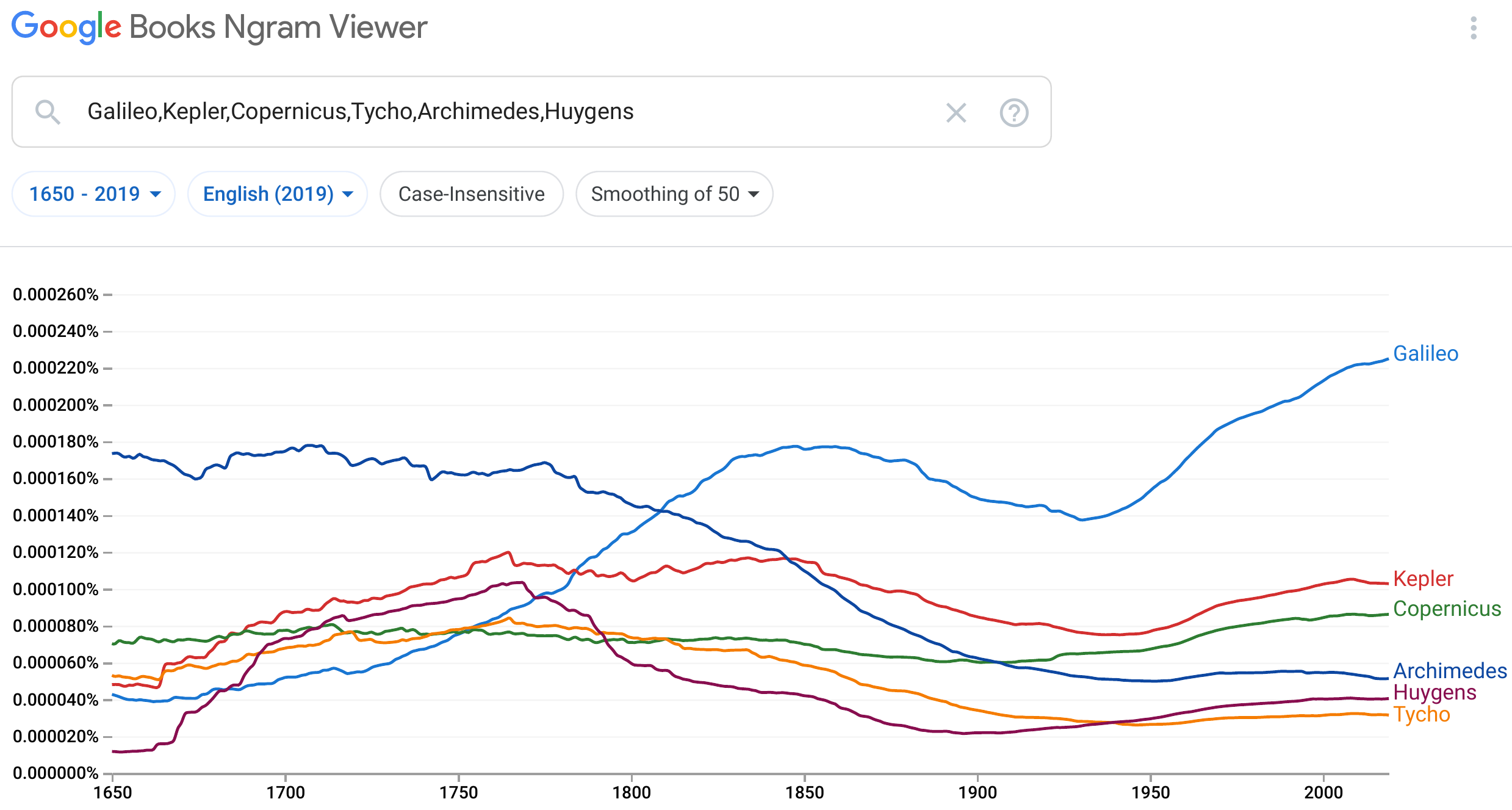}
\caption[Galileo Ngram.]{Mentions of Galileo compared to other scientists in the Google Books corpus.
}
\label{Gngramfig}
\end{figure}

But I challenge the notion that Galileo received ``endless praise'' from all quarters, especially if we focus on the judgements of mathematically competent people. I believe that, as Figure \ref{Gngramfig} suggests, Galileo's fame has snowballed over time, while contemporary scientists, who understood Galileo's work in the context of his day, were much less impressed. Indeed, it is my contention that Galileo erroneously gets credit for insights that were already well established in the Greek mathematical tradition and the works of people like Archimedes. This squares well with the statistical trend that Galileo's prominence stands in inverse proportion to that of Archimedes. Archimedes loomed large in scientific consciousness in the 17th century, but as his works moved from scientists' desks to dusty shelves of ancient history, people began to forget what had been obvious to those who had studied him with care. Meanwhile, Galileo's lively prose was a lot more accessible than abstruse Archimedes. And as knowledge of Archimedes and the true state of science in Galileo's time faded from living memory, readers more readily bought Galileo's self-aggrandising narrative that everyone but him was a foolish Aristotelian.

\subsection{Descartes}
\label{Descartes}

I have already quoted Descartes's opinion that Galileo ``is eloquent to refute Aristotle, but that is not hard.''\footnote{Descartes to Mersenne, 11 October 1638, \cite[390]{drakeGatwork}. \S\ref{mathvsidiots}.} Descartes continued in the same vein: ``I see nothing in his books to make me envious, and hardly anything I should wish to avow mine.''\footnote{Descartes to Mersenne, 11 October 1638, \cite[392]{drakeGatwork}.} Galileo's mathematical demonstrations in particular did not impress Descartes: ``he did not need to be a great geometer to discover them.''\footnote{Descartes to Mersenne, 11 October 1638, \cite[391]{drakeGatwork}.} These were not empty words. Descartes backed up his assessment with an extensive laundry list of specific critiques, a number of which we have had occasion to mention above.\footnote{\S\S\ref{bookofnature}, \ref{projectilemotion}, \ref{ballistics}, \ref{Babylonianeggs}.}

Defenders of the standard view of the great Galileo have been puzzled by this and often tried to write it off on the grounds that Descartes was ``perhaps too much [Galileo's] rival to judge his merits quite impartially.''\footnote{\cite[114]{HallRevSci}.} But ``can Descartes' critique of Galileo's last work on motion really be discarded as the envious comment of a stubborn philosopher, unable or unwilling to acknowledge that this work inaugurated a new physics?''\footnote{\cite[264]{LimitsPreclasMech}.}
More recent scholarship says: no. ``Descartes was not only right in asserting that the theory of motion presented in the {\it Discorsi} did not cohere but also in claiming that some of its foundational concepts were questionable.''\footnote{\cite[267]{LimitsPreclasMech}.}

So Descartes---one of Galileo's most mathematically competent contemporaries---thought his work was quite useless. And he supported this assessment with compelling and perceptive arguments. “Descartes’ judgement of Galileo’s mechanics deserves more sympathy than it usually receives.”\footnote{\cite[25]{moodygolino}.}

\subsection{Kepler}

Kepler was the best mathematical astronomer in Galileo's day. What was his relation to Galileo? Certainly not as substantive as one might expect.
\quote{One wonders why these two great men, who were both present and actual participants at the very birth of some of the most world-shaking scientific events, and who apparently were very much in accord in their astronomical views, did not engage in a more on-going correspondence over these years.\footnote{\cite[329]{PostlCorrKG}.}}
This is a puzzle and a paradox if one accepts the standard view of Galileo. But of course it becomes perfectly understandable as soon as one realises that Kepler, who was a brilliant mathematician, had very little to learn from a dilettante such as Galileo.

Their correspondence began when, ``in 1597, as a lowly high-school teacher of mathematics and a fledgling author, Kepler …\ vainly implored Galileo, the established university professor, to give him the benefit of a judgment of his first major work.''\footnote{\cite[263]{RosenCorrGK}. The work in question is Kepler's {\it Mysterium Cosmographicum} (1596).} Galileo replied briefly, declaring himself in agreement with the Copernican standpoint of Kepler's book, although ``I have preferred not to publish, intimidated by the fortune of our teacher Copernicus, who though he will be of immortal fame to some, is yet by an infinite number (for such is the multitude of fools) laughed at and rejected.''\footnote{Galileo to Kepler, 4 August 1597, \cite[41]{drakeGatwork}.} Kepler is happy to hear that Galileo, ``like so many learned mathematicians,'' has joined in supporting ``the Copernican heresy.''\footnote{Kepler, \cite[40--41]{BaumgardtKepler}.} Galileo should grow a spine though, ``for it is not only you Italians who do not believe that they move unless they feel it, but we in Germany, too, in no way make ourselves popular with this idea.''\footnote{Kepler to Galileo, 13 October 1597, \cite[41]{BaumgardtKepler}.} Kepler urges Galileo to focus on compelling mathematics instead of on the number of fools: ``Not many good mathematicians in Europe will want to differ from us; such is the power of truth.''\footnote{Kepler to Galileo, 13 October 1597, \cite[113]{heilbron}.}

At this point, Kepler naively mistook Galileo for a serious scientific interlocutor. In connection with their discussion of Copernicanism, Kepler noted the importance of parallax and asked Galileo if he could help him with observations for this, adding detailed instructions regarding the exact nature and timing of the requisite measurements.\footnote{Kepler to Galileo, 13 October 1597, \cite[42]{BaumgardtKepler}.} Kepler also sent additional copies of his book, as Galileo had requested, and ``asked only for a long letter of response as payment---which was, however, never forthcoming.''\footnote{\cite[326]{PostlCorrKG}. OGG.X.71.} Galileo stopped replying, presumably since this kind of actual, substantive mathematical astronomy was beyond his abilities.

Kepler's was not the only scientific correspondence Galileo shrunk from. He also neglected to reply to all three letters he received from Mersenne, offering only ``the rather limp excuse that he found Mersenne's handwriting too hard to read.''\footnote{\cite[744]{lewismersenne}, ``scrivere in carattere intelligibile,'' OGG.17.370.} It seems he had a point, for others complained similarly of Mersenne's letters that ``his hande is an Arabicke character to me.''\footnote{Charles Cavendish, \cite[72]{Halliwell}.} Nevertheless these are further instances of Galileo failing to reply to a serious scientific interlocutor.

The tables were turned in 1610. While Galileo had not seen the greatness in Kepler's book, more mathematically competent people had, and consequently Kepler had succeeded Tycho Brahe as the Imperial Mathematician of the Holy Roman Emperor. ``In that capacity Kepler's help was sorely needed by Galileo in 1610, when his momentous telescopic discoveries were being received on all sides with skepticism and hostility.''\footnote{\cite[263]{RosenCorrGK}.} ``To Kepler's credit …\ he manfully swallowed his justifiable resentment'' and ``ungrudgingly gave Galileo the authoritative support he could find nowhere else.''\footnote{\cite[263, 264]{RosenCorrGK}.} ``In spite of Galileo's earlier silence after his own request in 1597, Kepler quickly and enthusiastically responded to Galileo's findings, within 11 days.''\footnote{\cite[327]{PostlCorrKG}.} Galileo surely had this in mind when, in reply, he praised Kepler for ``your uprightness and loftiness of mind''---``you were the first one, and practically the only one, to have complete faith in my assertions'' regarding the telescopic discoveries.\footnote{OGG.X.421. \cite[264]{RosenCorrGK}.} Kepler's support was indeed crucial, and Galileo keenly flaunted it to his advantage.\footnote{\cite[60]{GalileoDiscOp}.}

Galileo did not take the occasion to revive their scientific discussion or comment on Kepler's brilliant new book, the {\it Astronomia Nova} (1609). Instead he only wanted to make fun of dumb philosophers:
\quote{Oh, my dear Kepler, how I wish that we could have one hearty laugh together! Here at Padua is the principal professor of philosophy, whom I have repeatedly and urgently requested to look at the moon and planets thorough my glass, which he pertinaciously refuses to do. Why are you not here? What shouts of laughter we should have at this glorious folly!\footnote{Galileo to Kepler, 1610. \cite[66--67]{BurttMetaph}.}}
Kepler wasn't there because he was busy doing real science. He ignored idiotic philosophers, as all mathematically competent people had done for thousands of years. Galileo, however, had nothing better to do than to sit around and laugh at idiots. To him, it seems, the most desirable application of science is a clever put-down and the last laugh.

Kepler eventually grew weary of Galileo's dilettantism. When, in later years, he found himself having to correct errors in Galileo's superficial writings, he fully justifiably took a patronising tone:
\quote{Galileo rejects Tycho's argument that there are no celestial orbs with definite surfaces because there are no refractions of the stars. …\ Rays reach the earth perpendicular to the spheres, says Galileo, and perpendicular rays are not refracted. But oh, Galileo, if there are orbs, it is necessary that they be eccentric. Therefore, no rays perpendicular to the spheres reach to the earth except at apogee and perigee. Hence, Tycho's argument is a strong one, if you are willing to listen.\footnote{Kepler, appendix to {\it Hyperaspistes} (1625), \cite[350]{DrakeOMalley}.}}
\quote{Galileo denies that the Ptolemaic hypothesis could be refuted by Tycho, Copernicus, or others, and says that it was refuted only by Galileo through the use of the telescope for observation of the variation of the discs of Mars and Venus. …\ Nothing is more valuable than that observation of yours Galileo; nothing is more advantageous for the advancement of astronomy. Yet, with your indulgence, if I may state what I believe, it seems to me that you would be well advised to collect those thoughts of yours that go wandering from the course of reason and memory in that vastness of many interrelated things. This observation of yours …\ does not refute the very distinguished system of Ptolemy nor add to it. Indeed, this observation of yours refutes not the Ptolemaic system but rather, I say, it refutes the traditions of the Ptolemaics regarding the least difference of planetary diameters. …\ Your own observation of the discs confirms the proportion for the eccentric to the epicycle in Ptolemy, as it does the orbit of the sun in Tycho or of the {\it orbis magnus} in Copernicus.\footnote{Kepler, appendix to {\it Hyperaspistes} (1625), \cite[344--345]{DrakeOMalley}. The passage criticised is Galileo, {\it Assayer} (1623), \cite[257--258]{DrakeOMalley}.}}
In both cases, Kepler exposes Galileo's true colours. Galileo doesn't treat the matter as a serious mathematical astronomer, but rather as a superficial and unscrupulous rhetorician. Kepler is right to scold him as he does.

\subsection{Huygens}

Huygens was perhaps the greatest physicist of the generation between Galileo and Newton. He is often portrayed as continuing the scientific program of Galileo. But, in the 22 thick volumes of his collected works, one searches in vain for any strong praise of Galileo, let alone anything remotely like calling him a ``father of science.'' The closest Huygens ever gets to mentioning Galileo favourably is in the context of a critique of Cartesianism.\footnote{Huygens, 1693, {\it Oeuvres} X.403--406.} In the late 17th century, the teachings of Descartes had attracted a strong following. In the eyes of many mathematicians, the way Cartesianism had become an entrenched belief system was uncomfortably similar to how Aristotelianism had been an all too dominant dogma a century before. Huygens makes the parallel explicit:
\quote{Descartes …\ had a great desire to be regarded as the author of a new philosophy …\ [and] it appears …\ that he wished to have it taught in the academies in place of Aristotle. …\ [Descartes] should have proposed his system of physics as an essay on what can be said with probability. …\ That would have been admirable. But in wishing to be thought to have found the truth, …\ he has done something which is a great detriment to the progress of philosophy. For those who believe him and who have become his disciples imagine themselves to possess an understanding of the causes of everything that it is possible to know; in this way, they often lose time in supporting the teaching of their master and not studying enough to fathom the true reasons of this great number of phenomena of which Des Cartes has only spread idle fancies.\footnote{Ibid. \cite[98]{westmanHuygens}.}}
It is in direct contrast with this that Huygens slips in a few kind words for Galileo:
\quote{[Galileo] had neither the audacity nor the vanity to wish to be the head of a sect. He was modest and loved the truth too much.\footnote{Ibid. \cite[97]{westmanHuygens}.}}
Historians have observed that Huygens in all likelihood quite consciously intended this passage to apply to himself as much as to Galileo.\footnote{\cite[32]{DugasH}, \cite[97]{westmanHuygens}, \cite[247--248]{FJDLenses}, \cite[xi]{AndriesseH}.} Perhaps this is why Huygens is surely too generous in praising Galileo's alleged ``modesty.''

In any case, it is very interesting to see what Huygens says about Galileo's actual science in this passage. Let us read it, and keep in mind that this is as close as Huygens ever gets to praising Galileo, and that the context of the passage---a scathing condemnation of Cartesianism---gives Huygens a notable incentive to put Galileo's scientific achievements in the most positive terms for the sake of contrast. In light of this, Huygens's ostensible praise for Galileo is most remarkable, I think, for how  qualified and restrained it is. Here is how I read it:
\quote{Galileo had, in spirit and awareness of mathematics, all that is needed to make progress in physics $\ldots$} Meaning: He said all the right things about about mathematics and scientific method, but he didn't actually carry through on it. Given his rhetoric, he ought to have been able to do it, but be didn't. \quote{$\ldots$ and one has to admit that he was the first to make very beautiful discoveries concerning the nature of motion $\ldots$} He wasn't the first, as we now know, but even though Huygens is overly generous his formulation is still very restrained: ``one has to admit'' ({\em il faut avouer})---a phrase that suggests reluctance to concede the point. Who speaks of their greatest hero in such terms? One ``has to admit'' that he made some discoveries? That seems more like the kind of phrasing you use to describe the work of someone who is overrated, not someone you esteem as the founder of science. \quote{$\ldots$ although he left very considerable things to be done.\footnote{Ibid.}} Exactly. What is most striking and remarkable about the work of Galileo is not the few discoveries he ``admittedly'' made, but how very little he actually accomplished despite all his posturing about mathematics and scientific method. It seems to me that Huygens and I agree on this.

\subsection{Newton}
\label{Newton}

Newton famously said that ``if I have seen further it is by standing on the shoulders of giants.''\footnote{Newton to Hooke, 5 February 1676, {\it The Correspondence of Isaac Newton} I.416.} Many have erroneously assumed that Galileo was one of these ``giants.''\footnote{\cite[106]{BrakeRev}.} One scholar even proposes to explain that ``when Newton credits Galileo with being one of the giants on whose shoulders he stood, he means …''\footnote{\cite[5]{pitt}.} We do not need to listen to what this philosopher thinks Newton meant, because the first part of the sentence is false already. The assumption that Galileo was one of the scientific giants in question has no basis in fact.

The closest Newton gets to praising Galileo is in the {\it Principia}. After introducing his laws of motion, Newton adds some notes on their history.
\quote{The principles I have set forth are accepted by mathematicians and confirmed by experiments of many kinds. …\ By means of the first two laws and the first two corollaries Galileo found that the descent of heavy bodies is in the squared ratio of the time and that the motion of projectiles occurs in a parabola.\footnote{\cite[424]{newtonprincipiacohened}. Note that Newton doesn’t say Galileo was the discoverer of these laws. Indeed, “Newton’s Latin contains some ambiguity” for it “can have two very different meanings: that the two laws were completely accepted by Galileo before he found that projectiles follow a parabolic path, or that these two laws were already generally accepted by scientists at the time that Galileo made his discovery of the parabolic path.” \cite[xxxviii]{cohenNlaws}.}}
The laws and corollaries in question are: the law of inertia, which Galileo did not know;\footnote{\S\ref{inertia}.} Newton's force law $F=ma$, which Galileo also did not know; and the composition of forces and motions, which was established in antiquity.\footnote{The Pseudo-Aristotelian {\it Mechanics} has a law for the composition of motions that ``has all the important features of the Parallelogram Rule'' of Newtonian physics \cite[162]{Miller2017}.}

Of course, once you are looking at the world though Newtonian mechanics it is natural to think that surely Galileo must have had these laws, because that is so obviously the right way to think about parabolic motion. Hence, ``[according to] the myth in which he appears as the founder of classical dynamics, …\ [Galileo] must surely have known the proportionality of force and acceleration. …\ But to those who have become acquainted with Galileo through his own works, not at second hand, there can be no doubt that he never possessed this insight.''\footnote{\cite[344]{DijksterhuisMech}.} Indeed, “Newton almost certainly did not read [Galileo’s] {\it Discorsi}---if, indeed, he ever did---until some considerable time after he had published the {\it Principia}.”\footnote{\cite[xxvi]{cohenNlaws}.}

``Hence Newton (rather too generously, for once!) allowed to Galileo the discovery of the first two laws of motion.''\footnote{\cite[108]{HallRevSci}.} The reason for Newton's excessive charity is not hard to divine. Newton’s {\it Principia} is marked by an obvious and vehement “anti-Cartesian bias.”\footnote{\cite[46]{newtonprincipiacohened}.} “Because of his strongly anti-Cartesian position, Newton might …\ have preferred to think of Galileo rather than Descartes as the originator of the First Law.”\footnote{\cite[xli]{cohenNlaws}.} Whereas, “in point of fact, the {\it Prima Lex} [i.e., the law of inertia] of Newton’s {\it Principia} was derived directly from the {\it Prima Lex} of Descartes’s {\it Principia}.”\footnote{\cite[xxvii]{cohenNlaws}.}

Clearly, then, Newton's attribution of these laws to Galileo means next to nothing. In fact, there is further evidence that it was not meant as high praise in any case. For when Newton continues his historical discussion he says on the very same page: ``Sir Christopher Wren, Dr.\ John Wallis, and Mr.\ Christiaan Huygens, {\em easily the foremost geometers of the previous generation}, independently found the rules of the collisions and reflections of hard bodies.''\footnote{\cite[424]{newtonprincipiacohened}. Emphasis added.} So evidently Newton was of a mind to point out who ``the foremost geometers'' of the past were, yet he had no such words for Galileo---a telling omission. Altogether there is no evidence that Newton regarded Galileo particularly highly, let alone considered him anywhere near a ``father of modern science.''

\section{Before Galileo}

My polemic against Galileo is over. I now wish to consider the broader implications of the story I have told. There are much more important matters at stake than deciding whether Galileo was smart or dumb. The traditional picture of Galileo is the linchpin of an entire historical worldview. When we pull the rug underneath him, a cascade of misconceptions come crashing down.

Galileo’s status stands and falls with our willingness to accept radical relativism. His discoveries are so basic and obvious that the only way to consider them profound is to maintain that they were once {\em not} basic and obvious. In other words, that they are fundamentally different in character from anything the Greeks were doing, for example. Believing in Galileo's greatness means believing that the history of science is a story of ``conceptual'' revolutions that made previously unimaginable things suddenly obvious. It means believing that basic principles of scientific method that seem so obvious today were in fact once completely outside the cognitive universe of even extremely sophisticated mathematical scientists like Archimedes. It means believing that first-rate mathematicians wasted enormous efforts on specialised technical work when vastly greater advances were to be had by simply postulating a few basic philosophical ideas.

The idea that modern science was born in a Galilean revolution is thus based on seeing history as soaked in cultural relativism and replete with dramatic Gestalt shifts. This is the worldview and historiographical approach of many who are far removed from mathematics. Mathematicians find this hard to stomach. They are more ready to say: There is a spiritual unity of scientific thought from ancient to modern times. Great minds think alike. What is obvious to us was obvious to the Greeks. It is ludicrous to think that generations of Greek mathematical geniuses of the first order, with their extensively documented interest science, all somehow failed to conceive basic principles of scientific method. Such is the historiographical outlook of mathematicians.

How did modern science grow out of mathematical and philosophical tradition? The humanistic perspective is that science needed both: it was born through the unification of the technical but insular know-how of the mathematicians with the conceptual depth and holistic vision of the philosophers. The mathematical perspective is that science is what the mathematicians were doing all along. Science did not need philosophy to be its eye-opener and better half; it merely needed the philosophers to step out of the way and let the mathematicians do their thing.

This is why Galileo is the idol of the humanists and the bane of the mathematicians. The philosophers say he invented modern science; the mathematicians that he's a poor man's Archimedes. The issue cuts much deeper than merely allotting credit to one century rather than another. Much more than a question of the detailed chronology of obscure scientific facts, it is a question of worldview and how one should approach and understand history.

If we think there is only one common sense, and that mathematical truth and thought is the same for everyone, then we are strongly inclined to see Galileo’s achievements as trifling. On the other hand, if we reject the very notion of a universal scientific common sense, then we are primed to think that Galileo opened up an entirely new world with his style of science, and that the Greeks couldn’t even think such things, because the way they approached the world was just inherently and profoundly different from ours.

So studying Galileo is a mirror to much larger questions. Either you are a cultural relativist and you think Galileo was a revolutionary, or you think mathematical thought is the same for you, me and everybody who ever lived, and then you think Galileo was just doing common-sense things. Those are the two possibilities. You have to pick sides. You can’t mix and match. You can’t have both mathematical universalism {\em and} Galileo being a revolutionary. The two contradict one another.

\subsection{Plato}
\label{Plato}

Plato is sometimes seen as the mathematician's philosopher. According to legend, an inscription above the entrance to his Academy admonished: ``Let no one enter here who is ignorant of geometry.'' Indeed, Plato founded his entire epistemology on the example of geometry,\footnote{{\it Meno}, 82--85.} and speculated at some length about the mathematical design of the universe.\footnote{{\it Timaeus}, 31--32, 55--61.}

But Plato was no mathematician. Many misconceptions about Greek mathematics and science stem from mistakenly assuming that he speaks on behalf of these fields. The image of mathematics conveyed in Plato's works has aptly been called ``Mathematics: The Movie''\footnote{\cite[215]{NetzGroupP}.}---that is, a dumbed-down, vulgarised page-turner, rich in grandiose sentimentality and sanctimonious moralising that strikes anyone familiar with the real thing as irresponsibly simplistic.

In Plato, mathematics is purer than snow. To apply it to the physical world is to defile it. Many modern scholars assume this was the view of Greek geometers generally, but in fact even Plato's own words show that this was clearly not the case:
\quote{No one with even a little experience of geometry will dispute that this science is entirely the opposite of what is said about it in the accounts of its practitioners. …\ They give ridiculous accounts of it, …\ for they speak like practical men, and all their accounts refer to doing things. They talk of ``squaring,'' ``applying,'' ``adding,'' and the like, whereas the entire subject is pursued for the sake of knowledge …\ [and] for the sake of knowing what always is, not what comes into being and passes away.\footnote{Plato, {\it Republic}, VII, 527, \cite[1143]{PlatoWorks}.}}
So Plato's view of mathematics is in fact, by his own admission, in direct and explicit opposition to the view of actual mathematicians.
\quote{It is evident that Plato's role [in the development of mathematics] has been widely exaggerated. His own direct contributions to mathematical knowledge were obviously nil. …\ The exceedingly elementary character of the examples of mathematical procedures quoted by Plato and Aristotle give no support to the hypothesis that [mathematicians] had anything to learn from [them].\footnote{\cite[152]{NeugebauerExactSci}.}}
Mathematicians do not pay any attention to philosophers trying to tell them about their own subject using one or two basic examples---not then, not in the 17th century, and not today.\footnote{\S\ref{mathvsidiots}.}

In keeping with his praise for the abstraction of mathematics, Plato looked down on empirical science and satirised it with no little scorn: ``Birds …\ descended from …\ simpleminded men, men who studied the heavenly bodies but in their naiveté believed that the most reliable proofs concerning them could be based upon visual observation,''\footnote{Plato, {\it Timaeus}, 91d, \cite[1290]{PlatoWorks}.} leading them to grow wings in order to be able to look at the heavens more closely. Focussing on empirical data is for unphilosophical beasts, in other words. Such an attitude is quite an obstacle to science. If you glorify pure and abstract thought as the only worthwhile pursuit of rational beings, and deride empiricism as fit only for brutes, then you’re not going to get a whole lot of science done.

But again there is very little evidence that mathematically competent people and working scientists ever shared this view, and plenty of evidence to the contrary. This would be evident to all ``if modern scholars had devoted as much attention to Galen or Ptolemy as they did to Plato and his followers.''\footnote{\cite[152]{NeugebauerExactSci}. See also \cite[194]{RussoSciRev}.} Plato is advancing his personal ideology, not describing Greek science. Indeed, ``the retarding effect that Platonism could thus exert on science …\ is encountered particularly during periods in which disparagement of the empirical study of nature on philosophical grounds was supported by a contempt for the material world from religious motives.''\footnote{\cite[15]{DijksterhuisMech}.} In other words, it is a view that has often appealed to philosophers and ideologues for reasons external to science itself.

\subsection{Aristotle}
\label{Aristotle}

\quote{Dante called Aristotle ‘the Master of those who know’. Aristotle was so regarded by learned men from the time of Aquinas to that of Galileo. If one wished to know, the way to go about it was to read the texts of Aristotle with care, to study commentaries on Aristotle in order to grasp his meaning in difficult passages, and to explore questions that had been raised and debated arising from Aristotle’s books.\footnote{\cite[1]{VShortGalileo}.}}
Yes, but note well the crucial qualifier: {\em from the time of Aquinas}, not from the time of Aristotle himself. The intellectual quality of the European Middle Ages was indeed so low that subservience to Aristotle was what passed for erudition. But in the far more advanced intellectual culture of ancient Greece people were not so foolish. Theophrastus---who ``was head of the Lyceum for some 36 years after Aristotle's death''---was ``highly critical of Aristotle, both of his specific physical theories and of his general doctrine of causation.''\footnote{\cite[8--9]{lloydGSAA}.} Theophrastus was succeeded by Strato, whose ``position is [even] further from Aristotle'': he ``rejected many of Aristotle's ideas'' and ``broke new ground in his attempts to investigate problems in dynamics and pneumatics experimentally.''\footnote{\cite[19]{lloydGSAA}.} And these were Aristotle's immediate successors at his own Lyceum. Not even they were particularly attached to his ideas. Mathematicians---who were used to progress in their field, and who were used to accepting propositions based on proof rather than philosophical authority---would have had no reason to adhere to Aristotle's teachings on mechanics, and indeed there is virtually no evidence that they ever did.

It is therefore a mistake to ask the question: ``How is it that the scientific enterprise undertaken by the Greeks, with their unique interest in the rational interpretation of nature, nevertheless culminated in the radically wrong natural philosophy of Aristotle?''\footnote{\cite[243]{CohenSciRev}.} It didn't. Aristotle was not the ``culmination'' of Greek thought. He was one particular philosopher who lived well before the true flourishing of Greek science and who didn't know any mathematics to speak of. The notion of taking Aristotle as their master in mathematical sciences such as mechanics and astronomy was laughable to mathematically competent people in Hellenistic antiquity.
\quote{Aristotle, for whom the brain had a cooling function, could
not possibly have enjoyed excessive authority in the eyes of Herophilus and his disciples, the founders of neurophysiology, nor in those of Archimedes and Ctesibius, who had designed machines that could perform operations whose impossibility Aristotle had “demonstrated”. [There is] an analogous supersession of Aristotle in Aristarchus’ heliocentric theory [and many other domains.] …\ It is clear that the “excessive authority of Aristotle” applies only to later ages and is often backdated.\footnote{\cite[233]{RussoSciRev}.}}
Altogether the spectres of Plato and Aristotle have ruined historical understanding of Greek science---the former by falsely portraying it as adverse to empiricism, the latter by giving the false impression that it treated motion and mechanics only dogmatically and through a qualitative, philosophical lens. ``Negative assessments …\ of ancient Greek mechanics need to be reconsidered. …\ The common view …\ is based partly on a misguided Platonizing tendency …, and partly on a mistaken view that the tradition Galileo is rejecting goes back to ancient Greek thought.''\footnote{\cite[176]{Berryman2009}.}

\subsection{Golden age}
\label{ArchHydrostatics}

Plato and Aristotle were not the pinnacle of Greek thought, as people who read too much philosophy and not enough mathematics are inclined to believe. In the century or two after their death, the Greeks made enormous strides in mathematics and science. This is the age that gave us excellent mathematicians and scientists like Euclid, Aristarchus, Archimedes, Eratosthenes, Apollonius, and Hipparchus, all of whom were miles ahead of Galileo in mathematical sophistication. We shall have occasion to mention their work in suitable contexts below. For now, let us take just one example as emblematic of them all.

\begin{figure}[tp]\centering
\includegraphics[width=0.6\textwidth]{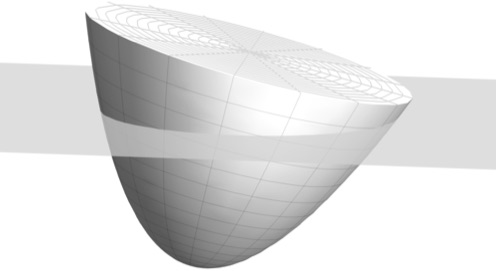}
\caption[Paraboloid submerged in water.]{A paraboloid submerged in water. Archimedes derived exact, detailed, quantitative, highly non-trivial, empirically verifiable results regarding the floatation behaviour and stability conditions of such paraboloids in terms of their shape, tilt, and density.
}
\label{ArchParaboloidFig}
\end{figure}

Archimedes is the quintessential Hellenistic scientist: a first-rate mathematician who was also deeply immersed in practical engineering. That this combination gave rise to brilliant science is proved perhaps most clearly by Archimedes’s work of floating bodies (Figure \ref{ArchParaboloidFig}). This work ``ranks with Newton’s {\it Principia Mathematica} as a work in which basic physical laws are both formulated and accompanied by superb applications,'' namely a detailed investigation of the floatation behaviour of paraboloids that was ``the standard starting point for scientists and naval architects examining the stability of ships'' still thousands of years later, and that can also be used to explain phenomena such as ``the sudden tumbling of a melting iceberg or the toppling of a tall structure due to liquefaction of the ground beneath it.''\footnote{\cite[32--33]{Rorres}. These applications were not discussed by Archimedes in his treatise.} What more could one ask for? It is quite simply an outstanding masterpiece of science by the standards of any age. Only the mathematically illiterate could fail to grasp its immense significance---as indeed they have.

Now, poor Archimedes, he is often misunderstood. Many people who don't care for mathematics hardly even know this work exists. But if they do look at it they say: What's this? It's just a bunch of technical geometry about parabolas and such. Archimedes says not a word about any experiment, not a word about any empirical data, nothing about testing his theory. It seems to have very little to do with the real world. It's not science; it's an exercise in pure geometry. The Greeks may have been excellent geometers but they didn't really do science, you see. They were speculative thinkers, philosophers. They were great with abstract stuff but they didn't have the sense to ground their fanciful theories in reality.

That attitude is completely wrong, in my opinion. Ask yourself: What are the odds that Archimedes got his detailed, quantitative theory of floating bodies absolutely spot-on right if he was ignorant of empiricism and experiment and scientific method? Was he just sitting around doings speculative armchair geometry and, whoops, it just happened to come out exactly equal to empirical facts in a range of far from obvious ways? Are we supposed to believe that was just dumb luck? It doesn't make any sense.

I much more plausible interpretation is that of course Archimedes knew about the scientific method. Of course Archimedes tested his theory by experiment. That is obvious from the accuracy of his results. His text doesn't say that because he was too good of a mathematician to think that kind of kid's stuff counted for much of anything. He only published the actual theory, not the obvious tests that any simpleton could do for themselves. Galileo was precisely that simpleton. He spent his whole life spelling out those parts that Archimedes thought were too trivial to mention. People ignorant of Archimedes are readily tricked into thinking that this was somehow profound. But mathematicians know better.

\subsection{Lost and ignored}

\quote{In the historiography of the Scientific Revolution there has been a notorious difficulty with …\ the idea of ``the Archimedean origins of early modern science'' …: If the 16th-century impact of Archimedes' work in mathematical physics was so revolutionary, why, then, was its original impact in the 2nd century B.C. so negligible?\footnote{\cite[277]{CohenSciRev}.}}
A simple rebuttal is: it wasn't. When we study Greek science and mathematics, we must remember that only a fraction of even the very best works of this era have survived. 17th-century mathematicians had boundless respect and admiration for the ancient Greeks, but they never doubted that there were many more ``long buried monuments of geometry in which so many great findings of the Ancients lie with the roaches and worms.''\footnote{Fermat, \cite[119]{MahoneyFermat}.} And they were right.

In the 20th century a few such masterpieces were recovered. In 1906, a work of Archimedes that had been lost since antiquity was rediscovered in a dusty Constantinople library. The valuable parchment on which it was written had been scrubbed and reused for some religious text, but the original could still just about be made out underneath it. ``Our admiration of the genius of the greatest mathematician of antiquity must surely be increased, if that were possible,'' by this ``astounding'' work,\footnote{\cite[10]{HeathMethod}.} which draws creative inspiration from the mechanical law of the lever to solve advanced geometrical problems. If this brilliant work by antiquity's greatest geometer only survived by the skin of its teeth and dumb luck, just imagine how many more works are lost forever.

Also in the 20th century, divers chanced upon an ancient shipwreck, which turned out to contain a complex machine (the so-called Antikythera mechanism). ``From all we know of science and technology in the Hellenistic age we should have felt that such a device could not exist.''\footnote{\cite[60]{SollaPriceSciAm}.} ``This singular artifact is now identified as an astronomical or calendrical calculating device involving a very sophisticated arrangement of more than thirty gear-wheels. It transcends all that we had previously known from textual and literary sources and may involve a completely new appraisal of the scientific technology of the Hellenistic period.''\footnote{\cite[5]{SollaPriceAntikythera}.}

Another area in which the Greeks appear to have been much further ahead than conventional sources would lead one to believe is combinatorics. Of this entire mathematical field little more survives than one stray remark mentioned parenthetically in a non-mathematical work:
\quote{Chrysippus said that the number of intertwinings obtainable from ten simple statements is over one million. Hipparchus contradicted him, showing that affirmatively there are 103,049 intertwinings.\footnote{Plutarch, \cite[281]{RussoSciRev}.}}
``This passage stumped commentators until 1994,'' when a mathematician realised that it corresponds to the correct solution of a complex combinatorial problem worked out in modern Europe in 1870, thereby forcing ``a reevaluation of our notions of what was known about combinatorics in Antiquity.''\footnote{\cite[281]{RussoSciRev}. See also \cite[282]{AcerbiOxHB}.} It is undeniable from this evidence that this entire field of mathematics must have reached an advanced stage, yet not one single treatise on it survives.

These are just a few striking examples illustrating an indisputable point: the Hellenistic age was extremely sophisticated mathematically and scientifically, and we don't even know the half of it. When ``human learning …\ suffered shipwreck,'' ``the systems of Aristotle and Plato, like planks of lighter and less solid material, floated on the waves of time and were preserved,''\footnote{Bacon, {\it Novum Organum} (1620), I.77, \cite[115]{BaconSelPW}.} while treasure troves of much more mathematically advanced works were lost forever.

With so many key works being lost, we are forced to rely on later commentators and compilers for accounts of the works of Hellenistic authors. This is not unlike forming an image of modern science and mathematics from popularisations in the Sunday paper. Such coverage is invariably oversimplified and dumbed-down, reducing the matter to one or two simplistic ideas while conveying nothing whatsoever of the often massive technical groundwork underlying it. Actually this is a misleading analogy. The situation is much worse than this.
\quote{Nearly all that we know on …\ observations and experiments …\ among the Greeks …\ comes from compilations and manuals composed centuries later, by men who were not themselves interested in science, and for readers who were even less so. …\ Even worse, these works were to a great extent inspired by the desire to discredit science by emphasizing the way in which men of science contradicted each other, and the paradoxical character of the conclusions at which they arrived. This being the object, it was obviously useless, and even out of place, to say much about the methods employed in arriving at the conclusions. It suited Epicurean and Sceptic, as also Christian, writers to represent them as arbitrary dogmas. We can get a slight idea of the situation by imagining, some centuries hence, contemporary science as represented by elementary manuals, second- and third-hand compilations, drawn up in a spirit hostile to science and scientific methods. Such being the nature of the evidence with which we have to deal, it is obvious that all the actual examples of the use of sound scientific methods that we can discover will carry much more weight than would otherwise be the case. If we can point to indubitable examples of the use of experiment and observation, we are justified in supposing that there were others of which we know nothing because they did not happen to interest the compilers on whom we are dependent. As a matter of fact, there are a fair number of such examples.\footnote{\cite[253--254]{BurnetEssays}. For further references to the same effect, see \cite[75]{lloydMP}.}}
Sadly, the lack of appreciation for science among these third-rate commentators continues among scholars today.
\quote{More of Greek ``science'' survives than does any other category of Greek literature, and yet much of that is obscure even to classicists.\footnote{\cite[xxi]{HellenisticEraSourcebook}.}}
\quote{The state of editions and translations of ancient scientific works as a whole remains scandalous by comparison with the torrent of modern works on anything unscientific -- about 100 papers per year on Homer, for example. And an embarrassingly large number of classicists are …\ ignorant of Greek scientific works.\footnote{\cite[137]{rihllGS}.}}
\quote{Classicists include many who have chosen Latin and Greek precisely to {\em escape} from science at the very early stage of specialisation that our schools' curricula permit: and often a very successful escape it is, to judge from the depth of ignorance of science ancient {\em and} modern that it often secures.\footnote{\cite[354]{lloydMP}.}}
\quote{Modern scholars, persistently regarding this era as somehow inferior to the Athens of Perikles and Demosthenes, have often disregarded Hellenistic science and technology in favor of later Roman achievements or earlier Greek work.\footnote{\cite[243]{CambridgeHellenisticWorld}.}}
Little wonder then that Greek science is systematically misunderstood and undervalued, and that simplistic ideas of philosophical authors and commentators are substituted for the real thing.

\subsection{Greek mechanics}
\label{Greekmech}

``The study of mechanics …\ is eagerly pursued by all those interested in mathematics,''\footnote{Pappus, {\it Collection}, VIII.1, \cite[91]{lloydGSAA}.} Greek sources credibly report. Perhaps even Aristotle himself wrote on technical mechanics. At any rate a work on this subject survived under his name, which Galileo and his contemporaries still regarded as a work by Aristotle. In more recent times, ``the attribution to Aristotle has been questioned mainly on the ground that the treatise's attention to practical problems is `quite un-Aristotelian', which is doubtful reasoning at best.'' In any case, it is ``agreed by those who would question its attribution to Aristotle that the treatise was composed by an Aristotelian shortly after Aristotle's time.''\footnote{\cite[4]{ClagettMA}.} The authors of the treatise, whoever they are, ``discuss the lever, the pulley, and the balance, and expound with considerable success some of the main principles of statics---the law of virtual velocities, the parallelogram of forces, and the law of inertia.''\footnote{\cite[12]{RossAristotle}.} Quite serious and quantitative mechanics, in other words. And this is near the time of Aristotle, still before the golden age of Greek science.

Doubtless the Hellenistic age saw great improvements on these respectable beginnings. In connection with the ``Galilean'' theory of falling bodies and projectiles, it is intriguing that both Strato and Hipparchus wrote treatises on falling bodies that are now lost.\footnote{\cite[70--71]{ClagettGS}.} It is squarely within the realm of possibility that these may have contained most or all of ``Galileo's'' discoveries.

Strato had ``an apparently deserved reputation as an experimenter.''\footnote{\cite[69]{ClagettGS}.} To prove that falling objects speed up, he reasoned as follows.\footnote{\cite[16]{lloydGSAA}.} Pour water slowly from a vessel. At first it flows in a continuous stream, but then further down its fall it breaks up into drops and trickles. This is because the water is speeding up. So the water spreads out, like cars let loose on a highway after a congested area. Another experiment proving the same point is stones dropped into a sand bed from various heights. The stone makes craters of different depths depending on height fallen. These are the kinds of basics preserved in the superficial commentary literature. Quite possibly the original treatise backed these things up with a more mathematical treatment.

Hipparchus was the greatest mathematical astronomer of the Hellenistic era and certainly more than capable of giving a mathematically sophisticated treatment of falling bodies. Later commentators tell us that “Hipparchus contradicts Aristotle regarding weight, as he says that the further something is, the heavier it is.”\footnote{Simplicius, \cite[291]{RussoSciRev}.} ``The only way to make [this statement] comprehensible is to suppose that Hipparchus meant the weight of bodies inside the earth, recognizing that it decreases as the body nears the center,''\footnote{\cite[293]{RussoSciRev}.} which is in accordance with modern gravitational theory. Another commentator evidently had this in mind when describing the question posed by “the folks who introduced the thrust toward the center” as to whether “boulders thrust through [a tunnel into] the depths of the earth, upon reaching the center, should stay still with nothing touching or supporting them; [or whether] if thrust down with impetus they should overshoot the center and turn back again and keep bobbing back and forth.”\footnote{Plutarch, \cite[293]{RussoSciRev}.}

These things are very much in line with 17th-century physics. Dropping stones into sand, thinking about how gravity varies on a super-terrestrial scale and inside the earth:  scientists spent a lot of serious effort on exactly that in the 17th century. All in all, it is certainly possible that Hipparchus and his contemporaries were familiar with principles of ``Galilean'' mechanics. As one modern scholar has observed: ``The ease of stumbling upon this discovery renders it highly improbable that natural philosophers had ever searched for the law of fall'' before Galileo.\footnote{\cite[2]{drakeHistFall}.} Perhaps a more natural conclusion from this logic is that Hellenistic scientists in fact did search for it and in fact did stumble upon it in those lost works on falling bodies that we know they wrote.

\subsection{Greek astronomy}
\label{Greekastro}

Ptolemy is the canonical source for Greek astronomy, and the target of Galileo's attacks. But it is quite plausible that Ptolemy, who lived hundreds of years after the golden age of Greek science, was not the pinnacle of Greek astronomy, but should rather be seen as a regressive later author. Ptolemy's big book possibly did ``more damage to astronomy than any other work ever written'' by displacing much better Hellenistic Greek astronomy, so that ``astronomy would be better off if it had never existed.''\footnote{\cite[379]{Pcrime}.} It is in any case certain that Ptolemy was merely the last in a very long tradition of mathematical astronomy.\footnote{The problem of creating geometrical models for planetary motions goes back at least to Aristotle's time, as is evident from e.g.\ {\it Metaphysics} XII.8.} Sadly, virtually none of the technical works on astronomy before Ptolemy have survived.

We know for a fact that heliocentrism was pursued very seriously in the Hellenistic era, long before Ptolemy. Aristarchus, who was a good mathematician, wrote a now lost treatise arguing that the earth revolves around the sun. Archimedes cites this work with tacit approval.\footnote{\cite[222]{HeathArchimedes}.} There are some indications that others kept pursuing heliocentrism.\footnote{\cite[80--82, 285--286, 294--297]{RussoSciRev}.} Indeed, “several technical elements of Ptolemaic astronomy can only be explained as derivatives of an earlier heliocentric model.”\footnote{\cite[317]{RussoSciRev}, referring to \cite{RawlinsHelio}. See also \cite{ThurstonIsis}, \cite{RawlinsIsis}.}

There are clear indications that Greek heliocentrists supported their theory with physical arguments. First of all it is quite evident that the spherical shape of the earth can be explained as a consequence of gravitational forces. In fact, Archimedes proves as a theorem in his hydrostatics that a spherical shape is the result or equilibrium outcome of basic gravitational assumptions. This idea, that gravitational forces are the cause of the spherical shape of the earth, is explicitly stated in ancient sources.\footnote{\cite[303]{RussoSciRev}.}

But of course other heavenly bodies are round too, such as the moon for example. So that very naturally suggests that they have their own gravity just like the earth. This conclusion too is explicitly spelled out in ancient sources. Thus Plutarch says: “The downward tendency of falling bodies is evidence not of the earth’s centrality but of the affinity and cohesion to earth of those bodies which when thrust away fall back again. …\ The way in which things here [fall] upon the earth suggests how in all probability things [on the moon] fall …\ upon the moon and remain there.”\footnote{Plutarch, \cite[72]{HellenisticEraSourcebook}.}

Now, from this way of thinking, it is a short step to the idea that the heavenly bodies pull not only on nearby objects but also on each other. This is again explicit in ancient sources. This is why Seneca, for example, says that “if ever [these bodies] stop, they will fall upon one another.”\footnote{Seneca, \cite[294]{RussoSciRev}.} That is correct, of course. The planets would “fall upon one another” if it wasn’t for their orbital speed. The Greeks were well aware that tides can be explained in terms of such attractions, as we have seen.\footnote{\S\ref{tides}.}

This point of view explains the motions of the planets in terms of physical forces. It’s not that the planets have circularity of motion as an inherent attribute imbedded in their essence, as Aristotle would have it. Rather, circularity is a secondary effect, the result of the interaction of two primary forces: a tangential force from motion and a radial force from gravity. There are clear indications that ancient astronomers worked out such a theory, including a mathematical treatment. Thus Vitruvius says: “the sun’s powerful force attracts to itself the planets by means of rays projected in the shape of triangles; as if braking their forward movement or holding them back, the sun does not allow them to go forth but [forces them] to return to it.”\footnote{Vitruvius, \cite[297]{RussoSciRev}.} Pliny says the same thing: planets are “prevented by a triangular solar ray from following a straight path.”\footnote{Pliny, \cite[298]{RussoSciRev}.}

All this talk of triangles, in both of these authors, certainly suggests an underlying mathematical treatment. Indeed, the Greeks knew very well the parallelogram law for the composition of forces or displacements, and in fact explicitly used this to explain circular motion as the net result of a tangential and a radial motion.\footnote{\cite[301--302]{RussoSciRev}.}

It is beautiful how coherently all of that fits together and how naturally we were led from one idea to the other. Just like the water of the oceans naturally seeks a spherical shape, so the spherical shape of the earth has been formed by the same forces. And just as gravity explains why the earth is round, so it must explain why other planets are round. Hence they have gravity. But just as they attract nearby objects, so they attract each other. So the heavens have a perpetual tendency to lump itself up, except this tendency is counterbalanced by the tendency of speeding objects to shoot off in a straight line.

Hence top mathematicians in the Hellenistic era advocated heliocentrism, and evidently integrated this view with mechanical considerations comparable in spirit to the works of people like Kepler in the 17th century. Many people refuse to believe this. It has recently been claimed, for example, that “pre-Copernican heliocentrisms (that of Aristarchus, for example) have all the disadvantages and none of the advantages of Copernican heliocentrism,” because they postulated only that the earth revolves around the sun, not, as has commonly been assumed, that all the other planets do so as well. This supposedly “explains why Copernicus’s heliocentrism was accepted …, while pre-Copernican heliocentrism” was not.\footnote{\cite[16]{CarmanCop}.} This is completely wrong, in my opinion. And for an obvious reason. Namely: Why would Aristarchus have affirmed and written a treatise on heliocentrism if it had nothing but disadvantages? What possible reason could he have had done for doing so? None, in fact. Yet this is exactly what this recent article proposes.

It is a fact that Aristarchus asserted the physical reality of his hypothesis. And it is a fact that he recognised the parallax argument against it (discussed in \S\ref{parallax}), as Archimedes implies when he mentions Aristarchus.\footnote{\cite[222]{HeathArchimedes}.} Even the recent article I cited admits this. So why, then, would Aristarchus write a treatise proposing this bold hypothesis, discuss a major argument against it (namely the parallax argument) and no arguments in favour of it, and then nevertheless conclude that the hypothesis is true? And, furthermore, why would Archimedes, perhaps the greatest mathematician of all time, cite this treatise with tacit approval as a viable description of physical reality? None of that makes any sense.

The only reasonable explanation is that Aristarchus did in fact recognise an advantage of placing the sun in the center. And the obvious guess for what this was is that he saw the same advantages as Copernicus did, including the more natural explanations of the retrograde motion of the outer planets and the bounded deviation from the sun of the inner planets.

As we have seen, an important argument for heliocentrism in the 17th century was that its is natural for a smaller body to orbit a larger one, but not the other way around.\footnote{\S\ref{againsttycho}.} Strikingly, Aristarchus wrote a treatise calculating the relative distances and sizes of the sun, the earth, and the moon.\footnote{His only surviving work, {\it On the Sizes and Distances of the Sun and Moon}.} This treatise shows that Aristarchus was at any rate a highly competent mathematician. But I think it shows much more than that. I think it feeds directly into his heliocentrism. Or are we supposed to believe that Aristarchus calculated the sizes of heavenly bodies just for kicks in one treatise and did not see any connection with the heliocentrism he advanced in another treatise even though the obvious connection was right under his nose? What is the probability that he suffered from such schizophrenia? Virtually zero, in my opinion.

In fact there are certain aspects Aristarchus’s treatise that make much more sense when you read it this way. On its own it is a weird treatise. On the one hand it calculates the sizes and distances of the sun, moon, and earth in a mathematically sophisticated manner. Very detailed, technical stuff, including the completely pointless complication that the sun does not quite illuminate half the moon but ever so slightly more than half, since the sun is larger than the moon. This is “pure mathematical pedantry.”\footnote{\cite[643]{NeugebauerHAMA}.} It makes the geometrical calculations ten times more intricate while having only the most minuscule and completely insignificant impact on the final results.

On the other hand, the observational data that Aristarchus uses for his calculations are extremely crude. He says that the angular distance between the sun and the moon at half moon is 87 degrees: a pretty terrible value. The real value is more like 89.9 degrees.\footnote{\cite[642]{NeugebauerHAMA}.} Because of this his results are way off. For instance, his calculated distance to the sun is off by a factor of 20 or so. So what’s going on? Why do such intricate mathematics with such worthless data? Did he just care about the mathematical ideas and not about the actual numbers? I think it would be a mistake to jump to that simplistic conclusion, even though many people have done so.

In fact, it is easy to see how Aristarchus had a purpose in underestimating the angle. His purpose with the treatise, I propose, is to support his heliocentric cosmology based on the principle that smaller bodies orbit bigger ones. This hypothesis fits very well with the structure of Aristarchus’s treatise. The treatise has 18 propositions. Proposition 16 says that the sun has a volume about 300 times greater than the earth, and Proposition 18, the very last proposition, says that the earth has a volume about 20 times greater than the moon. These are exactly the propositions you need to explain which body should orbit which. And that is exactly where Aristarchus chooses to end his treatise.

Many commentators have been puzzled by why Aristarchus ends in that strange place.\footnote{\cite[636]{NeugebauerHAMA}, \cite[221--222]{BerggrenSidoli}.} In particular, many have been baffled by why he does not give distances and sizes in terms of earth radii. This seems like the natural and obvious thing to do, and doing so would have been easily within his reach. Many modern commentators add the small extra steps along the same lines needed to fill this obvious “gap.”\footnote{\cite[222--223]{BerggrenSidoli}, \cite[8--9]{HeldenMeasuringBook}.} Except it’s not a gap at all and there is no need to be puzzled by Aristarchus’s choices. If we accept my hypothesis, everything he does makes perfect sense all of a sudden. He carries his calculations precisely as far as he needs for this purpose, and no further.

My hypothesis also explains why he chose such a poor value for the angular measurement. He has every reason to purposefully use a value that is much too small. Underestimating this angle means that the size of the sun will be underestimated. And his goal is to show that the sun is much bigger than the earth. So he has shown that even if we grossly underestimate the angle, the sun is still much bigger than the earth. So he has considered the worst case scenario for his desired conclusion, and he still comes out on top. That just makes his case all the stronger, of course.

Clearly my interpretation requires that Aristarchus knew that 87 degrees was an underestimate. The standard view in the literature is that he could not have known this. Aristarchus’s numerical data are “nothing but arithmetically convenient parameters, chosen without consideration for observational facts”. “It is obvious that [Aristarchus’s] fundamental idea …\ is totally impracticable. …\ 87 degrees is a purely fictitious number. …\ The actual value …\ must …\ [have] elude[d] direct determination by methods available to ancient observers.”\footnote{\cite[642--643]{NeugebauerHAMA}.} The argument for this is as follows. You are trying to measure the angle between the sun and the moon at half moon. But to do that you need to pinpoint the moment of half moon, which can only be done with an accuracy of maybe half a day. But in half a day the moon has moved six degrees, and therefore radically changed the angle you are trying to measure. So your observational value is going to have a margin of error of 6 degrees, which is enormous and makes the whole thing completely pointless.

But I’m not convinced that it’s as hopeless as all that. One way to work around the problem would be to use not one single observation, but the average of many observations. There is little evidence that the Greeks ever made use of averaging that way,\footnote{\cite[1--2]{Astrostats}.} but the idea is simple enough. I did a bit of statistics to see if this would be viable. Let’s assume that our angular measurements are normally distributed about the true value. And let's accept the assumption in the critique of Aristarchus that one would be lucky to get the moment of half moon correct to the day.\footnote{\cite[642]{NeugebauerHAMA}.} So we can tell it’s today rather than yesterday, but we can’t tell at what exact hour the moon is exactly half full. Let us roughly translate this into statistical terms by saying that the observations have a standard deviation of 12 hours, or six degrees.

Now, an astronomer active for, let’s say, two decades would have occasion to observe about 500 half moons. So say he makes 500 angular measurements and then average them. This would produce an estimate of the true value with a 95\% confidence interval of $\pm$ about half a degree. A margin of error of half a degree is easily enough to support my interpretation that 87 is a conscious underestimate.

Naturally, Aristarchus would not have reasoned in such terms exactly, but it is not necessary to know any statistical theory to get an intuitive sense of the order of magnitude of the error in such an estimate. As you keep adding observations, and keep averaging them, you will see the average stabilising over time, of course. So certainly it will become clear after a while that, whatever the true value is, it must be greater than 87 degrees at any rate.

It is certainly extremely speculative to imagine that Aristarchus might have had something like this in mind. But in any case my argument shows that it cannot be ruled out as out of the question that Aristarchus could in principle have had solid empirical evidence that his value of 87 degrees was certainly an underestimate.

So, in conclusion: Aristarchus was a good mathematician. He proposed a heliocentric theory that was taken very seriously by Archimedes. There was a long tradition in Greek thought of trying to account for the motions of the planets in terms of everyday physics. This is naturally connected to heliocentrism because of the natural idea that smaller bodies orbit bigger ones. Aristarchus in fact wrote a major treatise devoted specifically to comparing the sizes of the sun and the earth, and the earth and the moon. Several otherwise peculiar aspects of the treatise fit like a glove the idea that it was written precisely to lend credibility to heliocentrism.

The specific details of this speculative reconstruction are not important as far as our Galilean story is concerned. But it does indicate the scope of uncertainty that the ancient sources leave us with. For the purposes of a contextualised assessment of Galileo, the important question is whether we should think of the emergence of modern science as a drastic breakthrough and discontinuity. Was Galileo's mode of scientific thought fundamentally different in kind from that of even the best minds of previous eras? Some will have you believe that it is. I disagree. What I have said about the possible contents of the lost scientific works of antiquity is tentative and speculative. I have certainly not proved that the ancients made those discoveries. But I have proved something, namely that the assumption that ancient science was fundamentally different in kind from the ``Scientific Revolution'' of Galileo's age is by no means an established fact. My speculative reconstructions do not show that this assumption is false, but they do show that it is {\em possible} that it is false. This in itself is enough to call into question the routine assumptions often made about how ``everyone'' before Galileo allegedly lacked this or that scientific insight.

To sum up my discussion of Hellenistic science more broadly, it is my contention that there is a strong continuity between the mathematical tradition in ancient science and 17th-century science. It's not for nothing that Archimedes was widely regarded to have been a veritable God among men in the 17th century. Furthermore, a series of first-rate minds in Greek times, working in an extremely sophisticated mathematical culture, wrote numerous scientific treatises that are now lost. It is quite possible that these works could have included the better part of ``Galileo's'' science, for example. The knee-jerk reaction among historians is to object that this is speculation that cannot be directly proven from the sources. If these historians were consistent in applying this principle, however, they would have to reject with equal zeal and fervour any speculation that Aristotle and Plato had any influence whatsoever on the mathematical-scientific tradition in Greek antiquity---a widely accepted view for which there is even less source evidence. In any case we know for a fact that many first-rate works are lost, so stubbornly refusing to extrapolate beyond the available sources is an absolute guarantee that we will be wrong about the scope of Greek science. Finally, to understand the impact of Greek thought on the early modern period, we must realise that 17th-century mathematicians did not share the reluctance of present-day historians to make bold inferences about the store of lost wisdom once known to ``the ancients.''

\section{``Father of Science''}

Galileo ``may properly be regarded as the `father of modern science'.{''}\footnote{\cite[323]{OxCompHistModSci}.} This is still the accepted view among modern historians, as in this quote from the recent {\it Oxford Companion to the History of Modern Science}. This view is considered so unassailable that even the very Pope once conceded that Galileo ``is justly entitled the founder of modern physics.''\footnote{Pope John Paul II, 1979, \cite[198]{poupard}.} There is less agreement on what warrants this epithet. ``No one indeed is prepared to challenge [Galileo's] scientific greatness or to deny that he was perhaps the man who made the greatest contribution to the growth of classical science. But on the question of what precisely his contribution was and wherein his greatness essentially lay there seems to be no unanimity at all.''\footnote{\cite[333]{DijksterhuisMech}.} I shall go though all major attempts at capturing Galileo's alleged greatness, and offer a consistent rebuttal case against them.

\subsection{Mathematics and nature}
\label{mathandnature}

It is a common view that Galileo was the first to bring together mathematics and the study of the natural world.
\quote{The momentous change that Galileo …\ introduced into scientific ontology was to identify the substance of the real world with …\ mathematical entities. The important practical result of this identification was to open the physical world to an unrestricted use of mathematics.\footnote{\cite[310]{CrombieGross}.}}
\quote{Galileo is one of the Founders of Modern Science [primarily on account of] his formulation and defence of the idea of a mathematical method as the most appropriate means for uncovering the secrets of nature. [This was a great achievement because] explaining how mathematics could show anything about the properties of matter [posed] an enormous conceptual difficulty, [since] geometry was recognized to be an abstract discipline, not one appropriate for dealing with physical matters.\footnote{\cite[xii, 5, 6]{pitt}.}}
\quote{Galileo was the ‘first inventor’ of mathematical physics: I mean, of truly modern physics. …\ It was Galileo who, by consistently applying mathematics to physics and physics to astronomy, first brought mathematics, physics, and astronomy together in a truly significant and fruitful way. The three disciplines had always been looked upon as essentially separate: Galileo revealed their triply-paired relationships.\footnote{\cite[xxi, 64--65]{drakeessays1}.}}
\quote{The traditional historical view, which sees Galileo as the father of classical science, is not wrong. For it is in his work …\ that the idea of mathematical physics …\ was realised for the first time in the history of human thought.\footnote{\cite[201]{KoyreGS}.}}
\quote{What was fundamentally original and revolutionary in the conception of Galileo …\ was the assurance that, in principle, the potentialities of mathematical reasoning went far beyond the narrow limits allowed by traditional philosophy. …\ [Thus Galileo was] the founder of the mathematical philosophy of nature.\footnote{\cite[113, 285]{HallRevSci}.}}
This standard view obviously rests on the assumption that, before Galileo, mathematics and natural science were fundamentally disjoint. This assumption is plainly and unequivocally false. In Greek works by mathematically competent authors, there is zero evidence for this assumption and a mountain of evidence to the contrary. ``We attack mathematically everything in nature'' said Iamblichus of Greek science,\footnote{\cite[156]{lloydGSAA}.} and he was right. Far from being unable to conceive the unity of mathematics and science, ``Hellenistic natural philosophers often took mathematics as the paradigm of science and sought to mathematize their study, that is, to ground all its claims in mathematical theorems and procedures, a goal shared by modern scientists.''\footnote{\cite[245]{CambridgeHellenisticWorld}.}

How can so many historians get it exactly backwards? By ignoring the entire corpus of Greek mathematics and instead relying exclusively on philosophical authors. Thus we are told that, following ``the classification of philosophical knowledge deriving from Aristotle,'' a sharp division prevailed among ``the Greeks'' between ``natural science (or `physics'), which studied the causes of change in material things,'' and ``mathematics, which was the science of abstract quantity.''\footnote{\cite[1--2]{CrombieGross}.} This was perhaps a problem for philosophers who spent their time trying to classify scientific knowledge instead of contributing to it. But I challenge you to produce one single piece of evidence that this division had any impact whatsoever on any mathematically creative person in antiquity.

The alleged divide doesn't exist in Aristotle's own works either, for that matter. Aristotle lived well before the glory days of Greek science, and he was clearly no mathematician. But even Aristotle lists mechanics, optics, harmonics, and astronomy as fields based on mathematical demonstrations,\footnote{Aristotle, {\it Posterior Analytics}, 1.9: 76a, 1.13: 78b, {\it Metaphysics}, 13.3: 1078a.} and even explicitly calls them ``branches of mathematics.''\footnote{Aristotle, {\it Physics}, II, 194a8, \cite[239]{AristotleBW}.} How can anyone infer from this that Aristotle saw the very notion of mathematical science as a conceptual impossibility? Remarkably, historians in fact do so, by insisting that these fields are mere exceptions:
\quote{Previous assumptions [before Galileo], encouraged by Aristotle and scholastic philosophers, held that mathematics was only relevant to our understanding of very specific aspects of the natural world, such as astronomy, and the behaviour of light rays (optics), both of which could be reduced to exercises in geometry. Otherwise, mathematics was just too abstract to have any relevance to the physical world.\footnote{\cite[109]{HenryShortHist}.}}
The implausibility of this view is obvious. If, as Aristotle himself clearly states, mechanics, optics, harmonics, and astronomy are four entire fields of knowledge that successfully use mathematics to understand the natural world, who in their right mind would then categorically insist that, nevertheless, {\em other than that} mathematics surely has nothing to contribute to science. It makes no sense. If mathematics has already given you four entire branches of science, why close your mind to the possibility of any further success along similar lines? It is hard to think of any reason for taking such a stance, except perhaps for someone who themselves lack mathematical ability and want to justify their neglect of this field.

The strange habit of writing off the numerous branches of mathematical science in antiquity as so many exceptions is necessary to maintain triumphalist narratives of the great Galilean revolution. For example, we are told that ``it was Galileo who first subjected other natural phenomena to mathematical treatment than the …\ Alexandrian ones.''\footnote{\cite[114]{CohenRiseExpl}.} In other words, except mechanics, astronomy, optics, statics, and hydrostatics,\footnote{\S\ref{ArchHydrostatics}.} Galileo was {\em the very first} to take this step. That is to say, if you ignore all previous mathematicians who did this exact thing in great detail, Galileo's step was revolutionary.

Some have tried to defend the alleged originality of Galileo's step by arguing that, before Galileo, there were various successful applications of mathematics, but what was lacking was the vision of nature herself as inherently mathematical.\footnote{\cite{AlexanderProof}.} But many Greek achievements do not fit that characterisation. The Pythagorean insight that numbers govern musical harmonies can precisely {\em not} be interpreted as a man-made application of mathematics; rather it is clearly and unequivocally an indication that the world is fundamentally, inherently mathematical somehow, in a deep way that is not immediately obvious. The obvious conclusion is the one the Pythagoreans explicitly drew: that ``all is number,'' not that mathematics sometimes happens to be applicable to some limited aspects of nature. In astronomy, the Greeks insisted on explaining everything in terms of uniform circular motions. Aristotle and Plato are both very explicit about this.\footnote{Plato {\it Timaeus} 34a, Aristotle {\it De caelo} II.6.} What is this but an explicit commitment to the universe being inherently mathematical in its very metaphysical essence? In optics, light chooses the shortest path, whence light rays are straight and are reflected with the outgoing angle equal to the incoming angle. Again, not really an application of mathematics but rather an indication that nature herself is doing mathematics to determine the outcome of fundamental processes. Plato's {\it Timaeus} is explicit about nature being mathematically designed, down to its very elements, which are regular polyhedra. Others too saw mathematical design everywhere they looked. Honeycombs are hexagonal because it optimises the amount of area per perimeter among regular polygons, as Pappus says.\footnote{\cite[526--527]{BlasjoIsop}.} Ptolemy infers planetary distances by assuming that the universe was designed to fit the epicyclic models with no waste in between.\footnote{\S\ref{realism}.} And so on. The Greek universe is through and through the work of mathematical demiurges.

Another strategy for explaining away the obvious fact of extensive mathematical sciences in antiquity is to discount them as genuine science on the grounds that they were abstractions. Thus some claim that, despite ostensible applications of mathematics in numerous fields, ``mathematical theory and natural reality remain almost entirely separate entities'' due to the ``high level of abstraction'' of the mathematical theories, which mean that they were ``barely connected with the real world.''\footnote{\cite[19]{CohenRiseExpl}.} Supposedly, Galileo broke this spell---an absurd claim since this critique is all the more true for his science: even Galileo's supposedly ``best'' discoveries are often way out of touch with reality,\footnote{\S\S\ref{fallandweight}, \ref{ballistics}.} not to mention his many erroneous theories. Meanwhile, Greek scientific laws of statics, optics, hydrostatics, and harmonics concern everyday phenomena that can be verified by anyone in their own back yard using common household items. Indeed, they are still part of modern physics textbooks---and high school laboratory demonstrations---to this day. Take optics, for example. Heron of Alexandria proved the law of reflection, which anyone with a mirror can readily check, using the distance-minimisation argument still found in every textbook today.\footnote{\cite[80]{ClagettGS}, \cite[264--265]{CohenDrabkinSourceB}, \cite[195]{HellenisticEraSourcebook}.} Diocles demonstrated the reflective property of the parabola and used it to ``cause burning'' by concentrating the rays of the sun with a paraboloid mirror;\footnote{\cite[189--192]{HellenisticEraSourcebook}.} a principle still widely applied today, for example in satellite dishes and flashlights. Ptolemy demonstrated the magnifying property of concave mirrors, such as modern makeup mirrors.\footnote{\cite[199--200]{HellenisticEraSourcebook}.} These kinds of results, which are not atypical, are clearly not disconnected from reality by any means.

The false notion of a divide between mathematics and science also rests on a conception of mathematics itself as a purely abstract field.
\quote{Traditionally, geometry was taken to be an abstract inquiry into the properties of magnitudes that are not to be found in nature. Dimensionless points, breadthless lines, and depthless surfaces of Euclidean geometry were not traditionally taken to be the sort of thing one might encounter while walking down the street. Whether such items were characterized as Platonic objects inhabiting a separate realm of geometric forms, or as abstractions arising from experience, it was generally agreed that the objects of geometry and the space in which they are located could not be identified with material objects or the space of everyday experience.\footnote{\cite[205]{Jesseph2015}.}}
This is again a view expressed by philosophers only. Nothing of the sort is ever stated by any mathematically competent author in antiquity. On the contrary, mathematicians routinely take the exact opposite for granted. Allegedly ``abstract'' geometry is constantly applied to physical objects in Greek mathematical works without ado.

It is obvious that the long list of Greek mathematicians who studied the natural world always took for granted the identification of geometry with the space and material objects around us. And why shouldn't they? For thousands of years geometry had been used to delineate fields, draw up buildings, measure volumes of produce, and a thousand other practical purposes---exactly ``the sort of thing one might encounter while walking down the street.'' Every single theorem of Euclid's geometry can be verified by concrete measurements and constructions with physical tools and materials. So why would mathematicians suddenly insist that their field is completely divorced from reality? What could possibly be their motivation for doing so? It accomplishes nothing and creates tons of obvious problems when one wants to apply mathematics far and wide in numerous areas, as mathematicians always did. The only people with any motive to take such an extremist stance are philosophers with an axe to grind.

\begin{figure}[tp]
   \centering
   \begin{tabular}{ccc}
   \includegraphics[height=0.15\textwidth]{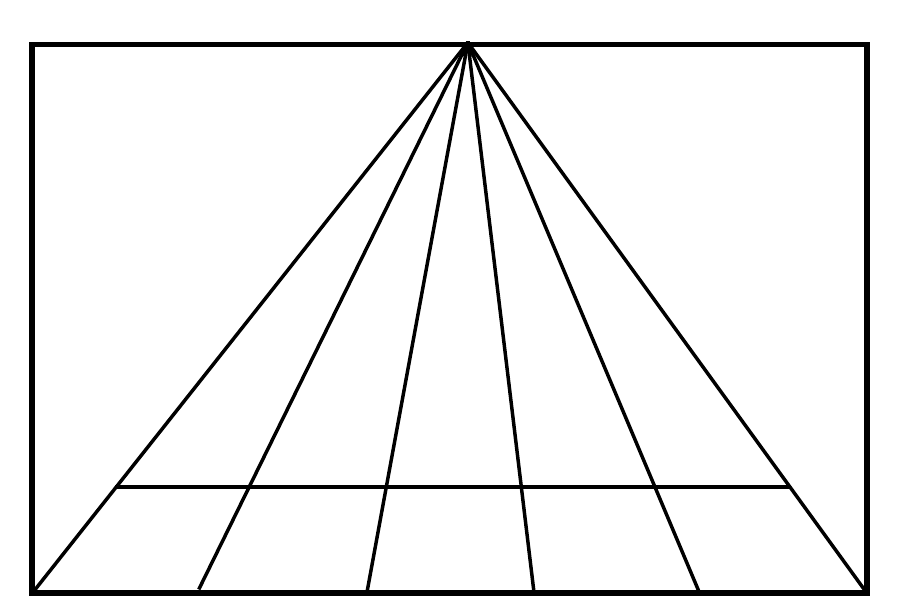} & \includegraphics[height=0.15\textwidth]{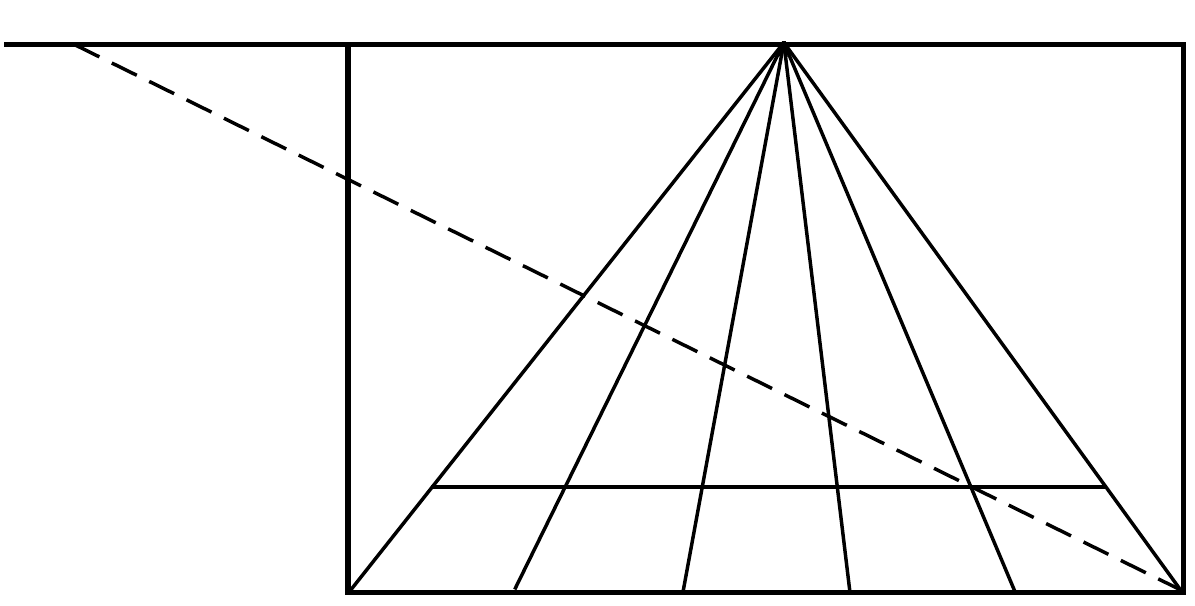} & \includegraphics[height=0.15\textwidth]{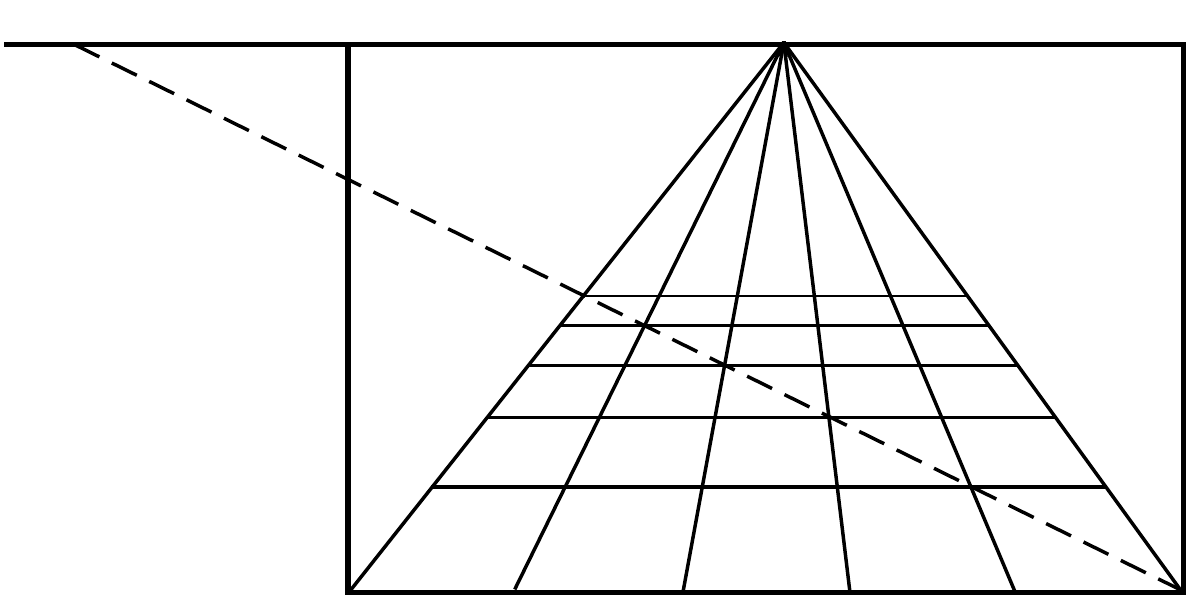} \\
    \end{tabular}\\
    \caption[Drawing a tiled floor.]{Steps for drawing a tiled floor given in Alberti's {\it De Pictura} (1435)---an example of the geometry involved in perspective painting. The construction is based on the elementary geometrical fact that the projection of a straight line onto the canvas is a straight line. From this fact it follows that the diagonal of the first tile is also the diagonal of successive tiles, whence the image of a tiled floor can be drawn according to the steps indicated.} \label{Alberti}
\end{figure}

Only those ignorant of the vast tradition of Greek mathematical science can maintain that the unity of mathematics and science in the 17th century was in any way revolutionary. However, even if one accepts this completely wrongheaded view, credit still should not go to Galileo. Some recent historians have begun to stress that ``the mathematization of the sublunary world begins not with Galileo but with Alberti,''\footnote{\cite[\S5.8]{woottoninvention}.} who wrote on the geometrical principles of perspective in painting in the 15th century (Figure \ref{Alberti}).
\quote{The invention of perspective by the Renaissance artists, …\ by demonstrating that mathematics could be usefully applied to physical space itself, [constituted] a momentous step …\ toward the general representation of physical phenomena in mathematical terms.\footnote{\cite[270]{ConnerPpls}.}}
These historians correctly challenge the narrative of Galileo as the heroic visionary who united mathematics and the physical world, but they retain the erroneous underlying assumption that this unification was revolutionary to begin with. Perspective painting is fine mathematics, but it wasn't a ``momentous step'' ``demonstrating'' that mathematics could be applied to the world, because that had already been demonstrated over and over again thousands of years before. Vitruvius, to take just one example, had pointed out the obvious: ``an architect should be …\ instructed in geometry,'' which ``is of much assistance in architecture.''\footnote{Vitruvius, {\it De Architectura} (1st century BC), I.3--4, \cite[5--6]{Vitruvius}. Indeed, Vitruvius is perfectly well aquatinted with perspective painting, which he describes in VII.11, \cite[198]{Vitruvius}.} Certainly a strange thing to say if the ``momentous'' insight that geometry is relevant to ``the space of everyday experience'' is still more than a thousand years in the future! No, the absurd notion that the application of geometry to physical space was somehow a Renaissance revolution can only occur to those who spend too much time reading philosophical authors pontificating about the divisions of knowledge instead of reading authors actually active in those fields.

The restriction to ``the sublunary world'' in the above quotation is also telling. The allegedly profound conceptual divide between heaven and earth in this period is a standard trope among historians, as we have seen.\footnote{\S\ref{moonmountains}.} Of course, the Greeks mathematised the sublunary world too, but you have to read specialised works to find out much about that. Astronomy, on the other hand, is such an obvious example of an extremely successful and detailed mathematisation of one aspect of reality that even philosophers and historians cannot ignore this elephant in the room. Hence they rely on the qualifier that the allegedly revolutionary step was ``the mathematization of the {\em sublunary} world.'' Aristotle did indeed make much of the difference between the earthly, sublunary world and the world of heavenly motions. But this is one particular dogma of one particular school of philosophy. There is no reason for any mathematician to accept it, nor is there any evidence that any mathematically competent person in the golden age of Greek science did so. The Aristotelian dichotomy is far from natural or necessary: ``Aristotle argues, {\em against his predecessors}, that the celestial world is radically different from the sublunary world,''\footnote{\cite[ix]{FalconA}. Emphasis added.} nota bene. For that matter, even if Aristotle's dogmatic and arbitrary dichotomy is accepted, it would still be madness to acknowledge the undeniable success of mathematics on one side of the divide, yet consider its application on the other side of the divide a conceptual impossibility.

Ptolemy, who was by all appearances much more philosophically conservative than the astronomers of the Hellenistic age, speaks in Aristotelian terms when he contrasts astronomy with physics. The subject matter of astronomy is ``eternal and unchanging,'' while physics ``investigates material and ever-moving nature …\ situated (for the most part) amongst corruptible bodies and below the lunar sphere.''\footnote{Ptolemy, {\it Almagest}, I.1, \cite[36]{AlmagestToomer}.} This is arguably more of a fact than a philosophical commitment: planetary motions are regular and periodic, whereas falling bodies, projectile motion, and other phenomena of terrestrial physics are inherently fleeting and limited to a short time span. But it is conceivable that someone might seize on this dichotomy to ``explain'' why mathematics is suitable for the heavens only, and not for the sublunary world. This, however, is definitely not Ptolemy's stance. He unequivocally expresses the exact opposite view: ``as for physics, mathematics can make a significant contribution'' there too.\footnote{Ptolemy, {\it Almagest}, I.1, \cite[36]{AlmagestToomer}.}

In sum, the Aristotelian dichotomy was never an obstacle to mathematicians. And this with good reason. The whole business of emphasising the dichotomy in the context of the mathematisation of the world is a figment of the imagination of historians, who find themselves having to somehow explain away astronomy as irrelevant when they want to claim that there was a mathematical revolution in early modern science. We do not need to resort to such fictions if we instead accept that the unity of mathematics and science had been obvious since time immemorial.

Another argument for Galileo as the unifier of physics and mathematics consists in stressing that other mathematicians of his day were often more concerned with pure geometry than with projectile motion and the like. For instance, in France there were highly capable ``new Archimedeans'' like Descartes, Roberval, and Fermat, but their focus differed from that of Galileo.
\quote{They were indeed good mathematicians, but they did not consider mathematics as a method for understanding physical things. Mathematical constructions were only abstractions to them, with which it was fun to play, but which were not to be confused with what really happened in nature. Moreover, …\ They were not interested in the ways in which motion intervened in natural processes.
\footnote{\cite[66]{PalmerinoThijssen}.}}
In my view, Galileo would have loved to have been this kind of ``new Archimedean'' too if only he had been capable of it.\footnote{\S\ref{cycloid}.} And it is not true that these Frenchmen ignored motion and the mathematisation of nature. We have already noted that Descartes studied the law of fall,\footnote{\S\ref{lawoffall}.} and that Fermat corrected Galileo on the path of a falling object in absolute space.\footnote{\S\ref{pathoffallabs}.} Both Descartes and Fermat also wrote on the law of refraction of optics, deriving it from physical considerations regarding the speed of light in different media. Also, Descartes explained the motion of the planets, and the fact that they all revolve in the same direction about the sun, by postulating that they were carried along by a vortex. So these mathematicians were clearly not ignorant of or averse to studying how ``motion intervened in natural processes.''

So it is not attention to motion per se, but the study of projectile motion specifically, that sets Galileo apart from these mathematical contemporaries. Does Galileo deserve great credit in this regard? I think not. Why is projectile motion important? With Newton, projectile motion took on a fundamental importance because he saw that planetary motion was governed by the same principles. Galileo had no inkling of this insight. With Newton, projectile motion is also fundamental as a paradigm illustration of the principles---such as inertia and Newton's force law---that govern all other mechanics. In Newtonian mechanics this is the basis for understanding phenomena such as pendulum motion. Galileo, however, got this wrong,\footnote{\S\S\ref{inertia}, \ref{pendulum}, \ref{Newton}.} so he cannot be celebrated for this insight either. Thus we see that praising Galileo for studying projectile motion is anachronistic. Galileo got lucky: the topic he studied later turned out to be very important for reasons he did not perceive, so that in retrospect his work seems much more prescient and groundbreaking than it really was. He himself in fact motivates the theory of projectile motion almost exclusively in terms of practical ballistics---a nonsensical application of zero practical value,\footnote{\S\ref{ballistics}.} which one cannot blame other mathematicians for ignoring.

\subsection{Empirical method}
\label{empiricalmethod}

``Galileo became (and still is) the model for the empiricist scientist who, unlike the natural philosophers of his day, sought to answer questions not by reading philosophical works, but rather through direct contact with nature.''\footnote{\cite[398]{camcomp}.} This is an image Galileo eagerly (but dishonestly) sought to promote, as we have seen.\footnote{\S\S\ref{fallandweight}, \ref{Babylonianeggs}.}

Praise for Galileo in this regard goes hand in hand with ``the verdict that Greek science suffered from an overdose of rash generalizations at the expense of a careful scrutiny, whether experimental or observational, of the relevant facts.''\footnote{\cite[245]{CohenSciRev}.}
In other words, ``Greek thinkers generally …\ overrated the power of unchecked, speculative thought in the natural sciences.''\footnote{\cite[I.92]{DijksterhuisMech}.}

In reality, an empirical approach to the study of nature is not a newfangled invention by Galileo but just common sense. It was obviously adopted by the Greeks, especially the mathematicians. 
In particular, ``ancient mechanics never lost [its] empirical, experiential thread.''\footnote{\cite[43]{IrbyComp}.} Even Aristotle, who practiced ``speculative thought in the natural sciences'' to a much greater extent than mathematicians, was a keen empiricist, and his followers insisted on this as one of the key principles of his philosophy. Aristotle's zoology largely follows a laudable empirical method quite modern in spirit. The same approach was applied by his immediate followers in botany and petrology, including for example cataloging extensive empirical data on how a wide variety of minerals react to heating.\footnote{\cite[11]{lloydGSAA}.} This was far from forgotten in Galileo's day, where one often encounters passages like these from committed Aristotelians:
\quote{We made use of a material instrument to establish by means of our senses what the demonstration had disclosed to our intellect. Such an experimental verification is very important according to [Aristotelian] doctrine.\footnote{Piccolomini, 16th century, \cite[147]{DuhemStatics}.}}
Not infrequently, Galileo's Aristotelian opponents attacked him for being too speculative while they saw themselves as representing the empirical approach. For example, one critic writes to Galileo:
\quote{At the beginning of your work, you often proclaim that you wish to follow the way of the senses so closely that Aristotle (who promised to follow this method and taught it to others) would have changed his opinion, having seen what you have observed. Nonetheless, in the progress of the book you have always been so much a stranger to this way of proceeding that …\ all your controversial conclusions go against our sense knowledge, as anyone can see by himself, and as you expressly say yourself, …\ speaking of the theory of Copernicus, which was rendered plausible and admirable to many by abstract reasoning although it was against all sensory experience.\footnote{Antonio Rocco, \cite[183--184]{sheaGrev}. OGG.VII.712.}}
It is true that there were also many spineless ``Aristotelians'' in Galileo's day who preferred hiding behind textual studies rather than engaging with actual science. But this was one perverse sect of scholasticism, not the overall state of human knowledge before Galileo. A contemporary colleague of Galileo put is well:
\quote{The Science of Nature has been already too long made only a work of the Brain and the Fancy: It is now high time that it should {\em return} to the plainness and soundness of Observations on material and obvious things.\footnote{Robert Hooke, {\it Micrographia} (1665), preface. Emphasis added.}}
``Return'' indeed. Not: Galileo invented this new thing, empiricism. Rather: empiricism is the natural and obvious way to study nature, and the departure from it in certain philosophical circles is a corrupt aberration.

The misconception that the Greeks were anti-empirical stems from a foolish reading of the mathematical tradition:
\quote{Archimedes never appealed to actual measurements in any of his proofs, or even in confirmation of his theorems. …\ The idea that actual measurement could contribute anything of real value …\ was absent from physics for two millennia.\footnote{Drake, \cite[xvii--xviii]{galileo2newsci2ndedWT}.}}
\quote{The mathematics of Euclid and the physics of Archimedes were …\ necessary, but not sufficient, for Galileo’s science. …\ They leave unexplained Galileo’s repeated appeals to sensate experience.\footnote{\cite[281--282]{drakeessays1}.}}
On a superficial reading this may indeed appear so. Open, say, Archimedes's treatise on floating bodies and you will find no mention of any measurement or experiment or data of any kind, only theorems and proofs. It may seem natural to infer from this that Archimedes was doing speculative mathematics divorced from reality, and that he had no understanding of the importance of empirical tests. This is what it looks like to historians who insist on an overly literal reading of the text and lack a sympathetic understanding of how the mathematical mind works. The fact of the matter is that Archimedes's theorems are empirically excellent. It makes no sense to imagine that Archimedes was reasoning about abstractions as an intellectual game, and that his extremely elaborate and detailed claims about the floatation behaviour of various bodies given their shapes and densities just happened to align exactly with reality by pure chance. Archimedes doesn't have to point out that he made very careful empirical investigations, because it is obvious from the accuracy of his results that he did.\footnote{\S\ref{ArchHydrostatics}.}

This is a better way of putting the relation between mathematics and empirical data:
\quote{Mixed mathematics were often presented in axiomatic fashion, following the Archimedean tradition …\ In this tradition, …\ experiments were often conceived of as inherently uncertain and therefore they could not be placed at the foundation of a science, lest that science too be tainted with that same degree of uncertainty. To be sure, experiments were still used as heuristic tools, for example, but their role often remained private, concealed from public presentations.\footnote{\cite[200]{OxHBHistPhys}.}}
So the point is not that empirical data is neglected, but that it is a mere preliminary step. Anyone can make measurements and collect data. Self-respecting mathematicians do not publish such trivialities. Instead they go on to the really challenging step of synthesising it into a coherent mathematical theory. Galileo did not have the ability to do the latter, so he had to stick with the basics, and pretend, nonsensically, that this was somehow an important innovation. Then as now, there were enough non-mathematicians in the world for his cheap charade to be successful.

\subsection{Experimental method}
\label{experimentalmethod}

Some draw a distinction between passive observation and active experiment, and take the latter to be the key Galilean innovation.
\quote{When Galileo caused balls …\ to roll down an inclined plane …\ a light broke upon all students of nature. …\ Reason …\ must approach nature in order to be taught by it. It must not, however, do so in the character of a pupil who listens to everything that the teacher chooses to say, but of an appointed judge who compels the witnesses to answer questions which he has himself formulated.\footnote{Kant, {\it Critique of Pure Reason}, B.xii--xiii, \cite[20]{kantcrittransl}.}}
\quote{The originality of Galileo's method lay precisely in his effective combination of mathematics with experiment. …\ The distinctive feature of scientific method in the seventeenth century, as compared with that in ancient Greece, was its conception of how to relate a theory to the observed facts …\ and submitting them to experimental tests. …\ [This feature] transformed the Greek geometrical method into the experimental science of the modern world.\footnote{\cite[303, 1]{CrombieGross}.}}
In reality, the use of experiment in Greek science is abundantly documented to anyone who bothers to read mathematical authors. 

Greek scientists knew perfectly well that ``it is not possible for everything to be grasped by reasoning …, many things are also discovered through experience.''\footnote{Philon, c.\ 200 B.C., {\it Belopoeica} 50.14--29, \cite[625]{SchiefskyPhilo}.} This quote refers to the precise numerical proportions needed for the spring in a stone-throwing engine. The same author also offered an experimental demonstration that air is corporeal.\footnote{Philon, \cite[217]{HellenisticEraSourcebook}.} Ptolemy experimented with balloons (``inflated skins'') to investigate whether air or water has weight in their own medium---indeed he ``performed the experiment with the greatest possible care.''\footnote{As reported by Simplicius. \cite[225]{HellenisticEraSourcebook}, \cite[248]{CohenDrabkinSourceB}.} Heron of Alexandria gives a detailed description of an experimental setup to prove the existence of a vacuum. He explicitly states that ``referring to the appearances and to what is accessible to sensation'' trumps abstract arguments that there can be no vacuum.\footnote{\cite[17]{lloydGSAA}.} Such arguments had been given by Aristotle. In optics, Ptolemy explicitly verified the law of reflection by experiment, and studied refraction experimentally, giving tables for the angle of refraction of a light ray for various incoming angles in increments of 10 degrees for passages between air, water, and glass.\footnote{\cite[80--82]{ClagettGS}. \cite[200--201]{HellenisticEraSourcebook}. \cite[271--281]{CohenDrabkinSourceB}. \cite[133--135]{lloydGSAA}.}

\begin{figure}
\centering
\includegraphics[width=0.75\textwidth]{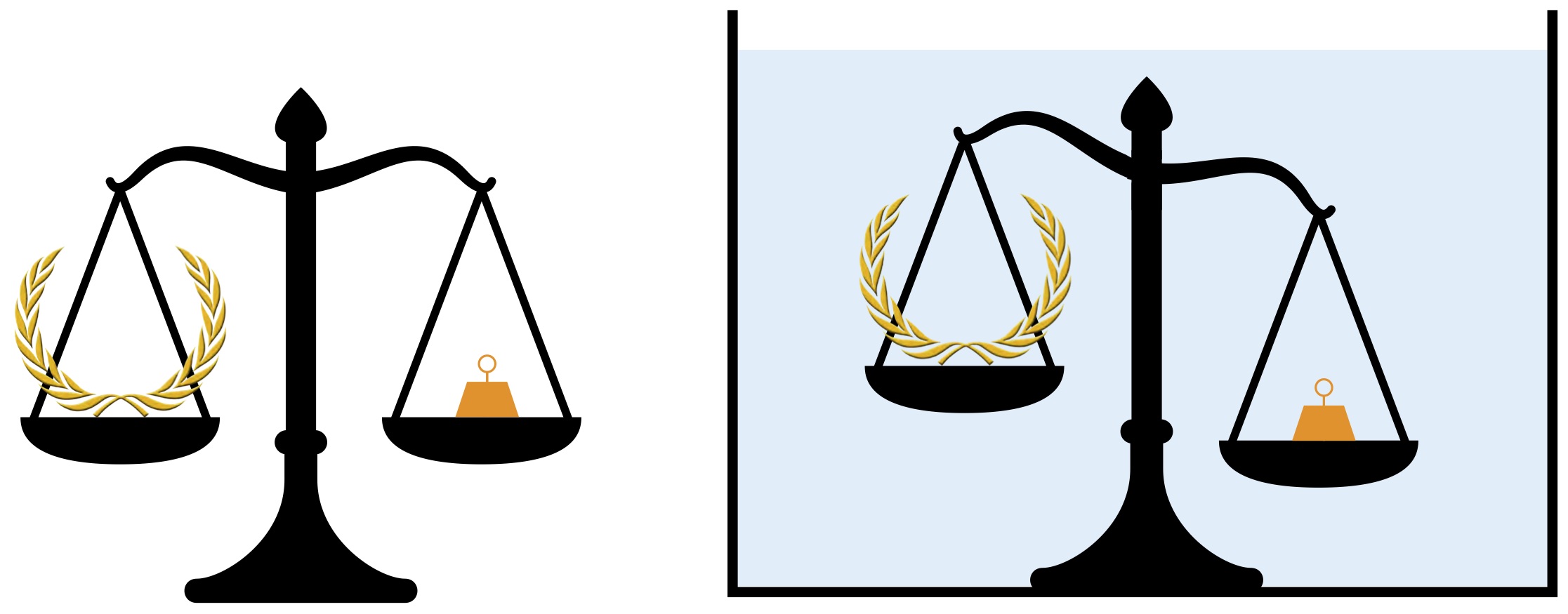}
\caption[Archimedes's gold wreath experiment.]{Archimedes used a scientific experiment to expose a forger who tried to pass off as pure gold a wreath that was actually gold-coated silver. Archimedes's principle of hydrostatics states that the upward buoyant force on a submerged object is equal to the volume of the displaced water. Therefore, if the wreath is lighter in water than an equal mass of gold, then it has greater volume, and hence lower density, and hence is not pure gold.} \label{Acrown}
\end{figure}

\begin{figure}
\centering
\includegraphics[width=0.75\textwidth]{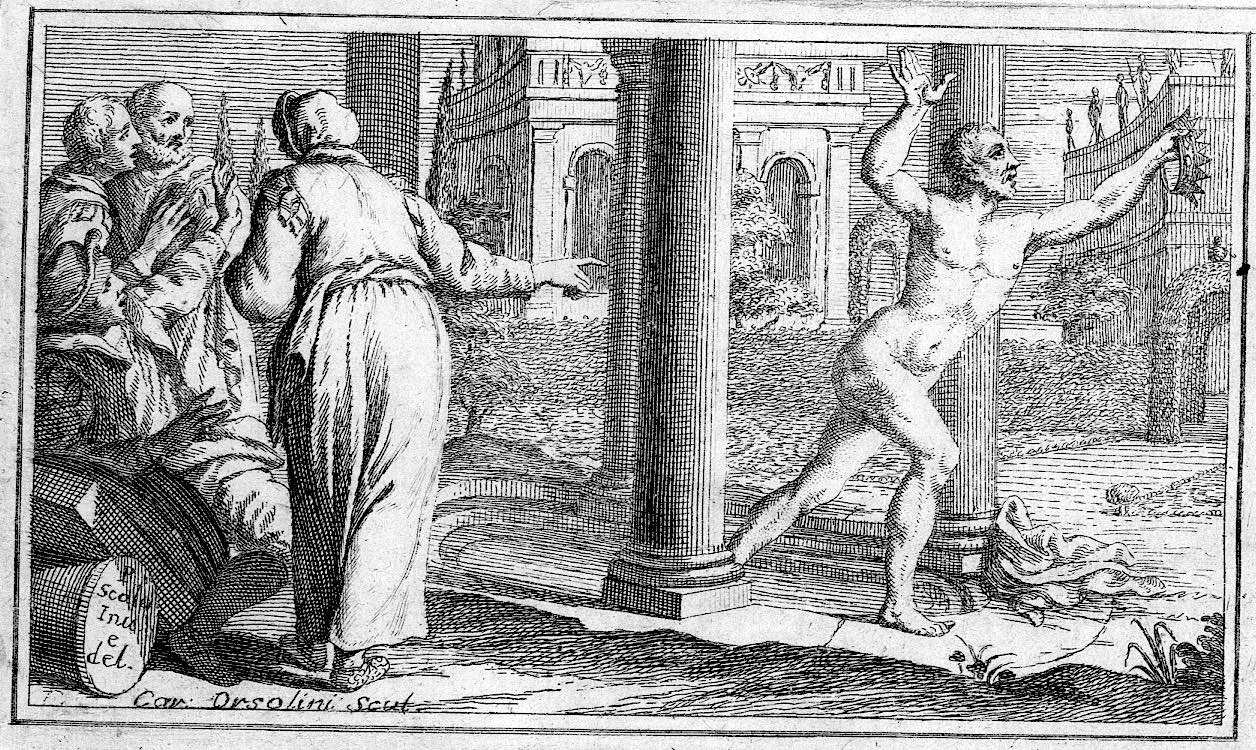}
\caption[Archimedes running.]{Archimedes's excitement at having had the idea for the experiment of Figure \ref{Acrown}. A testament to his love of experimental science. Frontispiece of Gian-Maria Mazzuchelli, {\it Notizie istoriche e critiche intorno alla vita, alle invenzioni, ed agli scritti di Archimede siracusano}, Brescia, 1737.} \label{Acrown2}
\end{figure}

Archimedes performed a scientific experiment (Figure \ref{Acrown}) that pleased him so much that he ran naked through the streets yelling ``eureka'' in excitement.\footnote{\cite[238--239]{CohenDrabkinSourceB}. Vitruvius, {\it De Architectura}, IX. \cite{RorresCrown}.} Such was his love of empirical, experimental science---yet many scholars keep insisting that, like a second Plato, all Archimedes really cared about was abstract geometry. Evidently, even running naked through the streets and screaming at the top of one's lungs is not enough for some people to open their eyes. It is hard to imagine what else one can do to draw their attention to the obvious, namely that Greek mathematicians embraced experimental method through and through.

\subsection{Technology and craft knowledge}

A prevalent view has it that Galileo was the first to bring together abstract mathematics and science with concrete technology and practical know-how of craftsmen and workers in mechanical fields.
\quote{Real science is born when, with the progress of technology, the experimental method of the craftsmen overcomes the prejudice against manual work and is adopted by rationally trained university-scholars. This is accomplished with Galileo.\footnote{\cite[5]{ZilselBoston}.}}
\quote{[Galileo was able] to bring together two once separate worlds that from his time on were destined to remain forever closely linked---the world of scientific research and that of technology.\footnote{Drake, \cite[78]{GalileoDiscOp}.}}
\quote{Galileo may fruitfully be seen as the culmination point of a tradition in Archimedean thought which, by itself, had run into a dead end. What enabled Galileo to overcome its limitations …\ seems easily explicable upon considering Galileo's background in the arts and crafts.\footnote{\cite[349]{CohenSciRev}.}}
\quote{The separation between theory and practice, imposed by university professors of natural philosophy, was repeatedly exposed as untenable. Of course the greatest figure in this movement is Galileo.\footnote{\cite[16]{HenrySciRev}.}}
Galileo himself eagerly cultivated this image. The very first words of his big book on mechanics are devoted to extolling the importance for science of observing ``every sort of instrument and machine'' in action at the ``famous arsenal'' of Venice, praising the experiential knowledge of the ``truly expert'' workmen there.\footnote{\cite[11]{galileo2newsci2ndedWT}, OGG.VIII.49.}

It is true that universities were filled with many blockheads who foolishly insisted on keeping intellectual work aloof from such connections to the real world. For example, when Wallis went to Oxford in 1632 there was no one at the university who could teach him mathematics. “For Mathematicks, (at that time, with us) were scarce looked upon as Accademical studies, but rather Mechanical; as the business of Traders, Merchants, Seamen, Carpenters, Surveyors of Lands, or the like.”\footnote{Wallis, in his autobiography. \cite[27]{AutobWallis}.}

But it would be mistake to infer from this that Galileo's step was an innovation. The stupidity of the university professors was the doing of one particular clique of mathematically ignorant people. Their attitude is not natural or representative of the state of human knowledge. Galileo is not a brilliant maverick thinking outside the box. Rather, he is merely doing what had, among mathematically competent people, been recognised as the natural and obviously right way to do science for thousands of years. Galileo is not taking a qualitative leap beyond limitations that had crippled all previous thinkers. Rather, he is merely reversing the obvious cardinal error of one particularly dumb philosophical movement that had happened to gain too much influence at the time, because people were too ignorant to recognise the evident superiority of more mathematical and scientific schools of thought that had already proven their worth in a large body of ancient works available to anyone who cared to read.

In order to defend the misconceived idea of Galileo the trailblazing innovator one must ignore the large body of obvious precedent for his view in antiquity, and project the foolish nonsense of medieval universities onto the Greeks. Indeed, historians have concocted a false narrative to this effect.
\quote{Greek technology and science were …\ rigidly separated.\footnote{\cite[2]{LindbergMidAg}.}}
\quote{The Greek hand worker was considered inferior to the brain worker or contemplative thinker. …\ So, despite the fact that the philosophers derived some of their conclusions as to how nature behaved from the work of the craftsmen, they rarely had experience of that work. What is more, they were seldom inclined to improve it, and so were powerless to pry apart its potential treasure of knowledge that was to lead to the scientific revolution in the Renaissance.\footnote{\cite[59]{BrakeRev}.}}
\quote{The fundamental brake upon the further progress of science in antiquity was slave labour [which precluded any] meaningful combination of theory and practice.\footnote{\cite[248]{CohenSciRev}, summarising the view of Farrington. A similar view is expressed in \cite[74]{DijksterhuisMech}.}}
More specialised scholarship knows better. The recent {\it Oxford Handbook of Engineering and Technology in the Classical World} is perfectly clear on the matter:
\quote{Many twentieth-century scholars hit upon …\ banausic prejudice [i.e., a snobbish contempt for manual labour] as an `explanation' for a perceived blockage of technological innovation in the Greco-Roman world. The presence of slave labor was felt to be a related, concomitant factor. …\ [But] this now discredited interpretation [should be rejected and we should] …\ put an end to the myth of a `technological blockage' in the classical cultures.\footnote{\cite[5--6]{OxHTechClassical}.}}
This is the view of experts on the matter, while the false narrative is promulgated by scholars who focus on Galileo, take it for granted that he is ``the Father of Modern Science,'' and postulate such nonsense about the Greeks because that's the only way to craft a narrative that fits with this false assumption.

\begin{figure}[tp]\centering
\includegraphics[width=0.5\textwidth]{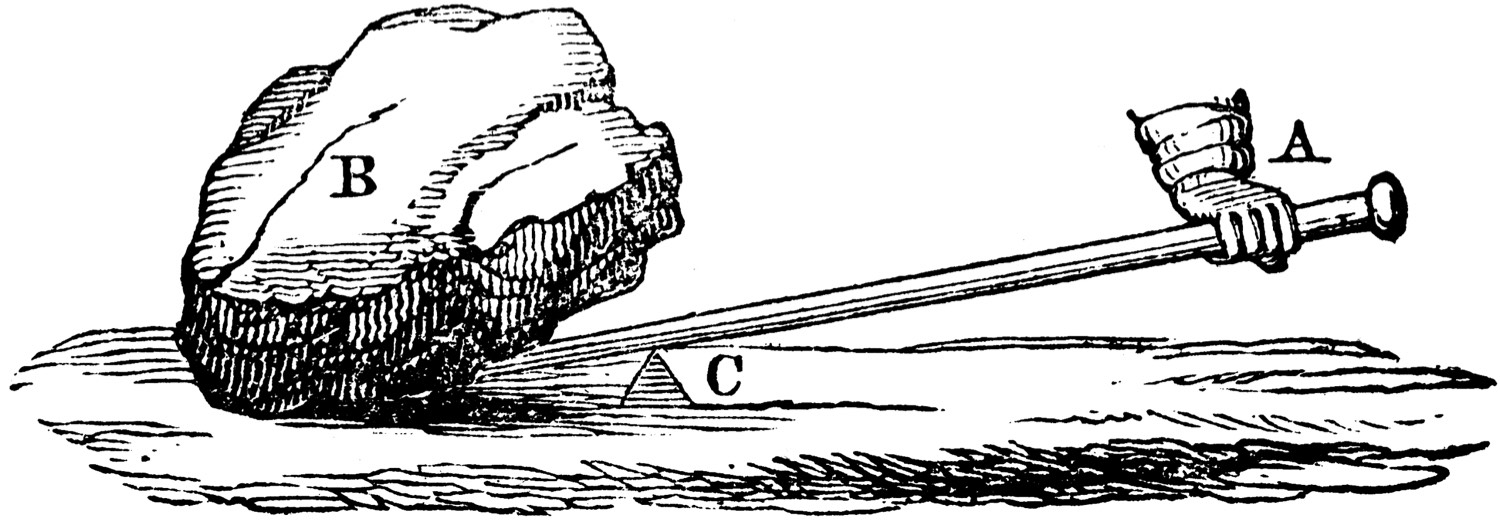}
\caption[The lever.]{The law of the lever states that if the lever arm $AC$ is $n$ times longer than $CB$, then only one $n^\textnormal{th}$ of the force is required to lift the stone with the lever compared to lifting it directly. Archimedes wrote an excellent mathematical treatise on the lever. Figure from \cite[69]{Comstock}.
}
\label{lever}
\end{figure}

\begin{figure}[tp]\centering
\includegraphics[width=0.2\textwidth]{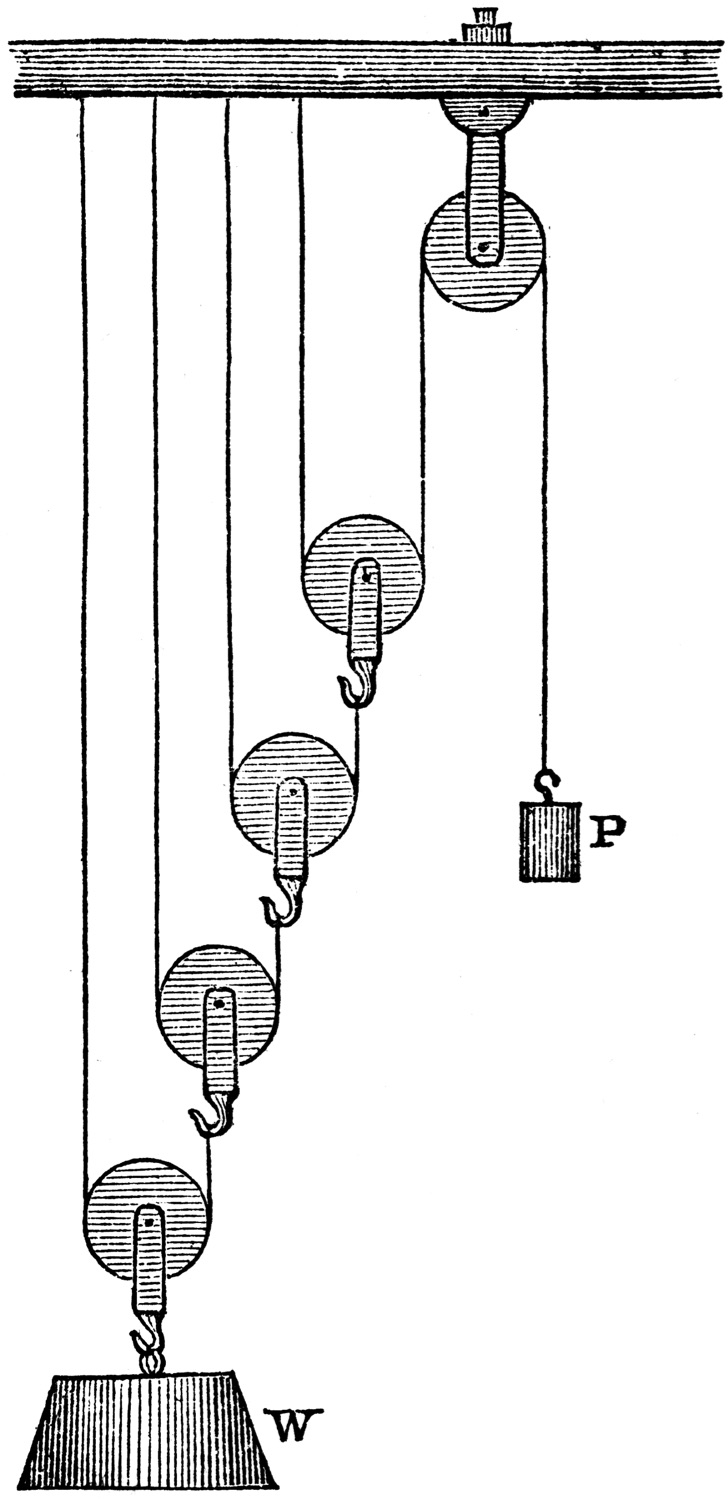}
\caption[The pulley.]{A compound pulley. Each pulley beyond the first doubles the mechanical advantage. In this case, since there are 4 pulleys beyond the first, a weight $P$ of $1$ balances a weight $W$ of $2^{4}=16$. Thus the heavy weight $W$ can be lifted by applying a small force at $P$. Greek mechanicians wrote extensively on the mathematical relations involved in this and other ``machines'' of this sort. Figure from \cite[87]{Comstock}.
}
\label{pulley}
\end{figure}

\begin{figure}[tp]\centering
\includegraphics[width=0.5\textwidth]{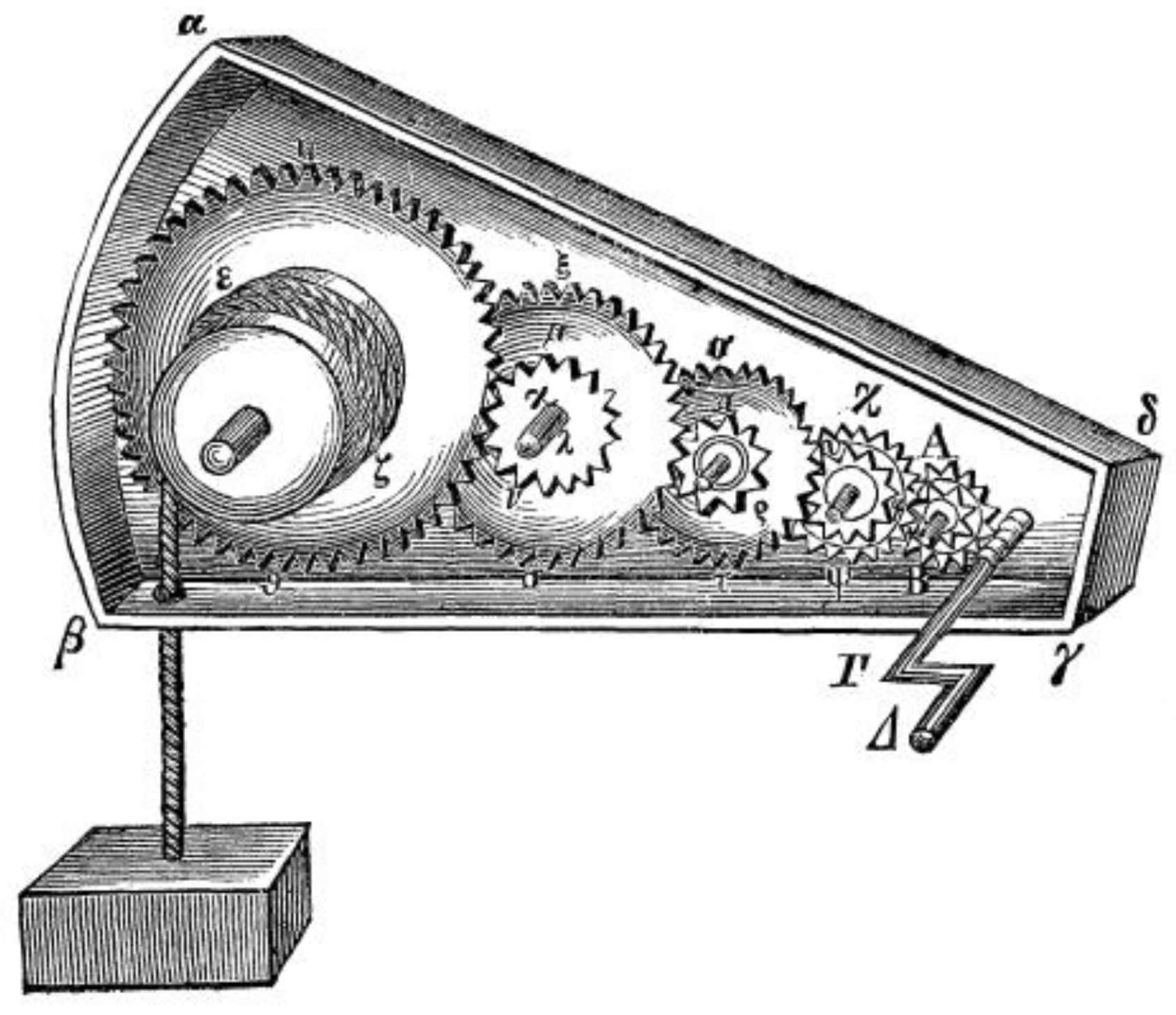}
\caption[Gear machine.]{How to lift a given weight by a given force, as explained by Pappus of Alexandria in the 4th century ({\it Collection}, VIII.X, \cite[333]{PappusMech}). In the case illustrated, a force of 4 units applied at the crank arm $\Delta$ lifts a weight that it would have taken 160 units of force to lift directly.}
\label{gears}
\end{figure}

Promulgators of the nonsense about practice-adverse Greeks have evidently not bothered to read mathematical authors. Pappus, for example, explains clearly that mathematicians enthusiastically embrace practical and manual skills:
\quote{The science of mechanics …\ has many important uses in practical life, …\ and is zealously studied by mathematicians. …\ Mechanics can be divided into a theoretical and a manual part; the theoretical part is composed of geometry, arithmetic, astronomy and physics, the manual of working in metals, architecture, carpentering and painting and anything involving skill with the hands.\footnote{\cite[II.615]{LoebMath}.}}
He praises the interaction of geometry with practical fields or ``arts'' as beneficial to both:
\quote{Geometry is in no way injured, but is capable of giving content to many arts by being associated with them, and, so far from being injured, it is obvious, while itself advancing those arts, appropriately honoured and adorned by them.\footnote{\cite[II.619--621]{LoebMath}.}}
These were no empty words. The Greeks had an extensive tradition of studying ``machines,'' meaning devices based on components such as the lever (Figure \ref{lever}), pulley (Figure \ref{pulley}), wheel and axle, winch, wedge, screw, and gear wheel (Figure \ref{gears}). The primary purpose of these machines was that of ``multiplying an effort to exert greater force than can human or animal muscle power alone.'' Such machines were ``used primarily in construction, water-lifting, mining, the processing of agricultural produce, and warfare.''\footnote{\cite[337]{OxHTechClassical}.} The Greeks also undertook advanced engineering projects, such as digging a tunnel of more than a kilometer through a mountain, the planning of which involved quite sophisticated geometry to enable the tunnel to be dug from both ends, with the diggers meeting in the middle.\footnote{\cite[324--325]{OxHTechClassical}.} In short, ``while it is crucial to distinguish between theoretical mechanics and practitioners' knowledge, there is substantial evidence of a two-way interaction between them in Antiquity.''\footnote{\cite[15]{LairdRoux}.}

Mathematicians were very much involved with such things. Indeed, ``the modern distinction between physical and mathematical sciences was alien to Hellenistic science, which was unitary.''\footnote{\cite[189]{RussoSciRev}.} There are many testimonies attributing to Archimedes various accomplishments in engineering, such as moving a ship singlehandedly by means of pulleys, destroying enemy ships using machines, building a screw for lifting water, and so on.\footnote{Plutarch, {\it Marcellus}, 14.8. \cite[51]{Proclus}. \cite[14--29]{DArchimedes}.} Apollonius wrote a very advanced and thorough treatise on conic sections, which is studiously abstract and undoubtedly {\it l'art pour l'art} pure mathematics if there ever was such a thing. Yet the same Apollonius ``besides writing on conic sections produced a now lost work on a flute-player driven by compressed air released by valves controlled by the operation of a water wheel.''\footnote{\cite[338]{OxHTechClassical}.} The title page of the Arabic manuscript that has preserved this work for us reads: ``by Apollonius, the carpenter, the geometer.''\footnote{\cite[X.334]{HillArabic}.} The cliche of Greek geometry as nothing but abstruse abstractions divorced from reality is a modern fiction. The sources tell a different story. It is not for nothing that one of the most refined mathematicians of antiquity went by the moniker ``the carpenter.''

Unfortunately,
\quote{Renaissance intellectuals were not in a position to understand Hellenistic scientific theories, but, like bright children whose lively curiosity is set astir by a first visit to the library, they found in the manuscripts many captivating topics, especially those that came with illustrations. …\ The most famous intellectual attracted by all these `novelties' was Leonardo da Vinci. …\ Leonardo’s `futuristic' technical drawings …\ was not a science-fiction voyage into the future so much as a plunge into a distant past. Leonardo’s drawings often show objects that could not have been built in his time because the relevant technology did not exist. This is not due to a special genius for divining the future, but to the mundane fact that behind those drawings …\ there were older drawings from a time when technology was far more advanced.\footnote{\cite[335--336]{RussoSciRev}.}}
The false narrative of the mechanically ignorant, anti-practical Greeks has obscured this fact, and led to an exaggerated evaluation of Renaissance technology, such as instruments for navigation, surveying, drawing, timekeeping, and so on.
\quote{Renaissance developments in practical mathematics predated the intellectual shifts in natural philosophy. …\ Historians of the early modern reform of natural philosophy have failed to appreciate the significance of the prior success of the practical mathematical programme, [which] must figure in an explanation of why the new dogma of the seventeenth century embraced mathematics, mechanism, experiment and instrumentation.\footnote{\cite[176, 189]{Bennett1991}.}}
This author proves at length that the practical mathematical tradition had much to commend it, which I do not dispute. But then he casually asserts with hardly any justification that there was nothing comparable in Greek times.\footnote{\cite[182]{Bennett1991}.} This is typical of much scholarship of this period. The deeply entrenched standard view of the Galilean revolution is basically taken for granted and subsequent work is presented as emendations to it. For instance, if you want to prove the importance of a Renaissance pre-revolution in practical mathematics, you need to prove two things: first that it was relevant to the scientific revolution, and second that it was not present long before. It is a typical pattern to see historians put all their efforts toward proving the first point, and glossing over the second point in sentence or two. They can get away with this since the alleged shortcomings of the Greeks is supposedly common knowledge, while the first point is the one that departs from the standard narrative. Hence, if the standard narrative is misconceived in the first place, so is all this more specialised research, which, although it ostensibly departs from the standard view, actually retains its most fundamental errors in the very framing of its argument.

It is right to emphasise that the practical mathematical tradition stood for a much more fruitful and progressive approach to nature than that dominant among the philosophy professors of the time. But it is a mistake to believe that these professors represented the considered opinion of the best minds, while the mathematical practitioners were oddball underdogs whose pioneering success eventually proved undeniable to the surprise of everyone. The mathematical practitioners stood for simple common sense, not renegade iconoclasm. They practiced the same common sense that their peers had in antiquity, with much the same results. The university professors, meanwhile, should not be mistaken for a neutral representation of the state of human knowledge at the time. Rather, they formed one particular philosophical sect which retained its domination of the universities not because of the preeminence of its teachings but because of the incestuous appointment practices and obsequiousness of academics.

\subsection{Instrumentalism versus realism}
\label{realism}

A standard view is that ``the Scientific Revolution saw the replacement of a predominantly instrumentalist attitude to mathematical analysis with a more realist outlook.''\footnote{\cite[8]{HenrySciRev}.} Instrumentalism means the following:
\quote{An explanation which conforms to the facts does not imply that the hypotheses are real and exist. …\ [Astronomers] have been unable to establish in what sense, exactly, the consequences entailed by these arrangements are merely fictive and not real at all. So they are satisfied to assert that it is possible, by means of circular and uniform movements, always in the same direction, to save the apparent movements of the wandering stars.\footnote{Simplicius, \cite[23]{DuhemSavePhen}.}}
Instrumentalism, as opposed to realism, was supposedly the accepted philosophy of science among ``the Greeks'':
\quote{[Ancient Greek astronomers] balked at the idea that the eccentrics and epicycles are bodies, really up there on the vaults of the heavens. For the Greeks they were simply geometrical fictions requisite to the subjection of celestial phenomena to calculation. If these calculations are in accord with the results of observation, if the `hypotheses' succeed in `saving the phenomena', the astronomer's problem is solved.\footnote{\cite[25]{DuhemSavePhen}.}}
\quote{An astronomer who understands the true purpose of science, as defined by men like Posidonius, Ptolemy, Proclus, and Simplicius, …\ would not require the hypotheses supporting his system to be {\em true}, that is, in conformity with things. For him it will be enough if the results of calculation agree with the results of observation---{\em if appearances are saved}.\footnote{\cite[31]{DuhemSavePhen}.}}
\quote{The Greek geometer in formulating his astronomical theories does not make any statements about physical nature at all. His theories are purely geometrical fictions. That means that to save the appearances became a purely mathematical task, it was an exercise in geometry, no more, but, of course, also no less.\footnote{\cite[57]{WassersteinGSci}.}}
Galileo, by contrast, brought ``a radically new mode of realist-mathematical nature knowledge.''\footnote{\cite[146]{CohenRiseExpl}.}
\quote{Galileo endorsed a view that was [contrary to] that of the Greeks but was also much more creative …\ It is a crippling restriction to hold that no theory about reality can be in mathematical form; the Renaissance rejected this restriction, holding that it was a worthwhile enterprise to search for mathematical theories which also---by metaphysical criteria---could be supposed `real'. …\ The most eloquent and full defence of this process was given by Galileo.\footnote{\cite[11--12]{HallRevSci}.}}
Hence the Scientific Revolution owes much to ``the novel quality of realism that the abstract-mathematical mode of nature-knowledge acquired in …\ Galileo's hands.''\footnote{\cite[120]{CohenRiseExpl}.}

In reality, no mathematically competent Greek author ever advocated instrumentalism. The notion that ``the Greeks'' were instrumentalists relies exclusively on passages by philosophical commentators. The notion that Ptolemy believed his planetary models were ``fictional …\ combinations of circles which could never exist in celestial reality''\footnote{\cite[17]{CohenRiseExpl}.} is demonstrably false. First of all Ptolemy opens his big book with physical arguments for why the earth is in the center of the universe---a blatantly realist justification for this aspect of his astronomical models. Furthermore, Ptolemy has a detailed discussion of the order and distances of the planets that obviously assumes that the planetary models, epicycles and all, are physically real. ``The distances of the …\ planets may be determined without difficulty from the nesting of the spheres, where the least distance of a sphere is considered equal to the greatest distance of a sphere below it.''\footnote{Ptolemy, {\it Planetary Hypotheses}, I, \cite[7]{PtolemyPlanHyp}.} That is to say, according to Ptolemy's epicyclic planetary models, each planet sways back and forth between a nearest and a furthest distance from the earth. The ``sphere'' of each planet must be just thick enough to contain these motions. Ptolemy assumes that ``there is no space between the greatest and least distances [of adjacent spheres],'' which ``is most plausible, for it is not conceivable that there be in Nature a vacuum, or any meaningless and useless thing.''\footnote{Ptolemy, {\it Planetary Hypotheses}, I, \cite[8]{PtolemyPlanHyp}.} Clearly this is based on taking planetary models to be very real indeed, and not at all mathematical fictions invented for calculation.

Nor was Ptolemy an exception in his realism. His colleague Geminos ``was a thoughtful realist'' as well.\footnote{\cite[58]{EBgeminos}.} Hipparchus too evidently chose models for planetary motion on realist grounds. His works are lost, but we know that he proved the mathematical equivalence of epicyclic and equant motion. In other words, he showed that two different geometrical models of planetary motion are observationally equivalent; they lead to the exact same visual impressions seen from earth, but they are brought about by different mechanisms. How should one choose between the two models in such a case? If Hipparchus was an instrumentalist, he wouldn't care one way or the other, or he would just pick whichever was more mathematically convenient. But if he was a realist he would be interested in which model could more plausibly correspond to actual physical reality. So what did he do? Here is what Theon says: ``Hipparchus, convinced that this is how the phenomena are brought about, adopted the epicyclic hypothesis as his own and says that it is likely that all the heavenly bodies are uniformly placed with respect to the center of the world and that they are united to it in a similar way.''\footnote{Theon, {\it Astronomia} 34, \cite[9]{DuhemSavePhen}.} So Hipparchus decided between equivalent models based on physical plausibility. This is quite clearly a realist argument.

Historians have brought up other ``evidence'' that ``the Greeks'' were instrumentalists. One thing they point to is the alleged compartmentalisation of Greek science.
\quote{Phenomena [such as] consonance, light, planetary trajectories and the two states of equilibrium [i.e., statics and hydrostatics] are investigated separately. There is no search for interconnections, let alone for an overarching unity.\footnote{\cite[19--20]{CohenRiseExpl}.}}
This would indeed make sense if mathematical science was just instrumental computation tools with no genuine anchoring in reality. The only problem is that the claim is false. Greek science is in fact full of interconnections, just as one would expect if they were committed realists. Ptolemy uses mechanics to justify geocentrism; Archimedean hydrostatics explains shapes of planets and ``casts light on the earth’s geological past'';\footnote{\cite[303]{RussoSciRev}. See \S\ref{Greekastro}.} Archimedes used statical principles to compute areas in geometry. Ptolemy applies ``consonance'' (that is, musical theory) to ``the human soul, the ecliptic, zodiac, fixed stars, and planets.''\footnote{Ptolemy, {\it Harmonics}, \cite[xxxiv]{PtolHarmonics}.} He also applies the law of refraction of optics to atmospheric refraction, noting its importance for astronomical observations.\footnote{\cite[281--283]{CohenDrabkinSourceB}.}

In Galileo's time, the same pattern as among the Greeks prevails: mathematically competent people are unabashed realists, while philosophers and theologians often find instrumentalism more appealing for reasons that have nothing to do with science. Copernicus's book, for example, is unequivocally realist. Spineless philosophers and theologians could not accept this. One even resorted to the ugly trick of inserting an unsigned foreword in the book without Copernicus's authorisation, in which they espoused instrumentalism:
\quote{It is the job of the astronomer to use painstaking and skilled observation in gathering together the history of the celestial movements, and then---since he cannot by any line of reasoning reach the true causes of these movements---to think up or construct whatever causes or hypotheses he pleases such that, by the assumption of these causes, those same movements can be calculated from the principles of geometry for the past and for the future too. …\ It is not necessary that these hypotheses should be true. …\ It is enough if they provide a calculus which fits the observations.\footnote{Osiander, foreword to Copernicus's book, \cite[3]{CopRevTransl}. For the attribution to Osiander, see \cite[152]{KeplerApologia}.
}}
This surely fooled no one who actually read the book, with all its blatant realism. Giordano Bruno, for one, thought ``there can be no question that Copernicus believed in this motion [of the earth],'' and hence concluded that the timid foreword must have been written ``by I know not what ignorant and presumptuous ass.''\footnote{Bruno, \cite[99]{DuhemSavePhen}.} But then again the mathematically incompetent people whom this foreword was designed to appease could not read the book anyway.

In medieval and renaissance philosophical texts it is not hard to find many assertions to the effect that ``real astronomy is nonexistent'' and what passes for astronomy ``is merely something suitable for computing the entries in astronomical almanacs.''\footnote{Achillini (16th century), \cite[47]{DuhemSavePhen}.} There were many instrumentalists at the time,\footnote{\cite[231--237]{KeplerApologiaJardine}.} to be sure, but the challenge is to find a single serious mathematical astronomer among them. They were exclusively theologians and philosophers.

All historians nowadays recognise that ``Copernicus clearly believed in the physical reality of his astronomical system,'' but their inference that he ``thus broke down the traditional disciplinary boundary between astronomy (a branch of mixed mathematics) and physics (or natural philosophy)''\footnote{\cite[52]{OslerReconfiguring}.} is dubious. This was ``the traditional'' view only in a very limited sense. It was traditional among the particular sect of Aristotelians that occupied the universities, but outside this narrow clique it had no credibility or standing whatsoever. Among mathematicians, Copernicus's view was exactly the traditional one.

All mathematically competent people continued in the same vein, long before Galileo entered the scene. Already in the 16th century, ``Tycho and Rothman, Maestlin, and even Ursus …\ openly deploy a wide range of physical arguments in debating the issue between the rival world-systems.''\footnote{\cite[244]{KeplerApologiaJardine}.} Kepler puts the matter very clearly:
\quote{One who predicts as accurately as possible the movements and positions of the stars performs the task of the astronomers well. But one who, in addition to this, also employs true opinions about the form of the universe performs it better and is held worthy of greater praise. The former, indeed, draws conclusions that are true as far as what is observed is concerned; the latter not only does justice in his conclusions to what is seen, but also …\ in drawing conclusions embracing the inmost form of nature.\footnote{\cite[145]{KeplerApologia}.}}
As Kepler notes, this was all obviously well-known and accepted since antiquity, for ``to predict the motions of the planets Ptolemy did not have to consider the order of the planetary spheres, and yet he certainly did so diligently.''\footnote{\cite[145]{KeplerApologia}.}

\subsection{Mechanical philosophy}

Some say that ``the mechanization of the world-picture'' was the defining ingredient of ``the transition from ancient to classical science.''\footnote{\cite[501]{DijksterhuisMech}.} A paradigm conception at the heart of the new science was that of the world as a machine: a ``clockwork universe'' in which everything is caused by bodies pushing one another according to basic mechanical laws, as opposed to a world governed by teleological purpose, divine will and intervention, anthropomorphised desires and sympathies ascribed to physical objects, or other supernatural forces. Galileo was supposedly a pioneer in how he always stuck to the right side in this divide.
\quote{Galileo possessed in a high degree one special faculty. …\ That is the faculty of thinking correctly about physical problems as such, and not confusing them with either mathematical or philosophical problems. It is a faculty rare enough still, but much more frequently encountered today than it was in Galileo's time, if only because nowadays we all cope with mechanical devices from childhood on.\footnote{\cite[603]{drakeinertia}.}}
Of course, this ``special faculty'' is precisely what led Galileo to reject as occult the correct explanation of the tides and propose his own embarrassing nonstarter of a tidal theory based on an analogy with ``mechanical devices.''\footnote{\S\ref{tides}.} But let us put that aside.

There is nothing modern about the mechanical philosophy. ``We all cope with mechanical devices from childhood on,'' but so did the Greeks, who built automata such as entirely mechanical puppet-theatres, self-opening temple doors, a coin-operated holy water dispenser, and so on.\footnote{\cite{GRTechnologySourcebook}, \cite[224--234]{CohenDrabkinSourceB}.} Pappus notes that ``the science of mechanics'' has many applications ``of practical utility,'' including machines for lifting weights, warfare machines such as catapults, water-lifting machines, and ``marvellous devices'' using ``ropes and cables to simulate the motions of living things.''\footnote{Pappus, {\it Collection}, VIII.1, \cite[183--184]{CohenDrabkinSourceB}.} Clearly, then, ``Ancient Greek mechanics offered working artifacts complex enough to suggest that …\ organisms, the cosmos as a whole, or we ourselves, might `work like that'.''\footnote{\cite[229--230]{Berryman2009}.} Thus we read in ancient sources that ``the universe is like a single mechanism'' governed by simple and deterministic laws that ultimately lead to ``all the varieties of tragic and comedic and other interactions of human affairs.''\footnote{Theodorus, as paraphrased by Proclus, {\it On providence} 2, \cite[243]{JonesAntiky}.} This line of reasoning soon lead to a secularisation of science. ``Bit by bit, Zeus was relieved of thunderbolt duty, Poseidon of earthquakes, Apollo of epidemic disease, Hera of births, and the rest of the pantheon of gods were pensioned off'' in the same manner.\footnote{\cite[17]{rihllGS}.}

Mechanical explanations are widespread in Greek science. The Aristotelian {\it Mechanics} uses the law of lever to explain ``why rowers who are …\ in the middle of the ship …\ move the ship the most,'' and ``how it is that dentists extract teeth more easily by …\ a tooth-extractor [forceps] than with the bare hand only.''\footnote{(Pseudo-)Aristotle, {\it Mechanics}, \cite[194]{CohenDrabkinSourceB}.} Greek scientists explained perfectly clearly that sound is a ``wave of air in motion,'' comparable to the rings forming on a pond when when one throws in a stone.\footnote{\cite[74]{ClagettGS}.} Atomism---a widely espoused conception of the world in Greek antiquity---is of course in effect a plan to ``make material principles the basis of all reality.''\footnote{\cite[41]{IrbyComp}.}

Greek astronomy went hand in hand with mechanical planetaria that directly reproduced a scale model of planetary motion. And not just basic toy models, but ``complex and scientifically ambitious instruments'' that could generate all heavenly motions mechanically from a single generating motion (the turn of a crank, as it were).\footnote{\cite[239--242]{JonesAntiky}.}

The possibility that even biological phenomena worked on the same principle immediately suggested itself and was eagerly pursued.
\quote{Just as people who imitate the revolutions of the wandering stars by means of certain instruments instill a principle of motion in them and then go away, while [the devices] operate just as if the craftsman was there and overseeing them in everything, I think in the same way each of the parts in the body operates by some succession and reception of motion from the first principle to every part, needing no overseer.\footnote{Galen, {\it On the Use of the Parts} 12.5, \cite[245]{JonesAntiky}.}}
Indeed, ancient medical research put this vision into practice. ``The use of what we should call mechanical ideas to explain organic processes''---such as digestion and other physiological functions---is ``the most prominent feature'' of the work of Erasistratus in medicine,\footnote{\cite[80]{lloydGSAA}.} who also tested his ideas experimentally.\footnote{\cite[85]{lloydGSAA}.}

In sum, the world did not need Galileo to tell them about the mechanical philosophy, since it had been widely regarded as common sense already thousands of years before.

\subsection{Continuity thesis}
\label{continuitythesis}

Many have tried to stress commonalities between Galileo and the Aristotelian philosophers who preceded him. That is to say, they argue for the ``continuity thesis'' which says that the so-called ``Scientific Revolution'' was not a radical or revolutionary break with previous thought.
\quote{Galileo essentially pursued a progressive Aristotelianism [during the first half of his life---the period of] positive growth that laid the foundation …\ for the new sciences.\footnote{\cite[I.350]{WallaceColl}.}}
\quote{A particular school of Renaissance Aristotelians, located at the University of Padua, constructed a very sophisticated methodology for experimental science; …\ Galileo knew this school of thought and built upon its results; …\ this goes a long way toward explaining the birth of early modern science.\footnote{\cite[279]{CohenSciRev}, summarising \cite{RandallThesis}.}}
\quote{The mechanical and physical science of which the present day is so proud comes to us through an uninterrupted sequence of almost imperceptible refinements from the doctrines professed within the Schools of the Middle Ages.\footnote{\cite[9]{DuhemStatics}.}}
\quote{Galileo was clearly the heir of the medieval kinematicists.\footnote{\cite[666]{ClagettMA}.}}
I agree with these authors that ``those great truths for which Galileo …\ received credit'' are not his.\footnote{\cite[17]{DuhemStatics}.} But the notion that they were first conceived in Aristotelian schools of philosophy is wrongheaded.

The argument of these historians is based on a simple logic. First they show that various concepts of ``Galilean'' science are prefigured in earlier sources. Then they want to infer from this that these sources marked the true beginning of the scientific revolution. But in order to draw this inference they need two assumptions: first, that Galileo was the father of modern science; and second, that the Greeks were nowhere near the same accomplishments. These two assumptions are simply taken for granted by these authors, as a matter of common knowledge. But in reality both assumptions are dead wrong, and therefore the inference to the significance of the Aristotelian sources is unwarranted.

The continuity thesis, then, devalues the contributions of Galileo, yet at the same time desperately needs to reassert the traditional view that ``Galileo has a clear and undisputed title as the `father of modern science',''\footnote{\cite[57]{wallace1987}.} since this is what gives them the one point of connection they are able to establish between medieval and modern science. Historians have been able to match up some aspects of Renaissance Aristotelian tradition with the scientific methodology and practice of Galileo. From this they are inclined to infer that Galileo was influenced by this tradition. But there is little direct evidence to this effect. Galileo is said to have found Aristotelian teachings "boring."\footnote{Vincenzo Viviani, {\it Racconto istorico della vita del sig.r Galileo Galileo} (1654), OGG.XIX.602, \cite[5]{GatteiGLife}.} Insofar as there are similarities between him and earlier thinkers who tried to tackle similar questions, this can very well be due largely to common sense rather than direct philosophical influence.

It would be much harder to find substantive links between medieval authors and mathematically competent people such as Kepler or Newton. Since Galileo mostly writes philosophical prose and rarely if ever substantive mathematics, it is much easier to try to construe his works as related to the preceding philosophical tradition. The strength of those alleged continuities in philosophical thought between Galileo and preceding Aristotelian tradition is highly debatable. But it is a moot point in any case unless Galileo is accepted as a founder of modern science. The entire argument stands and falls with this false premiss. Therefore, if one proves, as we have done above, that Galileo was a mediocre scientists of negligible importance to the mathematically competent people who actually achieved the scientific revolution, then the continuity thesis collapses like a house of cards.

The defenders of the continuity thesis are equally ineffectual in establishing the second false premiss of their argument, namely the alleged absence of these ``new'' ideas in Greek thought. In fact, even continuity thesis advocates make no secret of the fact that the medieval tradition was built on ``remnants of Alexandrian science.''\footnote{\cite[75]{DuhemStatics}.} For example, ``although we are left with few monuments from the profound research of the Ancients into the laws of equilibrium, those few are worthy of eternal admiration.''\footnote{\cite[11]{DuhemStatics}.} Obviously, ``masterpieces of Greek science …\ [such as the works of] Pappus, and especially Archimedes, …\ are proof that the deductive method can be applied with as much rigor to the field of mechanics as to the demonstrations of geometry.''\footnote{\cite[149]{DuhemStatics}.} Galileo himself embraced Archimedes as his role model in no uncertain terms. ``As far as genius is concerned, [Galileo] claimed Archimedes had exceeded everybody else, and called him his master.'' ``Galileo claimed it was possible to walk safely, without stumbling, on earth as well as in heaven, as long as we remain in Archimedes's footsteps.''\footnote{Niccolo Gherardini, {\it Vita di Galileo Galilei} (1654), OGG.XIX.645, 637, \cite[159, 145]{GatteiGLife}.}

How, then, can continuity thesis advocates acknowledge these ``masterpieces'' ``worthy of eternal admiration'' from antiquity, yet at the same time attribute the scientific revolution to medieval or renaissance philosophers? By writing off those ancient works as minor technical footnotes to an otherwise thoroughly Aristotelian paradigm. Only if this picture is accepted can any kind of greatness be ascribed to the pre-Galileans, as is evident from passages such as these:
\quote{Some philosophers in medieval universities were teaching ideas about motion and mechanics that were totally non-Aristotelian [and] were consciously based on criticisms of Aristotle's own pronouncements.\footnote{\cite[30]{HallRevSci}.}}
\quote{Admittedly, most of these …\ significant medieval mechanical doctrines …\ were formed within the Aristotelian framework of mechanics. But these medieval doctrines contained within them the seeds of a critical refutation of that mechanics.\footnote{\cite[682]{ClagettMA}.}}
\quote{The medieval mechanics occupied an important middle position between …\ Aristotelian and Newtonian mechanics. …\ [Hence it was] an important link in man's efforts to represent the laws that concern bodies at rest and in movement.\footnote{\cite[670--671]{ClagettMA}.}}
\quote{The impressive set of departures from Aristotelianism achieved by medieval science nevertheless failed to produce genuine efforts to reconstruct, or replace, the Aristotelian world picture.\footnote{\cite[266]{CohenSciRev}.}}
If Aristotle is taken as the baseline, this looks quite impressive indeed. But why should Aristotle be accepted as the default opinion? Aristotle was one particular philosopher who was a nobody in mathematics and lived well before the golden age of Greek science. Medieval and renaissance thinkers indeed mustered up the courage to challenge isolated claims of his teachings almost two thousand years later, while mostly retaining his overall outlook. This does not constitute great open-mindedness and progress. Rather it is a sign of small-mindedness that these people paid so much attention to Aristotle at all in the first place. In my view, it is not so much impressive that they deviated a bit from Aristotle as it is deplorable that they framed so much of what they did relative to Aristotle, even when they disagreed with him. This is very different from post-Aristotelian thought in Greek times, where there is no evidence that any mathematician paid any attention to Aristotle's mechanics.

In any case, ``extravagant claims for the modernity of medieval concepts'' suffer from ``serious defects.''\footnote{\cite[xxi]{ClagettMA}.}
\quote{There was no
such thing as a fourteenth-century science of mechanics
in the sense of a general theory of local motion applicable
throughout nature, and based on a few unified principles.
By searching the literature of late medieval physics for just
those ideas and those pieces of quantitative analysis that
turned out, three centuries later, to be important in
seventeenth-century mechanics, one can find them; and
one can construct …\ a “medieval science of
mechanics” that {\em appears} to form a coherent whole and to
be built on new foundations replacing those of Aristotle's
physics. But this is an illusion, and an anachronistic fiction,
which we are able to construct only because Galileo and
Newton gave us the pattern by which to select the right
pieces and put them together.\footnote{\cite[42--43]{moodygolino}.}}
The main piece of such precursorism is the so-called ``mean speed theorem.'' This trivial theorem states that, in terms of distance covered, a uniformly accelerated motion is equivalent to a constant-speed motion with the same average speed (Figure \ref{velocitytimeparadox}). Some people praise this as an ``impressive'' achievement\footnote{\cite[25]{Crowe}.}---``probably the most outstanding single medieval contribution to the history of physics,'' derived by ``admirable and ingenious'' reasoning\footnote{\cite[56]{GrantMA}.}---even though these authors did nothing with this trivial theorem and only deduced it to illustrate the notion of uniform change abstractly within Aristotelian philosophy. Later the theorem became central in ``Galilean'' mechanics since free fall is uniformly accelerated.\footnote{\S\ref{lawoffall}.} But it ``was, in fact, never applied to motion in fall from rest during the 14th, or even in the 15th century'' (only in the mid-16th century there is a passing remark to this effect within the Aristotelian tradition, ``without any accompanying evidence'').\footnote{\cite[17]{drakeHistFall}.} Let us not radically inflate our esteem for the Middle Ages by anachronistically praising them for pointing out a trivial thing that centuries later took on a significance of which they had no inkling.\footnote{For example, it played no significant role in the thought of the first person to prove it, Oresme. \cite[19]{LimitsPreclasMech}.} Let us instead recognise the theorem for the trifle that it is. Then we shall also not have any need to be surprised when it turns out that Babylonian astronomers assumed it without fanfare thousands of years earlier still.\footnote{\cite{Ossendrijver}.}

In a similar vein we are told that there are ``unmistakeable Jesuit influences in Galileo's work''\footnote{\cite[55]{wallace1987}.}: ``Above all Galileo was intent in following out Clavius's program of applying mathematics to the study of nature and to generating a mathematical physics.''\footnote{\cite[57]{wallace1987}.} The preposterous notion that this was ``Clavius's'' program can only enter one's mind if one only reads philosophy. It was obviously Archimedes's program, except, unlike Clavius, he proved his point by actually carrying it out instead of sermonising about what one ought to do in philosophical prose. Philosophers (ancient and modern alike) have a tendency to place disproportionate value on explaining something conceptually as opposed to actually doing it. After all, that is virtually the definition of philosophy. Hence they praise certain Aristotelians for explaining some supposedly profound principles of scientific method even when ``it is quite clear that [none of them] ever applied his advocated methods to actual scientific problems.''\footnote{\cite[283]{CohenSciRev}.} Descartes---a mathematically creative person---knew better: ``we ought not to believe an alchemist who boasts he has the technique of making gold, unless he is extremely wealthy; and by the same token we should not believe the learned writer who promises new sciences, unless he demonstrates that he has discovered many things that have been unknown up till now.''\footnote{Descartes to Van Hogelande, 1639/1640, \cite[384]{VenBos}.} Unfortunately, such basic common sense is often lacking among historians and philosophers assigning credit for basic principles of the scientific method.

There is a contradiction in the way modern historians try to trace many aspects of the scientific revolution to roots in the middle ages. On the one hand these historians like to claim that the traditional view of the scientific revolution is ahistorical and based on an anachronistic mindset, whereas their own account that sees continuity with the middle ages is more sensitive to how people actually thought at the time itself. Ironically, however, their view, which is supposed to be more true to the historical actors' way of thinking, is actually all the more blatantly at odds with how virtually all leaders of the scientific revolution thought of the middle ages. ``The scientific achievement of the Middle Ages was held in unanimous contempt from Galileo's time onward by those who adhered to the new science. Leibniz' scathing verdict `barbarismus physicus' neatly encapsulates the reigning sentiment.''\footnote{\cite[260]{CohenSciRev}.} This was not for nothing. Leibniz was an erudite scholar well versed in the philosophy of the schools. But he was also an excellent mathematician. The latter enabled him to pass a sound judgement on medieval science.

\subsection{Epilogue}

Galileo could not have asked for better co-conspirators than those modern academia have provided for him. He desperately needs his audience to view Aristotle as the default baseline against which all his works should be evaluated. He desperately needs his audience to be ignorant of mathematical authors. And he's in luck. The default training of historians of science is not higher mathematics and physics, but seminars based on non-mathematical authors such as Aristotle. So the people tasked with being Galileo experts are by design the people most inclined to accept Galileo’s deceit. Pretending that Archimedes doesn’t exist serves both their purposes and Galileo’s, since they share Galileo’s aversion to proper mathematics. For example, the annual bibliography of works on the history of science published by the flagship journal of the History of Science Society consistently lists well over ten times as many works on Aristotle as on Archimedes, and almost as many more on ``Aristotelianism.''\footnote{Number of entries in the subject index of the annual {\it Isis} comprehensive bibliography for the years 2003--2017: Archimedes 42, Aristotle 482, Aristotelianism 339. There has never been an entry on ``Archimedeanism,'' though I for one would welcome it.}

Such are the inclinations and predispositions of the philosophy-trained humanists who dominate the field today. But people steeped in mathematics see the world differently---a fundamental schism in the outlooks of modern historians of science and the historical figures they are trying to understand. If we read authors like Copernicus, Kepler, Galileo, Descartes, and virtually everyone else who made a contribution to the mathematical sciences, we find endless praise for Archimedes and bottomless contempt for Aristotle. With this as our baseline our view of Galileo is radically transformed.

Following Galileo's sentencing by the church, an author of a book to be printed at Florence was told by the Inquisition to change the phrase ``most distinguished Galileo'' into ``Galileo, man of noted name.''\footnote{\cite[414]{drakeGatwork}.} Though I am not generally on the side of the Inquisition, I have come to the conclusion that this particular decree is sound. Instead of ``Galileo, father of modern science,'' we would do better to make it ``Galileo, man of noted name.''

\end{document}